\input amstex
\expandafter\ifx\csname mathdefs.tex\endcsname\relax
  \expandafter\gdef\csname mathdefs.tex\endcsname{}
\else \message{Hey!  Apparently you were trying to
  \string\input{mathdefs.tex} twice.   This does not make sense.} 
\errmessage{Please edit your file (probably \jobname.tex) and remove
any duplicate ``\string\input'' lines}\endinput\fi




\catcode`\X=12\catcode`\@=11

\def\n@wcount{\alloc@0\count\countdef\insc@unt}
\def\n@wwrite{\alloc@7\write\chardef\sixt@@n}
\def\n@wread{\alloc@6\read\chardef\sixt@@n}
\def\r@s@t{\relax}\def\v@idline{\par}\def\@mputate#1/{#1}
\def\l@c@l#1X{\firstpart.#1}\def\gl@b@l#1X{#1}\def\t@d@l#1X{{}}

\def\crossrefs#1{\ifx\all#1\let\tr@ce=\all\else\def\tr@ce{#1,}\fi
   \n@wwrite\cit@tionsout\openout\cit@tionsout=\jobname.cit 
   \write\cit@tionsout{\tr@ce}\expandafter\setfl@gs\tr@ce,}
\def\setfl@gs#1,{\def\@{#1}\ifx\@\empty\let\next=\relax
   \else\let\next=\setfl@gs\expandafter\xdef
   \csname#1tr@cetrue\endcsname{}\fi\next}
\def\m@ketag#1#2{\expandafter\n@wcount\csname#2tagno\endcsname
     \csname#2tagno\endcsname=0\let\tail=\all\xdef\all{\tail#2,}
   \ifx#1\l@c@l\let\tail=\r@s@t\xdef\r@s@t{\csname#2tagno\endcsname=0\tail}\fi
   \expandafter\gdef\csname#2cite\endcsname##1{\expandafter
     \ifx\csname#2tag##1\endcsname\relax?\else\csname#2tag##1\endcsname\fi
     \expandafter\ifx\csname#2tr@cetrue\endcsname\relax\else
     \write\cit@tionsout{#2tag ##1 cited on page \folio.}\fi}
   \expandafter\gdef\csname#2page\endcsname##1{\expandafter
     \ifx\csname#2page##1\endcsname\relax?\else\csname#2page##1\endcsname\fi
     \expandafter\ifx\csname#2tr@cetrue\endcsname\relax\else
     \write\cit@tionsout{#2tag ##1 cited on page \folio.}\fi}
   \expandafter\gdef\csname#2tag\endcsname##1{\expandafter
      \ifx\csname#2check##1\endcsname\relax
      \expandafter\xdef\csname#2check##1\endcsname{}%
      \else\immediate\write16{Warning: #2tag ##1 used more than once.}\fi
      \multit@g{#1}{#2}##1/X%
      \write\t@gsout{#2tag ##1 assigned number \csname#2tag##1\endcsname\space
      on page \number\count0.}%
   \csname#2tag##1\endcsname}}

\def\multit@g#1#2#3/#4X{\def\t@mp{#4}\ifx\t@mp\empty%
      \global\advance\csname#2tagno\endcsname by 1 
      \expandafter\xdef\csname#2tag#3\endcsname
      {#1\number\csname#2tagno\endcsnameX}%
   \else\expandafter\ifx\csname#2last#3\endcsname\relax
      \expandafter\n@wcount\csname#2last#3\endcsname
      \global\advance\csname#2tagno\endcsname by 1 
      \expandafter\xdef\csname#2tag#3\endcsname
      {#1\number\csname#2tagno\endcsnameX}
      \write\t@gsout{#2tag #3 assigned number \csname#2tag#3\endcsname\space
      on page \number\count0.}\fi
   \global\advance\csname#2last#3\endcsname by 1
   \def\t@mp{\expandafter\xdef\csname#2tag#3/}%
   \expandafter\t@mp\@mputate#4\endcsname
   {\csname#2tag#3\endcsname\lastpart{\csname#2last#3\endcsname}}\fi}
\def\t@gs#1{\def\all{}\m@ketag#1e\m@ketag#1s\m@ketag\t@d@l p
\let\realscite\scite
\let\realstag\stag
   \m@ketag\gl@b@l r \n@wread\t@gsin
   \openin\t@gsin=\jobname.tgs \re@der \closein\t@gsin
   \n@wwrite\t@gsout\openout\t@gsout=\jobname.tgs }
\outer\def\localtags{\t@gs\l@c@l}
\outer\def\globaltags{\t@gs\gl@b@l}
\outer\def\newlocaltag#1{\m@ketag\l@c@l{#1}}
\outer\def\newglobaltag#1{\m@ketag\gl@b@l{#1}}

\newif\ifpr@ 
\def\m@kecs #1tag #2 assigned number #3 on page #4.%
   {\expandafter\gdef\csname#1tag#2\endcsname{#3}
   \expandafter\gdef\csname#1page#2\endcsname{#4}
   \ifpr@\expandafter\xdef\csname#1check#2\endcsname{}\fi}
\def\re@der{\ifeof\t@gsin\let\next=\relax\else
   \read\t@gsin to\t@gline\ifx\t@gline\v@idline\else
   \expandafter\m@kecs \t@gline\fi\let \next=\re@der\fi\next}
\def\pretags#1{\pr@true\pret@gs#1,,}
\def\pret@gs#1,{\def\@{#1}\ifx\@\empty\let\n@xtfile=\relax
   \else\let\n@xtfile=\pret@gs \openin\t@gsin=#1.tgs \message{#1} \re@der 
   \closein\t@gsin\fi \n@xtfile}

\newcount\sectno\sectno=0\newcount\subsectno\subsectno=0
\newif\ifultr@local \def\ultralocal{\ultr@localtrue}
\def\firstpart{\number\sectno}
\def\lastpart#1{\ifcase#1 \or a\or b\or c\or d\or e\or f\or g\or h\or 
   i\or k\or l\or m\or n\or o\or p\or q\or r\or s\or t\or u\or v\or w\or 
   x\or y\or z \fi}

\def\resetall{\global\advance\sectno by 1\subsectno=0
   \gdef\firstpart{\number\sectno}\r@s@t}
\def\resetsub{\global\advance\subsectno by 1
   \gdef\firstpart{\number\sectno.\number\subsectno}\r@s@t}
\def\newsection#1\par{\resetall\vskip0pt plus.3\vsize\penalty-250
   \vskip0pt plus-.3\vsize\bigskip\bigskip
   \message{#1}\leftline{\bf#1}\nobreak\bigskip}
\def\subsection#1\par{\ifultr@local\resetsub\fi
   \vskip0pt plus.2\vsize\penalty-250\vskip0pt plus-.2\vsize
   \bigskip\smallskip\message{#1}\leftline{\bf#1}\nobreak\medskip}


\newdimen\marginshift

\newdimen\margindelta
\newdimen\marginmax
\newdimen\marginmin

\def\margininit{       
\marginmax=3 true cm                  
				      
\margindelta=0.1 true cm              
\marginmin=0.1true cm                 
\marginshift=\marginmin
}    

\def\t@gsjj#1,{\def\@{#1}\ifx\@\empty\let\next=\relax\else\let\next=\t@gsjj
   \def\@@{p}\ifx\@\@@\else
   \expandafter\gdef\csname#1cite\endcsname##1{\citejj{##1}}
   \expandafter\gdef\csname#1page\endcsname##1{?}
   \expandafter\gdef\csname#1tag\endcsname##1{\tagjj{##1}}\fi\fi\next}
\newif\ifshowstuffinmargin
\showstuffinmarginfalse
\def\jjtags{\ifx\shlhetal\relax 
  \else
\ifx\shlhetal\undefinedcontrolseq
\else
\showstuffinmargintrue
\ifx\all\relax\else\expandafter\t@gsjj\all,\fi\fi \fi
}

\def\tagjj#1{\realstag{#1}\oldmginpar{\zeigen{#1}}}
\def\citejj#1{\rechnen{#1}\mginpar{\zeigen{#1}}}     

\def\rechnen#1{\expandafter\ifx\csname stag#1\endcsname\relax ??\else
                           \csname stag#1\endcsname\fi}

\newdimen\theight

\def\marginfont{\sevenrm}

\def\trymarginbox#1{\setbox0=\hbox{\marginfont\hskip\marginshift #1}%
		\global\marginshift\wd0 
		\global\advance\marginshift\margindelta}

\def \oldmginpar#1{%
\ifvmode\setbox0\hbox to \hsize{\hfill\rlap{\marginfont\quad#1}}%
\ht0 0cm
\dp0 0cm
\box0\vskip-\baselineskip
\else 
             \vadjust{\trymarginbox{#1}%
		\ifdim\marginshift>\marginmax \global\marginshift\marginmin
			\trymarginbox{#1}%
                \fi
             \theight=\ht0
             \advance\theight by \dp0    \advance\theight by \lineskip
             \kern -\theight \vbox to \theight{\rightline{\rlap{\box0}}%
\vss}}\fi}

\newdimen\upordown
\global\upordown=8pt
\font\tinyfont=cmtt8 
\def\mginpar#1{\smash{\hbox to 0cm{\kern-10pt\raise7pt\hbox{\tinyfont #1}\hss}}}
\def\mginpar#1{{\hbox to 0cm{\kern-10pt\raise\upordown\hbox{\tinyfont #1}\hss}}\global\upordown-\upordown}


\def\t@gsoff#1,{\def\@{#1}\ifx\@\empty\let\next=\relax\else\let\next=\t@gsoff
   \def\@@{p}\ifx\@\@@\else
   \expandafter\gdef\csname#1cite\endcsname##1{\zeigen{##1}}
   \expandafter\gdef\csname#1page\endcsname##1{?}
   \expandafter\gdef\csname#1tag\endcsname##1{\zeigen{##1}}\fi\fi\next}
\def\verbatimtags{\showstuffinmarginfalse
\ifx\all\relax\else\expandafter\t@gsoff\all,\fi}
\def\zeigen#1{\hbox{$\scriptstyle\langle$}#1\hbox{$\scriptstyle\rangle$}}


\def\margintag#1{\ifshowstuffinmargin\oldmginpar{\zeigen{#1}}\fi}

\def\marginplain#1{\ifshowstuffinmargin\mginpar{{#1}}\fi}
\def\marginbf#1{\marginplain{{\bf \ \ #1}}}

\def\(#1){\edef\dot@g{\ifmmode\ifinner(\hbox{\noexpand\etag{#1}})
   \else\noexpand\eqno(\hbox{\noexpand\etag{#1}})\fi
   \else(\noexpand\ecite{#1})\fi}\dot@g}

\newif\ifbr@ck
\def\eat#1{}
\def\[#1]{\br@cktrue[\br@cket#1'X]}
\def\br@cket#1'#2X{\def\temp{#2}\ifx\temp\empty\let\next\eat
   \else\let\next\br@cket\fi
   \ifbr@ck\br@ckfalse\br@ck@t#1,X\else\br@cktrue#1\fi\next#2X}
\def\br@ck@t#1,#2X{\def\temp{#2}\ifx\temp\empty\let\neext\eat
   \else\let\neext\br@ck@t\def\temp{,}\fi
   \def\teemp{#1}\ifx\teemp\empty\else\rcite{#1}\fi\temp\neext#2X}
\def\resetbr@cket{\gdef\[##1]{[\rtag{##1}]}}
\def\references{\resetbr@cket\newsection References\par}

\newtoks\symb@ls\newtoks\s@mb@ls\newtoks\p@gelist\n@wcount\ftn@mber
    \ftn@mber=1\newif\ifftn@mbers\ftn@mbersfalse\newif\ifbyp@ge\byp@gefalse
\def\defm@rk{\ifftn@mbers\n@mberm@rk\else\symb@lm@rk\fi}
\def\n@mberm@rk{\xdef\m@rk{{\the\ftn@mber}}%
    \global\advance\ftn@mber by 1 }
\def\rot@te#1{\let\temp=#1\global#1=\expandafter\r@t@te\the\temp,X}
\def\r@t@te#1,#2X{{#2#1}\xdef\m@rk{{#1}}}
\def\b@@st#1{{$^{#1}$}}\def\str@p#1{#1}
\def\symb@lm@rk{\ifbyp@ge\rot@te\p@gelist\ifnum\expandafter\str@p\m@rk=1 
    \s@mb@ls=\symb@ls\fi\write\f@nsout{\number\count0}\fi \rot@te\s@mb@ls}
\def\byp@ge{\byp@getrue\n@wwrite\f@nsin\openin\f@nsin=\jobname.fns 
    \n@wcount\currentp@ge\currentp@ge=0\p@gelist={0}
    \re@dfns\closein\f@nsin\rot@te\p@gelist
    \n@wread\f@nsout\openout\f@nsout=\jobname.fns }
\def\m@kelist#1X#2{{#1,#2}}
\def\re@dfns{\ifeof\f@nsin\let\next=\relax\else\read\f@nsin to \f@nline
    \ifx\f@nline\v@idline\else\let\t@mplist=\p@gelist
    \ifnum\currentp@ge=\f@nline
    \global\p@gelist=\expandafter\m@kelist\the\t@mplistX0
    \else\currentp@ge=\f@nline
    \global\p@gelist=\expandafter\m@kelist\the\t@mplistX1\fi\fi
    \let\next=\re@dfns\fi\next}
\def\symbols#1{\symb@ls={#1}\s@mb@ls=\symb@ls} 
\def\bigsymbol{\textstyle}
\symbols{\bigsymbol\ast,\dagger,\ddagger,\sharp,\flat,\natural,\star}
\def\ftnumbers{\ftn@mberstrue} \def\ftsymbols{\ftn@mbersfalse}
\def\paginal{\byp@ge} \def\resetftnumbers{\ftn@mber=1}
\def\ftnote#1{\defm@rk\expandafter\expandafter\expandafter\footnote
    \expandafter\b@@st\m@rk{#1}}

\long\def\jump#1\endjump{}
\def\ssum{\mathop{\lower .1em\hbox{$\textstyle\Sigma$}}\nolimits}

\def\qed{\nobreak\kern 1em \vrule height .5em width .5em depth 0em}
\def\newneq{\hbox{\rlap{\hbox to 1\wd9{\hss$=$\hss}}\raise .1em 
   \hbox to 1\wd9{\hss$\scriptscriptstyle/$\hss}}}
\def\subsetne{\setbox9 = \hbox{$\subset$}\mathrel{\hbox{\rlap
   {\lower .4em \newneq}\raise .13em \hbox{$\subset$}}}}
\def\supsetne{\setbox9 = \hbox{$\subset$}\mathrel{\hbox{\rlap
   {\lower .4em \newneq}\raise .13em \hbox{$\supset$}}}}

\def\vbar{\mathchoice{\vrule height6.3ptdepth-.5ptwidth.8pt\kern-.8pt}
   {\vrule height6.3ptdepth-.5ptwidth.8pt\kern-.8pt}
   {\vrule height4.1ptdepth-.35ptwidth.6pt\kern-.6pt}
   {\vrule height3.1ptdepth-.25ptwidth.5pt\kern-.5pt}}
\def\f@dge{\mathchoice{}{}{\mkern.5mu}{\mkern.8mu}}
\def\b@c#1#2{{\rm \mkern#2mu\vbar\mkern-#2mu#1}}
\def\b@b#1{{\rm I\mkern-3.5mu #1}}
\def\b@a#1#2{{\rm #1\mkern-#2mu\f@dge #1}}
\def\bb#1{{\count4=`#1 \advance\count4by-64 \ifcase\count4\or\b@a A{11.5}\or
   \b@b B\or\b@c C{5}\or\b@b D\or\b@b E\or\b@b F \or\b@c G{5}\or\b@b H\or
   \b@b I\or\b@c J{3}\or\b@b K\or\b@b L \or\b@b M\or\b@b N\or\b@c O{5} \or
   \b@b P\or\b@c Q{5}\or\b@b R\or\b@a S{8}\or\b@a T{10.5}\or\b@c U{5}\or
   \b@a V{12}\or\b@a W{16.5}\or\b@a X{11}\or\b@a Y{11.7}\or\b@a Z{7.5}\fi}}

\catcode`\X=11 \catcode`\@=12




\let\thischap\jobname

\def\partof#1{\csname returnthe#1part\endcsname}
\def\chapof#1{\csname returnthe#1chap\endcsname}

\def\setchapter#1,#2,#3;{%
  \expandafter\def\csname returnthe#1part\endcsname{#2}%
  \expandafter\def\csname returnthe#1chap\endcsname{#3}%
}

\setchapter 300a,A,II.A;
\setchapter 300b,A,II.B;
\setchapter 300c,A,II.C;
\setchapter 300d,A,II.D;
\setchapter 300e,A,II.E;
\setchapter 300f,A,II.F;
\setchapter 300g,A,II.G;
\setchapter  E53,B,N;
\setchapter  88r,B,I;
\setchapter  600,B,III;
\setchapter  705,B,IV;
\setchapter  734,B,V;

\def\cprefix#1{
\edef\theotherpart{\partof{#1}}\edef\theotherchap{\chapof{#1}}%
\ifx\theotherpart\thispart
   \ifx\theotherchap\thischap 
    \else 
     \theotherchap%
    \fi
   \else 
     \theotherchap\fi}

\def\sectioncite[#1]#2{%
     \cprefix{#2}#1}

\def\chaptercite#1{Chapter \cprefix{#1}}

\edef\thispart{\partof{\thischap}}
\edef\thischap{\chapof{\thischap}}

\def\lastpage of '#1' is #2.{\expandafter\def\csname lastpage#1\endcsname{#2}}


\def\spuriousreset{}


\expandafter\ifx\csname citeadd.tex\endcsname\relax
\expandafter\gdef\csname citeadd.tex\endcsname{}
\else \message{Hey!  Apparently you were trying to
\string\input{citeadd.tex} twice.   This does not make sense.} 
\errmessage{Please edit your file (probably \jobname.tex) and remove
any duplicate ``\string\input'' lines}\endinput\fi

\sectno=-1   
\localtags
\jjtags
\NoBlackBoxes
\define\mr{\medskip\roster}
\define\sn{\smallskip\noindent}
\define\mn{\medskip\noindent}
\define\bn{\bigskip\noindent}
\define\ub{\underbar}
\define\wilog{\text{without loss of generality}}
\define\ermn{\endroster\medskip\noindent}
\define\dbca{\dsize\bigcap}
\define\dbcu{\dsize\bigcup}
\define \nl{\newline}
\magnification=\magstep 1
\documentstyle{amsppt}

{    
\catcode`@11

\ifx\alicetwothousandloaded@\relax
  \endinput\else\global\let\alicetwothousandloaded@\relax\fi

\gdef\subjclass{\let\savedef@\subjclass
 \def\subjclass##1\endsubjclass{\let\subjclass\savedef@
   \toks@{\def\usualspace{{\rm\enspace}}\eightpoint}%
   \toks@@{##1\unskip.}%
   \edef\thesubjclass@{\the\toks@
     \frills@{{\noexpand\rm2000 {\noexpand\it Mathematics Subject
       Classification}.\noexpand\enspace}}%
     \the\toks@@}}%
  \nofrillscheck\subjclass}
} 


\expandafter\ifx\csname alice2jlem.tex\endcsname\relax
  \expandafter\xdef\csname alice2jlem.tex\endcsname{\the\catcode`@}
\else \message{Hey!  Apparently you were trying to
\string\input{alice2jlem.tex}  twice.   This does not make sense.}
\errmessage{Please edit your file (probably \jobname.tex) and remove
any duplicate ``\string\input'' lines}\endinput\fi

\expandafter\ifx\csname bib4plain.tex\endcsname\relax
  \expandafter\gdef\csname bib4plain.tex\endcsname{}
\else \message{Hey!  Apparently you were trying to \string\input
  bib4plain.tex twice.   This does not make sense.}
\errmessage{Please edit your file (probably \jobname.tex) and remove
any duplicate ``\string\input'' lines}\endinput\fi

\def\renewcommand{\newcommand}	       
\edef\cite{\the\catcode`@}%
\catcode`@ = 11
\let\@oldatcatcode = \cite
\chardef\@letter = 11
\chardef\@other = 12
%
%
%
%
\def\@innerdef#1#2{\edef#1{\expandafter\noexpand\csname #2\endcsname}}%
%
%
\@innerdef\@innernewcount{newcount}%
\@innerdef\@innernewdimen{newdimen}%
\@innerdef\@innernewif{newif}%
\@innerdef\@innernewwrite{newwrite}%
%
%
%
\def\@gobble#1{}%
%
%
%
\ifx\inputlineno\@undefined
   \let\@linenumber = \empty 
\else
   \def\@linenumber{\the\inputlineno:\space}%
\fi
%
%
%
\def\@futurenonspacelet#1{\def\cs{#1}%
   \afterassignment\@stepone\let\@nexttoken=
}%
\begingroup 
\def\\{\global\let\@stoken= }%
\\ 
\endgroup
\def\@stepone{\expandafter\futurelet\cs\@steptwo}%
\def\@steptwo{\expandafter\ifx\cs\@stoken\let\@@next=\@stepthree
   \else\let\@@next=\@nexttoken\fi \@@next}%
\def\@stepthree{\afterassignment\@stepone\let\@@next= }%
%
%
%
\def\@getoptionalarg#1{%
   \let\@optionaltemp = #1%
   \let\@optionalnext = \relax
   \@futurenonspacelet\@optionalnext\@bracketcheck
}%
%
%
\def\@bracketcheck{%
   \ifx [\@optionalnext
      \expandafter\@@getoptionalarg
   \else
      \let\@optionalarg = \empty
      \expandafter\@optionaltemp
   \fi
}%
\def\@@getoptionalarg[#1]{%
   \def\@optionalarg{#1}%
   \@optionaltemp
}%
%
%
%
\def\@nnil{\@nil}%
\def\@fornoop#1\@@#2#3{}%
\def\@for#1:=#2\do#3{%
   \edef\@fortmp{#2}%
   \ifx\@fortmp\empty \else
      \expandafter\@forloop#2,\@nil,\@nil\@@#1{#3}%
   \fi
}%
\def\@forloop#1,#2,#3\@@#4#5{\def#4{#1}\ifx #4\@nnil \else
       #5\def#4{#2}\ifx #4\@nnil \else#5\@iforloop #3\@@#4{#5}\fi\fi
}%
\def\@iforloop#1,#2\@@#3#4{\def#3{#1}\ifx #3\@nnil
       \let\@nextwhile=\@fornoop \else
      #4\relax\let\@nextwhile=\@iforloop\fi\@nextwhile#2\@@#3{#4}%
}%
%
%
%
\@innernewif\if@fileexists
\def\@testfileexistence{\@getoptionalarg\@finishtestfileexistence}%
\def\@finishtestfileexistence#1{%
   \begingroup
      \def\extension{#1}%
      \immediate\openin0 =
         \ifx\@optionalarg\empty\jobname\else\@optionalarg\fi
         \ifx\extension\empty \else .#1\fi
         \space
      \ifeof 0
         \global\@fileexistsfalse
      \else
         \global\@fileexiststrue
      \fi
      \immediate\closein0
   \endgroup
}%
%
%
%
%
\def\bibliographystyle#1{%
   \@readauxfile
   \@writeaux{\string\bibstyle{#1}}%
}%
\let\bibstyle = \@gobble
%
%
\let\bblfilebasename = \jobname
\def\bibliography#1{%
   \@readauxfile
   \@writeaux{\string\bibdata{#1}}%
   \@testfileexistence[\bblfilebasename]{bbl}%
   \if@fileexists
      \nobreak
      \@readbblfile
   \fi
}%
\let\bibdata = \@gobble
%
%
\def\nocite#1{%
   \@readauxfile
   \@writeaux{\string\citation{#1}}%
}%
\@innernewif\if@notfirstcitation
%
%
\def\cite{\@getoptionalarg\@cite}%
%
%
\def\@cite#1{%
   \let\@citenotetext = \@optionalarg
   \printcitestart
   \nocite{#1}%
   \@notfirstcitationfalse
   \@for \@citation :=#1\do
   {%
      \expandafter\@onecitation\@citation\@@
   }%
   \ifx\empty\@citenotetext\else
      \printcitenote{\@citenotetext}%
   \fi
   \printcitefinish
}%
\newif\ifweareinprivate
\weareinprivatetrue
\ifx\shlhetal\undefinedcontrolseq\weareinprivatefalse\fi
\ifx\shlhetal\relax\weareinprivatefalse\fi
\def\@onecitation#1\@@{%
   \if@notfirstcitation
      \printbetweencitations
   \fi
   \expandafter \ifx \csname\@citelabel{#1}\endcsname \relax
      \if@citewarning
         \message{\@linenumber Undefined citation `#1'.}%
      \fi
     \ifweareinprivate
      \expandafter\gdef\csname\@citelabel{#1}\endcsname{%
\strut 
\vadjust{\vskip-\dp\strutbox
\vbox to 0pt{\vss\parindent0cm \leftskip=\hsize 
\advance\leftskip3mm
\advance\hsize 4cm\strut\openup-4pt 
\rightskip 0cm plus 1cm minus 0.5cm ?  #1 ?\strut}}
         {\tt
            \escapechar = -1
            \nobreak\hskip0pt\pfeilsw
            \expandafter\string\csname#1\endcsname
             \pfeilso
            \nobreak\hskip0pt
         }%
      }%
     \else  
      \expandafter\gdef\csname\@citelabel{#1}\endcsname{%
            {\tt\expandafter\string\csname#1\endcsname}
      }%
     \fi  
   \fi
   \csname\@citelabel{#1}\endcsname
   \@notfirstcitationtrue
}%
%
%
\def\@citelabel#1{b@#1}%
%
%
\def\@citedef#1#2{\expandafter\gdef\csname\@citelabel{#1}\endcsname{#2}}%
%
%
%
\def\@readbblfile{%
   \ifx\@itemnum\@undefined
      \@innernewcount\@itemnum
   \fi
   \begingroup
      \def\begin##1##2{%
         \setbox0 = \hbox{\biblabelcontents{##2}}%
         \biblabelwidth = \wd0
      }%
      \def\end##1{}
      %
      %
      \@itemnum = 0
      \def\bibitem{\@getoptionalarg\@bibitem}%
      \def\@bibitem{%
         \ifx\@optionalarg\empty
            \expandafter\@numberedbibitem
         \else
            \expandafter\@alphabibitem
         \fi
      }%
      \def\@alphabibitem##1{%
         \expandafter \xdef\csname\@citelabel{##1}\endcsname {\@optionalarg}%
         \ifx\biblabelprecontents\@undefined
            \let\biblabelprecontents = \relax
         \fi
         \ifx\biblabelpostcontents\@undefined
            \let\biblabelpostcontents = \hss
         \fi
         \@finishbibitem{##1}%
      }%
      \def\@numberedbibitem##1{%
         \advance\@itemnum by 1
         \expandafter \xdef\csname\@citelabel{##1}\endcsname{\number\@itemnum}%
         \ifx\biblabelprecontents\@undefined
            \let\biblabelprecontents = \hss
         \fi
         \ifx\biblabelpostcontents\@undefined
            \let\biblabelpostcontents = \relax
         \fi
         \@finishbibitem{##1}%
      }%
      \def\@finishbibitem##1{%
         \biblabelprint{\csname\@citelabel{##1}\endcsname}%
         \@writeaux{\string\@citedef{##1}{\csname\@citelabel{##1}\endcsname}}%
         \ignorespaces
      }%
      %
      %
      \let\em = \bblem
      \let\newblock = \bblnewblock
      \let\sc = \bblsc
      \frenchspacing
      \clubpenalty = 4000 \widowpenalty = 4000
      \tolerance = 10000 \hfuzz = .5pt
      \everypar = {\hangindent = \biblabelwidth
                      \advance\hangindent by \biblabelextraspace}%
      \bblrm
      \parskip = 1.5ex plus .5ex minus .5ex
      \biblabelextraspace = .5em
      \bblhook
      \input \bblfilebasename.bbl
   \endgroup
}%
%
%
\@innernewdimen\biblabelwidth
\@innernewdimen\biblabelextraspace
%
%
%
\def\biblabelprint#1{%
   \noindent
   \hbox to \biblabelwidth{%
      \biblabelprecontents
      \biblabelcontents{#1}%
      \biblabelpostcontents
   }%
   \kern\biblabelextraspace
}%
%
%
%
\def\biblabelcontents#1{{\bblrm [#1]}}%
%
%
\def\bblrm{\rm}%
%
%
\def\bblem{\it}%
%
%
\def\bblsc{\ifx\@scfont\@undefined
              \font\@scfont = cmcsc10
           \fi
           \@scfont
}%
%
%
\def\bblnewblock{\hskip .11em plus .33em minus .07em }%
%
%
\let\bblhook = \empty
%
%
%
\def\printcitestart{[}
\def\printcitefinish{]}
\def\printbetweencitations{, }
\def\printcitenote#1{, #1}
%
%
%
\let\citation = \@gobble
%
%
%
\@innernewcount\@numparams
%
%
\def\newcommand#1{%
   \def\@commandname{#1}%
   \@getoptionalarg\@continuenewcommand
}%
%
%
\def\@continuenewcommand{%
   \@numparams = \ifx\@optionalarg\empty 0\else\@optionalarg \fi \relax
   \@newcommand
}%
%
%
\def\@newcommand#1{%
   \def\@startdef{\expandafter\edef\@commandname}%
   \ifnum\@numparams=0
      \let\@paramdef = \empty
   \else
      \ifnum\@numparams>9
         \errmessage{\the\@numparams\space is too many parameters}%
      \else
         \ifnum\@numparams<0
            \errmessage{\the\@numparams\space is too few parameters}%
         \else
            \edef\@paramdef{%
               \ifcase\@numparams
                  \empty  No arguments.
               \or ####1%
               \or ####1####2%
               \or ####1####2####3%
               \or ####1####2####3####4%
               \or ####1####2####3####4####5%
               \or ####1####2####3####4####5####6%
               \or ####1####2####3####4####5####6####7%
               \or ####1####2####3####4####5####6####7####8%
               \or ####1####2####3####4####5####6####7####8####9%
               \fi
            }%
         \fi
      \fi
   \fi
   \expandafter\@startdef\@paramdef{#1}%
}%
%
%
%
%
\def\@readauxfile{%
   \if@auxfiledone \else 
      \global\@auxfiledonetrue
      \@testfileexistence{aux}%
      \if@fileexists
         \begingroup
            \endlinechar = -1
            \catcode`@ = 11
            \input \jobname.aux
         \endgroup
      \else
         \message{\@undefinedmessage}%
         \global\@citewarningfalse
      \fi
      \immediate\openout\@auxfile = \jobname.aux
   \fi
}%
%
%
\newif\if@auxfiledone
\ifx\noauxfile\@undefined \else \@auxfiledonetrue\fi
%
%
%
%
\@innernewwrite\@auxfile
\def\@writeaux#1{\ifx\noauxfile\@undefined \write\@auxfile{#1}\fi}%
%
%
%
\ifx\@undefinedmessage\@undefined
   \def\@undefinedmessage{No .aux file; I won't give you warnings about
                          undefined citations.}%
\fi
%
%
\@innernewif\if@citewarning
\ifx\noauxfile\@undefined \@citewarningtrue\fi
%
%
%
\catcode`@ = \@oldatcatcode

\def\pfeilso{\leavevmode
            \vrule width 1pt height9pt depth 0pt\relax
           \vrule width 1pt height8.7pt depth 0pt\relax
           \vrule width 1pt height8.3pt depth 0pt\relax
           \vrule width 1pt height8.0pt depth 0pt\relax
           \vrule width 1pt height7.7pt depth 0pt\relax
            \vrule width 1pt height7.3pt depth 0pt\relax
            \vrule width 1pt height7.0pt depth 0pt\relax
            \vrule width 1pt height6.7pt depth 0pt\relax
            \vrule width 1pt height6.3pt depth 0pt\relax
            \vrule width 1pt height6.0pt depth 0pt\relax
            \vrule width 1pt height5.7pt depth 0pt\relax
            \vrule width 1pt height5.3pt depth 0pt\relax
            \vrule width 1pt height5.0pt depth 0pt\relax
            \vrule width 1pt height4.7pt depth 0pt\relax
            \vrule width 1pt height4.3pt depth 0pt\relax
            \vrule width 1pt height4.0pt depth 0pt\relax
            \vrule width 1pt height3.7pt depth 0pt\relax
            \vrule width 1pt height3.3pt depth 0pt\relax
            \vrule width 1pt height3.0pt depth 0pt\relax
            \vrule width 1pt height2.7pt depth 0pt\relax
            \vrule width 1pt height2.3pt depth 0pt\relax
            \vrule width 1pt height2.0pt depth 0pt\relax
            \vrule width 1pt height1.7pt depth 0pt\relax
            \vrule width 1pt height1.3pt depth 0pt\relax
            \vrule width 1pt height1.0pt depth 0pt\relax
            \vrule width 1pt height0.7pt depth 0pt\relax
            \vrule width 1pt height0.3pt depth 0pt\relax}

\def\pfeilsw{ \leavevmode 
            \vrule width 1pt height0.3pt depth 0pt\relax
            \vrule width 1pt height0.7pt depth 0pt\relax
            \vrule width 1pt height1.0pt depth 0pt\relax
            \vrule width 1pt height1.3pt depth 0pt\relax
            \vrule width 1pt height1.7pt depth 0pt\relax
            \vrule width 1pt height2.0pt depth 0pt\relax
            \vrule width 1pt height2.3pt depth 0pt\relax
            \vrule width 1pt height2.7pt depth 0pt\relax
            \vrule width 1pt height3.0pt depth 0pt\relax
            \vrule width 1pt height3.3pt depth 0pt\relax
            \vrule width 1pt height3.7pt depth 0pt\relax
            \vrule width 1pt height4.0pt depth 0pt\relax
            \vrule width 1pt height4.3pt depth 0pt\relax
            \vrule width 1pt height4.7pt depth 0pt\relax
            \vrule width 1pt height5.0pt depth 0pt\relax
            \vrule width 1pt height5.3pt depth 0pt\relax
            \vrule width 1pt height5.7pt depth 0pt\relax
            \vrule width 1pt height6.0pt depth 0pt\relax
            \vrule width 1pt height6.3pt depth 0pt\relax
            \vrule width 1pt height6.7pt depth 0pt\relax
            \vrule width 1pt height7.0pt depth 0pt\relax
            \vrule width 1pt height7.3pt depth 0pt\relax
            \vrule width 1pt height7.7pt depth 0pt\relax
            \vrule width 1pt height8.0pt depth 0pt\relax
            \vrule width 1pt height8.3pt depth 0pt\relax
            \vrule width 1pt height8.7pt depth 0pt\relax
            \vrule width 1pt height9pt depth 0pt\relax
      }


\def\widestnumber#1#2{}

\def\citewarning#1{\ifx\shlhetal\relax 
    \else
    \par{#1}\par
    \fi
}

\def\rm{\fam0 \tenrm}

\def\fakesubhead#1\endsubhead{\bigskip\noindent{\bf#1}\par}



%
%
%

%

\font\textrsfs=rsfs10
\font\scriptrsfs=rsfs7
\font\scriptscriptrsfs=rsfs5

\newfam\rsfsfam
\textfont\rsfsfam=\textrsfs
\scriptfont\rsfsfam=\scriptrsfs
\scriptscriptfont\rsfsfam=\scriptscriptrsfs

\edef\oldcatcodeofat{\the\catcode`\@}
\catcode`\@11

\def\Cal@@#1{\noaccents@ \fam \rsfsfam #1}

\catcode`\@\oldcatcodeofat


\expandafter\ifx \csname margininit\endcsname \relax\else\margininit\fi

\long\def\red#1\endred{}
\long\def\green#1\endgreen{}
\long\def\blue#1\endblue{}
\long\def\private#1\endprivate{}

\def\endred{ \unmatched endred! }
\def\endgreen{ \unmatched endgreen! }
\def\endblue{ \unmatched endblue! }
\def\endprivate{ \unmatched endprivate! }

\ifx\latexcolors\undefinedcs\def\latexcolors{}\fi

\def\emptycs{}
\def\evaluatelatexcolors{%
        \ifx\latexcolors\emptycs\else
        \expandafter\xxevaluate\latexcolors\xxfertig\evaluatelatexcolors\fi}
\def\xxevaluate#1,#2\xxfertig{\setupthiscolor{#1}%
        \def\latexcolors{#2}}


\font\smallfont=cmsl7
\def\rutgerscolor{\ifmmode\else\endgraf\fi\smallfont
\advance\leftskip0.5cm\relax}
\def\setupthiscolor#1{\edef\tmptmpcs{\noexpand\bgroup\noexpand\rutgerscolor
\noexpand\def\noexpand\currentcolor{#1}%
\noexpand}%
\expandafter\let\csname#1\endcsname\tmptmpcs
\def\tmptmpcs{\checkColorUnmatched{#1}\popthecolor}
\expandafter\let\csname end#1\endcsname\tmptmpcs}

\def\checkColorUnmatched#1{\def\expectcolor{#1}%
    \ifx\expectcolor\currentcolor   
    \else \edef\failhere{\noexpand\tryingToClose '\currentcolor' with end\expectcolor}\failhere\fi}

\def\currentcolor{???}

\def\popthecolor{\ifmmode\else\endgraf\fi\egroup}

\expandafter\def\csname#1\endcsname{}

\evaluatelatexcolors

 \let\outerhead\head
 \def\head{\innerhead}
 \let\innerhead\outerhead

 \let\outersubhead\subhead
 \def\subhead{\innersubhead}
 \let\innersubhead\outersubhead

 \let\outersubsubhead\subsubhead
 \def\subsubhead{\innersubsubhead}
 \let\innersubsubhead\outersubsubhead

 \let\outerproclaim\proclaim
 \def\proclaim{\innerproclaim}
 \let\innerproclaim\outerproclaim

 %
 %
 %
 %

\def\demo#1{\medskip\noindent{\it #1.\/}}
\def\enddemo{\smallskip}

\def\remark#1{\medskip\noindent{\it #1.\/}}
\def\endremark{\smallskip}

\pageheight{8.5truein}
\topmatter
\title{Abstract elementary classes near $\aleph_1$} \endtitle
\rightheadtext{ Classification of NE classes}
\author {Saharon Shelah \thanks {\null\newline I would like to thank 
Alice Leonhardt for the beautiful typing. \null\newline
 This research was partially supported by the United States Israel
Binational Science Foundation (BSF) and the NSF. Publication 88r.} 
\endthanks} \endauthor  

\affil{The Hebrew University of Jerusalem \\
Einstein Institute of Mathematics \\
Edmond J. Safra Campus, Givat Ram \\
Jerusalem 91904, Israel
 \medskip
 University of Wisconsin \\
 Madison, Wisconsin \\
 \medskip
 Institute of Advanced Studies \\
 Jerusalem, Israel
\medskip
Department of Mathematics \\ 
Hill Center-Busch Campus \\
Rutgers, The State University of New Jersey \\
110 Frelinghuysen Road \\
Piscataway, NJ 08854-8019  USA} \endaffil
\medskip

\abstract  We prove in ZFC, no $\psi \in L_{\omega_1,\omega}[\bold Q]$
have unique model of uncountable cardinality, this confirms the Baldwin
conjecture.  But we analyze this in more general terms.  We introduce
and investigate a.e.c. and also versions of limit models, and prove
some basic properties like representation by PC class, for any a.e.c.
For PC$_{\aleph_0}$-representable a.e.c. we investigate the
conclusion of having not too many non-isomorphic models in $\aleph_1$
and $\aleph_2$, but have to assume $2^{\aleph_0} < 2^{\aleph_1}$ and
even $2^{\aleph_1} < 2^{\aleph_2}$.
\endabstract
\endtopmatter
\document  
 
\pretags{300a,300b,300c,300d,300e,300f,300g,300x,300y,300z,600,705,E53,734}
\newpage

\head {\S0 Introduction} \endhead  \resetall \sectno=0
 \spuriousreset
\bigskip

In \cite{Sh:48}, proving a conjecture of Baldwin, we show
that ($\bold Q$ here stands for the quantifier $\bold
Q^{\text{car}}_{\ge \aleph_1}$, there are uncountably many)
\mr
\item "{$(*)_1$}"  no $\psi \in \Bbb L_{\omega_1,\omega}(\bold Q)$ has a
unique uncountable model up to isomorphism 
\ermn
by showing that
\mr
\item "{$(*)_2$}"  categoricity (of $\psi \in \Bbb L_{\omega_1,\omega}
(\bold Q))$ in $\aleph_1$ implies the existence of a model of $\psi$ of
cardinality $\aleph_2$ (so $\psi$ has $\ge 2$ non-isomorphism models).
\ermn
Unfortunately, this was not proved in ZFC because diamond on
$\aleph_1$ was assumed.  In \cite{Sh:87a} and \cite{Sh:87b} this set
theoretic assumption was weakened to
$2^{\aleph_0} < 2^{\aleph_1}$; here we shall prove it in ZFC (see \S3).
(However, for getting the conclusion from the weaker model theoretic assumption
$\dot I(\aleph_1,\psi) < 2^{\aleph_1}$ as there, we still need $2^{\aleph_0}
< 2^{\aleph_1}$).

The main result of \cite{Sh:87a}, \cite{Sh:87b} was:
\mr
\item "{$(*)_3$}"  if $n > 0,2^{\aleph_0} < 2^{\aleph_1} < \ldots <
2^{\aleph_n},\psi \in \Bbb L_{\omega_1,\omega},1 \le \dot I(\aleph_\ell,\psi)
< \mu_{\text{wd}}(\aleph_\ell)$ for $\ell \le n,\ell \ge 1$ (where
$\mu_{\text{wd}}(\aleph_\ell)$ is usually $2^{\aleph_\ell}$ and always $>
2^{\aleph_{\ell-1}}$, see \scite{88r-0.wD} below) 
\ub{then} $\psi$ has a model of cardinality $\aleph_{n+1}$
\sn
\item "{$(*)_4$}"  if $2^{\aleph_0} < 2^{\aleph_1} < \ldots <
2^{\aleph_n} < 2^{\aleph_{n+1}} < \ldots$ and
$\psi \in \Bbb L_{\omega_1,\omega},1 \le \dot I(\aleph_\ell,\psi) <
\mu_{\text{wd}}(\aleph_\ell)$ for 
$\ell < \omega$ \ub{then} $\psi$ has a model in
every infinite cardinal (and satisfies \L os Conjecture), (note that
$(*)_3$ for $n=1$, assuming $\diamondsuit_{\aleph_1}$ was proved in
\cite{Sh:48}).
\ermn
In $(*)_4$, it is proved that \wilog \, ${\frak K}$ is excellent; this
means in particular that $K$ is the class of atomic models of some
countable first order $T$.  The point is that an excellent 
class ${\frak K}$ is similar to the class of
models of an $\aleph_0$-stable first order $T$.  
In particular the set of relevant
types, ${\bold S}_{\frak K}(A,M)$ is defined as $\{p(x):p(x)$ a complete type
over $A$ in $M$ in the first order sense such that $p \restriction B$ is
isolated for every finite $B \subseteq A\}$.  \ub{But} we better
restrict ourselves to ``nice $A$", that is $A$ which is 
the universe of some $N \prec
M$ or $A = N_1 \cup N_2$ where $N_0,N_1,N_2$ are in stable amalgamation
or $\cup\{N_u:u \in {\Cal P} \subseteq {\Cal P}(n)\}$ for some (so
called) stable
system $\langle N_u:u \in {\Cal P} \rangle$: (for stable such systems
in the stable first order 
case see \cite[XII,\S5]{Sh:c}).  So types are quite like
the first order case.  In particular we say $M \in {\frak K}$ is
$\lambda$-full if $p \in {\bold S}_{\frak K}(A,M),A$ as above, $|A| < \lambda$
implies $p$ is realized in $M$; this is the replacement of
$\lambda$-saturated for that context.

Why in \cite{Sh:87a} and \cite{Sh:87b}, $\psi$ was assumed to be just in
$\Bbb L_{\omega_1,\omega}$ and not more generally in
$\Bbb L_{\omega_1,\omega}(\bold Q)$?  Mainly because we feel that in
\cite{Sh:48}, the logic $\Bbb L_{\omega_1,\omega}(\bold Q)$ was incidental.
We delay the search for the right context to this sequel.  So here we
are working in a.e.c., ``abstract elementary class" (so no logic is present in
the context) whose main feature is the absence of amalgamation, it is
${\frak K} = (K,\le_{\frak K})$ where 
$\le_{\frak K}$ the ``abstract" notion of elementary submodel.  So if
${\Cal L}$ is a fragment of $\Bbb L_{\infty,\omega}(\tau)$ (for a fixed
vocabulary), $T \subseteq {\Cal L}$ a theory included in ${\Cal L}$,
and we let $K = \{M:M \models T\},M
\le_{\frak K} N$ if and only if $M \prec_{\Cal L} N$, we get such a
class; if ${\Cal L}$ is countable then ${\frak K}$ has L.S. number $\aleph_0$.
So the class of models of
$\psi \in \Bbb L_{\omega_1,\omega}(\bold Q)$ is not represented
directly, but can be with minor adaptation; see \scite{88r-3.9}(2).  
Surprisingly (and by not so hard proof), every a.e.c. ${\frak K}$ can be 
represented as a pseudo elementary
class \ub{if} we allow omitting types, (see \scite{88r-1.8}).  We introduce a
replacement for saturated models (for stable first order $T$)
and full models (for excellent classes, see \cite{Sh:87a} and \cite{Sh:87b}): 
limit models; really several variants of this notion.  
See Definition \scite{88r-3.1}.  The strongest and most
important variant is ``$M \in K_\lambda$ superlimit" which means: $M$
is universal (under $\le_{\frak K}$)
$(\exists N)(M \le_{\frak K} N \wedge M \ne N)$ and if $M_i \cong
M$ for $i < \delta \le \|M\|$ and $M_i$ is $\le_{\frak K}$-increasing then
$\dbcu_{i < \delta} M_i \cong M$.  If we restrict ourselves to
$\delta$'s of cofinality $\kappa$ we get $(\lambda,\kappa)$-superlimit.
Such $M$ exists for a first order $T$ for some pairs
$\lambda,\kappa$.  In particular (see more in \cite{Sh:868})
\mr
\item "{$(*)_5$}"  for every $\lambda \ge 2^{|T|} + \beth_\omega$, 
a superlimit model of
$T$ of cardinality $\lambda$ exists if and only if $T$ is superstable
(by \cite[3.1]{Sh:868}).
\ermn
Moreover
\mr
\item "{$(*)_6$}"   ``almost always"; for 
$\lambda \ge 2^{|T|} + \kappa,\kappa = 
\text{ cf}(\kappa)$ (for simplicity) we have: \nl
a $(\lambda,\kappa)$-superlimit model exists iff $T$
is stable in $\lambda \and \kappa \ge \kappa(T)$ or $\lambda =
\lambda^{< \kappa}$.
\ermn
But we can prove something under those circumstances: if
$K$ is categorical in $\lambda$ or just have a superlimit model $M^*$
in $\lambda$, but the $\lambda$-amalgamation property fails for $M^*$
and $2^\lambda < 2^{\lambda^+}$ \ub{then} $\dot I(\lambda^+,K) =
2^{\lambda^+}$ (see \scite{88r-3.5}).  With some reasonable
restrictions on $\lambda$
and $K$, we can prove e.g. $\dot I(\lambda,K) = \dot I(\lambda^+,K) = 1
\Rightarrow \dot I(\lambda^{++},K) \ge 1$, (see \scite{88r-3.7}, \scite{88r-3.8}).

However, our long term main aim was to do the parallel of \cite{Sh:87a} and
\cite{Sh:87b} in the present context, i.e., for an a.e.c. ${\frak K}$ and 
it is natural to assume ${\frak K}$ is PC$_{\aleph_0}$, here we
prepare the ground.

Sections 4,5 present work toward this goal (\S5 assuming $2^{\aleph_0}
< 2^{\aleph_1}$; \S4 without it).  We should note that dealing with
superlimit models rather than full ones make problems, as well as the
fact that the class is not necessarily elementary in some reasonable
logics.  Because of the second we were driven to use formulas which
hold ``generically", are ``forced" instead of are satisfied, and
``the type $\bar a$ materialize" instead of realize and 
${\text{\rm gtp\/}}(\bar a,N,M)$ instead of tp$(\bar a,N,M)$.  We also
(necessarily) encounter the case $|D| = \aleph_1$.  Because
of the first, the scenario for getting a full model in $\aleph_1$
(which can be adapted to $(\aleph_1,\{\aleph_1\})$-superlimit - see
\scite{88r-5.9}) does not seem to be enough for getting superlimit models
in $\aleph_1$ (see \scite{88r-5.24}).

We had felt that arriving at enough conclusions on the models of
cardinality $\aleph_1$ to start dealing with models of cardinality
$\aleph_2$, will be a strong indication that we can complete the
generalization of \cite{Sh:87a} and \cite{Sh:87b}, so 
getting superlimits in $\aleph_1$ is
the culmination of this paper and a natural stopping point.  Trying to
do the rest
(of the parallel to \cite{Sh:87a} and \cite{Sh:87b}) was delayed.
\nl
Much remains to be done,
\bn
\margintag{88r-0.0}\ub{\stag{88r-0.0} Tasks}:
\sn
1) Proving $(*)_3,(*)_4$ in our context. \nl
2) Parallel results in ZFC; e.g. prove $(*)_3$ for $n=1,2^{\aleph_0} =
2^{\aleph_1}$. \nl
Note that if $2^{\aleph_0} = 2^{\aleph_1}$, assuming $1 \le
\dot I(\aleph_1,K) < 2^{\aleph_1}$ give really less model theoretic
consequences, as new phenomena arise (see \S6).  See \S4 (and 
its concluding remarks). \nl
3) Construct examples; e.g. $K$ (or $\psi \in \Bbb L_{\omega_1,\omega})$,
categorical in $\aleph_0,\aleph_1,\dotsc,\aleph_n$ but not in
$\aleph_{n+1}$.
\nl
4) If ${\frak K}$ is PC$_\lambda$, categorical in $\lambda,\lambda^+$, does it
necessarily have a model in $\lambda^{++}$?

See the book's introduction on the progress on those problems.
This is a revised version of \cite{Sh:88} which continues
\cite{Sh:87a}, \cite{Sh:87b} but do not use them.  The paper
\cite{Sh:88} and the present chapter relies on
\cite{Sh:48} only when deducing results on $\psi \in \Bbb
L_{\omega_1,\omega}(\bold Q)$; it improves 
some of its early results and extends the context.
The work on \cite{Sh:88} was done in 1977, and a preprint 
was circulated.  Before the paper had appeared, a 
user-friendly expository article of 
Makowsky \cite{Mw85a} represent, give background
and explain the easy parts of the paper and the author have corrected
and replaced some proofs and added mainly \S6.

We thank Rami Grossberg for lots of work in the early eighties
on previous versions, which
improved this paper, and the writing up of an earlier version of \S6
and Assaf Hasson on helpful comments in 2002 and Alex Usvyatsov for very
careful reading, corrections and comments and Adi Jarden and Alon
Siton on help in the final stages.
\bn
\centerline {$* \qquad * \qquad *$}
\bn
On history and background on
$\Bbb L_{\omega_1,\omega},\Bbb L_{\infty,\omega}$ and the quantifier
$\bold Q$ see \cite{Ke71}.  On $(D,\lambda)$-sequence-homogeneous 
(which \scite{88r-2.1} - \scite{88r-2.4} here
generalized) see Keisler and Morley \cite{KM67}, this is defined in
\scite{88r-2.2}(5), and \scite{88r-2.4} is from there.  
Theorem \scite{88r-3.5} is similar to \cite[2.7]{Sh:87a}
and \cite[6.3]{Sh:87b}.  
\bigskip

\remark{Remark}  On non-splitting used here in \scite{88r-5.4} see \cite{Sh:3},
\cite[Ch.I,Def.2.6,p.11]{Sh:c} or \cite{Sh:48}.
\endremark
\bigskip

By \cite{Ke70} and \cite{Mo70},
\proclaim{\stag{88r-0.1} Claim}  1) Assume that $\psi \in \Bbb
L_{\omega_1,\omega}(\bold Q)$ has a model $M$ in which
$\{{\text{\rm tp\/}}_\Delta(\bar a,\emptyset,M):\bar a \in M\}$ is uncountable
where $\Delta \subseteq \Bbb L_{\omega_1,\omega}(\bold Q)$ is countable,
\ub{then} $\psi$ has $2^{\aleph_1}$ pairwise non-isomorphic models 
of cardinality $\aleph_1$, in fact we can find models
$M_\alpha$ of $\psi$ of cardinality $\aleph_1$ for 
$\alpha < 2^{\aleph_1}$ such that
$\{{\text{\rm tp\/}}_\Delta(a;\emptyset,M_\alpha):a \in M_\alpha\}$ are
pairwise distinct where {\rm tp}$_\Delta(\bar a,A,M) = \{\varphi(\bar x,
\bar b):\varphi(\bar x,\bar y) \in \Delta$ and $M \models
\varphi[\bar a,\bar b]$ and $\bar b \in {}^{\omega >} A\}$. 
\nl
2) If $\psi \in \Bbb L_{\omega_1,\omega}(\bold Q),
\Delta \subseteq \Bbb L_{\omega_1,\omega}(\bold Q)$ is countable 
and $\{{\text{\rm tp\/}}_\Delta(\bar a,\emptyset,M):
\bar a \in {}^{\omega >}M$ and $M$ is a model of $\psi\}$ is 
uncountable, \ub{then} it has cardinality $2^{\aleph_0}$.
\endproclaim
\bigskip

\demo{\stag{88r-0.9} Observation}  Assume ($\tau$ is a vocabulary and)
\mr
\item "{$(a)$}"  $K$ is a family of $\tau$-models of cardinality
$\lambda$ 
\sn
\item "{$(b)$}"  $\mu > \lambda^\kappa$
\sn
\item "{$(c)$}"  $\{(M,\bar a):M \in K$ and $\bar a \in
{}^\kappa M\}$ has $\ge \mu$ members up to isomorphism.
\ermn
\ub{Then} $K$ has $\ge \mu$ models up to isomorphisms (similarly for $=\mu$).
\enddemo
\bigskip

\demo{Proof}  See \cite[VIII,1.3]{Sh:a} or just check.
\enddemo
\bigskip

\proclaim{\stag{88r-0.2} Claim}  1) Assume $\lambda$ is regular
uncountable, $M_0$ is a model with countable vocabulary 
and $T = { \text{\rm Th\/}}_{\Bbb L}(M_0)$, $<$ a binary predicate from
$\tau(T)$ and $(P^{M_0},<^{M_0}) = (\lambda,<)$.  \ub{Then}
every countable model $M$ of $T$ has an end extension, i.e., $M \prec
N$ and $P^M \ne P^N$ and $a \in P^N \wedge b \in P^M \wedge a <^N b
\Rightarrow a \in M$. \nl
2) Moreover, we can further demand $(P^N,<^N)$ is non-well ordered and
we can demand $|P^N| = \aleph_1,(P^N,<^N)$ is $\aleph_1$-like (which
means that it has cardinality $\aleph_1$ but every (proper) initial
segment has cardinality $< \aleph_1$); and we can demand $N$ is countable.
\nl
3) Moreover, we can add the demand that 
in $(P^N,<^N)$ there is a first element in
$P^N \backslash P^M$ and we can add the demand: in
$(P^N,<^N)$, there is no first element in $P^N \backslash P^M$.
\endproclaim
\bigskip

\demo{Proof}  1),2)  Keisler \cite{Ke70}. 
\nl
3) By \cite{Sh:43} and independently Schmerl \cite{Sc76}.
\hfill$\square_{\scite{88r-0.2}}$
\enddemo
\bn
By Devlin-Shelah \cite{DvSh:65}, and 
\cite[Ap,\S1]{Sh:f} (the so-called weak diamond).
\proclaim{\stag{88r-0.wD} Theorem}  Assume that $2^\lambda <
2^{\lambda^+}$. \nl
1)  There is a normal ideal {\rm WDmId}$_{\lambda^+}$
on $\lambda^+$ and $\lambda^+ \notin { \text{\rm WDmId\/}}_{\lambda^+}$, of
course, (the members are called small set) 
such that: if $S \in (\text{\rm WDmId}_{\lambda^+})^+$ 
(e.g., $S = \lambda^+$) and $\bold c:{}^{\lambda^+ >}(\lambda^+) 
\rightarrow \{0,1\}$, \ub{then} there is $\bar \ell
= \langle \ell_\alpha:\alpha < \lambda^+ \rangle \in {}^{\lambda^+} 2$ 
such that for every $\eta \in {}^{\lambda^+}(\lambda^+)$ the 
set $\{\delta \in S:\bold c
(\eta \restriction \delta) = \ell_\alpha\}$ is stationary; we call
$\bar \ell$ a weak diamond sequence (for the colouring $\bold c$ and
the stationary set $S$). \nl
2) $\mu_* = \mu_{\text{wd}}(\lambda^+)$, the cardinal defined by $(*)$
below, is $> 2^\lambda$ (we do not say $\ge 2^{\lambda^+}$!)
\mr
\item "{$(*)$}"   $(\alpha) \quad$  if $\mu < \mu_*$ and
$\bold c_\varepsilon$ for $\varepsilon < \mu$ is as above then we can 
find $\bar \ell$ as in \nl

\hskip25pt part (1) for all the $\bold c_\varepsilon$'s simultaneously
\sn
\item "{${{}}$}"  $(\beta) \quad \mu_*$ is maximal such that clause
$(\alpha)$ holds.
\ermn
3) $\mu_* = \mu_{\text{unif}}(\lambda^+,2^\lambda)$ satisfies
$\mu^{\aleph_0}_* = 2^{\lambda^+}$ and moreover $\lambda \ge
\beth_\omega \Rightarrow \mu_* = 2^\lambda$ where
$\mu_{\text{unif}}(\lambda^+,\chi)$ is the first cardinal $\mu$
such that we can find $\langle \bold c_\alpha:\alpha < \mu\rangle$
such that:
\mr
\item "{$(a)$}"  $\bold c_\alpha$ is a function from
${}^{\lambda^+>}(\lambda^+)$ to $\chi$
\sn
\item "{$(b)$}"  there is no $\rho \in {}^{\lambda^+}\chi$ such that
for every $\alpha < \mu$ for some $\eta \in {}^{\lambda^+}(\lambda^+)$
the set $\{\delta < \lambda:\bold c_\alpha(\eta \restriction \delta)
\ne \rho(\delta)\}$ is stationary (so $\mu_{\text{wd}}(\lambda^+) =
\mu_{\text{unif}}(\lambda^+,2))$. 
\ermn
See more in \cite{Sh:E45}.
\endproclaim
\bn
The following are used in \S2.
\definition{\stag{88r-0.5} Definition}  1) For a regular uncountable
cardinal $\lambda$ let $\check I[\lambda] = \{S \subseteq \lambda$: some
pair $(E,\bar a)$ witnesses $S \in \check I(\lambda)$, see below$\}$. \nl
2) We say that $(E,u)$ is a witness for $S \in \check I[\lambda]$
\ub{if}:
\mr
\item "{$(a)$}"  $E$ is a club of the regular cardinal $\lambda$
\sn
\item "{$(b)$}"  $u = \langle u_\alpha:\alpha < \lambda
\rangle,a_\alpha \subseteq \alpha$ and $\beta \in a_\alpha \Rightarrow
a_\beta = \beta \cap a_\alpha$
\sn
\item "{$(c)$}"  for every $\delta \in E \cap S,u_\delta$ is an
unbounded subset of $\delta$ of order-type $< \delta$ (and $\delta$ is
a limit ordinal).
\endroster
\enddefinition
\bn
By \cite{Sh:420} and \cite{Sh:E12}
\proclaim{\stag{88r-0.6} Claim}  Let $\lambda$ be regular uncountable. \nl
1) If $S \in \check I[\lambda]$ \ub{then} 
we can find a witness $(E,\bar a)$ for $S \in \check I[\lambda]$ such that:
\mr
\item "{$(a)$}"  $\delta \in S \cap E \Rightarrow \,{\text{\rm otp\/}}
(a_\delta) = \,{\text{\rm cf\/}}(\delta)$
\sn
\item "{$(b)$}"  if $\alpha \notin S$ then ${\text{\rm otp\/}}
(a_\alpha) < \,{\text{\rm cf\/}}(\delta)$ for some $\delta \in S \cap E$.
\ermn
2) $S \in \check I[\lambda]$ \ub{iff} there is a pair $(E,\bar{\Cal P})$
such that:
\mr
\item "{$(a)$}"  $E$ is a club of the regular uncountable $\lambda$
\sn
\item "{$(b)$}"  $\bar{\Cal P} = \langle {\Cal P}_\alpha:\alpha <
\lambda\rangle$, where ${\Cal P}_\alpha \subseteq \{u:u \subseteq
\alpha\}$ has cardinality $< \lambda$
\sn
\item "{$(c)$}"  if $\alpha < \beta < \lambda$ and $\alpha \in u \in 
{\Cal P}_\beta$ then $u \cap \alpha \in {\Cal P}_\alpha$
\sn
\item "{$(d)$}"  if $\delta \in E \cap S$ then some $u \in 
{\Cal P}_\delta$ is an unbounded subset of $\delta$ (and $\delta$ is a limit
ordinal).
\endroster
\endproclaim
\newpage

\head {\S1 Axioms and simple properties for classes of models}
\endhead  \resetall 
 \spuriousreset
\bigskip

\demo{\stag{88r-1.1} Context}  Here in \S1-\S5, $\tau$ is a vocabulary, 
$K$ will be a class of $\tau$-models
and $\le_{\frak K}$ a two-place relation on the models in $K$. 
We do not always strictly distinguish between $K$ and ${\frak K} =
(K,\le_{\frak K})$.  We shall assume that $K,\le_{\frak K}$ are fixed;
and usually we assume that it is an a.e.c. (abstract elementary class) which
means that the following axioms hold.  We may write $\le$ or $\le_K$
instead of $\le_{\frak K}$.  For a logic ${\Cal L}$ let 
$M \prec_{\Cal L} N$ mean $M$ is an
elementary submodel of $N$ for the language ${\Cal L}(\tau_M)$ and
$\tau_M \subseteq \tau_N$, i.e., if $\varphi(\bar x) \in {\Cal
L}(\tau_M)$ and $\bar a \in {}^{\ell g(\bar x)} M$ then $M \models
\varphi[\bar a] \Leftrightarrow N \models \varphi[\bar a]$; similarly
$M \prec_L N$ for $L$ a language.   So $M \prec N$ in the usual sense means 
$M \prec_{\Bbb L} N$ and $M \subseteq N$ means $M$ is a
submodel of $N$. 
\enddemo
\bigskip

\definition{\stag{88r-1.2} Definition}  1) We say ${\frak K}$ is a a.e.c. with
L.S. number $\lambda({\frak K}) = \text{ LS}({\frak K})$ \ub{if}: \nl
\ub{Ax 0}:  The holding of $M \in K,N \le_{\frak K} M$ 
depend on $N,M$ only up to isomorphism, i.e. $[M \in K,M
\cong N \Rightarrow N \in K]$ and 
[if $N \le_{\frak K} M$ and $f$ is an isomorphism
from $M$ onto the $\tau$-model $M',f \restriction N$ is an isomorphism
from $N$ onto $N'$ \ub{then} $N' \le_{\frak K} M'$].
\mn
\ub{Ax I}:  if $M \le_{\frak K} N$ 
then $M \subseteq N$ (i.e. $M$ is a submodel of $N$).
\mn
\ub{Ax II}:  $M_0 \le_{\frak K} M_1 \le_{\frak K} M_2$ 
implies $M_0 \le_{\frak K} M_2$ and $M \le_{\frak K} M$ for $M \in K$.
\mn
\ub{Ax III}:  If $\lambda$ is a regular cardinal, $M_i(i < \lambda)$
is a $\le_{\frak K}$-increasing (i.e. $i < j < \lambda$ implies 
$M_i \le_{\frak K} M_j$) and
continuous (i.e. for $\delta < \lambda,M_\delta = \dbcu_{i < \delta}
M_i$) \ub{then} $M_0 \le_{\frak K} \dbcu_{i < \lambda} M_i$.
\mn
\ub{Ax IV}:  If $\lambda$ is a regular cardinal and $M_i$ 
(for $i < \lambda)$ is $\le_{\frak K}$-increasing continuous and
$M_i \le_{\frak K} N$ for $i < \lambda$ \ub{then} 
$\dbcu_{i < \lambda} M_i \le_{\frak K} N$.
\mn
\ub{Ax V}:  If $N_0 \subseteq N_1 \le_{\frak K} M$ and $N_0 \le_{\frak K} M$ 
\ub{then} $N_0 \le_{\frak K} N_1$.
\mn
\ub{Ax VI}:  If $A \subseteq N \in K$ and $|A| \le \text{ LS}({\frak K})$ 
then for some $M \le_{\frak K} N,A \subseteq |M|$ and $\|M\| \le 
\text{ LS}({\frak K})$ (and $\text{\rm LS}({\frak K})$ is the minimal
infinite cardinal satisfying this axiom which is $\ge |\tau|$; 
the $\ge |\tau|$ is for notational simplicity). \nl
2) We say ${\frak K}$ is a weak \footnote{this is not really
investigated here}  a.e.c. \ub{if} above we omit clause IV.
\enddefinition
\bigskip

\remark{Remark}  Note that AxV holds for $\prec_{\Cal L}$ for any logic
${\Cal L}$.
\endremark
\bn
\ub{Notation}:  Let $K_\lambda = \{M \in K:\|M\| = \lambda\}$ and 
$K_{< \lambda} = \dbcu_{\mu < \lambda} K_\mu$ and ${\frak K}_\lambda =
(K_\lambda,\le_{\frak K} \restriction K_\lambda)$.
\bigskip

\definition{\stag{88r-1.3} Definition}  The embedding $f:N \rightarrow M$
is called a $\le_{\frak K}$-embedding \ub{if} the range of $f$ is the 
universe of a model $N' \le_{\frak K} M$ 
(so $f:N \rightarrow N'$ is an isomorphism onto).
\enddefinition
\bigskip

\definition{\stag{88r-1.4} Definition}  Let $T_1$ be a theory in
${\Cal L}(\tau_1),\Gamma$ a set of types in ${\Cal L}(\tau_1)$ for 
some logic ${\Cal L}$, usually first order. \nl
1) EC$(T_1,\Gamma) = \{M:M$ an $\tau_1$-model of $T_1$ which omits every
$p \in \Gamma\}$. \nl
We implicitly use that $\tau_1$ is reconstructible from $T_1,\Gamma$.
A problem may arise only if some symbols from $\tau_1$ are not mentioned
in $T_1$ and in $\Gamma$, so we may write EC$(T_1,\Gamma,\tau_1)$, but
usually we ignore this point. \nl
2) For $\tau \subseteq \tau_1$ we let 
PC$(T_1,\Gamma,\tau) = \text{ PC}_\tau(T_1,\Gamma) =
\{M:M$ is a $\tau$-reduct of some $M_1 \in \text{ EC}(T_1,\Gamma)\}$. \nl
3) We say that $K$, a class of $\tau$-models, is a PC$^\mu_\lambda$ or
PC$_{\lambda,\mu}$ 
class if for some $T_1,\Gamma_1,\tau_1$ we have $\tau \subseteq
\tau_1,T_1$ a first order theory in the vocabulary $\tau_1,\Gamma_1$ a
set of types in $\Bbb L(\tau_1),K = \text{ PC}_\tau(T_1,\Gamma_1)$ and
$|T_1| \le \lambda,|\Gamma_1| \le \mu$. \nl
4) We say ${\frak K}$ is PC$^\mu_\lambda$ or PC$_{\lambda,\mu}$ 
\ub{if} for some $T_1,T_2,\Gamma_1,\Gamma_2$ as in part (3) and 
$\tau_1,\tau_2$ we have
$K = \text{ PC}(T_1,\Gamma_1,\tau)$ and $\{(M,N):M \le_{\frak K} N$ hence
$M,N \in K\} = \text{ PC}
(T_2,\Gamma_2,\tau')$ where $\tau' = \tau \cup \{P\},P$ a
new one-place predicate, $|T_\ell| \le \lambda,|\Gamma_\ell| \le \mu$
for $\ell = 1,2$.
\nl
If $\mu = \lambda$ we may omit $\mu$. \nl
5) In (4) we may say ``${\frak K}$ is $(\lambda,\mu)$-presentable" and
if $\lambda = \mu$ we may say ``${\frak K}$ is $\lambda$-presentable".
\enddefinition
\bn
\margintag{88r-1.5}\ub{\stag{88r-1.5} Example}:  If $T \subseteq \Bbb L(\tau),\Gamma$ a set of
types in $\Bbb L(\tau)$, then 
$K =: \text{ EC}(T,\Gamma),\le_{\frak K} =: \prec_{{\Bbb L}_{\omega,\omega}}$
form an a.e.c. with LS-number $\le |T| + |\tau| + \aleph_0$, 
that is, satisfy the Axioms from 
\scite{88r-1.2} (for LS$({\frak K}) =: |\tau| + \aleph_0$).
\bigskip

\demo{\stag{88r-1.6} Observation}  Let $I$ be a directed set (i.e. partially
ordered by $\le$, such that any two elements have a common upper
bound). \nl
1) If $M_t$ is defined for $t \in I$ and $t \le s \in I$ implies $M_t
\le_{\frak K} M_s$ \ub{then} for every $t \in I$ we have
$M_t \le_{\frak K} \dbcu_{s \in I} M_s$. \nl
2) If in addition $t \in I$ implies $M_t \le_{\frak K} N$ 
\ub{then} $\dbcu_{s \in I} M_s \le_{\frak K} N$.
\enddemo
\bigskip

\demo{Proof}  By induction on $|I|$ (simultaneously for (1) and (2)).

If $I$ is finite, then $I$ has a maximal element $t(0)$, hence
$\dbcu_{t \in I} M_t = M_{t(0)}$, so there is nothing to prove. \nl
So suppose $|I| = \mu$ and we have proved the assertion when $|I| <
\mu$.  Let $\lambda = \text{ cf}(\mu)$ so $\lambda$ is a regular
cardinal; hence we can find $I_\alpha$ (for $\alpha < \lambda$) such that
$|I_\alpha| < |I|,\alpha < \beta < \lambda$ implies $I_\alpha
\subseteq I_\beta \subseteq I,\dbcu_{\alpha < \lambda} I_\alpha = I$,
for limit $\delta < \lambda,I_\delta = \dbcu_{\alpha < \delta}
I_\alpha$ and each $I_\alpha$ is directed and non-empty.  Let $M^\alpha = 
\dbcu_{t \in I_\alpha} M_t$; so by the induction hypothesis on (1) we
know that $t \in I_\alpha$ implies $M_t \le_{\frak K} M^\alpha$.
If $\alpha < \beta$
then $t \in I_\alpha$ implies $t \in I_\beta$ hence
$M_t \le_{\frak K} M^\beta$; 
hence by the induction hypothesis on (2) applied to $\langle M_t:t \in
I_\alpha \rangle$ we have $M^\alpha = \dbcu_{t \in I_\alpha} M_t 
\le_{\frak K} M^\beta$.
So by Ax III, applied to $\langle M^\alpha:\alpha < \lambda \rangle$
we have $M^\alpha \le_{\frak K} 
\dbcu_{\beta < \lambda} M^\beta = \dbcu_{t \in I} M_t$, 
and as $t \in I_\alpha$ implies $M_t \le_{\frak K} M^\alpha$, by Ax II, 
$t \in I$ implies $M_t \le_{\frak K} \dbcu_{s \in I} M_s$.  
So we have finished proving part (1) for the case $|I|=\mu$.  To prove
(2) in this case note that for each $\alpha < \lambda,
\langle M_t:t \in I_\alpha
\rangle$ is $\le_{\frak K}$-directed and $t \in I_\alpha \Rightarrow
M_t \le_{\frak K} N$, so clearly by the induction hypothesis for (2) we have
$M^\alpha =: \cup\{M_t:t \in I_\alpha\}$ is $\le_{\frak K} N$.  So
$\alpha < \lambda \Rightarrow M^\alpha \le_{\frak K} N$ and as proved
above $\langle M^\alpha:\alpha < \lambda \rangle$ is $\le_{\frak
K}$-increasing, hence by Ax IV, $\dbcu_{s \in I} M_s = 
\dbcu_{\alpha < \lambda} M^\alpha \le_{\frak K} N$.
\hfill$\square_{\scite{88r-1.6}}$
\enddemo
\bigskip

\proclaim{\stag{88r-1.7} Lemma}  Let $\tau_1 = \tau \cup \{F^n_i:i <
{ \text{\rm LS\/}}({\frak K}),n < \omega\},F^n_i$ an $n$-place function symbol 
(assuming, of course, $F^n_i \notin \tau$).

Every model $M$ (in $K$) can be expanded to an $\tau_1$-model $M_1$ such
that:
\mr
\item "{$(A)$}"    $M_{\bar a} \le_{\frak K} M$ 
when $\bar a \in {}^n|M|$ and where $M_{\bar a}$ is the 
submodel of $M$ with universe $\{F^n_i(\bar a):i < { \text{\rm LS\/}}
({\frak K})\}$
\sn
\item "{$(B)$}"  if $\bar a \in {}^n|M|$ then $\|M_{\bar a}\| \le
{ \text{\rm LS\/}}({\frak K})$
\sn
\item "{$(C)$}"  if $\bar b$ is a subsequence of a permutation of
$\bar a$, then $M_{\bar b} \le_{\frak K} M_{\bar a}$
\sn
\item "{$(D)$}"  for every $N_1 \subseteq M_1$ we have $N_1 \restriction \tau 
\le_{\frak K} M$.
\endroster
\endproclaim
\bigskip

\demo{Proof}  We define by induction on $n$, the values of $M_{\bar
a}$ and of $F^n_i(\bar a)$ for every $i < \text{ LS}({\frak K}),\bar a \in 
{}^n|M|$.  Arriving to $n$, for each $\bar a \in {}^n M$ by Ax VI there is
an $M_{\bar a} \le_{\frak K} M$ such that 
$\|M_{\bar a}\| \le \text{ LS}({\frak K}),|M_{\bar a}|$ include
$\cup\{M_{\bar b}:\bar b$ a subsequence of $\bar a$ of length $<n\}
\cup \bar a$ and $M_{\bar a}$ does not depend on the order of $\bar
a$.  Let $|M_{\bar a}| = \{c_i:i < i_0 \le \text{ LS}({\frak K})\}$ and define
$F^n_i(\bar a) = c_i$ for $i < i_0$ and $c_0$ for $i_0 \le i < 
\text{ LS}({\frak K})$.

Clearly our conditions are satisfied; in particular, if 
$\bar b$ is a subsequence of $\bar a,M_{\bar b} 
\le_{\frak K} M_{\bar a}$ by Ax V and clause (D) holds by \scite{88r-1.6}
and Ax IV.  \hfill$\square_{\scite{88r-1.7}}$
\enddemo
\bigskip

\remark{\stag{88r-1.7A} Remark}  1) This is the ``main" place we use Ax V,VI;
it seems that we use it rarely, e.g., in \scite{88r-2.7} which is not 
used later.  It is clear that we can omit
Ax V if we strengthen somewhat Ax VI for the proofs above. \nl
2) Note that in \scite{88r-1.7}, we do not require that $M_{\bar a}$ is closed
under the functions $(F^n_i)^{M_1}$.
By a different bookkeeping we can have it: renaming
$\tau_{1,\varepsilon} = \tau \cup\{F^n_i:i < \text{ LS}({\frak K})
\times \varepsilon,n < \omega\}$ for $\varepsilon \le \omega$ and we
choose a $\tau_{1,n}$-expansion $M_{1,n}$ of $M$ such that $m < n
\Rightarrow M_{1,n} \restriction \tau_{1,m} = M_{1,m}$.  Let $M_{1,0} = M$,
and if $M_{1,n}$ is defined, choose for every $\bar a \in {}^{\omega
>}(M_{1,n})$ a (non-empty) 
subset $A^{1,n}_{\bar a}$ of $M_{1,n}$ of cardinality
$\le$ LS$({\frak K})$ such that $A^{1,n}_{\bar a}$ is closed under the
functions of $M_{1,n}$ and $M \restriction A^{1,n}_{\bar a} \le_{\frak
K} M$, let $A^{1,n}_{\bar a} = \{c_{\bar a,i}:i \in [\text{LS}({\frak
K}) \times n$, LS$({\frak K}) \times (n+1))$ and define $M_{1,n+1}$ by
letting $(F^m_i)^{M_{1,n+1}}(\bar a) = c_{\bar a,i}$.  Let $M_1 =
M_{1,\omega}$ be the $\tau_\omega$-model with the universe of $M$ such that
$n < \omega \Rightarrow M_1 \restriction \tau_{1,n} =
M_{1,n}$. 
\nl
3) Actually $M_{1,1}$ suffices if we expand it by making every term
$\tau(\bar x)$ equal to some function $F(\bar x)$.
\endremark
\bigskip

\proclaim{\stag{88r-1.8} Lemma}  1) ${\frak K}$ is $({\text{\rm LS\/}}
({\frak K}),2^{\text{LS}({\frak K})})$-presentable. \nl
2) There is a set $\Gamma$ of types in $\Bbb L(\tau_1)$ (where
$\tau_1$ is from Lemma \scite{88r-1.7}) such that $K = 
{ \text{\rm PC\/}}_\tau(\emptyset,\Gamma)$. \nl
3) For the $\Gamma$ from part (2), if 
$M_1 \subseteq N_1 \in { \text{\rm EC\/}}(\emptyset,\Gamma)$ and $M,N$ are 
the $\tau$-reducts of $M_1,N_1$ respectively
then $M \le_{\frak K} N$. 
\nl
4) For the $\Gamma$ from part (2), we have 
$\{(M,N):M \le_{\frak K} N;N,M \in K\} = \{(M_1 \restriction
\tau,N_1 \restriction \tau):M_1 \subseteq N_1$ are both from
{\rm PC}$_\Gamma(\emptyset,\Gamma)\}$.
\endproclaim
\bigskip

\demo{Proof}  1) By part (2) the first half of ``${\frak K}$ is
$(\text{LS}({\frak K}),2^{\text{LS}({\frak K})})$-presentable
holds. The second part will be proved with part (4). 
\nl
2) Let $\Gamma_n$ be the set of complete quantifier free
$n$-types $p(x_0,\dotsc,x_{n-1})$ in $\Bbb L(\tau_1)$ such that: if $M_1$ is an
$\tau_1$-model, $\bar a$ realizes $p$ in $M_1$ and $M$ is the $\tau$-reduct
of $M_1$, then $M_{\bar a} \in K$ and
$M_{\bar b} \le_{\frak K} M_{\bar a}$ for any subsequence $\bar b$ of
any permutation of $\bar a$; where $M_{\bar c}(\bar c \in {}^m|M_1|)$ is the
submodel of $M$ whose universe is $\{F^m_i(\bar c):i < 
\text{ LS}({\frak K})\}$ and there are such submodels.

Let $\Gamma$ be the set of $p$ which, for some $n$, are 
complete quantifier free
$n$-types (in $\Bbb L(\tau_1)$) which do not belong to $\Gamma_n$.  
By \scite{88r-1.6} we have PC$_\tau(\emptyset,\Gamma) 
\subseteq K$ and by \scite{88r-1.7} $K \subseteq
\text{ PC}_\tau(\emptyset,\Gamma)$. \nl
3) Similar to the proof of (2) using \scite{88r-1.6}. \nl
4) The inclusion $\supseteq$ holds by part (3); so let us prove the
other direction.  
Given $N \le_{\frak K} M$ we apply the proof of \scite{88r-1.7} to $M$,
but demand further  $\bar a \in {}^n N \Rightarrow M_{\bar a} \subseteq N$;
simply add this demand to the choice of the $M_{\bar a}$'s (hence of
the $F^n_i$'s).  We still have a debt from part (1).

We let $\Gamma'_n$ be the set of complete quantifier free $n$-types in
$\tau'_1 =: \tau_1 \cup \{P\}$ ($P$ a new unary predicate),
$p(x_0,\dotsc,x_{n-1})$ such that:
\mr
\item "{$(*)$}"  if $M_1$ is an $\tau'_1$-model, $\bar a$ realizes
$p$ in $M_1,M$ the $\tau$-reduct of $M_1$, \ub{then} 
{\roster
\itemitem{ $(\alpha)$ }  $M_{\bar b} \le_{\frak K}
M_{\bar a}$ for any subsequence $\bar b$ of $\bar a$ where 
$M_{\bar c}$ (for $\bar c \in |M_1|$) is the submodel of $M$ whose universe is
$\{F^{M_1}(\bar c):i < \text{ LS}({\frak K})\}$, (and there are such models),
\sn
\itemitem{ $(\beta)$ }  $\bar b \subseteq P^{M_1} \Rightarrow M_{\bar
b} \subseteq P^{M_1}$ for $\bar b \subseteq \bar a$.
\endroster}
\ermn
We leave the rest to the reader (alternatively, use
PC$_{\tau'_1}(T',\Gamma),T'$ saying ``$P$ is closed under all the
functions $F^n_i$).   \hfill$\square_{\scite{88r-1.8}}$
\enddemo
\bn
By the proof of \scite{88r-1.8}(4).
\demo{\stag{88r-1.9} Conclusion}  The $\tau_1$ and $\Gamma$ from
\scite{88r-1.8} (so $|\tau_1| \le \text{ LS}({\frak K}))$ satisfy: 
for any $M \in K$ and any $\tau_1$-expansion $M_1$
of $M$ which is in EC$_{\tau_1}(\emptyset,\Gamma)$
\mr
\item "{$(a)$}"   $N_1 \prec_{\Bbb L} M_1 \Rightarrow N_1 \subseteq M_1
\Rightarrow N_1 \restriction \tau \le_{\frak K} M$
\sn
\item "{$(b)$}"  $N_1 \prec_{\Bbb L} N_2 \prec_{\Bbb L}
M_1 \Rightarrow N_1 \subseteq N_2 \subseteq M_1 \Rightarrow
N_1 \restriction \tau \le_{\frak K} N_2 \restriction \tau$
\sn
\item "{$(c)$}"  if $M \le_{\frak K} N$ then there is a
$\tau_1$-expansion $N_1$ of $N$ from
EC$_{\tau_1}(\emptyset,\Gamma)$ which extends $M_1$.
\endroster
\enddemo
\bn
\margintag{88r-1.10}\ub{\stag{88r-1.10} Conclusion}  If for every $\alpha <
(2^{\text{LS}({\frak K})})^+,
{\frak K}$ has a model of cardinality $\ge \beth_\alpha$
\ub{then} $K$ has a model in every cardinality $\ge \text{ LS}({\frak K})$.
\bigskip

\demo{Proof}  Use \scite{88r-1.8} and the classical upper bound on
value of the Hanf number for:
first order theory and omitting any set of types, for languages of
cardinality LS$({\frak K})$ (see, e.g., \cite[VII,5.3,5.5]{Sh:c}). 
 \hfill$\square_{\scite{88r-1.10}}$
\enddemo
\bigskip

\remark{\stag{88r-1.11} Remark}  1) Clearly $\{\mu:\mu \ge \text{
LS}({\frak K})$ and $K_\mu \ne 0\}$ is an initial segment of the class
of cardinals $\ge \text{ LS}({\frak K})$.
\nl
2) For every cardinal $\kappa(\ge \aleph_0)$ and ordinal $\alpha <
(2^\kappa)^+$ there is an a.e.c. ${\frak K}$ such that: LS$({\frak K})
= \kappa = |\tau_{\frak K}|$ and ${\frak K}$ has a model of
cardinality $\lambda$ iff $\lambda \in
[\kappa,\beth_\alpha(\kappa))$.  This follows by
\cite[VII,\S5,p.432]{Sh:c} in particular \cite[VII,5.5]{Sh:c}(6),
because
\mr
\item "{$(a)$}"  if a vocabulary of cardinality $\le \kappa$ and $T
\subseteq \Bbb L(\tau)$ and $\Gamma$ a set of $(\Bbb
L(\tau),<\omega)$-types then $K = \{M:M$ a $\tau$-model of $T$ omitting
every $\in \Gamma\}$ and $\le_{\frak K} = \prec \restriction K$ form
an a.e.c. (we can use $\Gamma$ a set of quantifier free types, $T =
\emptyset$), with LS$(({\frak K},\le_{\frak K}) \le \kappa$
\sn
\item "{$(b)$}"  if $\{c_i \ne c_j:i < j < \kappa\} \subseteq T$ then
$K$ above has no model of cardinality $< \kappa$.
\ermn
3) More on such theorems see \cite{Sh:394}.
\endremark
\newpage

\head {\S2 Amalgamation properties and homogeneity} \endhead  \resetall \sectno=2
 \spuriousreset
\bigskip

\demo{\stag{88r-2.0} Context}  ${\frak K}$ is an a.e.c.

The main theorem \scite{88r-2.6}, the existence and uniqueness of the
model homogeneous models, is from Jonsson \cite{Jo56}, \cite{Jo60}.
The result on the number of $D$-homogeneous universal models is 
from Keisler and Morley
\cite{KM67}.  Probably a new result is \scite{88r-2.4}(2) (and \scite{88r-2.11}).
\enddemo
\bigskip

\definition{\stag{88r-2.1} Definition}  ${\Bbb D}(M) := \{N/ \cong: N \le_{\frak K}
M,\|N\| \le \text{ LS}({\frak K})\}$

$$
{\Bbb D}({\frak K}) := \{N/ \cong: N \in K,\|N\| \le \text{ LS}({\frak K})\}.
$$
\enddefinition
\bigskip

\definition{\stag{88r-2.2} Definition}  Let $\lambda > \text{ LS}({\frak K})$. \nl
1) A model $M$ is $\lambda$-model homogeneous when: 
\ub{if} $N_0 \le_{\frak K} N_1 \le_{\frak K} M,
\|N_1\| < \lambda,f$ an $\le_{\frak K}$-embedding of $N_0$ into $M$, 
\ub{then} some $\le_{\frak K}$-embedding 
$f':N_1 \rightarrow M$ extends $f$. \nl
1A) A model $M$ is $({\Bbb D},\lambda)$-model homogeneous \ub{if}
${\Bbb D} = {\Bbb D}(M)$ and $M$ is a $\lambda$-model homogeneous. \nl
2) $M$ is $\lambda$-strongly model homogeneous if: for every $N \in
K_{< \lambda}$ such that $N \le_{\frak K} M$ and 
a $\le_{\frak K}$-embedding $f:N \rightarrow M$ there exists an 
automorphism $g$ of $M$ extending $f$. \nl
3) $M$ is $\lambda$-model universal homogeneous when: every $N \in
K_{\le \lambda}$ is $\le_{\frak K}$-embeddable into $M$ and 
for every $N_\ell \in
K_{< \lambda}$ (for $\ell=0,1$) 
such that $N_0 \le_{\frak K} N_1$ and
$\le_{\frak K}$-embedding of $f:N_0 \rightarrow M$ there exists a 
$\le_{\frak K}$-embedding $g:N_1 \rightarrow M$ extending $f$ (unlike
(1), we do not demand that $N_1$ is $\le_{\frak K}$-embeddable into
$M$; the universal is related to $\lambda$, it does not imply $M$ is
universal). 
\nl
4) For each of the above three properties and the one below, if 
$M$ has cardinality $\lambda$ and has the $\lambda$-property then we 
may say for short that $M$ has the
property (i.e. omitting $\lambda$).
\nl
5) $M$ is $(D,\lambda)$-sequence homogeneous \ub{if}: 
\mr
\item "{$(a)$}"  $D = D(M) = \{\text{tp}_{\Bbb L(\tau_M)}(\bar a,
\emptyset,M):\bar a \in |M|$, i.e., $\bar a$ a finite sequence from $M\}$ and
\sn
\item "{$(b)$}"   if $a_i \in M$ for $i \le \alpha
< \lambda,b_j \in M$ for $j < \alpha$ and tp$_{\Bbb L(\tau M)}
(\langle a_i:i < \alpha \rangle,\emptyset,M) = 
\text{ tp}_{\Bbb L(\tau_M)}(\langle b_i:i < \alpha
\rangle,\emptyset,M)$, \ub{then} for some $b_\alpha \in M$,
tp$_{\Bbb L(\tau_M)}(\langle a_i:i \le \alpha \rangle,
\emptyset,M) = \text{ tp}_{\Bbb L(\tau_M)}(\langle b_i:i \le \alpha 
\rangle,\emptyset,M)$.
\ermn
6) We omit the ``model/sequence", when which one is clear from the context,
i.e., if $D$ is as in \scite{88r-2.2}(5)(a), $(D,\lambda)$-homogeneous
means $(D,\lambda)$-sequence-homogeneous: if ${\Bbb D}$ is as in
Definition \scite{88r-2.1}, $({\Bbb D},\lambda)$-homogeneous means 
$({\Bbb D},\lambda)$-model-homogeneous, if not obvious we mean 
the model version.
\enddefinition
\bigskip

\proclaim{\stag{88r-2.3} Claim}  Assume $N$ is $\lambda$-model homogeneous and
${\Bbb D}(M) \subseteq {\Bbb D}(N)$, (and {\rm LS}$({\frak K}) 
< \lambda$, of course). \nl
1) If $M_0 \le_{\frak K} M_1 \le_{\frak K} M,\|M_0\| < \lambda,
\|M_1\| \le \lambda$ and $f$ is a $\le_{\frak K}$-embedding of $M_0$ into $N$, 
\ub{then} we can extend $f$ to a $\le_{\frak K}$-embedding of $M_1$ into
$N$.
\nl
2) If $M_1 \le_{\frak K} M,\|M_1\| \le \lambda$ 
\ub{then} there is a $\le_{\frak K}$-embedding of $M_1$ into $N$.
\endproclaim
\bigskip

\demo{Proof}  We prove by induction on $\mu \le \lambda$ simultaneously that
\mr
\item "{$(i)_\mu$}"   for every $M_1 \le_{\frak K} M,\|M_1\| \le \mu$ 
(yes! not $< \mu$) there is
a $\le_{\frak K}$-embedding of $M_1$ into $N$
\sn
\item "{$(ii)_\mu$}"  if $M_0 \le_{\frak K} M_1 \le_{\frak K} M,
\|M_1\| \le \mu,\|M_0\| < \lambda$ then any $\le_{\frak K}$-embedding $f$
of $M_0$ into $N$ can be extended to a $\le_{\frak K}$-embedding 
of $M_1$ into $N$.
\ermn
Clearly $(i)_\lambda$ is part (2) and $(ii)_\lambda$ is part (1) so
this is enough.
\enddemo
\bigskip

\demo{Proof of $(i)_\mu$}  If $\mu \le \text{ LS}({\frak K})$, this follows
by ${\Bbb D}(M) \subseteq {\Bbb D}(N)$. \nl
If $\mu > \text{ LS}({\frak K})$, then by \scite{88r-1.9} we can find
$\bar M_1 = \langle M^\alpha_1:\alpha < \mu \rangle$ such that
$M_1 = \dbcu_{\alpha < \mu} M^\alpha_1$ and $\alpha < \mu \Rightarrow
M^\alpha_1 \le_{\frak K} M_1$ and $M^\alpha_1$ is $\le_{\frak K}$-increasing
continuous with $\alpha$ and $\alpha < \mu \Rightarrow \|M^\alpha_1\| < \mu$.  
We define by induction on
$\alpha$, a $\le_{\frak K}$-embedding $f_\alpha:M^\alpha_1 \rightarrow N$, such
that for $\beta < \alpha,f_\alpha$ extend $f_\beta$.  For $\alpha = 0$
we can define $f_\alpha$ by 
$(i)_{\chi(0)}$ which holds as by the induction hypothesis,
where $\chi(\beta) = \|M^\beta_1\|$.  We next 
define $f_\alpha$ for $\alpha = \gamma+1$: by
$(ii)_{\chi(\alpha)}$ which holds by the induction hypothesis
there is a $\le_{\frak K}$-embedding $f_\alpha$ of
$M^\alpha_1$ into $N$ extending $f_\gamma$.

Lastly, for limit $\alpha$ we let $f_\alpha = \dbcu_{\beta < \alpha}
f_\beta$, it is a $\le_{\frak K}$-embedding into $N$ by \scite{88r-1.6}.  So we
finish the induction and $\dbcu_{\alpha < \mu} f_\alpha$ is as required.
\enddemo
\bigskip

\demo{Proof of $(ii)_\mu$}  First, assume that $\mu = \lambda$ so
we have proved $(ii)_\theta$ for $\theta < \lambda$ and $\|M_1\| = 
\lambda > \|M_0\|$, so LS$({\frak K}) < \mu = \lambda$ hence 
we can find $\langle M^\alpha_1:\alpha < \mu \rangle$ as in the
proof of $(i)_\mu$ such that $M^0_1 = M_0$.  Now we define $f_\beta$
by induction on $\beta \le \mu$ such that $f_\beta$ is a $\le_{\frak
K}$-embedding of $M^1_\beta$ into $N$ and $f_\beta$ is increasing
continuous in $\beta$ and $f_0 = f$.  We can do this as in the proof
of $(i)_\mu$ by $(ii)_{\chi(\alpha)}$ for $\alpha < \mu$.  

Second, assume $\|M_1\| < \lambda$.  Let $g$ be a 
$\le_{\frak K}$-embedding of $M_1$ into
$N$, it exists by $(i)_\mu$ which we have just proved.  Let $g$ be onto
$N'_1 \le_{\frak K} N$, and let $g \restriction M_0$ be onto 
$N'_0 \le_{\frak K} N'_1$, and
let $f$ be onto $N_0 \le_{\frak K} N$.  So clearly 
$h:N'_0 \rightarrow N_0$ defined by $h(g(a)) = f(a)$ for 
$a \in |M_0|$ is an isomorphism from $N'_0$
onto $N_0$.  So $N_0,N'_0,N'_1 \le_{\frak K} N$.  As $\|M_1\| <
\lambda$ clearly $\|N'_1\| < \lambda$ so 
(by the assumption ``$N$ is $\lambda$-model homogeneous")
we can extend
$h$ to an isomorphism $h'$ from $N'_1$ onto some $N_1 \le_{\frak K} N$, 
so $h' \circ g:M_1 \rightarrow N$ is as required.
\hfill$\square_{\scite{88r-2.3}}$
\enddemo
\bn
\margintag{88r-2.4}\ub{\stag{88r-2.4} Conclusion}  1) If $M,N$ are model-homogeneous, of the
same cardinality ($> \text{ LS}({\frak K}$))
and ${\Bbb D}(M) = {\Bbb D}(N)$ \ub{then} $M,N$ are
isomorphic.  Moreover, if $M_0 \le_{\frak K} M,\|M_0\| < \|M\|$, \ub{then} any
$\le_{\frak K}$-embedding of $M_0$ into $N$ can be extended to an isomorphism from
$M$ onto $N$. \nl
2) The number of model homogeneous models from ${\frak K}$ of 
cardinality $\lambda$ is
$\le 2^{2^{\text{LS}({\frak K})}}$. \nl
3) If $M$ is $\lambda$-model-homogeneous, ${\Bbb D}(M) = {\Bbb D}({\frak K})$
\ub{then} $M$ is $\lambda$-universal, i.e. every model $N$ (in $K$) of
cardinality $\le \lambda$, has a $\le_{\frak K}$-embedding into $M$.  
So if $\Bbb D(M) = \Bbb D({\frak K})$ \ub{then}: $M$ is
$\lambda$-model universal homogeneous (see Definition \scite{88r-2.2}(3))
iff $M$ is a $\lambda$-model homogeneous iff $M$ is $(\lambda,\Bbb
D({\frak K}))$-homogeneous.
\nl
4) If $M$ is $\lambda$-model-homogeneous \ub{then} it is
$\lambda$-universal for $\{N \in K_\lambda:{\Bbb D}(N) \subseteq 
{\Bbb D}(M)\}$. \nl
5) If $M$ is $(D,\lambda)$-sequence homogeneous, $(\lambda > \text{ LS}
({\frak K}))$ \ub{then} $M$ is a $\lambda$-model homogeneous.
\bigskip

\demo{Proof}  1) Immediate by \scite{88r-2.3}(1), using the standard hence
and forth argument. \nl
2) The number of models (in $K$) of power $\le \text{ LS}
({\frak K})$ is, up to
isomorphism, $\le 2^{\text{LS}({\frak K})}$ (recalling that we are assuming
$|\tau({\frak K})| \le \text{ LS}({\frak K}))$.  Hence the number of possible
${\Bbb D}(M)$ is $\le 2^{2^{\text{LS}
({\frak K})}}$.  So by \scite{88r-2.4}(1) we are done.
\nl
3),4),5)  Immediate.   \hfill$\square_{\scite{88r-2.4}}$
\enddemo
\bigskip

\proclaim{\stag{88r-2.4A} Claim}  The results parallel to \scite{88r-2.4}(1)-(4)
 for $\lambda$-sequence homogeneous holds.
\endproclaim
\bigskip

\definition{\stag{88r-2.5} Definition}  1) A model $M$ has the
$(\lambda,\mu)$-amalgamation property (= am.p., in ${\frak K}$, of
course) \ub{if}: for every $M_1,M_2$
such that $\|M_1\| = \lambda,\|M_2\| = \mu,M \le_{\frak K} M_1$ and
$M \le_{\frak K} M_2$, there is
a model $N$ and $\le_{\frak K}$-embeddings $f_1:M_1 \rightarrow N$ and $f_2:M_2
\rightarrow N$ such that $f_1 \restriction |M| = f_2 \restriction
|M|$.  Now the meaning of e.g. the $(\le \lambda,< \mu)$-amalgamation
property should be clear.  Always $\lambda,\mu \ge \text{ LS}
({\frak K})$ (and, of course, if we use $< \mu,\mu > \text{\rm
LS}({\frak K}))$. \nl
2) ${\frak K}$ has the 
$(\kappa,\lambda,\mu)$-amalgamation property \ub{if} every
model $M$ (in $K$) of cardinality $\kappa$ has the
$(\lambda,\mu)$-amalgamation property.  The
$(\kappa,\lambda)$-amalgamation property for ${\frak K}$ means just the
$(\kappa,\kappa,\lambda)$-amalgamation property.  The
$\kappa$-amalgamation property for ${\frak K}$ is just the
$(\kappa,\kappa,\kappa)$-amalgamation property. \nl
3) ${\frak K}$ has the 
$(\lambda,\mu)$-JEP (joint embedding property) if for any
$M_1 \in K,M_2 \in K$ of cardinality $\lambda,\mu$ respectively there
is $N \in K$ into which $M_1$ and $M_2$ are $\le_{\frak K}$-embeddable. \nl
4) The $\lambda$-JEP is the $(\lambda,\lambda)$-JEP. \nl
5) The amalgamation property means the
$(\kappa,\lambda,\mu)$-amalgamation property for every $\lambda,\mu
\ge \kappa (\ge \text{ LS}({\frak K}))$. \nl
6) The JEP means the $(\lambda,\mu)$-JEP for every $\lambda,\mu \ge
\text{ LS}({\frak K})$.
\enddefinition
\bigskip

\remark{Remark}  Clearly in \scite{88r-2.5}, parts (1), (2) first 
sentence, (3),(5),  the roles of $\lambda,\mu$ are symmetric.
\endremark
\bigskip

\proclaim{\stag{88r-2.6} Theorem}  1) If {\rm LS}$({\frak K}) < \kappa \le
\lambda,\lambda = \lambda^{< \kappa},K_\lambda \ne \emptyset$ and
${\frak K}$ has the $(< \kappa,\lambda)$-amalgamation property 
\ub{then} for every
model $M$ of cardinality $\lambda$, there is a $\kappa$-model
homogeneous model $N$ of cardinality $\lambda$ satisfying $M
\le_{\frak K} N$.  If $\kappa = \lambda$, alternatively
the $(< \kappa,< \lambda)$-amalgamation property suffices.
\nl
2) So in (1) if $\kappa = \lambda$, there is a universal, model homogeneous
model of cardinality $\lambda$, provided that for some $M \in K_{\le
\lambda},{\Bbb D}(M) = {\Bbb D}({\frak K})$ or just ${\frak K}$ has the
{\rm LS}$({\frak K})$-{\rm JEP}. \nl
3) If ${\frak K}$ has the amalgamation property and the {\rm LS}$({\frak
K})$-{\rm JEP}, \ub{then} ${\frak K}$ has the {\rm JEP}.
\endproclaim
\bigskip

\remark{\stag{88r-2.6A} Remark}  1) The last assumption of \scite{88r-2.6}(2)
holds, e.g., if $(\le \text{ LS}({\frak K}),
< 2^{\text{LS}({\frak K})})$-JEP holds and
$|{\Bbb D}({\frak K})| \le \lambda$. \nl
2) If for some $M \in K,{\Bbb D}(M) = {\Bbb D}({\frak K})$ \ub{then} we can
have such $M$ of cardinality $\le 2^{\text{LS}({\frak K})}$. \nl
3) We can in \scite{88r-2.6} replace the assumption 
``$(< \kappa,\lambda)$-amalgamation property" by
``$(< \kappa,< \lambda)$-amalgamation property" if, e.g., no 
$M \in K_{< \lambda}$ is maximal.
\endremark
\bigskip

\demo{Proof}  Immediate; in (1) note that if $\kappa$ is singular
then necessarily $\lambda > \kappa \cap \lambda = \lambda^\kappa =
\lambda^{< \kappa^+}$ so we can replace $\kappa$ by $\kappa^+$.
\enddemo
\bigskip

\remark{\stag{88r-2.6B} Remark}  Also the corresponding converses hold.
\endremark
\bigskip

\proclaim{\stag{88r-2.11.1} Claim}  Assume that $\lambda =
\lambda^{<\lambda}$ and ${\frak K}$ has
$(<\lambda,<\lambda,<\lambda)$-amalgamation property but in the
stronger version demanding that the resulting model has cardinality $<
\lambda$.  \ub{Then} 
${\frak K}_\lambda$ has a smooth model homogeneous members.
\endproclaim
\bigskip

\proclaim{\stag{88r-2.7} Lemma}  1) If ${\frak K}$ has the $\kappa$-amalgamation
property \ub{then} ${\frak K}$ 
has the $(\kappa,\kappa^+)$-amalgamation property
and even the $(\kappa,\kappa^+,\kappa^+)$-amalgamation property. \nl
2) If $\kappa \le \mu \le \lambda$ and ${\frak K}$ has the
$(\kappa,\mu)$-amalgamation property and the
$(\mu,\lambda)$-amalgamation property \ub{then} ${\frak K}$ has the
$(\kappa,\lambda)$-amalgamation property.  If ${\frak K}$ has the
$(\kappa,\mu,\mu)$ and the $(\mu,\lambda)$-amalgamation property,
\ub{then} ${\frak K}$ has the 
$(\kappa,\lambda,\mu)$-amalgamation property. 
\nl
3) If $\lambda_i(i \le \alpha)$ is increasing and continuous,
{\rm LS}$({\frak K}) 
\le \lambda_0$ and for every $i < \alpha,{\frak K}$ has the
$(\lambda_i,\mu + \lambda_i,\lambda_{i+1})$-amalgamation property \ub{then}
${\frak K}$ has the $(\lambda_0,\mu + \lambda_0,
\lambda_\alpha)$-amalgamation property. \nl
4) If $\kappa \le \mu_1 \le \mu$ and for every $M,\|M\| = \mu_1$,
there is $N,M \le_{\frak K} N,\|N\| = \mu$, \ub{then} the
$(\kappa,\mu,\lambda)$-amalgamation property (for ${\frak K}$) implies the
$(\kappa,\mu_1,\lambda)$-amalgamation property (for ${\frak K}$).
\endproclaim
\bigskip

\demo{Proof}  Straightforward.
\enddemo
\bigskip

\demo{\stag{88r-2.8} Conclusion}  If LS$({\frak K}) \le \chi_1 < \chi_2$
and ${\frak K}$ has the $\kappa$-amalgamation property whenever 
$\chi_1 \le \kappa < \chi_2$ \ub{then} $K$ has the
$(\kappa,\lambda,\mu)$-amalgamation property  whenever $\chi_1 \le
\kappa \le \lambda \le \chi_2,\kappa 
\le \mu \le \chi_2$ and $\kappa < \chi_2$.
\enddemo
\bn
\centerline {$* \qquad * \qquad *$}
\bn
It may be interesting to note that even waiving AX IV we can say
something.
\bn
\margintag{88r-2.9}\ub{\stag{88r-2.9} Context}:  For the remainder of this section ${\frak K}$
is just a weak a.e.c., i.e.,  Ax IV is not assumed.
\bigskip

\definition{\stag{88r-2.10} Definition}  Let $M \in K$ have cardinality
$\lambda$, a regular uncountable cardinal $> \text{ LS}({\frak K})$.  
We say $M$ is \ub{smooth} if there is a sequence 
$\langle M_i:i < \lambda \rangle$ with $M_i$ being 
$\le_{\frak K}$-increasing continuous, $M_i \le_{\frak K} M$ and
$\|M_i\| < \lambda$ for $i < \lambda$ and 
$M = \dbcu_{i < \lambda} M_i$.
\enddefinition
\bigskip

\remark{\stag{88r-2.10A} Remark}  We can define $S/{\Cal D}$-smooth, for $S$ a
subset of ${\Cal P}(\lambda),{\Cal D}$ a filter on ${\Cal
P}(\lambda)$, that is: $M \in K_\lambda$ is $(S/{\Cal D})$-smooth when
for every one-to-one function $f$ from $|M|$ onto $\lambda$ the set
$\{u \in {\Cal P} (\lambda):M \restriction \{a:f(a) \in u\} \le_{\frak
K} M\} \in D$.  Usually we demand that for every permutation $f$ on
$\lambda,\{u \subseteq \lambda$: $u$ is closed under $f\} \in {\Cal D}$, 
and usually we demand that ${\Cal D}$ is a 
normal LS$({\frak K})^+$-complete filter).
\endremark
\bigskip

\proclaim{\stag{88r-2.11} Lemma}  If $M,N \in K_\lambda(\lambda > 
{ \text{\rm LS\/}}({\frak K}))$ 
are smooth, model homogeneous and ${\Bbb D}(M) = {\Bbb D}(N)$ 
\ub{then} $M \cong N$.
\endproclaim
\bigskip

\demo{Proof}  By the hence and forth argument left to the reader (the
set of approximations is $\{f:f$ isomorphism from some $M' \le_{\frak K} M$
of cardinality $< \lambda$ onto some $N' \le_{\frak K} N\}$ but note that
not for any increasing continuous sequence of approximations is the
union an approximation).   \hfill$\square_{\scite{88r-2.11}}$ 
\enddemo
\bigskip

\remark{\stag{88r-2.11A} Remark}  It is reasonable to consider
\mr
\item "{$(*)$}"  if $M \in K_\lambda,(\lambda > 
\text{ LS}({\frak K}))$ is smooth
and model homogeneous and $N \in K_\lambda$ is smooth, ${\Bbb D}(N)
\subseteq {\Bbb D}(M)$ \ub{then} $N$ can be $\le_{\frak K}$-embedded into $M$.
\endroster
\endremark
\bn
This can be proved in the context of universal classes (e.g. $Ax
Fr_1$ from \cite{Sh:300b}).
\bn
\margintag{88r-2.12}\ub{\stag{88r-2.12} Fact}:  1) If $(K_i,<_i)$ satisfies the axioms with
$\lambda_i = \text{ LS}(K_i,\le_i)$ where $\lambda_i \ge \aleph_0$ 
for $i < \alpha,i < \alpha \Rightarrow \tau_{K_i} = \tau$ and
$K = \dbca_{i < \alpha} K_i$ and $\le$ is defined by $M \le N$ if and only
if for every $i < \alpha,
M \le_i N$ \ub{then} $(K,\le)$ satisfies the axioms with
LS$(K,\le) \le \dsize \sum_{i < \alpha} \lambda_i$. \nl
2) Concerning AxI-V, we can omit some of them in the assumption and
still get the rest in the conclusion.  But for AxVI we need in
addition to assume AxV + AxIV$_\theta$ for at least one $\theta =
\text{ cf}(\theta) \le \dsize \sum_{i < \alpha} \lambda_i$.
\bigskip

\demo{Proof}  Easy.
\enddemo
\bn
\margintag{88r-2.13}\ub{\stag{88r-2.13} Example}  Consider the class $K$ of norm spaces over
the reals with $M \le_{\frak K} N$ iff $M \subseteq N$ and $M$ is
complete inside $N$.  Now ${\frak K} = (K,\le_{\frak K})$ is a weak
a.e.c. with LS$({\frak K}) = 2^{\aleph_0}$ and it is as required in
\scite{88r-2.11.1}. 
\newpage

\head {\S3 Limit models and other results} \endhead  \resetall \sectno=3
 \spuriousreset
\bigskip

In this section we introduce various variants of limit models (the
most important are the superlimit ones).  We
prove that if ${\frak K}$ has a superlimit model $M^*$ of cardinality
$\lambda$ for
which the $\lambda$-amalgamation property fails and $2^\lambda <
2^{\lambda^+}$ \ub{then} $\dot I(\lambda,K) = 2^\lambda$ (see
\scite{88r-3.5}).  We later prove that if $\psi \in
\Bbb L_{\omega_1,\omega}(\bold Q)$ is categorical in $\aleph_1$ then it
has model in $\aleph_2$ see \scite{88r-3.9}(2).  This finally solves
Baldwin's problem (see \S0).  In fact we prove an
essentially more general result on a.e.c. and $\lambda$ 
(see \scite{88r-3.7}, \scite{88r-3.8}).  

The reader can read \scite{88r-3.1}(1), ignore the other definitions, and
continue with \scite{88r-3.4}(2),(5) and everything from \scite{88r-3.5}
(interpreting all variants as superlimits).
\nl
We may wonder can we prove the parallel to Baldwin conjecture in
$\lambda^+$ if $\lambda > \aleph_0$; it is
\mr
\item "{$\circledast_\lambda$}"  if ${\frak K}$ is
$\lambda$-presentable a.e.c. with LS$({\frak K}) = \lambda$,
categorical in $\lambda^+$ then $K_{\lambda^{++}} \ne \emptyset$.
\ermn
This is easily false when cf$(\lambda) > \aleph_0$.
\bigskip

\demo{\stag{88r-3.0} Context}  ${\frak K}$ is an a.e.c.
\enddemo
\bn
\margintag{88r-3.0.2}\ub{\stag{88r-3.0.2} Example}:  Let $\lambda$ be given and
${\frak K} = (K,\le_{\frak K})$ be defined by

$$
K = \{(A,<):(A,<) \text{ a well order of order type } \le \lambda^+\}
$$

$$
\le_{\frak K} = \{(M,N):M,N \in K \text { and }  N 
\text{ is an end extension of } M\}.
$$
\mn
Now
\mr
\item "{$(a)$}"   ${\frak K}$ 
is an abstract elementary class with LS$({\frak K}) = \lambda$ and
${\frak K}$ categorical in $\lambda^+$
\sn
\item "{$(b)$}"  if $\lambda$ has cofinality $\ge \aleph_1$ \ub{then}
${\frak K}$ is $\lambda$-presentable (see, e.g., \cite[VII,\S5]{Sh:c} and
history there); by clause (a) it is always $(\lambda,2^\lambda)$-presentable, 
\sn
\item "{$(c)$}"  ${\frak K}$ has no model of cardinality $> \lambda^+$.
\ermn
Note that if we are dealing with classes which are categorical (or just
simple in some sense), we have a good chance to find limit models and
they are useful in constructions.
\bigskip

\definition{\stag{88r-3.1} Definition}  Let $\lambda$ be a cardinal $\ge
\text{ LS}({\frak K})$.  For parts 3) - 7) but not 8), for simplifying
the presentation we assume the axiom of global choice (alternatively,
we restrict ourselves to models with universe an ordinal $< \lambda^+$).
\nl
1) $M \in K_\lambda$ is locally superlimit (for ${\frak K}$) \ub{if}:
\mr
\item "{$(a)$}"  for every $N \in K_\lambda$ such that $M \le_{\frak
K} N$ there is $M' \in K_\lambda$ isomorphic to $M$ such that
$N \le_{\frak K} M'$ and $N \ne M'$
\sn
\item "{$(b)$}"  if $\delta < \lambda^+$ is a limit ordinal and 
$\langle M_i:i < \delta \rangle$ is $\le_{\frak K}$-increasing
sequence and $M_i \cong M$
for $i < \delta$ then $\dbcu_{i < \delta} M_i \cong M$.
\ermn
1A) $M \in K_\lambda$ is globally superlimit if (a) +(b) and 
\mr
\item "{$(c)$}"   $M$ is universal in ${\frak K}_\lambda$, i.e., 
any $N \in K_\lambda$ can be $\le_{\frak K}$-embedded into $M$.
\ermn
1B) Just superlimit means globally.  Similarly with the other notions
below we define the global version as adding clause (c) from (1A) and
the default version is the global one.  (Note that in the local 
version we can restrict our class to $\{N \in
K_\lambda:M$ can be $\le_{\frak K}$-embedded into $N\}$ and get the
global one).
\nl
2) For $\Theta \subseteq \{\mu:\aleph_0 \le \mu < \lambda,\mu$ regular$\},M
\in K_\lambda$ is locally $(\lambda,\Theta)$-superlimit \ub{if}:
\mr
\item "{$(a)$}"  from above holds and
\sn
\item "{$(b)$}"  if $\langle M_i:i \le \mu \rangle$ 
is $\le_{\frak K}$-increasing, $M_i \cong M$  for $i < \mu$ and 
$\mu \in \Theta$ then $\cup \{M_i:i < \mu\} \cong M$.
\ermn
2A) If $\Theta$ is a singleton, say $\Theta = \{\theta\}$, we may say that
$M$ is locally $(\lambda,\theta)$-superlimit.
\nl
3) Let $S \subseteq \lambda^+$ be stationary.   $M \in K_\lambda$
is called locally $S$-strongly limit or locally
$(\lambda,S)$-strongly limit \ub{when} for some function: 
$\bold F:K_\lambda
\rightarrow K_\lambda$ we have:
\mr
\item "{$(a)$}"  for $N \in K_\lambda$ we have $N \le_{\frak K} \bold F(N)$
\sn
\item "{$(b)$}"  if $\delta \in S$ is a limit ordinal
and $\langle M_i:i < \delta
\rangle$ is a $\le_{\frak K}$-increasing continuous sequence 
\footnote{no loss if we add $M_{i+1} \cong M$, so this simplifies the
demand on $\bold F$, i.e., only $\bold F(M)$ is required}
in $K_\lambda$ and $M_0 \cong M$ and $i < \delta \Rightarrow \bold F(M_{i+1})
\le_{\frak K} M_{i+2}$, \ub{then} $M \cong \cup\{M_i:i < \delta\}$
\sn
\item "{$(c)$}"  if $M \le_{\frak K} M_1 \in K_\lambda$ then there is $N$ such
that $M_1 <_{\frak K} N \in K_\lambda$. 
\ermn
4) Let $S \subseteq \lambda^+$ be stationary.  $M \in
K_\lambda$ is called locally $S$-limit or locally 
$(\lambda,S)$-limit \ub{if} for some function 
$\bold F:K_\lambda \rightarrow K_\lambda$ we have:
\mr
\item "{$(\alpha)$}"  for every $N \in K_\lambda$ we have 
$N \le_{\frak K} \bold F(N)$ 
\sn
\item "{$(\beta)$}"   if $\langle M_i:i < \lambda^+ \rangle$ is
a $\le_{\frak K}$-increasing continuous sequence of members of 
$K_\lambda,M_0 \cong M,\bold F(M_{i+1}) \le_{\frak K} M_{i+2}$ \ub{then} 
for some closed unbounded \footnote{we can use a filter as a parameter}
subset $C$ of $\lambda^+$,

$$
[\delta \in S \cap C \Rightarrow M_\delta \cong M].
$$
\sn
\item "{$(\gamma)$}"   if $M \le_{\frak K} M_1 \in K_\lambda$ then
there is $N,M_1 <_{\frak K} N \in K_\lambda$.
\ermn
5) We define ``locally $S$-weakly limit", ``locally $S$-medium limit" 
like ``locally $S$-limit", ``locally $S$-strongly limit" 
respectively by demanding that the domain of
$\bold F$ is the family of $\le_{\frak K}$-increasing continuous
sequence of members of ${\frak K}_{< \lambda}$ of length $< \lambda$ and
replacing ``$\bold F(M_{i+1}) 
\le_{\frak K} M_{i+2}$" by 
``$M_{i+1} \le_{\frak K} \bold F(\langle M_j:j \le i+1 \rangle) 
\le_{\frak K} M_{i+2}"$.  We replace ``limit" by ``limit$^-$" if
$``\bold F(M_{i+1}) \le_{\frak K} M_{i+2}",``M_{i+1} \le_{\frak K} 
\bold F(\langle M_j:j \le i+1 \rangle) \le_{\frak K} M_{i+2}"$ are replaced 
by $``\bold F(M_i) \le_{\frak K} M_{i+1}",``M_i \le_{\frak K} \bold F
(\langle M_j:j \le i \rangle) \le_{\frak K} M_{i+1}"$ respectively.
\nl
6) If $S = \lambda^+$ then we omit $S$ (in parts (3), (4), (5)).  
\nl
7) For $\Theta \subseteq \{\mu:\aleph_0 \le \mu \le \lambda$ and $\mu$ is 
regular$\},M$ is locally $(\lambda,\Theta)$-strongly limit if $M$ is locally
$\{\delta < \lambda^+:\text{cf}(\delta) \in \Theta\}$-strongly limit.  
Similarly for the other notions (where 
$\Theta \subseteq \{\mu:\mu \text{ regular } \le
\lambda\},S_1 \subseteq \{\delta < \lambda^+:\text{cf}(\delta) \in \Theta\}$
is a stationary subset of $\lambda^+$).  If we do not write $\lambda$
we mean $\lambda = \|M\|$.  Let locally $(\lambda,\theta)$-strongly limit mean
locally $(\lambda,\theta)$-strongly limit. 
\nl
8) We say that $M \in K_\lambda$ is invariantly strong limit when in
part (3) we demand that $\bold F$ is just a subset of $\{(M,N)/\cong:M
\le_{\frak K} N$ are from $K_\lambda\}$ and in clause (b) of part (3)
we replace ``$\bold F(M_{i+1}) \le_{\frak K} M_{i+2}"$ by ``$(\exists
N)(M_{i+1} \le_{\frak K} N \le_{\frak K} M_{i+2} \wedge 
((M_{i+1},N)/\cong) \in \bold F)$" but abusing notation we still 
write $N = \bold F(M)$ instead $((M,N)/ \cong) \in \bold F$.  
Similarly with the other notions, so if $\bold F$ acts on suitable
$\le_{\frak K}$-increasing sequence of models then 
we use the isomorphic type of $\bar M \char 94 \langle N\rangle$.  
\bn
\margintag{88r-3.1A}\ub{\stag{88r-3.1A} Obvious implication diagram}:  For 
$\Theta,S_1$ as in \scite{88r-3.1}(7):
\bn

$$
\text{superlimit } = (\lambda,\{\mu:\mu \le \lambda \text{
regular}\})\text{-superlimit}
$$
\centerline {$\downarrow$}
$$
(\lambda,\Theta)\text{-superlimit}
$$
\centerline {$\downarrow$}
$$
S_1 \text{-strongly limit}
$$
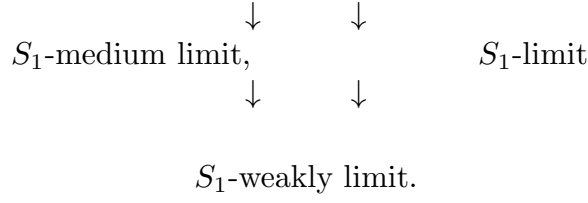
\centerline {$\downarrow$ \hskip25pt $\downarrow$}

\hskip75pt $S_1 \text{-medium limit}, \qquad \qquad \qquad \quad 
S_1\text{-limit}$

\centerline {$\downarrow$ \hskip25pt $\downarrow$}
$$
S_1 \text{-weakly limit}.
$$
\enddefinition
\bigskip

\proclaim{\stag{88r-3.2} Lemma}  0) All the properties are preserved if $S$
is replaced by a subset and if ${\frak K}$ has the $\lambda$-JEP, the
local and global version in Definition \scite{88r-3.1} are equivalent. \nl
1) If $S_i \subseteq \lambda^+$, for $i <
\lambda^+,S = \{\alpha < \lambda^+:(\exists i < \alpha) \alpha \in
S_i\}$ and $S_i \cap i = \emptyset$ for $i < \lambda$ \ub{then}: $M$ is
$S_i$-strongly limit for each $i < \lambda$ if and only if $M$ is
$S$-strongly limit. \nl
2) Suppose $\kappa \le \lambda$ is regular and $S \subseteq \{\delta <
\lambda^+:{\text{\rm cf\/}}(\delta) = \kappa\}$ is a stationary set and $M \in
K_\lambda$ \ub{then} the following are equivalent:
\mr
\item "{$(a)$}"  $M$ is $S$-strongly limit
\sn
\item "{$(b)$}"  $M$ is $(\lambda,\{\kappa\})$-strongly limit
\sn
\item "{$(c)$}"  $M \in {\frak K}_\lambda$ is $\le_{\frak
K}$-universal not $<_{\frak K}$-maximal and
there is a function $\bold F:K_\lambda \rightarrow
K_\lambda$ satisfying
$(\forall N \in K_\lambda)[N \le_{\frak K} \bold F(N)]$
such that if $M_i \in K_\lambda$ for $i < \kappa,[i < j \Rightarrow
M_i \le_{\frak K} M_j],\bold F(M_{i+1}) \le_{\frak K} 
M_{i+2}$ and $M_0 \cong M$ \ub{then} $\dbcu_{i < \kappa} M_i \cong M$. 
\ermn
2A) If $S \subseteq \lambda^+,\Theta = \{{\text{\rm cf\/}}(\delta):\delta \in
S\}$ \ub{then} $M$ is $S$-strongly limit iff clause (c) in part (2) 
above holds for every $\kappa \in \Theta$. 
\nl
3) In part (1) we can replace ``strongly limit" by ``limit", ``medium
limit" and ``weakly limit".  
\nl
4) Suppose $\kappa \le \lambda$ is regular, $S \subseteq \{\delta <
\lambda^+:{\text{\rm cf\/}}(\delta) = \kappa\}$ is a stationary set which
belongs to $\check I[\lambda]$ (see \scite{88r-0.5}, \scite{88r-0.6} above) and $M \in
K_\lambda$.
\sn
The following are equivalent
\mr
\item "{$(a)$}"  $M$ is $S$-medium limit in ${\frak K}_\lambda$
\sn
\item "{$(b)$}"  $M \in K_\lambda$ is $\le_{\frak K}$-universal not
maximal and there is a function $\bold F$ 
from ${\underset{\alpha < \kappa}
{}\to \bigcup^\alpha} (K_\lambda)$ to $K$ such that
{\roster
\itemitem{ $(\alpha)$ }   for any $\le_{\frak K}$-increasing 
$\langle M_i:i \le \alpha
\rangle$ if $M_0 = M,\alpha < \kappa,M_i$ is 
$\le_{\frak K}$-increasing,  $M_i \in K_\lambda$, then 
$M_\alpha \le_{\frak K} \bold F(\langle M_i:i \le \alpha \rangle)$
\sn
\itemitem{ $(\beta)$ }  if $\langle M_i:i < \kappa \rangle$ is
$\le_{\frak K}$-increasing, $M_0 = M,M_i \in K_\lambda$ and for $i <
\kappa$ we have $M_{i+1} \le_{\frak K} \bold F(\langle M_j:j 
\le i+1 \rangle) \le_{\frak K} M_{i+2}$ \ub{then}
$\dbcu_{i < \kappa} M_i \cong M$.
\endroster}
\endroster
\endproclaim
\bigskip

\demo{Proof}  Straightforward.
\nl
0) Trivial.
\nl
1) Recall that in Definition \scite{88r-3.1}(3), clause (b) 
we use $\bold F$ only on $M_{i+1}$; (see the proof of (2A) below,
second part).
\nl
2) For (c) $\Rightarrow$ (a) note that the demands on the sequence are
``local", $M_{i+1} \le_{\frak K} \bold F(M_{i+1}) \le_{\frak K} M_{i+2}$,
(whereas in (4) they are ``global").
\nl
2A) First assume that $M$ is $S$-strongly limit and let $\bold F$
witness it.  Suppose $\kappa \in \Theta$, so we choose $\delta_\kappa \in
S$ with cf$(\delta_\kappa) = \kappa$ and let $\langle \alpha_i:i <
\kappa \rangle$ be increasing continuous with limit $\delta,\alpha_0 =
0,\alpha_{i+1}$ a successor of a successor ordinal for each 
$i < \kappa$.  
We now define $\bold F_\kappa$ as follows: to define $\bold
F_\kappa(M)$ we define $\bold F_{\kappa,\alpha}$ for $\alpha \le
\delta$ by induction on $\alpha \le \delta$.  Let:
\mr
\item "{$(a)$}"  if $\alpha = 0$ then $\bold F_{\kappa,0}(M) = M$
\sn
\item "{$(b)$}"  if $\alpha = \beta +1$ then $\bold
F_{\kappa,\alpha}(M) = \bold F(\bold F_{\kappa,\beta}(M))$
\sn
\item "{$(c)$}"  if $\alpha \le \delta$ a limit ordinal then $\bold
F_{\kappa,\alpha}(M) = \cup\{\bold F_{\kappa,\beta}(M):\beta <
\alpha\}$.
\ermn
Lastly, let $\bold F_\kappa(M)$ be $\bold F_{\kappa,\delta}(M)$.

Now suppose $\langle N_i:i \le \kappa \rangle$ is 
$\le_{\frak K}$-increasing continuous, $N_i \in K_\lambda$ and
$\bold F_\kappa(N_{i+1}) \le_{\frak K} N_{i+2}$ for $i < \kappa$ and
we should prove $N_\kappa \cong M$.  Now we can find $\langle M_j:j <
\lambda^+ \rangle$ such that it obeys $\bold F$ and $M_{\alpha_i} =
N_i$ for $i < \kappa$< so clearly we are done.

Second, assume that for each $\kappa \in \Theta$, clause (c) of
\scite{88r-3.2}(2) holds and let $\bold F_\kappa$ exemplify this.  Let
$\langle \kappa_\varepsilon:\varepsilon < \varepsilon(*) \rangle$ list
$\Theta$ so $\varepsilon(*) < \lambda^+$ 
and define $\bold F$ as follows.  For any $M \in {\frak K}$
choose $M_{[\varepsilon]}$ by induction on $\varepsilon \le
\varepsilon(*)$ as follows: $M_{[0]} = M,M_{[\varepsilon+1]} = \bold
F_{\kappa_\varepsilon}(M_{[\varepsilon]})$ and for $\varepsilon$ limit
ordinal or $\varepsilon = \varepsilon(*)$ 
let $M_{[\varepsilon]} = \cup\{M_{[\zeta]}:\zeta
<\varepsilon\}$.  Lastly, let $\bold F[M] = M_{[\varepsilon(*)]}$.
Now check.
\nl
3) No new point.
\nl
4) First note that $(a) \Rightarrow (b)$ should be clear.  Second,
we prove that $(b) \Rightarrow (a)$ so let $\bold F$ witness
that clause (b) holds.  Let $E,\langle u_\alpha:\alpha < \lambda
\rangle$ witness that $S \in \check I[\lambda]$, i.e.
\mr
\item "{$(*)_1$}"  $(a) \quad E$ a club of $\lambda$
\sn
\item "{${{}}$}"  $(b) \quad a_\alpha \subseteq \alpha$ and 
otp$(a_\alpha) \le \kappa$ for $\alpha < \lambda$
\sn
\item "{${{}}$}"  $(c) \quad$ if $\alpha \in S \cap E$ then $\alpha =
\sup(u_\alpha)$ and otp$(u_\alpha) = \kappa$
\sn
\item "{${{}}$}"  $(d) \quad$ if $\alpha \in \lambda \backslash S \cap
E$ then otp$(u_\alpha) < \kappa$
\sn
\item "{${{}}$}"  $(e) \quad$ if $\alpha \in u_\beta$ then $u_\alpha =
u_\beta \cap \alpha$.
\ermn
We can add
\mr
\item "{$(*)_2$}"  $(f) \quad$ if $\beta \in u_\alpha$ then $\beta$ has the
form $3 \gamma+1$.
\ermn
Let $\langle \alpha_\varepsilon:\varepsilon < \lambda \rangle$ list
$E$ in increasing order and \wilog \, $\alpha_0 = 0,\alpha_{1 +
\varepsilon}$ is a limit ordinal (only the limit ordinals of $S$
count).

To define $\bold F'$ as required we shall deal with the requirement
according to whether $\delta \in S$ is ``easy", i.e. $\delta \notin E$
so $\delta \in (\alpha_\varepsilon,\alpha_{\varepsilon +1})$ so after
$\alpha_\varepsilon$ we can ``take care of it", and the ``hard"
$\delta$, i.e. $\delta \in E$ so we use the $\alpha \in u_\delta$.

We choose $\langle e_\delta:\delta \in S \backslash E \rangle$ such
that $\delta \in (\alpha_\varepsilon,\alpha_{\varepsilon +1}) \cap S$ implies
$e_\delta \subseteq \delta = \sup(e_\delta)$ and min$(e_\delta) >
\alpha_\varepsilon$, otp$(e_\delta) = \kappa,e_\delta$ is 
closed and $\alpha \in e_\delta
\Rightarrow \alpha = \sup(e_\delta \cap \alpha) \vee (\alpha$
successor).  Let $e_\delta = \langle \gamma_{\delta,\zeta}:\zeta <
\kappa \rangle$ increasing.

We now define a function $\bold F'$ so let $\langle M_j:j \le i+1
\rangle$ be given and let $\alpha_\varepsilon \le i <
\alpha_{\varepsilon +1}$.  We fix $\varepsilon$ so
$(\alpha_\varepsilon,\alpha_{\varepsilon +1})$ and now define 
$\bold F'(\langle M_j:j \le i+1 \rangle)$ by induction on $i \in
[\alpha_\varepsilon,\alpha_{\varepsilon +1})$ assuming that if
$\alpha_\varepsilon \le j' +1 < i+1$ then $\bold F'(\langle M_j:j \le
j'+1 \rangle) \le_{\frak K} M_{j'+2}$ and further there is $\bar
N^{j'+1} = \langle N_{j'+1,\xi}:\xi < \alpha_{\varepsilon +1}\rangle$ such that
the following holds:
\mr
\item "{$(*)_3$}"  $\bar N^{j'+1}$ is $\le_{{\frak
K}_\lambda}$-increasing continuous, $M_{j'+1} \le_{\frak K} N_{j'+1,0}$
and $N_{j'+1,\xi} \le_{{\frak K}_\lambda} M_{j'+2}$
\sn
\item "{$(*)_4$}"  if $\delta \in (S \backslash E) \cap
(\alpha_{\varepsilon +1} \backslash \alpha_\varepsilon),j' +1 =
\gamma_{\delta,\zeta}$ (note that necessarily $\zeta$ is a successor ordinal)
then let $\bar N^*_{\delta,j'} = \langle
N^*_{\delta,j',\zeta'}:\zeta' \le \zeta \rangle$ be the following sequence of
length $\zeta +1,N^*_{\delta,j',\zeta'}$ is $N_{\gamma_{\delta,\zeta'},\zeta'}$
if $\xi$ is a successor ordinal and is $M_{\gamma_{\delta,\zeta'}}$ if 
$\zeta'$ is limit, and we demand $\bold F(\langle
N^*_{\delta,j',\zeta'}:
\zeta' \le \zeta \rangle) \le_{\frak K} N_{j'+1,\zeta+1}$
\sn
\item "{$(*)_5$}"  if $\bold F(\langle
M_{\alpha_{f_\varepsilon(\zeta')}}:\zeta' \le \zeta\rangle)
\le_{\frak K} M_{\alpha_\varepsilon +1}$ when $\zeta' = \text{\rm
otp}(a_\varepsilon) < \kappa$ and $f_\varepsilon$ is the one-to-one
order preserving function from $\zeta' +1$ onto $c \ell(u_\zeta \cup
\{\zeta\})$ and $\zeta'$ is a successor.
\ermn
This implicitly defines $\bold F'$.  Now $\bold F'$ is as required:
$M_i \cong M$ when $i < \lambda$, cf$(i)=\kappa$ by $(*)_4$ when
$(\exists \varepsilon)(\alpha_\varepsilon < i < \alpha_{\varepsilon
+1})$ and by $(*)_5$ when $(\exists \varepsilon)(i = \alpha_\varepsilon)$.
\hfill$\square_{\scite{88r-3.2}}$
\enddemo
\bigskip

\proclaim{\stag{88r-3.3} Lemma}  Let $T$ be a first order complete theory,
$K$ its class of models and $\le_{\frak K} = \prec_{\Bbb L}$. \nl
1) If $\lambda$ is regular, $M$ a saturated model of $T$ of cardinality
$\lambda$, \ub{then} $M$ is $(\lambda,\{\lambda\})$-superlimit. 
\nl
2) If $T$ is stable, and $M$ a saturated model of $T$ of cardinality
$\lambda$ \ub{then} $M$ is $(\lambda,\{\mu:\kappa(T) \le \mu \le
\lambda$ and $\mu$ is regular$\})$-superlimit (on
$\kappa(T)$-see \cite[III,\S3]{Sh:c}).  (Note that by \cite{Sh:c} if
$\lambda$ is singular and $T$ has a saturated model of cardinality
$\lambda$ \ub{then} $T$ is stable and ${\text{\rm cf\/}}(\lambda) 
\ge \kappa(T))$. \nl
3) If $T$ is stable, $\lambda$ singular $> \kappa(T),M$ a special model
of $T$ of cardinality $\lambda,S \subseteq \{\delta <
\lambda^+:{\text{\rm cf\/}}(\delta) = \text{\rm cf}(\lambda)\}$ 
is stationary and $S \in \check I[\lambda]$ 
(see above \scite{88r-0.5}, \scite{88r-0.6}) \ub{then} 
$M$ is $(\lambda,S)$-medium limit.
\endproclaim
\bigskip

\remark{Remark}  See more in \cite{Sh:868}.
\endremark
\bigskip

\demo{Proof}  1) Because if $M_i$ is a $\lambda$-saturated model of $T$ for $i
< \delta$, cf$(\delta) \ge \lambda$, \ub{then} $\dbcu_{i < \delta}
M_i$ is $\lambda$-saturated.  Remembering the uniqueness of a
$\lambda$-saturated model of $T$ of cardinality $\lambda$ we finish. 
\nl
2) Use \cite[III,3.11]{Sh:c}: if $M_i$ is a $\lambda$-saturated model
of $T,\langle M_i:i < \delta \rangle$ increasing cf$(\delta) \ge
\kappa(T)$ \ub{then} $\dbcu_{i < \delta} M_i$ is
$\lambda$-saturated. 
\nl
3) Should be clear by now.   \hfill$\square_{\scite{88r-3.3}}$
\enddemo
\bigskip

\proclaim{\stag{88r-3.4} Claim}  1) If $M_\ell \in K_\lambda$ are 
$S_\ell$-weakly limit and $S_0 \cap S_1$ is stationary, \ub{then} 
$M_0 \cong M_1$, provided $\kappa$ has $(\lambda,\lambda)$-{\rm JEP}.
\nl
2) $K$ has at most one locally weakly limit model of cardinality $\lambda$
provided $K$ has $(\lambda,\lambda)$-{\rm JEP}. 
\nl
3) If $M \in K_\lambda$ \ub{then} $\{S \subseteq \lambda^+:M$ is $S$-weakly
limit or $S$ not stationary$\}$ is a normal ideal over
$\lambda^+$. \nl
Instead ``$S$-weakly limit", ``$S$-medium limit",``$S$-limit",
``$S$-strongly limit" can be used. 
\nl
4) In Definition \scite{88r-3.1} \wilog \, 
$\bold F(N) \cong M$ or $\bold F(\bar M)\cong M$ according to the case
(and we can add $N <_{\frak K} \bold F(N)$, etc.)
\nl
5) If $K$ is categorical in $\lambda$, \ub{then} the $M \in K_\lambda$
is superlimit provided that $K_{\lambda^+} \ne \emptyset$ (or, what is
equivalent, $M$ has a proper $\le_{\frak K}$-extension).
\endproclaim
\bigskip

\demo{Proof}  Easy.
\nl
1) E.g., let $\bold F_\ell$ witness that $M_\ell$ is $S_\ell$-weakly
limit.  We can choose $(M^0_\alpha,M^1_\alpha)$ by induction on
$\alpha$ such that: $\langle M^\ell_\beta:\beta \le \alpha \rangle$ is
$\le_{\frak K}$-increasing continuous for $\ell=0,1,M^0_\alpha
\le_{\frak K} M^1_{\alpha +1},M^1_\alpha \le_{\frak K} M^0_{\alpha
+1}$ and $\bold F_\ell(\langle M^\ell_\beta:\beta \le \alpha +1\rangle) \le
M^\ell_{\alpha +2}$.  So for some club $E_\ell$ of $\lambda^+,\delta
\in E_\ell \Rightarrow M^\ell_\delta \cong M_\ell$ for $\ell=0,1$.
But $S_0 \cap S_1$ is stationary hence there is a limit ordinal
$\delta \in S_0 \cap
S_1 \cap E_0 \cap E_1$, hence $M_0 \cong M^0_\delta = M^1_\delta \cong
M_1$ as required.  \hfill$\square_{\scite{88r-3.4}}$ 
\enddemo
\bigskip

\proclaim{\stag{88r-3.5} Theorem}  If $2^\lambda < 2^{\lambda^+},
M \in K_\lambda$ is $S$-weakly limit, $S$ is not small (see Definition
\scite{88r-0.wD}) and $M$ does not have the
$\lambda$-amalgamation property (in ${\frak K}$) \ub{then} 
$\dot I(\lambda^+,K) =
2^{\lambda^+}$, moreover there is no universal member in
${\frak K}_{\lambda^+}$ and $(2^\lambda)^+ < 2^{\lambda^+} \Rightarrow
\dot I \dot E(\lambda^+,K) = 2^{\lambda^+}$, that is 
there are $2^{\lambda^+}$ models $M \in K_{\lambda^+}$ no one
$\le_{\frak K}$-embeddable into another.
\endproclaim
\bigskip

\remark{\stag{88r-3.5A} Remark}  1) We can define a superlimit for a
family of models, i.e., when $\bold N  = 
\{N_t:t \in I\} \subseteq {\frak K}_\lambda$ is superlimit
(i.e., if $\langle M_i:i < \delta \rangle$ is $\le_{\frak
K}$-increasing, $i < \delta \Rightarrow M_i \in {\frak
K}_\lambda,\delta$ a limit ordinal $< \lambda^+,M_\delta = \cup\{M_i:i <
\delta\}$ \ub{then} $\dsize \bigwedge_{i < \delta} \, \dsize
\bigvee_{t \in I} M_i \cong N_t \Rightarrow \dsize \bigvee_{t \in I}
M_\delta \cong N_t$  (and the other variants). 
Of course, the family is $\subseteq K_\lambda$ and is not empty.
Essentially everything generalizes \ub{but} in \scite{88r-3.5} the
hypothesis should be stronger: the family should satisfy that any member does
not have the amalgamation property.  But this complicates the
situation, and the gain is not clear, so we abandon this.
\nl
2) We can many times (and in particular in \scite{88r-3.5})
strengthen ``there is no $\le_{\frak K}$-universal 
$M \in K_{\lambda^+}$" to ``there is
no $M \in K_\mu$ into which every $N \in K_{\lambda^+}$ can be
$\le_{\frak K}$-embedded" for $\mu$ not too large.  
We need $\neg \text{ Unif}(\lambda^+,S,2,\mu)$ (see \cite[AP,\S1]{Sh:f}).
\endremark
\bigskip

\demo{Proof}  Let $\bold F$ be as in Definition \scite{88r-3.1}(5) for $M$.
We now choose by induction on $\alpha < \lambda^+$, models
$M_\eta$ for $\eta \in {}^\alpha 2$ such that:
\mr
\widestnumber\item{$(iii)$}
\item "{$\circledast_1$}"  $(i) \quad M_\eta \in K_\lambda,M_{<>} = M$,
\sn 
\item "{${{}}$}"  $(ii) \quad$ for $\beta < \alpha,\eta \in {}^\alpha 2,
M_{\eta \restriction \beta} \le_{\frak K} M_\eta$
\sn
\item "{${{}}$}"  $(iii) \quad$ for $i +2 \le \alpha,\eta \in
{}^\alpha 2,(\bold F(\langle M_{\eta \restriction j}:j \le i+1\rangle)) 
\le_{\frak K} M_{\eta \restriction (i+2)}$ 
\sn
\item "{${{}}$}"  $(iv) \quad$ if 
$\alpha = \beta +1$ and $\beta$ non-limit, $\eta \in
{}^\alpha 2$, then $M_{\eta \restriction \beta} \ne M_\eta$
\sn
\item "{${{}}$}"  $(v) \quad$ if $\alpha < \lambda$ is a limit ordinal and
$\eta \in {}^\alpha 2$ \ub{then}:
{\roster
\itemitem{ ${{}}$ }  $(a) \quad M_\eta = 
\cup\{M_{\eta \restriction \beta}:\beta <\ell g(\eta)\}$ and
\sn
\itemitem{ ${{}}$ }   $(b) \quad$ if $M_\eta$ fails the 
$\lambda$-amalgamation property \ub{then}
$M_{\eta \char 94 <0>},M_{\eta \char 94 <1>}$
\nl

\hskip25pt  cannot be amalgamated
over $M_\eta$, i.e. for no $N$ do we have: 
\nl

\hskip25pt $M_\eta \le_{\frak K} N \in K$ and 
$M_{\eta \char 94 <0>},M_{<\eta \char 94 <1>}$ can be 
$\le_{\frak K}$-embedded 
\nl

\hskip25pt into $N$ over $M_\eta$.
\endroster}
\ermn
For $\alpha =0,\alpha$ limit, we have no problem, for $\alpha
+1,\alpha$-limit: if $M_\eta$ fails the $\lambda$-amalgamation
property - use its definition, otherwise let 
$M_{\eta \char 94 <1>} = M_\eta = M_{\eta \char 94 <0>}$; 
for $\alpha +1,\alpha$ non-limit - use $\bold F$ to guaranteee clause
(iii), and then for clause (iv) use clause $(\gamma)$ of 
Definition \scite{88r-3.1}(5), i.e., \scite{88r-3.1}(4).

Let for $\eta \in {}^{\lambda^+} 2,M_\eta = \dbcu_{\alpha < \lambda^+}
M_{\eta \restriction \alpha}$.  By changing names we can assume that
\mr
\item "{$\circledast_1$}"  $(vi) \quad$ for 
$\eta \in {}^\alpha 2(\alpha < \lambda^+)$ the universe of $M_\eta$
is an ordinal $< \lambda^+$ (or even
\nl

\hskip25pt  $\subseteq \lambda \times (1 + \ell g(\eta))$ and we could even demand equality).
\ermn
So (by  clause (iv)) for
$\eta \in {}^{\lambda^+} 2,M_\eta$ has universe $\lambda^+$.

First, why is there no universal member in ${\frak K}_{\lambda^+}$?  If $N \in
K_{\lambda^+}$ is universal (by $\le_{\frak K}$, of course), 
without loss of generality 
its universe is $\lambda^+$.  For
$\eta \in {}^{\lambda^+} 2$ as $M_\eta \in K_{\lambda^+}$, there is a
$\le_{\frak K}$-embedding $f_\eta$ of $M_\eta$ into $N$.  So $f_\eta$ is a
function from $\lambda^+$ to $\lambda^+$.  Let $\eta \in
{}^{\lambda^+} 2$, by the choice of $\bold F$ and of $\langle M_{\eta
\restriction \alpha}:\alpha < \lambda^+ \rangle$ there is a closed
unbounded $C_\eta \subseteq \lambda^+$ such that $\alpha \in S \cap
C_\eta \Rightarrow M_{\eta \restriction \alpha} \cong M$, hence 
$M_{\eta \restriction \alpha}$ fails the $\lambda$-amalgamation
property.  Without loss of generality for $\delta \in C_\eta,M_{\eta
\restriction \delta}$ has universe $\delta$.  Now by \scite{88r-0.wD}, if
$\langle (f_\rho,C_\rho):\rho \in {}^{\lambda^+} 2 \rangle$ satisfies
that for each $\rho \in {}^{\lambda^+}2,f_\rho:\lambda^+ \rightarrow
\lambda^+,C_\rho \subseteq \lambda^+$ is closed unbounded \ub{then} for
some $\eta \ne \nu \in {}^{\lambda^+}2$ and $\delta \in C_\eta \cap S$
we have $\eta \restriction \delta = \nu \restriction \delta,\eta(\delta) \ne
\nu(\delta)$ and $f_\eta \restriction \delta = f_\nu \restriction
\delta$.
\nl
[Why?  For every $\delta < \lambda^+,\rho \in {}^\delta 2$ and
$f:\delta \rightarrow \lambda^+$ we define $\bold c(\rho,f) \in 2$ as
follows: it is 1 iff there is $\nu \in {}^{\lambda^+}2$ such that $\rho = \nu
\restriction \delta \and f = f_\nu \restriction \delta \and
\nu(\delta) = 0$ and is 0 otherwise.  So some $\eta \in
{}^{\lambda^+}2$ is a weak diamond sequence for the colouring $\bold
c$ and the stationary set $S$.  
Now $C_\eta,f_\eta$ are well defined and $S' = \{\delta \in
S:\delta$ limit and $\eta(\delta) = \bold c(\eta \restriction \delta,f
\restriction \delta)\}$ is a stationary subset of $\lambda^+$, so we
can choose $\delta \in S' \cap C_\eta$.  If $\eta(\delta) = 0$, then
$\bold c(\eta \restriction \delta,f \restriction \delta) = 0$ but
$\eta$ witness that $\bold c(\eta \restriction \delta,f \restriction
\delta)$ is 1.  If $\eta(\delta)=1$ there is $\nu$ witnessing $\bold
c(\eta \restriction \delta,f_\eta \restriction \delta)=1$, in
particular $\nu(\delta)=0$,  so $\eta,\nu,\eta
\restriction \delta$, are as required.]

Now as $\delta \in S \cap C_\eta$ it follows that
$M_{\eta \restriction \delta} \cong M$ hence $M_{\eta \restriction
\delta}$  fails the $\lambda$-amalgamation property.  Also $M_{\eta
\restriction \delta}$ has universe $\delta$ as $\delta \in C_\eta$ and
$M_{\eta \restriction \delta} = M_{\nu \restriction \delta}$ as $\eta
\restriction \delta = \nu \restriction \delta$.

So $f_\eta \restriction M_{\eta \restriction \delta} = f_\eta
\restriction \delta = f_\nu \restriction \delta = f_\nu \restriction
M_{\nu \restriction \delta}$.  So $f_\eta \restriction
M_{\eta \restriction(\delta +1)},f_\nu \restriction M_{\nu
\restriction (\delta +1)}$ show that $M_{\eta \restriction (\delta
+1)},M_{\nu \restriction (\delta +1)}$, can be amalgamated over
$M_{\eta \restriction \delta}$ contradicting clause (v) of the
construction.  So there is no $\le_{\frak K}$-universal $N \in {\frak
K}_{\lambda^+}$. 

It takes some more effort to get $2^{\lambda^+}$ pairwise 
non-isomorphic models (rather than just quite many).
\enddemo
\bn
\ub{Case A} \footnote{we can make it a separate claim}:  
There is $M^* \in K_\lambda,M \le_{\frak K} M^*$ such that for
every $N$ satisfying $M^* \le_{\frak K} N \in K_\lambda$ 
there are $N^1,N^2 \in K_\lambda$ such that $N \le_{\frak K} N^1,
N \le_{\frak K} N^2$ and $N^2,N^1$ cannot be 
$\le_{\frak K}$-amalgamated over $M^*$ (not just $N$).
In this case we do not need ``$M$ is $S$-weakly limit".

We redefine $M_\eta,\eta \in {}^\alpha 2,\alpha < \lambda^+$ such
that:
\mr
\item "{$\circledast_2$}"  $(a) \quad \nu \triangleleft \eta \in {}^\alpha 2
\Rightarrow M_\nu \le_{\frak K} M_\eta \in K_\lambda$:
\sn
\item "{${{}}$}"  $(b) \quad$ if $\alpha =0,M_{<>} = M^*$; \nl
\sn
\item "{${{}}$}"  $(c) \quad$ if 
$\alpha$ limit and $\eta \in {}^\alpha 2:M_\eta 
= \dbcu_{\beta < \alpha} M_{\eta \restriction \beta}$; 
\sn
\item "{${{}}$}"  $(d) \quad$ if 
$\eta \in {}^\beta 2,\alpha = \beta +1$, use the assumption for $N =
M_\eta$, now 
\nl

\hskip25pt obviously the $(N^1,N^2)$ there satisfies
$N^1 \ne N$ and $N^2 \ne N$, so we 
\nl

\hskip25pt can have
$M_\eta <_{\frak K} M_{\eta \char 94 <1>} \in K_\lambda,
M_\eta <_{\frak K} M_{\eta \char 94 <0>} \in K_\lambda$, such that 
\nl

\hskip25pt $M_{\eta \char 94 <0>},M_{\eta \char 94 <1>}$
cannot be amalgamated over $M^*$.
\ermn
Obviously, the models $M_\eta =
\dbcu_{\alpha < \lambda^+} M_{\eta \restriction \alpha}$, for $\eta
\in {}^{\lambda^+} 2$ are pairwise non-isomorphic over $M^*$ and by
\scite{88r-0.9} as $2^\lambda < 2^{\lambda^+}$ we 
finish proving $\dot I(\lambda^+,{\frak K}) = 2^{\lambda^+}$.

Note also that for each $\eta \in {}^{\lambda^+} 2$ the set $\{\nu \in
{}^{\lambda^+} 2:M_\nu$ can be $\le_{\frak K}$-embedded into
$M_\eta\}$ has cardinality $\le |\{f:f$ a $\le_{\frak K}$-embedding of
$M^*$ into $M_\eta\}| \le 2^\lambda$.
So if $(2^\lambda)^+ < 2^{\lambda^+}$, then by Hajnal free subset
theorem (\cite{Ha61}), there are $2^{\lambda^+}$ models $M_\eta \in
K_{\lambda^+}(\eta \in {}^{\lambda^+} 2)$ no one 
$\le_{\frak K}$-embeddable into another.  
\bn
\ub{Case B}:  Not Case A.

Now we return to the first construction, but we can add
\mr
\item "{$(vii)$}"  if $\eta \in {}^{(\alpha +1)} 2$, \ub{then} if
$M_\eta \le_{\frak K} N^1,N^2$ both in $K_\lambda$, \ub{then} $N^1,N^2$ can be
$\le_{\frak K}$-amalgamated over $M_{\eta \restriction \alpha}$.
\ermn
As $\{W \subseteq \lambda^+:W$ is small$\}$ is a normal ideal (see
\scite{88r-0.wD})(and it is on a successor cardinal) it is well known that
we can find $\lambda^+$ pairwise disjoint non-small $S_\zeta \subseteq
S$ for $\zeta < \lambda^+$.  We define a colouring (= function) 
$\bold c$:
\mr
\item "{$\circledast_3$}"  $(a) \quad \bold c(\eta,\nu,f)$ will be 
defined \ub{iff} for some limit ordinal 
$\delta < \lambda^+,\eta \in {}^\delta 2,
\nu \in {}^\delta 2$ 
\nl

\hskip25pt and $f$ is a function from $\delta$ to $\lambda^+$
\sn
\item "{${{}}$}"  $(b) \quad \bold c(\eta,\nu,f)=1$ 
\ub{iff} the triple $(\eta,\nu,f)$ belongs to the domain of $\bold c$
(i.e., is 
\nl

\hskip25pt as in (a)) and $M_\eta,M_\nu$ have universe 
$\delta,f$ is a $\le_{\frak K}$-embedding
of $M_\eta$ 
\nl

\hskip25pt into $M_\nu$ and for some $\rho,\nu \char
94 <0> \triangleleft \rho \in {}^{\lambda^+}2$ the function 
$f$ can be 
\nl

\hskip25pt extended to a 
$\le_{\frak K}$-embedding of $M_{\eta \char 94 <0>}$ into $M_\rho$
\sn
\item "{${{}}$}"  $(c) \quad \bold c(\eta,\nu,f)$ is 
zero \ub{iff} it is defined but is $\ne 1$.
\ermn
For each $\zeta$, as $S_\zeta$ is not small, by simple coding, for
every $\zeta < \lambda^+$ there is $h_\zeta:S_\zeta 
\rightarrow \{0,1\}$ such that:
\mr
\item "{$(*)_\zeta$}"   for every $\eta \in {}^{\lambda^+}2,\nu \in
{}^{\lambda^+}2$ and $f:\lambda^+ \rightarrow \lambda^+$, for a stationary
set of $\delta \in S_\zeta$

$$
\bold c(\eta \restriction \delta,\nu \restriction \delta,f \restriction
\delta) = h_\zeta(\delta).
$$
\ermn
Now for every $W \subseteq \lambda^+$ we define $\eta_W \in
{}^{\lambda^+}2$ as follows:

$\eta_W(\alpha)$ is $h_\zeta(\alpha)$, if $\zeta \in W$ and
$\alpha \in S_\zeta$ (note that there is at most one $\zeta$) 

$\eta_W(\alpha)$ is zero if there is no such $\zeta$.

Now we can show (chasing the definitions) that
\mr
\item "{$\circledast_4$}"  if $W(1),W(2)
\subseteq \lambda^+,W(1) \nsubseteq W(2)$, \ub{then} 
$M_{\eta_{W(1)}}$ cannot be
$\le_{\frak K}$-embedded into $M_{\eta_{W(2)}}$.
\ermn
This clearly suffices. 
\nl
Why is $\circledast_4$ true?  Suppose $W(1) \nsubseteq W(2)$, let
$\zeta \in W(1) \backslash W(2)$ and toward contradiction let $f$ be a
$\le_{\frak K}$-embedding of $M_{\eta_{W(1)}}$ into 
$M_{\eta_{W(2)}}$, so $E = \{\delta:
M_{\eta_{W(1)} \restriction \delta},M_{\eta_{W(2)} \restriction
\delta}$ have universe $\delta$ and $f \restriction \delta$
is a $\le_{\frak K}$-embedding of $M_{\eta_{W(1)} \restriction
\delta}$ into $M_{\eta_{W(2)} \restriction \delta}\}$ is a club of
$\lambda^+$.  Hence by the choice of $\bold c$ and $h_\zeta$ there 
is $\delta \in E \cap S_\zeta$ such that
\mr
\item "{$\boxtimes$}"  $\bold c(\eta_{W(1)} \restriction
\delta,\eta_{W(2)} \restriction \delta,f \restriction \delta) =
h_\zeta(\delta)$ and $M_{\eta_{w(1)} \restriction \delta}$ is not an
amalgamation base.
\ermn
Now the proof splits to two cases.
\bn
\ub{Case 1}:  $h_\zeta(\delta) =0$.

So $\eta_{W(1)}(\delta) =0 = \eta_{W(2)}(\delta)$ and by clause (b) of
$\circledast_3$ above, i.e., the definition of $\bold c$ we have the objects
$\eta_{W(1)},\eta_{W(2)},f \restriction
M_{\eta_{W(1)} {}^\frown <0>} = f \restriction M_{\eta_{W(1)}
\restriction (\delta +1)}$ witness that $\bold c(\eta_{W(1)}
\restriction \delta,\eta_{W(2)} \restriction \delta,f \restriction
\delta)=1$, contradiction.
\bn
\ub{Case 2}:  $h_\zeta(\delta) =1$.

So $\eta_{W(1)}(\delta) =1,\eta_{W(2)}(\delta)=0,\bold c(\eta_{W(1)}
\restriction \delta,\eta_{W(2)} \restriction \delta,f \restriction
\delta)=1$.  By the definition of $\bold c$, we can find $\nu$ such that
$(\eta_{W(2)} \restriction \delta) {}^\frown
<0> \trianglelefteq \nu \in {}^{\lambda^+}2$ and a $\le_{\frak
K}$-embedding $g$ of $M_{(\eta_{W(1)} \restriction \delta) {}^\frown
<0>}$ into $M_\nu$.

For some $\alpha \in (\delta,\lambda^+),f$ embeds $M_{\eta_{W(1)} 
\restriction (\delta +1)} = M_{(\eta_{W(1)}  \restriction \delta)
{}^\frown <1>}$ into $M_{\eta_{W(2)} \restriction \alpha}$ and 
$g$ embeds $M_{(\eta_{W(1)} \restriction \delta) {}^\frown <0>}$ into 
$M_{\nu \restriction \alpha}$.

As $\eta_{W(2)} \restriction \delta 
{}^\frown <0> \triangleleft \nu \restriction
\alpha$ and $\eta_{W(2)} \restriction \delta
{}^\frown <0> \triangleleft \eta_{W(2)} \restriction
\alpha$ by clause (vii) above there are $f_1,g_1$ and 
$N \in K_\lambda$ such that
\mr
\item "{$(a)$}"  $M_{\eta_{W(2)} \restriction \delta} \le_{\frak K} N$
\sn
\item "{$(b)$}"  $f_1$ is a $\le_{\frak K}$-embedding of
$M_{\eta_{W(2)} \restriction \alpha}$ into $N$ over $M_{\eta_{W(2)} 
\restriction \delta}$
\sn
\item "{$(c)$}"  $g_1$ is a $\le_{\frak K}$-embedding of $M_{\nu
\restriction \alpha}$ into $N$ over $M_{\eta_{W(2)}  \restriction
\delta}$.
\ermn
So
\mr
\item "{$(b)^*$}"  $f_1 \circ f$ is a $\le_{\frak K}$-embedding of
$M_{(\eta_{W(1)} \restriction \delta) {}^\frown <1>}$ into $N$
\sn
\item "{$(c)^*$}"   $g_1 \circ g$ is a $\le_{\frak K}$-embedding of
$M_{(\eta_{W(1)} \restriction \delta) {}^\frown <0>}$ into $N$
\sn
\item "{$(d)^*$}"  $f_1 \circ f,g_1 \circ g$ extend
$f \restriction \delta:M_{\eta_{W(1)} \restriction \delta}
\rightarrow  N$ (both).
\ermn
So together we get a contradiction to assumption $(*)_1(d)$.  
\hfill$\square_{\scite{88r-3.5}}$
\bigskip

\proclaim{\stag{88r-3.6} Theorem}  1) Assume one of the following cases occurs:
\mr
\item "{$(a)_1$}"   ${\frak K}$ is {\rm PC}$_{\aleph_0}$ 
(hence {\rm LS}$({\frak K}) = \aleph_0)$ and $1 \le 
\dot I(\aleph_1,{\frak K}) < 2^{\aleph_1}$
\nl
or
\sn
\item "{$(a)_2$}"  ${\frak K}$ has models of arbitrarily large 
cardinality, {\rm LS}$({\frak K}) = \aleph_0$ and $\dot
I(\aleph_1,{\frak K}) < 2^{\aleph_1}$.
\ermn
\ub{Then} there is an a.e.c. ${\frak K}_1$ such that
\mr
\item "{$(A)$}"  $M \in K_1 
\Rightarrow M \in K$ and $M \le_{{\frak K}_1} N
\Rightarrow M \le_{\frak K} N$ and {\rm LS}$({\frak K}_1) 
= \text{ LS}({\frak K}) (= \aleph_0)$
\sn
\item "{$(B)$}"  if $K$ has models of arbitrarily large cardinality
\ub{then} so does $K_1$
\sn
\item "{$(C)$}"  ${\frak K}_1$ is {\rm PC}$_{\aleph_0}$
\sn
\item "{$(D)$}"  $(K_1)_{\aleph_1} \ne \emptyset$
\sn
\item "{$(E)$}"  all models of $K_1$ are
$\Bbb L_{\infty,\omega}$-equivalent and $M \le_{{\frak K}_1} N \Rightarrow M
\prec_{{\Bbb L}_{\infty,\omega}} N$ and $K_1$ is categorical in $\aleph_0$
\sn
\item "{$(F)$}"  if ${\frak K}$ is categorical in $\aleph_1$ \ub{then}
$({\frak K}_1)_\lambda = {\frak K}_\lambda$ for every $\lambda > \aleph_0$.
\ermn
2) If in (1) we add {\rm LS}$({\frak K})$ names to formulas in
$\Bbb L_{\infty,\omega}$ (i.e. to a set of representatios up to
equivalence) \ub{then} we can assume each member of $K$ is
$\aleph_0$-sequence-homogeneous.  The vocabulary remains countable, in fact,
for some countable first order theory $T$, the models of $K$ are the
atomic models of $T$ (in the first order sense) and
$\le_{\frak K}$ becomes $\subseteq$ (being a submodel).
\endproclaim
\bigskip

\demo{Proof}  Like \cite[2.3,2.5]{Sh:48} (using \scite{88r-2.12} here for
$\alpha =2$).
\enddemo
\bn
We arrive to the main theorem of this section.
\proclaim{\stag{88r-3.7} Theorem}  Suppose ${\frak K}$ and $\lambda$
satisfy the following conditions:
\mr
\item "{$(A)$}"  ${\frak K}$ 
has a superlimit member $M^*$ of cardinality
$\lambda,\lambda \ge { \text{\rm LS\/}}({\frak K})$, (if $K$ is categorical in
$\lambda$, then by assumption $(B)$ below there is such $M^*$; really
invariantly 
$\lambda^+$-strongly limit suffice if (d) of $(*)$ of \scite{88r-3.7A}(2)
below holds, see Definition \scite{88r-3.1})
\sn
\item "{$(B)$}"  ${\frak K}$ is categorical in $\lambda^+$
\sn
\item "{$(C)$}"  $(\alpha) \quad {\frak K}$ is {\rm PC}$_{\aleph_0},
\lambda = \aleph_0$ or
\sn
\item "{${{}}$}"  $(\beta) \quad {\frak K} 
= { \text{\rm PC\/}}_\lambda,\lambda = \beth_\delta,
{\text{\rm cf\/}}(\delta) = \aleph_0$ or
\sn
\item "{${{}}$}"   $(\gamma) \quad \lambda = \aleph_1,{\frak K}$ is
{\rm PC}$_{\aleph_0}$ or
\sn
\item "{${{}}$}"  $(\delta) \quad 
{\frak K}$ is {\rm PC}$_\mu,\lambda \ge \beth_{(2^\mu)^+}$.
\ermn
\ub{Then} $K$ has a model of cardinality $\lambda^{++}$.
\endproclaim 
\bigskip

\remark{\stag{88r-3.7A} Remark}  1) If 
$\lambda = \aleph_0$ we can wave hypothesis (A)
by the previous theorem \scite{88r-3.6}. \nl
2) Hypothesis (C) can be replaced by (giving a stronger theorem):
\mr
\item "{$(*)_{\lambda,\mu}(a)$}"  ${\frak K}$ is PC$_\mu$ and
\sn
\item "{$(b)$}"   any $\psi \in \Bbb L_{\mu^+,\omega}$
which has a model $M$ of order-type $\lambda^+,|P^M| = \lambda$, has a
non-well-ordered model $N$ of cardinality $\lambda$ 
\sn
\item "{$(c)$}"  $\{M \in K_\lambda:M \cong M^*\}$ is 
PC$_\mu$ (among models in $K_\lambda$) and 
\sn
\item "{$(d)$}"  for some $\bold F$ witnessing ``$M^*$ is invariantly
$\lambda$-strongly limit" the class  $\{(M,\bold F(M)):
M \in K_\lambda\}$ is PC$_\mu$ (if $M^*$ is superlimit this clause is not
required as $\bold F =$ the identity on $K_\lambda$ is O.K.)  
\ermn
3) It is well known, see e.g. \cite[VII,\S5]{Sh:c} that hypothesis (C)
implies $(*)_{\lambda,\mu}$ from part (2), see more \cite{GrSh:259}.
\endremark
\bigskip

\demo{Proof}  By \scite{88r-3.7A}(3) we can assume $(*)_{\lambda,\mu}$ from
\scite{88r-3.7A}(2).
\enddemo
\bn
\ub{Stage a}:  It suffices to find $N_0 \le_{\frak K} N_1,\|N_0\| =
\lambda^+,N_0 \ne N_1$.

Why?  We define by induction on $\alpha < \lambda^{++}$ a model $N_\alpha
\in K_{\lambda^+}$ such that $\beta < \alpha$ implies $N_\beta \le_{\frak K}
N_\alpha$ and $N_\beta \ne N_\alpha$.  Clearly $N_0,N_1$ are defined 
(without loss of generality
$\|N_1\| = \lambda^+$ as $\lambda \ge \text{ LS}({\frak K})$, also otherwise
we already have the desired conclusion), for limit $\delta
< \lambda^{++}$ the model 
$\dbcu_{\alpha < \delta} N_\alpha$ is as required.  For
$\alpha = \beta +1$, by the $\lambda^+$-categoricity, $N_0$ is
isomorphic to $N_\beta$, say by $f$ and we define $N_{\beta +1}$ such
that $f$ can be extended to an isomorphism from $N_1$ onto $N_{\beta
+1}$, so clearly $N_{\beta +1}$ is as required.  Now $\dbcu_{\alpha <
\lambda^{++}} N_\alpha \in K_{\lambda^{++}}$ is as required.
\bn
Hence the following theorem completes the proof of \scite{88r-3.7} 
(use $\bold F =$ the identity for the superlimit case).
\proclaim{\stag{88r-3.8} Theorem}  Suppose the following clauses:
\mr
\item "{$(A)$}"  ${\frak K}$ has an invariantly $\lambda$-strongly limit member
$M^*$ of cardinality $\lambda$, as exemplified by 
$\bold F:K_\lambda \rightarrow
K_\lambda$ and ${\frak K}_\lambda$ has the {\rm JEP} (see Definition
\scite{88r-3.1}) 
\sn
\item "{$(B)$}"  $\dot I(\lambda^+,K_{\lambda^+}) < 2^{\lambda^+}$ or even
just $\dot I(\lambda^+,K^{\bold F}_{\lambda^+}) < 2^{\lambda^+}$ 
(or just $\dot I \dot E(\lambda^+,K^{\bold F}_{\lambda^+}) < 2^{\lambda^+}$
(see below))
\sn
\item "{$(C)$}"  ${\frak K}$ is a {\rm PC}$_\mu$ class, as well as 
$\bold F$, i.e., $K'$ is {\rm PC}$_\mu$ where $K'$ is a class closed under an
isomorphism  of $(\tau_{\frak K} \cup \{P\})$-models, $P$ a unary
predicate such that $K' = \{(N,M):N = \bold F(M)\}$
\sn
\item "{$(D)$}"  $\mu = \lambda = \aleph_0$ or $\mu =
\lambda = \beth_\delta,{\text{\rm cf\/}}(\delta) = \aleph_0$ or $\mu =
\aleph_0,\lambda = \aleph_1$ or just $(*)_{\lambda,\mu}$ from \scite{88r-3.7A}(2)
\sn
\item "{$(E)$}"   $K$ categorical in $\lambda$ or at
least there is $\psi \in \Bbb L_{\omega_1,\omega}(\tau^+)$ such that 
$(M^* /\cong) = \{M \restriction \tau_{\frak K}:M \models \psi,\|M\| =
\lambda\}$. 
\endroster
\endproclaim
\bn
\ub{Then} we can find $N_0 \le_{\frak K} N_1,
N_0 \ne N_1$ such that $N_0,N_1 \in K^{\bold F}_{\lambda^+}$, \nl
where 
\definition{\stag{88r-3.8.1} Definition}  Assume $\bold F:K_\lambda
\rightarrow K_\lambda$ satisfies $M \le_{\frak K} \bold F(M)$ for $M \in
K_\lambda$ or more generally $\bold F \subseteq \{(M,N):M
\le_{\frak K} N$ are from $K_\lambda\}$ satisfies $(\forall M \in
K_\lambda)(\exists N)((M,N) \in  \bold F)$ or just $(\forall M \in
K_\lambda)(\exists N_0,N_1)[(N_0,N_1) \in \bold F \wedge M \le_{\frak K} 
N_0 \le_{\frak K} N_1]$.  Then we let
$K^{\bold F}_{\lambda^+} =: \{\dbcu_{i < \lambda^+} M_i:M_i \in
K_\lambda,\langle M_i:i < \lambda^+ \rangle$ is $\le_{\frak
K}$-increasing continuous not eventually constant and 
$\bold F(M_{i+1}) \le_{\frak K} M_{i+2}$
or $(M_{i+1},M_{i+2}) \in \bold F\}$.
\enddefinition
\bigskip

\remark{\stag{88r-3.8.3} Remark}  1) As the sequence in the definition of $K^{\bold
F}_{\lambda^+}$ is $\le_{\frak K}$-increasing and $(M,N) \in \bold F
 \Rightarrow M \ne N$  necessarily $K^{\bold F}_{\lambda^+}
\subseteq {\frak K}_{\lambda^+}$.
\nl
2) Theorem \scite{88r-3.8} is good for classes which are not
exactly a.e.c., see, e.g., \scite{88r-3.9}.
\endremark
\bn
Considering $K^{\bold F}_{\lambda^+}$ we may note that the proofs of
some earlier claims give more.  In particular, similarly to
\scite{88r-3.5}
\proclaim{\stag{88r-3.8.4} Claim}  Assume that
\mr
\item "{$(a)$}"  $2^\lambda < 2^{\lambda^+}$
\sn
\item "{$(b)$}"  ${\frak K}$ is an a.e.c. and {\rm LS}$({\frak K}) \le
\lambda$
\sn
\item "{$(c)$}"  $M \in K_\lambda$ is $S$-weakly limit, $S$ not small
(see Definition \scite{88r-0.wD})
\sn
\item "{$(d)$}"  $M$ does not have the amalgamation property in
${\frak K}$ (= is an amalgamation base)
\sn
\item "{$(e)$}"  $\bold F$ is as in \scite{88r-3.8.1}.
\ermn
\ub{Then} $\dot I(\lambda^+,K^{\bold F}_{\lambda^+}) = 2^{\lambda^+}$.
\endproclaim
\bigskip

\demo{Proof}  To avoid confusion rename $\bold F$ of clause (e) as
$\bold F_1$, and choose $\bold F_2$ which 
exemplifies ``$M$ is $S$-weakly limit",
i.e., as in Definition \scite{88r-3.1}(5).  Now we define $\bold F'$ with
the same domain as $\bold F_2$ by $\bold F'(\langle M_j:j \le i
\rangle) = \bold F_1(\bold F_2(\langle M_j:j \le i \rangle)$, and
continue as in the proof of \scite{88r-3.5} noting that $\bold F'$ works
as well there.

The sequence of models $\langle M_\eta:\eta \in {}^{\lambda^+}2\rangle$ we got
there are from $K^{\bold F_1}_{\lambda^+}$ (so witness that $\dot
I(\lambda^+,K^{\bold F_1}_{\lambda^+}) = 2^{\lambda^+}$) because:
\mr
\item "{$(*)$}"  if the sequence $\langle M_\alpha:\alpha <
\lambda^+\rangle,M_\alpha \in {\frak K}_\lambda$ for $\alpha <
\lambda^+$ is $\le_{\frak K}$-increasing continuous and $\bold
F'(\langle M_j:j \le i +1 \rangle) \le_{\frak K} M_{i+2}$ then
$\cup\{M_\alpha:\alpha < \lambda^+\}\in K^{\bold F_1}_{\lambda^+}$.
\endroster
\hfill$\square_{\scite{88r-3.8.4}}$
\enddemo
\bn
Also similarly to \scite{88r-3.6} we can prove:
\proclaim{\stag{88r-3.8.5} Claim}  Assume ${\frak K}$ is a
{\rm PC}$_{\aleph_0}$ and $\bold F$ a {\rm PC}$_{\aleph_0}$ is 
as in \scite{88r-3.8.1}. If $1 \le \dot I(\aleph_1,
K^{\bold F}_{\aleph_1}) < 2^{\aleph_1}$ then the
conclusion of \scite{88r-3.6} above holds.
\endproclaim
\bigskip

\demo{Proof of \scite{88r-3.8}}  (Hence of \scite{88r-3.7}).  The reader may
do well to read it with $\bold F =$ the identity in mind.
\enddemo
\bn
\ub{Stage b}:  We now try to find $N_0,N_1$ as mentioned in stage (a) above by
approximations of cardinality $\lambda$.  A triple will denote here
$(M,N,a)$ satisfying 
$M,N \cong M^*$ (see hypothesis $(A)$), $M \le_{\frak K} N$ and $a \in
N \backslash M$.  Let $<$ be the following partial 
order among this family of triples: $(M,N,a) <
(M',N',a')$ if $a = a',N \le_{\frak K} N',M \le_{\frak K} M',
M \ne M'$ and moreover $(\exists N'')[N \le_{\frak K} N'' \and \bold F(N'') 
\le_{\frak K} N']$ and $(\exists M'')[M \le_{\frak K} M'' \and 
\bold F(M'') \le_{\frak K} M']$.  
(It is tempting to omit $a$ and require $M = M' \cap N$, but
this apparently does not work as we do know if disjoint amalgamation
${\frak K}_{\aleph_0}$ exist).

We first note that there is at least one triple (as $M^*$ has a proper
elementary extension which is isomorphic to it, because it is a limit
model).
\bn
\ub{Stage c}:  We show that if there is no maximal triple, our
conclusions follows.  

We choose by induction on $\alpha$ a triple $(M_\alpha,N_\alpha,a)$ 
increasing by $<$.  For
$\alpha=0$ see the end of previous stage, for $\alpha = \beta +1$, we can
define $(M_\alpha,N_\alpha,a)$ by the hypothesis of this stage.  For
limit $\delta < \lambda^+,(M_\delta,N_\delta,a)$ will be
$(\dbcu_{\alpha < \delta} M_\alpha,\dbcu_{\alpha < \delta}
N_\alpha,a)$ 
(notice $M_\delta \le_{\frak K} N_\delta$ by AxIV and $M_\delta,N_\delta$ are
isomorphic to $M^*$ by the choice of $\bold F$ and the definition of
order on the family of triples).  Now similarly 
$M = \dbcu_{\alpha < \lambda^+} M_\alpha \le_{\frak K} 
N = \dbcu_{\alpha < \lambda^+} N_\alpha$ are both from 
${\frak K}^{\bold F}_{\lambda^+}$ and the element $a$ exemplifies
$M \ne N$, so by Stage (a) we finish.  
\nl
Recall
\mr
\item "{$\circledast$}"  if $(M,N,a)$ is a maximal triple \ub{then}
there is no triple $(M',N',a)$ such that $M' \le_{\frak K} N',
M <_{\frak K} M',N  \le_{\frak K} N',a \in N' \backslash M'$ and
$(\exists M'')(M \le_{\frak K} M'' \le_{\frak K} \bold F(M'')
\le_{\frak K} M')$ and $(\exists N'')(N \le_{\frak K} N'' \le_{\frak
K} \bold F(N'') \le_{\frak K} N')$.
\endroster
\bn
\ub{Stage d}:  There are $M_i \cong M^*$ for $i \le \omega$ such that
$[i < j \le \omega \Rightarrow M_j <_{\frak K} M_i],i < \omega \Rightarrow
\bold F(M_{i+1}) \le_{\frak K} M_i$ and
$|M_\omega| = \dbca_{n < \omega} |M_n|$ and note that $M_i$ is
$\lambda^+$-strongly limit.
\nl
This stage is dedicated to proving this statement.
As $M^*$ is superlimit (or just strongly limit), there is an 
$\le_{\frak K}$-increasing
continuous sequence $\langle M_i:i < \lambda^+ \rangle,M_i \cong
M^*$ and $\bold F(M_{i+1}) \le_{\frak K} M_{i+2}$.  (Note that this is
true also for limit models as we can restrict
ourselves to a club of $i$'s).
So without loss of generality $\dbcu_{i < \lambda^+} M_i$ has universe 
$\lambda^+,M_0$ has universe $\lambda$.
\nl
Define a model ${\frak B}$.

Its universe is $\lambda^+$.
\bn
\ub{Relations and Functions}:
\mr
\item "{$(a)$}"  those of $\dbcu_{i < \lambda^+} M_i$
\sn
\item "{$(b)$}"  $R$-two place: $aRi$ if and only if $a \in M_i$
\sn
\item "{$(c)$}"  $P$ (monadic relation) $P = \lambda$ which is the
universe of $M_0$
\sn
\item "{$(d)$}"  $g$, a two-place function such that for each $i,g(i,-)$
is an isomorphism from $M_0$ onto $M_i$
\sn
\item "{$(e)$}"  $<$ (two-place relation) - the usual ordering (on the
ordinals $< \lambda^+$) 
\sn
\item "{$(f)$}"  relations with parameter $i$ witnessing $M_i
\le_{\frak K} \dbcu_{j < \lambda^+} M_j$ (we can instead make
functions witnessing $M \in K$ as in \scite{88r-1.8} (the strong version) and have:
each $M_i$ is closed under them))
\sn
\item "{$(g)$}"  relations with parameter $i$ witnessing each
$\bold F(M_{i+1}) \le_{\frak K} M_{i+2}$ and $M_{i+1} \ne M_{i+2}$
(including $(M_{i+1},\bold F(M_{i+1})) \in \bold F$)
\sn
\item "{$(h)$}"  if $\mu = \lambda$, also individual constant for each
$a \in M_0$.
\ermn
Let $\psi \in \Bbb L_{\mu^+,\omega}$ describe this, in particular for
clauses (f), (g) use (C) of the assumptions.  So $\psi$ has a
non-well ordered model ${\frak B}^*,|P^{{\frak B}^*}| = \lambda$ (by
clause (D) of the assumption see \scite{88r-3.7A}(2)+(3)).  So
let

$$
{\frak B}^* \models ``a_{n+1} < a_n" \text{ for } n < \omega.
$$
\mn
Let for $a \in {\frak B}^*,A_a = \{x \in {\frak B}^*:{\frak B}^* \models
xRa\}$

$$
M_a = ({\frak B}^* \restriction \tau_{\frak K}) \restriction A_a.
$$
\mn
Easily $M_a \le_{\frak K} {\frak B}^* \restriction \tau_{\frak K}$ 
(use clause (f)) and $\|M_a\| =
\lambda$.  In fact $M_a$ is superlimit or just isomorphic to $M^*$
if $\mu = \lambda$, as $\psi$ includes the diagram of $M_0 = M^*$,
having names for all members, and if 
$\mu < \lambda$ see assumption (E).
So $M_{a_n} \le_{\frak K} {\frak B}^* \restriction \tau_{\frak K},
M_{a_{n+1}} \subseteq M_{a_n}$ hence $M_{a_{n+1}}
\le_{\frak K} M_{a_n}$ by Ax V.  Let $M_n =: M_{a_n}$.  Let $I = \{b \in
{\frak B}^*:\dsize \bigwedge_{n < \omega} [{\frak B}^* \models b < a_n]\}$.

Also as for $b \in I,M_b <_{\frak K} {\frak B}^* \restriction \tau_{\frak K}$ 
and $M_{b_1} <_{\frak K} M_{b_2}$ for $b_1 <^{{\frak B}^*} b_2$, by Ax IV 
clearly $M_\omega =: ({\frak B}^* \restriction 
(\tau_{\frak K})) \restriction \dbcu_{b\in I} A_b$ satisfies
$M_\omega \le_{\frak K} {\frak B}^* \restriction \tau_{\frak K}$ 
hence $M_\omega \le_{\frak K} M_n$ for
$n < \omega$.  Obviously $M_\omega \subseteq \dbca_{n < \omega} M_n$
and equality holds as $\psi$ guarantee
\mr
\item "{$(*)$}"  for every $y \in {\frak B}^*$ there is a minimal $x
\in {\frak B}^*$ such that $y \in M_x$.
\ermn
As each $M_b$ is isomorphic to $M^*$, of cardinality $\lambda$, also $M_\omega$
is.
\bn
\ub{Stage e}:  Suppose that there is a maximal triple, then we shall show
$\dot I(\lambda^+,K) = 2^{\lambda^+}$ and moreover
$\dot I(\lambda^+,K^{\bold F}_{\lambda^+}) = 
2^{\lambda^+}$, and so we shall get a
contradiction to assumption (B).

So there is a maximal triple $(M^0,N^0,a)$.  Hence by the uniqueness
of the limit model
for each $M \in K_\lambda$ which is isomorphic to $M^*$ hence to
$M^0$ there are $N,a$ satisfying $M \le_{\frak K} N \cong M^* \in
K_\lambda,a \in N \backslash M$ such that: if
$M <_{\frak K} M' \le_{\frak K} N' \in {\frak K}_\lambda,N
<_{\frak K} N',(\exists M'')(M \le_{\frak K} M'' \le_{\frak K} \bold F(M'')
\le_{\frak K} M' \cong M^*)$ and 
$(\exists N'')(N \le_{\frak K} N'' \le_{\frak K} \bold F(N'')
\le_{\frak K} N' \cong M^*)$ then $a \in M'$.  (That
is, in some sense $a$ is algebraic over $M$).
We can waive $(\exists N'')(N \le_{\frak K} N'' \le_{\frak K} \bold F(N'')
\le_{\frak K} N' \cong M^*)$ as by the definition of strongly limit there
is $N'_* \cong M^*$ such that $\bold F(N') \le_{\frak K} N'_*$.
On the other hand by Stage d
\mr
\item "{$(*)_1$}"  for each $M \in K_\lambda$ isomorphic to $M^*$
there are $M'_n(n < \omega)$ such that $M \le_{\frak K} M'_{n+1} <_{\frak K} 
M'_n \in K_\lambda,M'_n \cong M^*$ and $\bold F(M'_{n+1}) \le_{\frak
K} M'_n$ and $\dbca_{n < \omega} M'_n = M$.
\ermn
For notational simplicity: $M \in K_\lambda,|M|$ an ordinal $\Rightarrow
|\bold F(M)|$ an ordinal.

Now for each $S \subseteq \lambda^+$ we define by induction on $\alpha
\le \lambda^+,M^S_\alpha$, increasing (by $<_{\frak K}$) and continuous with
universe an ordinal $< \lambda^+$ such that $M^S_\alpha \cong M^*$ and
if $\beta + 2 \le \alpha$ then
$\bold F(M_{\beta +1}) \le_{\frak K} M_{\beta +1}$.  
Let $M^S_0 = M^*$ and for limit
$\delta < \lambda^+$ and let $M^S_\delta = \dbcu_{\alpha < \delta} M^S_\alpha$;
by the induction assumption and the choice of $M^*,\bold F$ clearly
$M^S_\delta$ is isomorphic to
$M^*$.  For $\alpha = \beta +1,\beta$ successor let $M^S_\alpha$ be
such that $\bold F(M^S_\beta) <_{\frak K} M^S_\alpha \cong M^*$.  
So we are left with the case $\alpha = \delta+1,\delta$ limit or zero.

Now if $\delta \in S$ hence $M^S_\delta \cong M^*$, choose 
$M_{\delta +1},a^S_\delta$ such that
$(M^S_{\delta +1},M^S_\delta,a^S_\delta)$ is a maximal triple (possible
as by the hypothesis of this case there is a maximal triple, and there
is a unique strong limit model).  
If $\delta \notin S$ we choose $M^{S,n}_\delta
\in K_\lambda$ for $n < \omega$ (not used) such that $M^S_\delta 
<_{\frak K} M^{S,n+1}_\delta \le_{\frak K} M^{S,n}_\delta$ and $\bold
F(M^{S,n+1}_\delta) \le_{\frak K} M^{S,n}_\delta$  for $n <
\omega$ and $M^S_\delta = \dbca_{n < \omega} M^{S,n}_\delta$ and
$M^{S,n}_\delta \cong M^*$; and let
$M^S_{\delta +1} = M^{S,0}_\delta$ (again possible as $M_\delta
\cong M^*$ and an $(*)_1$ above).
\nl
Lastly, let $M^S = \dbcu_\alpha M^S_\alpha$.

Now clearly it suffices to prove that if $S^0,S^1 \subseteq
\lambda^+,S^1 \backslash S^0$, is stationary then $M^{S^1} 
\ncong M^{S^0}$.  So suppose $f$ is a $\le_{\frak K}$-embedding from 
$M^{S^1}$ onto $M^{S^0}$ or just into $M^{S^0}$.  Then $E^2 = 
\{\delta < \lambda^+:M^{S^1}_\delta,M^{S^0}_\delta$ each has
universe $\delta$ and $i < \lambda^+$ implies $[i < \delta \Leftrightarrow
f(i) < \delta]\}$ is a closed unbounded subset of $\lambda^+$, hence
there is a limit ordinal 
$\delta \in (S^1 \backslash S^0) \cap E^2$.  Let us look at
$f(a^{S^1}_\delta)$; as $\delta \in S^1,a^{S^1}_\delta$ is well
defined, also $a^{S_1}_\delta \in
M^{S^1}_{\delta +1} \backslash M^{S^1}_\delta$, as $\delta \in E^2$ it
follows that $f(a^{S^1}_\delta) \nless \delta$ hence
$f(a^{S^1}_\delta)$ belongs to
$M^{S^0} \backslash M^{S^0}_\delta$ but $M^{S^0}_\delta = \dbca_{n < \omega}
M^{S^0,n}_\delta$ (as $\delta \notin S^0$). \nl
Hence for some $n,f(a^{S^1}_\delta) \notin M^{S^0,n}_\delta$.  
Let $\beta \in (\delta,\lambda^+)$ be large enough such that
$f(M^{S^1}_{\delta +1}) \subseteq M^{S^0}_\beta$.  But then 
$f(M^{S^1}_\delta) \le_{\frak K} M^{S^0,n}_\delta 
\le_{\frak K} M^{S^0}_\beta$ and $f(M^{S^1}_{\delta +1}) \le_{\frak K}
M^{S^0}_\beta$ and $a^{S^1}_\delta \notin f^{-1}(M^{S^0,n}_\delta)$.
Now $(f(M^{S^1}_\delta)),f(M^{S^1}_{\delta +1}),f(a^{S^1}_\delta))$
has the same properties as $(M^{S^1}_\delta,M^{S^1}_{\delta
+1},a^{S^1}_\delta)$ because if $f$ is an isomorphism from $M'$ onto
$M'' \in K_\lambda$ then we can extend $f$ to an isomorphism from
$\bold F(M')$ onto $\bold F(M'')$ (i.e., the ``invariant").  But
$(f(M^{S^1}_\delta),f(M^{S^1}_{\delta +1}),f(a^{S^1}_\delta)) < 
(M^{S^0,n}_\delta,M^{S^0}_\beta,f(a^{S^1}_\delta))$, contradiction.
So we are done.   \hfill$\square_{\scite{88r-3.8}}$ 
\bigskip

\demo{\stag{88r-3.9} Conclusion}  1) If LS$({\frak K}) = \aleph_0,K$ is
PC$_{\aleph_0}$ and $\dot I(\aleph_1,K) =1$, \ub{then} $K$ has a model of
cardinality $\aleph_2$. \nl
2) If $\psi \in \Bbb L_{\omega_1,\omega}(\bold Q)$ ($\bold Q$ is the
quantifier ``there are uncountably many") has one and only one model
of cardinality $\aleph_1$ up to isomorphism \ub{then} $\psi$ has a model in
$\aleph_2$.
\enddemo
\bigskip

\demo{Proof}  1) By \scite{88r-3.6} we get suitable ${\frak K}_1$ (as in its
conclusion) and by \scite{88r-3.7} the class
${\frak K}_1$ has a model in $\aleph_2$, hence
${\frak K}$ has a model in $\aleph_2$. 
\nl
2) We can replace $\psi$ by a countable theory $T \subseteq
\Bbb L_{\omega_1,\omega}(\bold Q)$. \nl
Let $L$ be a fragment of $\Bbb L_{\omega_1,\omega}(\bold Q)(\tau)$ in
which $T$ is included (e.g., $L$ is the closure of $T \cup$(the atomic
formulas) under subformulas, $\neg,\wedge,(\exists x),(\bold Q x)$; in
particular $L$ includes, of course, first order logic).  
By \cite{Sh:48}, without loss of generality $T$ ``says" that
every formula $\varphi(x_0,\dotsc,x_{n-1})$ of $L$ is
equivalent to an atomic formula (i.e., $P(x_0,\dotsc,x_{n-1}),P$ a predicate)
and every type realized in model of $T$ is isolated (i.e., every model is
atomic), and $T$ is complete in $L$.
Let

$$
\align
K = \{M:&M \text{ an atomic } \tau(T) \text{-model of } T \cap \Bbb L
\text{ and if } M \models P[\bar a] \\
  &\text{ and } (\forall \bar x)[P(\bar x) \equiv \neg(\bold Q
  y)R(y,\bar x)] \in T \\
  &\text{ then } \{b:M \models R[b,\bar a]\} \text{ is countable}\}
\endalign
$$

$$
\align
M \le_{\frak K} N \text{ iff } M \le^* N \text{ which means}: &(a)
\quad M \prec_{\Bbb L} N \\
  &(b) \quad \text{if } M \models P(\bar a) \text{ and } \forall \bar
x[P(\bar x) \equiv \neg \bold Q y R(y,\bar x)] \in T \\
  &\quad \qquad \text{ \ub{then} for no } b \in N \backslash M \text{
do we have } N \models R[b,\bar a].
\endalign 
$$
\mn
So ${\frak K} = (K,\le_{\frak K})$ is categorical in 
$\aleph_0$, is an a.e.c. and is
PC$_{\aleph_0}$.  Let $\bold F$ be (see \scite{88r-3.1}(8))
such that for $M \in K_{\aleph_0},N = \bold F(M)$ iff: 
$M <^{**} N$ which says $M \le_{\frak K} N \in K_{\aleph_0}$ and
if $\bar a \in M,M
\models P[\bar a],\forall \bar x[P(\bar x) \equiv {\bold Q} 
yR(y,\bar x)] \in T$, \ub{then} for some $b \in N \backslash M$ we have
$N \models R[b,\bar a]$.  So $\bold F$ is invariant.

Note that every $M \in K^{\bold F}_{\aleph_1}$ is a model of $\psi$.  So
\scite{88r-3.8} gives that some $M \in K^{\bold F}_{\aleph_1}$ has a proper
extension in $K^{\bold F}_{\aleph_1}$. \nl
The rest should be easy.  \hfill$\square_{\scite{88r-3.9}}$
\enddemo
\bn
\margintag{88r-3.9A}\ub{\stag{88r-3.9A} Question}  1) Under the assumptions of \scite{88r-3.9}(2), can
we get $M \in K_{\aleph_2}$, such that: if $M \models P[\bar a],
\forall \bar x[P(\bar x) \equiv (\bold Q y)R(y,\bar x)] \in T$ 
then $\{b \in M:M \models R[b,\bar a]\}$ has cardinality $\aleph_2$?
Note that in the proof of \scite{88r-3.8} we show that no triple is maximal.
\bigskip

\remark{\stag{88r-3.9.1} Remark}  1)  
We could have used multi-valued $\bold F$ then in the proof above
$N = \bold F(M)$ just means the demand there.
\nl
2) To answer \scite{88r-3.9A}, i.e., to prove the existence of $M \in
K_{\aleph_2}$ as above we have to prove:
\mr
\item "{$(*)_1$}"   there are $N,N_i \in K^{\bold
F}_{\aleph_1}$ for $i < \omega_1$ and $N \le_{\frak K} N_i$ such
that if $N \models P[\bar a]$ and the sentence $(\forall \bar
x)(P(\bar x) \equiv (\bold Q y)R(y,\bar x)]$ belongs to $T$, \ub{then}
for some $i < \omega_1$ there is $b_* \in N_i \backslash N$ such that
$N_i \models R[b,\bar a]$.
\nl
Clearly
\sn
\item "{$(*)_2$}"  the existence of $N,N_i$ as in $(*)_1$ is
equivalent to ``$\psi^*$ has a model" for some $\psi^* \in \Bbb
L_{\omega_1,\omega}(\bold Q)$ which is defined from $T,\le_{\frak K}$.
\ermn
Hence
\mr
\item "{$(*)_3$}"  it is enough to prove that for some forcing notion
$\Bbb P$ in $\bold V^{\Bbb P}$ there are $N,N_i$ as in $(*)_1$.
\ermn
There are some natural c.c.c. forcing notions tailor-made for this
\mr
\item "{$(*)_4$}"  consider the class of triples $(M,N,a)$ such that $M
\le_{\frak K} N \in K_{\aleph_0},\bar a \in {}^{\omega >} N,\ell <
\ell g(\bar a) \Rightarrow a_\ell \notin M$, order as in the proof of
\scite{88r-3.8}.  By the same proof there is no maximal triple.
\ermn
3) We can restrict ourselves in $(*)_2$ to

$$
\{R(y,\bar a):\bar a \in {}^{\ell g(\bar x)}N \text{ and } \bar a
\text{ realizes a type } p(\bar x)\}.
$$
\mn
Also we may demand $i < \omega_1 \Rightarrow N_i = N_0$ and we may try
to force such a sequence of models (or pairs) and there is a natural
forcing.  By absoluteness it is enough to prove that it satisfies the c.c.c.

\endremark
\bn
\margintag{88r-3.9B}\ub{\stag{88r-3.9B} Problem}:  If ${\frak K}$ is 
PC$_\lambda,K$ categorical in $\lambda$ and
$\lambda^+$, does it necessarily have a model in $\lambda^{++}$?
\bigskip

\remark{Remark}  The problem is proving $(*)$ of \scite{88r-3.7A}.
\endremark
\bn
\margintag{88r-3.9C}\ub{\stag{88r-3.9C} Question}:  Assume 
$\psi \in \Bbb L_{\omega_1,\omega}(\bold Q)(\tau)$ is
complete in $\Bbb L_{\omega_1,\omega}(\bold Q)(\tau)$, is categorical in
$\aleph_1$, has an uncountable model $M,\bar a \in {}^n M$
and $\varphi \in \Bbb L_{\omega_1,\omega}(\bold Q)(\tau)$ axiomatizes
the $\Bbb L_{\omega_1,\omega}(\bold Q)(\tau)$-theory of 
$(M,\bar a)$.  Is $\varphi$ categorical in $\aleph_1$?   
\bn
\margintag{88r-3.9D}\ub{\stag{88r-3.9D} Question}:  Can we weaken the demand on $M^*$ in
\scite{88r-3.8} to ``$M^*$ is a $\lambda^+$-limit model"?
\newpage

\head {\S4 Forcing and categoricity} \endhead  \resetall \sectno=4
 \spuriousreset
\bn

The main aim in this section is, for ${\frak K}$ as in \S1 with
LS$({\frak K}) = \aleph_0$, to find
what we can deduce from $1 \le \dot I(\aleph_1,K) < 2^{\aleph_1}$, first \
without assuming $2^{\aleph_0} < 2^{\aleph_1}$.

We can build a model of cardinality $\aleph_1$ by an $\omega_1$-sequence of
countable approximations.  Among those, there are models which are the
union of a quite generic $<_{\frak K}$-increasing 
sequence $\langle N_i:i < \omega_1 \rangle$
of countable models, so it is natural to look at them
(e.g. if ${\frak K}$ is categorical in $\aleph_1$, every model in
$K_{\aleph_1}$ is like that).  We say on such models that they are
quite generic.  More exactly, we look at countable
models and figure out properties of the quite generic models in
${\frak K}_{\aleph_1}$.  The main results are \scite{88r-4.8}(a),(f).
Note that $2^{\aleph_0} = 2^{\aleph_1}$, though in general 
making our work harder, can be utilized positively - see \scite{88r-4.7}.

A central notion is (e.g.) ``the type which $\bar a \in {}^{\omega
>}(N_1)$ materializes in $(N_1,N_0)$", $N_0 \le_{\frak K} N_1 \in
K_{\aleph_0}$.  This is as the name indicates,  the type materialized in
$N^+_1$, which is $N_1$ expanded by $P^{N^+_1} = N_0$; it consists of
the set of formulas forced (in the model theoretic sense started by
Robinson) to satisfy; here forced is defined thinking on
$(K_{\aleph_0},\le_{\aleph_0})$ so models in $K_{\aleph_1}$ can be
constructed as the union of quite generic $<_{\frak K}$-increasing
$\omega_1$-sequence.  As we would like to build models of cardinality
$\aleph_1$ by such sequence, the ``materialize" in $(N_1,N_0)$ becomes
realized in the (quite generic) 
$N \in K_{\aleph_1}$; but most of our work is in
$K_{\aleph_0}$.  This is also a way to express $\bold Q$ speaking on
countable models.

By the hypothesis \scite{88r-4.5} justified by \S3, the $\Bbb
L_{\infty,\omega}(\tau_{\frak K})$-theory of $M \in K$ is clear, in
particular has elimination of quantifiers hence $M \le_{\frak K} N
\Rightarrow M \prec_{\Bbb L_{\infty,\omega}} N$, but for $\bar N =
\langle N_\alpha:\alpha < \omega_1\rangle$ as above we would like to
understand $(N_\beta,N_\alpha)$ for $\alpha < \beta$ (from the point
of view of $N,\bar N$ is not reconstructible, but its behaviour on a
club is).  Toward a parallel analysis of such pairs we again analyze
them by $\langle L^0_\alpha:\alpha < \omega_1\rangle$ (similarly to
\cite{Mo70}). 
\bigskip

\demo{\stag{88r-4.0} Convention}  We fix $\lambda > \text{ LS}
({\frak K})$ as well as the a.e.c. ${\frak K}$.
\enddemo
\bn
The main case below is here $\lambda = \aleph_1,\kappa = \aleph_0$.
\definition{\stag{88r-4.1} Definition}  For $\lambda > \text{ LS}({\frak K})$ 
and $N_* \in K_{< \lambda}$ and $\mu,\kappa$ satisfying $\lambda \ge
\kappa \ge \aleph_0,\mu \ge \kappa$ and let
\nl
1) $\Bbb L^0_{\mu,\kappa}$ be first order logic enriched by conjunctions
(and disjunctions) of length $< \mu$, homogeneous strings of
existential quantifiers or of universal quantifiers of length $<
\kappa$, and the cardinality quantifier $\bold Q$ interpreted as
$\exists^{\ge \lambda}$.  But we apply those operations such that
any formula has $< \kappa$ free variables, and the non-logical symbols
are from $\tau({\frak K})$ so actually we should write $\Bbb
L^0_{\mu,\kappa}(\tau_{\frak K})$ but we may ``forget" to say this 
when clear; the syntax does not depend on $\lambda$ but we
shall mention it in the definition of satisfaction. \nl
2) For a logic ${\Cal L}$ and $A_i,A \subseteq N_*$ for $i <
\alpha,\alpha < \lambda$ let 
${\Cal L}(N_*,A_i;A)_{i < \alpha}$ be the language,
with the logic ${\Cal L}$, and with the vocabulary $\tau_{N_*,\bar
A,A}$ where $\bar A = \langle A_i:i < \alpha
\rangle$ and $\tau_{N_*,\bar A;A}$ consists of $\tau(K)$, the
predicates $x \in N_*$ and $x \in A_i$ for $i < \alpha$ and the
individual constants $c$ for $c \in A$.  (If $A = \emptyset$, we may
omit the $A$; if we omit $N_*$ then ``$x \in N_*$" is omitted, if the
sequence of the $A_i$ is omitted then the ``$x \in A_i"$ are omitted,
so ${\Cal L}()$ means having the vocabulary $\tau(K))$. 
So ${\Cal L}(N_*,A_i,A)_{i < \alpha}$ formally should have been
written ${\Cal L}(\tau_{N_*,\bar A;A})$.
\nl
3) $\Bbb L^1_{\mu,\kappa}$ is defined is as in part (1), but we have
also variables (and quantification) over relations of cardinality $<
\lambda$.  Let $\Bbb L^{-1}_{\mu,\kappa}$ be as in part (1) but not
allowing the cardinality quantifier $\bold Q$; this is the classical
logic $\Bbb L_{\mu,\kappa}$. 
\nl
4) $(N,N_*,A_i;A)_{i < \alpha}$ is the model $N$ expanded to a
$\tau_{N_*,\bar A;A}$-model by monadic
predicates for $N_*,A_i(i < \alpha)$ and individual constants for
every $c \in A$.  \nl
5) For $``x \in N_*",``x \in A_i"$ we use the predicates $P,P_i$ respectively, 
so we may write ${\Cal L}(\tau + P)$ instead ${\Cal L}(N_*)$,
but writing ${\Cal L}(N_*)$ we fix the interpretation of $P$. \nl
Let $\tau^{+ \alpha} = \tau \cup \{P,P_\beta:\beta < \alpha\}$ and if
$L = {\Cal L}(\tau^{+0})$, i.e., for $\alpha = 0$ then 
$L(N)$ means $L$ but we fix the interpretation
of $P$ as $N$, i.e., $|N|$, the set of elements of $N$. \nl
Let $L(N_*,N_i)_{i \in u}$ where $u$ a set of $< \kappa$ ordinals
means the language $L$ when we fix the interpretation of 
$P$ as $N_*$ and of $P_{\text{otp}(u \cap \alpha)}$ as $N_\alpha$.
\enddefinition
\bigskip

\definition{\stag{88r-4.2} Definition}  1) For $N_* \in K_{< \lambda}$ and
$\varphi(x_0,\ldots) \in \Bbb L^1_{\mu,\kappa}(N_*,\bar A;A)$ we define by
induction on $\varphi$ when $N_0 \Vdash^\lambda_{\frak K}
\varphi[a_0,\ldots]$ (where $N_* \le_{\frak K} N_0  \in K_{<
\lambda},a_0,\ldots$ are elements of $N_0$ or appropriate 
relations over it, depending on the kind of $x_i$).
Pedentically we should write $(N_0;N_*,\bar A;A) \Vdash^\lambda_{\frak
K} \varphi[a_0,\ldots]$ (and we may do it when not clear from the content).

For $\varphi$ atomic this means $N_0 \models \varphi[a_0,\ldots]$.
For $\varphi = \dsize \bigwedge_i \varphi_i$ this means

$$
N_0 \Vdash^\lambda_{\frak K} \varphi_i[a_0,\ldots] \text{ for each } i.
$$
\mn
For $\varphi = \exists \bar x \psi(\bar x,a_0,\ldots)$ this means that for
every $N_1$ satisfying $N_0 \le_{\frak K} N_1 \in K_{< \lambda}$ 
there is $N_2,N_1 \le_{\frak K} N_2 \in K_{< \lambda}$ 
and $\bar b$ from $N_2$ of the appropriate length (and kind) such
that $N_2 \Vdash^\lambda_{\frak K} \psi[\bar b,a]$.

For $\varphi = \neg \psi$ this means that for no $N_1$ do we have
$N_0 \le_{\frak K} N_1 \in K_{< \lambda}$ and $N_1
\Vdash^\lambda_{\frak K} \psi [a_0,\ldots]$.

For $\varphi(x_0,\ldots) = (\bold Q y)\psi(y,x_0,\ldots)$ this means for
every $N_1$ satisfying $N_0 \le_{\frak K} N_1 \in K_{< \lambda}$ there
is $N_2$ satisfying $N_0 \le_{\frak K} N_2 \in K_{< \lambda}$ and 
$a \in N_2 \backslash N_1$ such that 
$N_2 \Vdash^\lambda_{\frak K} \psi[a,a_0,\ldots]$. 
\nl
2) In part (1) if $\varphi \in \Bbb L^1_{\mu,\kappa}(N_*)$ we can omit
   the demand ``$N_* \le_{\frak K} N$" similarly below.
\nl
3) For a language $L \subseteq \Bbb L^1_{\mu,\kappa}(N_*,\bar
A;A)$ and a model $N$ satisfying $N_* \le_{\frak K} N \in
K_{< \lambda}$ and a sequence $\bar a \in {}^{\lambda >} N$
the $L$-generic type of $\bar a$ in $N$ is gtp$(\bar a;N_*,\bar A;A;N) =
\{\varphi(\bar x) \in L:N \Vdash^\lambda_{\frak K} \varphi[\bar a]\}$.   
\nl
4) Let gtp$^\lambda_L(\bar a;N_*,\bar A;A;N)$ where $N_* 
\le_{\frak K} N \in K_\lambda$ and $L \subseteq {\Cal L}(N_*,\bar A;A)$
be $\{\varphi(\bar x):\varphi \in {\Cal L}(N_*,\bar A;A)$ and for
some $N' \in K_{< \lambda}$ we have $N \le_{\frak K} N' \le_{\frak K}
N$ and $N' \Vdash^\lambda_{\frak K} \varphi[\bar a]\}$; we may 
omit $\bar A,A$ (and omit $\lambda$ if clear from the context) and may
write ${\Cal L}$ instead of $L = {\Cal L}(N_*,\bar A;A)$; but note
Definition \scite{88r-5.2.7}.
\nl
5) We say ``$\bar a$ materializes $p$ (or
$\varphi$)" if $p$ (or $\{\varphi\}$) is a subset of the
$L$-generic type of $\bar a$ in $N$. 
\enddefinition
\bigskip

\definition{\stag{88r-4.3} Definition}  Let $N_i(i < \lambda)$ be an
increasing (by $\le_{\frak K}$) continuous sequence, $N = \dbcu_{i < \lambda}
N_i,\|N_i\| < \lambda$ and $L^* \subseteq \dbcu_{\alpha < \kappa} \Bbb
L^1_{\infty,\kappa}(\tau^{+\alpha})$.
\nl
1) $N$ is $L^*$-generic, if for any formula
$\varphi(x_0,\ldots) \in L^* \cap \Bbb L^1_{\infty,\kappa}(\tau_{\frak
K})$ and $a_0,\ldots \in N$ we have: \nl
$N \models \varphi[a_0,\ldots] \Leftrightarrow$ for some $\alpha <
\lambda,N_\alpha \Vdash^\lambda_{\frak K} \varphi[a_0,\ldots]$. \nl
2) The $\le_{\frak K}$-presentation $\langle N_i:i < \lambda \rangle$ of $N$ is
$L^*$-generic \ub{if} for any $\alpha < \lambda$ of 
cofinality $\ge \kappa$ and $\psi(x_0,\ldots) \in L^*(N_\alpha,N_i)_{i
\in I}$ satisfying $I \subseteq \alpha,|I| < \kappa$
and $a_0,\ldots \in N$ we have:

$$
N \models \psi[a_0,\ldots] \Leftrightarrow \text{ for some } \gamma
< \lambda,N_\gamma \Vdash^\lambda_{\frak K} \psi[a_0,\ldots]
$$
\mn
and for each $\beta \ge \alpha$, with cofinality $\ge \kappa,N_\beta$
is almost $L^*(N_\alpha,N_i)_{i \in I}$-generic (see part (5)). \nl
3) $N$ is strongly $L^*$-generic \ub{if} it has an
$L^*$-generic presentation (in this case, if
$\lambda$ is regular, then for any presentation $\langle N_i:i <
\lambda \rangle$ of $N$ there is a closed unbounded $E \subseteq
\lambda$ such that $\langle N_i:i \in E \rangle$ is an
$L^*$-generic presentation). \nl
4) We say that $N \in K_{< \lambda}$ is pseudo $L^*$-generic \ub{if}
\mr
\item "{$(a)$}"  for every $\varphi(\bar x) = \exists \bar y \psi(\bar
x,\bar y) \in L^*$, if $N \Vdash^\lambda_{\frak K} \varphi(\bar a)$ 
then for some $\bar b,N \Vdash^\lambda_{\frak K} \psi(\bar a,\bar b)$
\sn 
\item "{$(b)$}"  for every $\bar a \in N,\bar a$ materializes in $N$ some
complete $L^*$-type.
\ermn
5) We add ``almost" to all the above defined notions if for
$\Vdash^\lambda_{\frak K}$,
the inductive definitions of satisfaction works except possibly for
$\bold Q$ (e.g., $N \Vdash^\lambda_{\frak K} \exists x 
\varphi(x,\ldots)$ iff for
some $a \in N,N \Vdash^\lambda_{\frak K} \varphi(a,\ldots))$. 
\enddefinition
\bigskip

\remark{\stag{88r-4.3A} Remark}  1) Notice we can choose $N_i = N_0 = N$,
so $\|N\| < \lambda$.  In particular almost (and pseudo) 
$L^*$-generic models of
cardinality $< \lambda$ may well exist. \nl
2) Here we concentrate on $\lambda = \aleph_1$ and fragments of
$\Bbb L^0_{\infty,\omega}$ (mainly $\Bbb L^0_{\omega_1,\omega}$ and its
countable fragments). \nl
3) There are obvious implications, and forcing is preserved by
isomorphism and replacing $N(\in K_{< \lambda})$ by $N',N \le_{\frak
K} N' \in K_{< \lambda}$.
\endremark
\bn
There are obvious theorems on the existence of generic models, e.g.,
\proclaim{\stag{88r-4.4} Theorem}  1) Assume $N_0 \in K_{< \lambda},\lambda =
\mu^+,\mu^{< \kappa} = \mu,L \subseteq \dbcu_{\alpha < \kappa} \Bbb
L_{\infty,\kappa}(\tau^{+ \alpha})$ and $L$ is 
closed under subformulas and $|L| < \lambda$.
\ub{Then} there are $N_i(i < \lambda)$ such that $\langle N_i:i <
\lambda \rangle$ is an $L$-generic representation of $N = 
\dbcu_{i < \lambda} N_i$, (hence $N$ is strongly $L$-generic). \nl
2)  In part (1), $N \in K_\lambda$ if no $N',N_0 \le_{\frak K} N' \in 
K_{< \lambda}$ is $\le_{\frak K}$-maximal.
\endproclaim
\bigskip

\demo{Proof}  Straightforward.
\enddemo
\bigskip

\remark{\stag{88r-4.4A} Remark}  1) If $L = \dbcu_{i < \lambda} L_i,
|L_i| < \lambda$, then we can get ``$\langle N_i:j < i 
< \lambda \rangle$ is an $L_j$-generic representation of $N$ for
each $j < \lambda$". \nl
2) When we speak on ``complete $L$-type $p$" we mean $p = 
p(x_0,\dotsc,x_{n-1})$ for some $n$. 
\endremark
\bn
\relax From time to time we add some hypothesis and prove a series of claims;
such that the hypothesis holds, at least \wilog, \, in the case we are
interested in.  We are mainly interested in the case 
$\dot I(\aleph_1,{\frak K}) < 2^{\aleph_1}$, etc., so 
by \scite{88r-3.6}, \scite{88r-3.8.5} it is reasonable to make:
\demo{\stag{88r-4.5} Hypothesis}  ${\frak K}$ is PC$_{\aleph_0},
\le_{\frak K}$ refines 
$\Bbb L_{\infty,\omega}$ and ${\frak K}$ is categorical in $\aleph_0$ and $1
\le \dot I(\aleph_1,K)$ and $\dot I(\aleph_1,K^{\bold
F}_{\aleph_1}) < 2^{\aleph_1}$ where $K^{\bold F}_{\aleph_1}$ is
as in \scite{88r-3.8.1} and is PC$_{\aleph_0}$ or just $\bold K^{\bold
F}_{\aleph_1} = \{M \restriction
\tau_{\frak K}:M \models \psi\}$ for some $\psi \in \Bbb
L_{\omega_1,\omega}(\bold Q)$ (if $\bold F$ is invariant, this follows).
\enddemo
\bigskip

\remark{\stag{88r-4.5.3} Remark}  1) Usually below we ignore the case 
$\dot I(\aleph_1,{\frak K}) < 2^{\aleph_0}$ as the proof is the same.
\nl
2) We can deal similarly with the case $1 \le \dot I(\aleph_1,K') 
< 2^{\aleph_0}$
where ${\frak K}_{\aleph_1} \subseteq K'_{\aleph_1} \subseteq \{M \in
{\frak K}_{\aleph_1}:M$ is strongly $L_*$-generic$\}$ and $K'$ is
PC$_{\aleph_0}$ (or less: $\{M \restriction \tau_{\frak K}:M$ a model
of $\psi \in \Bbb L_{\omega_1,\omega}(\bold Q)(\tau^*)\}$).
\nl
3) Can we use $\bold F$ a function with domain $K_{\aleph_0}$ such
that $M \le_{\frak K} \bold F(M_0) \in K_{\aleph_0}$ for $M \in
K_{\aleph_0}$ without the extra assumptions or even
$\bold F:\{\bar M = \langle M_i:i \le \alpha \rangle$ is $\le_{{\frak
K}_{\aleph_0}}$-increasing continous$\} \rightarrow {\frak
K}_{\aleph_0}$ such that $M_\alpha \le_{\frak K} \bold F(M_i:i \le 
\alpha \rangle)$?  We cannot use the non-definability of well ordering
(see \scite{88r-3.6}(3)); (as in the proof of (f) of \scite{88r-4.8}).
\endremark
\bigskip

\proclaim{\stag{88r-4.6} Claim}  1) If $\bar a \in N \in K_{\aleph_0}$
and $\varphi(\bar x) \in \Bbb L^0_{\infty,\omega}(\tau^{+0})$ 
(so $\bar a$ is a
finite sequence) \ub{then} $(N,N) \Vdash^{\aleph_1}_{\frak K} 
\varphi[\bar a]$ or $(N,N) \Vdash^{\aleph_1}_{\frak K} \neg \varphi[\bar a]$ 
(i.e. $P$ is interpreted as $N$). \nl
2) If $(N,N) \Vdash^{\aleph_1}_{\frak K} \exists \bar x \wedge p(\bar
x)$,  where $p(\bar x)$ is a not necessarily complete $n$-type
$(n = \ell g(\bar x))$ in $L$ where $L \subseteq \Bbb
L^0_{\omega_1,\omega}(\tau^{+0})$ is countable, \ub{then} for some
complete $n$-type $q$ in $L$ extending $p$ we have $(N,N)
\Vdash^{\aleph_1}_{\frak K} \exists \bar x \wedge q(\bar x)$. 
\endproclaim
\bigskip

\demo{Proof}  1) Suppose not, for each $S \subseteq \omega_1$, we
define by induction on $\alpha,N^S_\alpha \in K_{\aleph_0}(\alpha <
\omega_1)$, increasing (by $\le_{\frak K}$) and continuous, $N^S_0 = N$ and
for limit $\alpha,N^S_\alpha = \dbcu_{\beta < \alpha} N^S_\beta$.  For
$\alpha = 2 \beta +1$ remember that $(N^S_\beta,\bar a) \cong
(N,\bar a)$ because $N = N_0 \le_{\frak K} N^S_\beta$ hence $N_0
\prec_{{\Bbb L}_{\infty,\omega}} N^S_\beta \in K_{\aleph_0}$ hence
$(N^S_\beta,\bar a) \equiv_{\Bbb L_{\infty,\omega}} (N,\bar a)$ hence they
are isomorphic.  So $(N^S_\beta,N^S_\beta)$
forces $(\Vdash^{\aleph_1}_{\frak K})$ neither $\varphi[\bar a]$ nor 
$\neg \varphi[\bar a]$.
So there are $M_\ell$ (for $\ell = 0,1$) such that $N^S_\beta
\le_{\frak K} M_\ell \in K_{\aleph_0}$ and 
$(M_0,N^S_\beta) \Vdash^{\aleph_1}_{\frak K} \varphi[\bar a]$ 
but $(M_1,N^S_\beta) \Vdash^{\aleph_1}_{\frak K} \neg \varphi[\bar
a]$.  Now if $\beta \in S$ we let 
$N^S_\alpha = M_0$ and if $\beta \notin S$ we let $N^S_\alpha = M_1$.  

Lastly, $M_{2 \beta +2} = \bold F(M_{2 \beta+1})$.  
Let $N^S = \dbcu_{\alpha < \omega_1}
N^S_\alpha$.  Now if $S(0) \backslash S(1)$ is stationary then 
$(N^{S(0)},\bar a) \ncong (N^{S(1)},\bar a)$,  Why?  Because if 
$f:N^{S(0)} \rightarrow N^{S(1)}$ is
an isomorphism from $N^{S(0)}$ onto $N^{S(1)}$ mapping $\bar a$ to
$\bar a$ then for some closed unbounded set $E \subseteq \omega_1$, we have:
if $\alpha \in E$ then $f$ maps $N^{S(0)}_\alpha$ onto 
$N^{S(1)}_\alpha$, so choose some $\alpha \in E \cap S(0) \backslash
S(1)$ and choose $\beta \in E \backslash (\alpha +1)$.  Now
$(N^{S(0)}_{\alpha +1},N^{S(0)}_\alpha) \Vdash^{\aleph_1}_{\frak K}
\varphi[\bar a]$, hence
$(N^{S(0)}_\beta,N^{S(0)}_\alpha) \Vdash^{\aleph_1}_{\frak K} 
\varphi[\bar a]$, and
similarly $(N^{S(1)}_\beta,N^{S(1)}_\alpha) \Vdash^{\aleph_1}_{\frak
K} \neg \varphi(\bar a)$, 
but $f \restriction N^{S(0)}_\beta$ is an isomorphism from
$N^{S(0)}_\beta$ onto $N^{S(1)}_\beta$ mapping $N^{S(0)}_\alpha$ onto
$N^{S(1)}_\alpha$ and $\bar a$ to itself and we get a contradiction.  By
\scite{88r-0.9}, we get $\dot I(\aleph_1,K) = 2^{\aleph_1}$,
contradiction.  \nl
2) Easy by \scite{88r-4.4} and part (1). 
In detail, if $N \le_{\frak K} M_1 \in {\frak K}_{\aleph_0}$ then by
the definition of $\Vdash^{\aleph_1}_{\frak K}$ and the assumption we
can find $(M_2,\bar a)$ satisfying $M_1 \le_{\frak K} M_2 \in {\frak
K}_{\aleph_0}$ and $\bar a \in M_2$ such that $(M_2,N) \Vdash^{\aleph_1}_{\frak
K} \wedge p(\bar a)$.  As $L$ is countable and the definition of
$\Vdash^{\aleph_1}_{\frak K}$ \wilog \, for every formula
$\varphi(\bar x) \in L,(M_2,N) \Vdash^{\aleph_1}_{\frak K} \varphi[\bar
a]$ or $(M_2,N) \Vdash^{\aleph_1}_{\frak K} \neg \varphi[\bar a]$.
(Why?  Simply let $\langle \varphi_n(\bar x):n < \omega \rangle$ list
the formulas $\varphi(\bar x) \in L$ and choose $M_{2,n} \in {\frak
K}_{\aleph_0}$ by induction on $n$ such that $M_{2,0} = M_2,M_{2,n}
\le_{\frak K} M_{2,n+1}$ such that $(M_{2,n+1},N) \Vdash^{\aleph_1}_{\frak K}
\varphi_n(\bar x)$ or $(M_{2,n+1},N) \Vdash^{\aleph_1}_{\frak K} \neg
\varphi_n(\bar x)$; now replace $M_2$ by $\cup\{M_{2,n}:n <
\omega\})$.  Let $q = \text{ gtp}_{L(N)}(\bar a,N,M_2)$, 
it is a complete $(L(N),n)$-type.  So clearly
$(M_2,N) \Vdash^{\aleph_1}_{\frak K} (\exists \bar x) \wedge q(\bar x)$.   Now
apply part (1) to the formula $(\exists \bar x) \wedge 
q(\bar x)$ so we are done.  \hfill$\square_{\scite{88r-4.6}}$
\enddemo
\bigskip

\proclaim{\stag{88r-4.7} Claim}  For each countable 
$L \subseteq \Bbb L^0_{\omega_1,\omega}(\tau^{+0})$ and 
$N \in K_{\aleph_0}$ the number of
complete $L(N)$-types $p$ (with no parameters)
such that $N \Vdash^{\aleph_1}_{\frak K} (\exists \bar x) \wedge 
p(\bar x)$, is countable.
\endproclaim
\bigskip

\demo{Proof}  At first glance it seemed that \scite{88r-0.1} will
imply this trivially.  However, here we need the parameter $N$ as an
interpretation of the predicate $P$ and if $2^{\aleph_0} =
2^{\aleph_1}$ there are too many choices.  So we shall deal with
``every $N_\alpha$" in some presentation.  
Suppose the conclusion fails.  First we choose by
induction $N_\alpha$ (for $\alpha < \omega_1$) such that
\mr
\widestnumber\item{$(iii)$}
\item "{$(i)$}"  $N_\alpha \in K_{\aleph_0}$ is $\le_{\frak
K}$-increasing and $\langle N_\alpha:\alpha < \omega_1 \rangle$ is 
$L$-generic
\sn
\item "{$(ii)$}"  for each $\beta < \alpha$, there is $a^\beta_\alpha
\in N_{\alpha +1} \backslash N_\alpha$ materializing an 
$L(N_\beta)$-type not materialized in $N_\alpha$, (i.e. in
$(N_\alpha,N_\beta)$; see Definition \scite{88r-4.2}(2) on materialize)
(possible by \scite{88r-4.6} and our assumption toward contradiction)
\sn
\item "{$(iii)$}"  $|N_\alpha| = \omega \alpha$
\sn
\item "{$(iv)$}"  for $\alpha < \beta,N_\beta$ is pseudo
$L(N_\alpha)$-generic and $\bold F(N_{2 \beta +1}) \le_{\frak K} N_{2
\beta +2}$.
\ermn
Now let $N = \cup\{N_\alpha:\alpha < \omega_1\}$ and
we expand $N$ by all relevant information: the order $<$ on the countable
ordinals, $c(c \in N_0)$, enough ``set theory", ``witness" for 
$N_\beta \le_{\frak K} N_\alpha$ for 
$\beta < \alpha$ and $F,F(\beta,\alpha) = a^\beta_\alpha$ and
names for all formulas in $L(N_\alpha)$ (with
$\alpha$ as a parameter), i.e., the relations $R_{\varphi(\bar
x)} = \{\langle \alpha \rangle \char 94 \bar a:\alpha < \omega_1,\bar
a \in {}^{\ell g(x)} N$ and for every $\beta < \omega_1$ large enough
$(N_\beta,N_\alpha) \Vdash^{\aleph_1}_{\frak K}
``\varphi(\bar a)"\}$ for $\varphi(\bar x) \in L$.  
We get a model ${\frak B}$.  By \scite{88r-0.1}(1) applied to the case
$\Delta = L$, there are models ${\frak B}_i$ (for $i < 2^{\aleph_1})$ of
cardinality $\aleph_1$ (note $N_0 \le_{\frak K} {\frak B} \restriction
\tau_{\frak K}$), so that the set of $L(N_0)$-types
realizes in $N^i$ (the $\tau(K)$-reduct of ${\frak B}_i$) are distinct for
distinct $i$'s.  So $(N^i,c)_{c \in N_0}$ are pairwise non-isomorphic.
If $2^{\aleph_0} < 2^{\aleph_1}$ we finish by \scite{88r-0.9}.

So we can assume $2^{\aleph_0} = 2^{\aleph_1}$.
In $N$, uncountably many complete $L(N_0)$-types are realized hence by
\scite{88r-0.1}(2) the set 
$\{p:p$ a complete $L(N_0)-n$-type is realized in some
$N',N_0 \le_{\frak K} N' \in {\frak K}_{\aleph_1}\}$ has cardinality
continuum, hence by \scite{88r-4.6} the set of complete 
$L(N_0)$-types $p=p(x)$ such that $(N_0,N_0) 
\Vdash^{\aleph_1}_{\frak K} \exists \bar x \wedge p(\bar x)$ has cardinality
$2^{\aleph_0}$.  So we choose by induction on $\alpha < 2^{\aleph_0}$
a sequence $\langle N^\alpha_i,a^\alpha_i:i < \omega_1 \rangle$ such
that:
\mr
\item "{$(a)$}"   $N^\alpha_i \in {\frak K}_{\aleph_0}$
\sn
\item "{$(b)$}"  $N^\alpha_{i_0} \le_{\frak K} N^\alpha_i$ for $i_0 < i <
\omega_1$
\sn
\item "{$(c)$}"  $a^\alpha_i \in N^\alpha_{i+1} \backslash N^\alpha_i$
materialize a complete $L(N^\alpha_i)$-type $p^\alpha_i$
\sn
\item "{$(d)$}"  if $j < \omega_1$ is a limit ordinal then
$N^\alpha_j = \cup\{N^\alpha_i:i < j\}$
\sn
\item "{$(e)$}"  $p^\alpha_i \notin \{\text{gtp}(\bar a,N^\beta_{j_1},
N^\beta_{j_2}):j_1 < j_2 < \omega_1,\bar a \in 
{}^{\omega>}(N^\beta_{j_2})$ and $\beta < \alpha\}$ (see Definition
\scite{88r-4.2}(4))
\sn
\item "{$(f)$}"  $\bold F(N_{2 \beta +1}) \le_{\frak K} N_{2 \beta +2}$.
\ermn
As $\aleph_1 < 2^{\aleph_1} = 2^{\aleph_0}$ this is possible, i.e.,
in clause (e) we should find a type which is not in a set of $\le \aleph_1
\times |\alpha| < 2^{\aleph_0}$ types, as the number of possibilities
is $2^{\aleph_0}$;  let $N_\alpha =
\cup\{N^\alpha_i:i < \omega_1\}$ for $\alpha < 2^{\aleph_0}$, clearly
$N_\alpha \in K_{\aleph_1}$.  Now toward contradiction if
$\beta < \alpha < 2^{\aleph_0}$ and $N_\alpha \cong N_\beta$ then there is
an isomorphism $f$ from $N_\alpha$ onto $N_\beta$; necessarily $f$
maps $N^\alpha_i$
onto $N^\beta_i$ for a club of $i$.  For any such $i,p^\alpha_i \in
\text{ gtp}_L(f(\bar a^\alpha_i),N^\beta_i,N^\beta_j)$ for $j$ large
enough, contradiction. \nl
${{}}$   \hfill$\square_{\scite{88r-4.7}}$ 
\enddemo
\bigskip

\proclaim{\stag{88r-4.8} Lemma}  1) There are countable 
$L^0_\alpha \subseteq \Bbb L^0_{\omega_1,\omega}(\tau^{+0})$ for
$\alpha <\omega_1$ increasing continuous in $\alpha$, closed
under finitary operations and subformulas such that, 
letting $L^0_{< \omega_1} = \cup\{L^0_\alpha:\alpha < \omega_1\}$ we
have (some clauses do not metion the $L^0_\alpha$'s): 
\mr
\item "{$(a)$}"  for each $N \in K_{\aleph_0}$ and every complete
$L^0_\alpha(N)$-type $p(\bar x)$ we have $N \Vdash^{\aleph_1}_{\frak K}
(\exists \bar x) \wedge p(\bar x) \Rightarrow \wedge p \in 
L^0_{\alpha +1}(N)$.  Hence for every 
$\Bbb L^0_{\omega_1,\omega}(\tau^{+0})$-formula $\psi(\bar x)$ 
there are formulas $\varphi_n(\bar x) \in L^0_{< \omega_1}$ 
for $n < \omega$ such that $(N,N) \Vdash^{\aleph_1}_{\frak K} (\forall \bar x)
[\psi(\bar x) \equiv \dsize \bigvee_n \varphi_n(\bar x)]$
\sn
\item "{$(b)$}"  for every $N_0 \le_{\frak K} N_1 \in K_{\aleph_0}$ there is
$N_2,N_1 \le_{\frak K} N_2 \in K_{\aleph_0}$, such that for every
$\bar a \in N_2$ and $\varphi(\bar x) \in 
\Bbb L^0_{\omega_1,\omega}(N_0)$ we have 
$(N_2,N_0) \Vdash^{\aleph_1}_{\frak K} \varphi[\bar a]$ or $(N_2,N_0)
\Vdash^{\aleph_1}_{\frak K} \neg \varphi[\bar a]$
\sn
\item "{$(c)$}"  If $N \le_{\frak K} N_\ell \in K_{\aleph_0}(\ell = 1,2),\bar
a_\ell \in N_\ell$ and the 
$L^0_{< \omega_1}(N)$-generic types of $\bar a_\ell$ in $N_\ell$ are 
equal (though they are not necessarily complete; i.e., for 
every $\varphi(\bar x) \in L^0_{< \omega_1}(N)$ we have 
$N_1 \Vdash^{\aleph_1}_{\frak K} \varphi(\bar a_1)$ iff 
$N_2 \Vdash^{\aleph_1}_{\frak K} \varphi[\bar a_2]$), \ub{then} so 
are the $\Bbb L^0_{\infty,\omega}
(N)$-generic types.  In fact, there is $M,N
\le_{\frak K} M$ and $\le_{\frak K}$-embeddings $f_\ell:N_\ell \rightarrow M$ 
such that $f_\ell$ maps $N$ onto itself and $f_1(\bar a_1) = f_2(\bar a_2)$
though we do not claim $f_1 \restriction N = f_2 \restriction N$
\sn
\item "{$(d)$}"  for each $N \in K_{\aleph_0}$ and complete
$\Bbb L^0_{\omega_1,\omega}(N)$-type $p(\bar x)$, 
the class $K^1 =: \{(N,M,\bar a):M \in K_{\aleph_0},N \le_{\frak K} M,\bar a$ 
materialize $p$ in $(M,N)\}$ is a PC$_{\aleph_0}$-class
\sn
\item "{$(e)$}"  if 
$N \le_{\frak K} M \in K_{\aleph_0},\bar a \in M$, and for some complete
$\Bbb L^{-1}_{\omega_1,\omega}(N)$-type $p(\bar x),\bar a$ materialize $p$ in
$(M,N)$, \ub{then} for some complete $\Bbb L^0_{\omega_1,\omega}(N)$-type
$q_p,\bar a$ materialize $q_p$ in $(M,N)$; on $\Bbb L^0,\Bbb L^{-1}$
see Definition \scite{88r-4.1}(1),(3)
\sn
\item "{$(f)$}"  the number of complete
$\Bbb L^0_{\omega_1,\omega}(N)$-types $p$ which for some $\bar a,M$ we have 
$\bar a \in {}^{\omega >} M,M \in K_{\aleph_0},N \le_{\frak K} M$ and
$\bar a$ materialize in $(M,N)$ is $\le \aleph_1$
\sn
\item "{$(g)$}"  if in clause (f) we get that there are 
$\aleph_1$ such types then $\dot I(\aleph_1,K) \ge \aleph_1$
\sn
\item "{$(h)$}"   let $L^{-1}_\alpha =: L^0_\alpha \cap 
\Bbb L^{-1}_{\omega_1,\omega}(\tau^{+0})$ then the parallel clauses to (a)-(g)
holds.
\ermn
2) Clause (e) means that
\mr
\item "{$(i)$}"  assume further that $N_0 \le_{\frak K} N_\ell \in
K_{\aleph_0}$ for $\ell =1,2$ and $\bar a_\ell \in N_\ell$ and the 
$L^{-1}_{< \omega_1}(N)$-type which $\bar a_1$
materializes in $N_1$ is equal to the $L^{-1}_{< \omega_1}(N)$-type 
which $\bar a_2$ materializes in $N_2$.
\ub{Then} we can find $N^+_1,N^+_2$ such that $N_\ell \le_{\frak K}
N^+_\ell \in K_{\aleph_0}$ for $\ell = 1,2$ and isomorphism $f$ from
$N^+_1$ onto $N^+_2$ mapping $N$ onto itself and $\bar a_1$ to $\bar a_2$.
\endroster
\endproclaim
\bigskip

\remark{\stag{88r-4.8A} Remark}   1) We cannot get rid of the case of
$\aleph_1$ types (but see \scite{88r-5.16}, \scite{88r-5.17}) by the following
well known example.  For let $K =
\{(A,E,<):E$ an equivalence relation on $A$, each $E$-equivalence
class is countable, $x < y \Rightarrow xEy$ and on each
$E$-equivalence class $<$ is a 1-transitive linear order, and $xEy
\Rightarrow (x/E,<,x) \cong (y/E,$ \nl
$<,y)\}$ and $M \le_{\frak K} N$ if 
$M \subseteq N$ and $[x \in M \wedge y \in N \wedge xEy \Rightarrow y
\in M]$.
\nl
2) In clauses (c),(i) the mapping are not necessarily the identity on
$N$.  In clause (i) the assumption is apparently weaker (those by its
conclusion the assumption of (c) holds).
\nl
3) Note that clause (f) does not follow from clause (a) as there may
be $\aleph_1$-Kurepa trees.
\nl
4) In clause (c) of \scite{88r-4.8} for the second sentence we can weaken
the assumption: if $\varphi(\bar x) \in L^0_{< \omega_1}(N)$ and
$(N_1,N) \Vdash^{\aleph_1}_{\frak K} \varphi(\bar a_1)$ then
$(N_2,N) \nVdash^{\aleph_1}_{\frak K} \neg \varphi(\bar a_2)$.
This is enough to get the $M_{1,\alpha},M_{2,\alpha}$ from the proof.
(Why?  For each $\alpha < \omega_1$, there are $M_{1,\alpha}$ such
that $N_1 \le_{\frak K} M_{1,\alpha} \in K_{\aleph_0}$ and a complete
$\ell g(\bar a_i)-L^0_\alpha$-type $p_*(\bar x)$ such that
$(M_{1,\alpha},N) \Vdash \wedge p_1[\bar a_1]$.  But $\wedge p_1(\bar
x) \in L_{\alpha +1}$ hence $(N_2,N) \nVdash^{\aleph_1}_{\frak K} \neg
\wedge p_*(\bar a_2)$ hence there is $M_{2,\alpha}$ such that $N_2
\le_{\frak K} M_{2,\alpha} \in K_{\aleph_0}$ and $(M_{2,\alpha},N)
\Vdash^{\aleph_1}_{\frak K} \wedge p[\bar a_2]$.  Now continus as in
the proof below).
\endremark
\bigskip

\demo{Proof}  Note that proving clause (d) we say ``repeat the proof
of clause (a),(b),(c) for $L^{-1}_{\omega,\omega}$".
\sn
\ub{Clause $(a)$}:   We choose $L^0_\alpha$ by
induction on $\alpha$ using \scite{88r-4.7}.  The second phrase is proved
by induction on the depth of the formula using \scite{88r-4.6}.
\mn
\ub{Clause $(b)$}:   By iterating $\omega$ times, it suffices to prove
this for each $\bar a \in N_1$, so again by iterating $\omega$ times it
suffices to prove this for a fix $\bar a \in N_1$.
\sn
If the conclusion fails we can define by induction on $n < \omega$ for
every $\eta \in {}^n 2$, a model $M_\eta$ and $\varphi_\eta(\bar x) \in
\Bbb L^0_{\omega_1,\omega}(N)$ such that:
\mr
\widestnumber\item{$(iii)$}
\item "{$(i)$}"   $M_{<>} = N_1$
\sn
\item "{$(ii)$}"   $M_\eta \le_{\frak K} M_{\eta \char 94 <\ell>} \in
K_{\aleph_0}$ for $\ell = 0,1$
\sn
\item "{$(iii)$}"   $(M_\eta,N) \Vdash^{\aleph_1}_{\frak K} 
\varphi_\eta(\bar a)$
\sn
\item "{$(iv)$}"   $\varphi_{\eta \char 94 <1>} (\bar x) = \neg
\varphi_{\eta \char 94 <0>}(\bar x)$.
\ermn
Now for $\eta \in {}^\omega 2$, let $M_\eta = \dbcu_{n < \omega}
M_{\eta \restriction n}$.  Clearly for $\eta \in {}^\omega 2$ we have
$M_\eta \Vdash^{\aleph_1}_{\frak K} 
(\exists \bar x)[\dsize \bigwedge_{n < \omega}
\varphi_{\eta \restriction n}(\bar x)]$ and, after slight work, we get
contradiction to \scite{88r-4.7} + \scite{88r-4.6}.
\mn
\ub{Clause $(c)$}:   In general by clause (a) for each $\alpha < \omega_1$ 
we can find $M^\alpha_\ell \in K_{\aleph_1}$ for $\ell =1,2$ such that
$N_\ell \le_{\frak K} M^\alpha_\ell$ and 
$(M^\alpha_1,\bar a_1),(M^\alpha_2,\bar a_2)$ are 
$L^0_\alpha(N)$-equivalent and \wilog \, $N,N_\ell,M^\alpha_\ell$ have
universe an ordinal $< \omega_1$.  Let ${\frak A} = 
({\Cal H}(\aleph_2),N,N_1,N_2,\langle M^\alpha_1:\alpha < \omega_1
\rangle,\langle M^\alpha_2:\alpha < \omega_1 \rangle)$ let ${\frak
A}_1 \prec {\frak A}$ be countable  and find a non-well ordered 
countable model ${\frak A}_2$, which is an end extension of ${\frak
A}_1$ for $\omega^{{\frak A}_1}_1$, hence $\omega^{{\frak A}_2} =
\omega$ so $N^{{\frak A}_2} = N,N^{{\frak A}_2}_\ell = N_\ell$ for
$\ell = 1,2$.  Let $x \in (\omega_1)^{{\frak A}_2} \backslash
{\frak A}_1$ and $M^x_\ell = (M^x_\ell)^{{\frak A}_2}$ so $N_\ell
\le_{\frak K} M^x_\ell \in K_{\aleph_0}$.  Now 
there are $x_n$ such that ${\frak A}_2 \models
``x_{n+1} < x_n$ are countable ordinals"; \wilog \, $x_0 = x$ so 
using the hence and forth argument $(M^{x_0}_1,\bar a_1,N) 
\cong (M^{x_0}_2,\bar a_2,N)$.
\nl
[Why?  Let ${\Cal F}_n = \{(\bar b^1,\bar b^2):\bar b^\ell \in
{}^n(M^{x_0}_\ell)$ and ${\frak A}_2 \models \text{\rm
gtp}_{L^0_{x_n}}(\bar a^1 \char 94 \bar b^1,N;M^{x_0}_1) = \text{\rm
gtp}_{L^0_{x_n}}(\bar a^2 \char 94 \bar b^2;N;M^{x_0}_2)\}$.  Clearly
$(<>,<>) \in {\Cal F}_0$ and if $(\bar b^1,\bar b^2) \in {\Cal
F}_2,\ell \in \{1,2\}$ and $b^\ell_n \in M^{x_0}_\ell$ then there is
$b^{3 - \ell}_n \in M^{x_0}_{3 - \ell}$ such that $(\bar b^1 \char 94
\langle b^1_n\rangle,\bar b^2 \char 94 \langle b^n_2\rangle) \in {\Cal
F}_{n+1}$.  As $M^{x_0}_1,M^{x_0}_2$ are countable we can find an
isomorphism.]
\nl
But this is as required in the second phrase of (c).

We still have to prove the first phrase.  For this we prove by
induction on the ordinal $\alpha$ that
\mr
\item "{$\circledast_\alpha$}"  if for $\ell=1,2,\bar a_\ell \in
{}^{\omega >}(N_\ell)$ materialize in $(N_\ell,N_*)$ a complete
$L^0_{< \omega_1}$-type $p(\bar x)$ not depending on $\ell$ and
$\varphi(\bar x) \in \Bbb L^0_{\infty,\omega}(N_*)$ the quantifier
depth $< \alpha$ then: $\ell \in \{1,2\} \Rightarrow (N_\ell,N_*)
\Vdash^{\aleph_1}_{\frak K} \varphi(\bar a_\ell)$ \ub{or} $\ell \in
\{1,2\} \Rightarrow (N_\ell,N_*) \Vdash^{\aleph_1}_{\frak K} \neg
\varphi(\bar a_\ell)$.
\ermn
For $\alpha < \omega_1$ we already know it and for $\alpha$ limit
there is nothing to do.  For $\alpha = \beta +1$ it is enough to
consider $\varphi(\bar x)$ of the form $(\exists y) \dsize \bigwedge_i
\varphi_i(y,\bar x)$ where each $\varphi_i$ has quantifier depth $<
\beta$.

Assume this fails; recall if $(N_1,N_*)$ does not force
$(\Vdash^{\aleph_1}_{\frak K})$ the formula $\varphi(\bar a_\ell)$
then there is $N'_1$ such that $N_1 \le_{\frak K} N'_1 \in K_{\aleph_0},N'_1
\Vdash^{\aleph_1}_{\frak K} ``\neg \varphi(\bar a_1)"$.  Hence there are
$N^+_\ell$ such that $N_\ell \le_{\frak K} N^+_\ell \in
K_{\aleph_0},N^+_\ell$ forces $\varphi(\bar a_\ell) \equiv \bold t_\ell$
for some $\bold t_\ell \in \{\text{true,false}\}$ but $\bold t_1 \ne
\bold t_2$, so \wilog \, $N^+_1 \Vdash^{\aleph_1}_{\frak K}
\varphi(\bar a_1),N^+_2 \Vdash^{\aleph_1}_{\frak K} \neg \varphi(\bar
a_1)$, so \wilog \, for some $b \in N^+_1,N^+_1
\Vdash^{\aleph_1}_{\frak K} \dsize \bigwedge_i \varphi_i(b,\bar
a_\ell)$.  By the second sentence clause (c) which we have already
proved \wilog \, $(N^+_1,N_*,\bar a_1) \cong
(N^+_2,N_*,\bar a_2)$, and we get a contradiction.
\mn
\ub{Clause $(d)$}:   Let $N_0 \le_{\frak K} M_0 \in K_{\aleph_0}$ 
and $\bar a_0 \in M_0$ be such that 
$(N_0,M_0) \Vdash^{\aleph_1}_{\frak K} 
\dsize \bigwedge_{\varphi(\bar x) \in p} \varphi[\bar a_0]$ (if it
does not exist, the set of triples is empty).  Let $K'' =:
\{(N,M,\bar a):M \in K_{\aleph_0},N \in K_{\aleph_0},N \le_{\frak K} M$, 
and there are $M'' \in K_{\aleph_0},M \le_{\frak K} M''$ and 
$\le_{\frak K}$-embedding $f:M_0 \rightarrow
M''$, such that $f(N_0) = N,g(\bar a_0) = \bar a\}$.  Clearly it is 
a PC$_{\aleph_0}$ class.  

Now first if $(N,M,\bar a) \in K''$ let $(M'',f)$ witness this so by
applying clause (b) of \scite{88r-4.8} gtp$_{\Bbb L^0_{\omega_1,\omega}}(\bar
a;N;M) = \text{\rm gtp}_{\Bbb L^0_{\omega_1,\omega}}(\bar a,N,f(M_0)) =
\text{\rm gtp}_{\Bbb L^0_{\omega_1,\omega}}(a_0;N_0;M_0) = p$ so $(N,M,\bar
a) \in K''$.

Second, if $(N,M,\bar a) \in K''$ let $f_0$ be an isomorphism from
$M_0$ onto $M_0$.  Let $(M_1,f_1)$ be such that $N_0 \le_{\frak K} M_1
\in K_{\aleph_0},f_1 \supseteq f_0$ is an isomorphism from $M_1$ onto $M$ and
$\bar a_1 = f^{-1}_i(\bar a)$ hence $p =
\text{\rm gtp}_{\Bbb L^0_{\omega_1,\omega}}
(\bar a_1;N_0;M_1)$ and we apply clause
(c) of \scite{88r-4.8} with $N_0,M_0,\bar a_0,M_1,\bar a_1$ here standing
for $N,M_1,\bar a_1,M_2,\bar a_2$ there and can finish easily.
\mn
\ub{Clause $(e)$}:   We can define $\langle L^{-1}_\alpha:
\alpha < \omega_1 \rangle$ satisfying the parallel of Clause $(a)$ and
repeat the proofs of clauses (b),(c) and we are done. 
\mn
\ub{Clause $(f)$}:   Suppose this fails.

The proof splits to two cases.
\enddemo
\bn
\ub{Case A}:  $2^{\aleph_0} = 2^{\aleph_1}$.

We shall prove $\dot I(\aleph_1,K) \ge 2^{\aleph_0}$, thus, (as
$2^{\aleph_0} = 2^{\aleph_1}$) contradicting Hypothesis \scite{88r-4.5}
(this will be the only use of the hypothesis). \nl
Let $p_i$ (for $i < \omega_2)$ be distinct complete
$\Bbb L^0_{\omega_1,\omega}(\tau^{+0})$-types such that for each $i,p_i$ is 
materialized in some pair 
$(M,N)$, so $N \le_{\frak K} M \in K_{\aleph_0}$ 
(they exist by the assumption that (f) fails).  For
each $i < \omega_2$ we define 
$N_{i,\alpha},\xi_{i,\alpha}$ (for $\alpha < \omega_1$) and
$\bar a_{i,\alpha}$ such that:
\mr
\item "{$\boxtimes_1$}"  $(i) \quad N_{i,\alpha} \in K_{\aleph_0}$ has universe
$\omega(1 + \alpha),N_{0,0} = N$
\sn
\item "{${{}}$}"  $(ii) \quad \langle N_{i,\alpha}:
\alpha < \omega_1 \rangle$ is $\le_{\frak K}$-increasing continuous
\sn
\item "{${{}}$}"   $(iii) \quad \bar a_{i,\alpha} \in N_{i,\alpha +1},\bar
a_{i,\alpha}$ materialize $p_i$ in $(N_{i,\alpha +1},N_{i,\alpha})$
\sn
\item "{${{}}$}"  $(iv) \quad$ for every 
$\alpha < \beta < \omega_1$ and $\bar a
\in {}^{\omega >}(N_{i,\beta})$, the sequence $\bar a$ 
materialize \nl

\hskip20pt in $(N_{i,\beta},N_{i,\alpha})$ a
complete $\Bbb L^0_{\omega_1,\omega}(\tau^{+0})$-type
\sn
\item "{${{}}$}"  $(v) \quad \xi_{i,\alpha} < \omega_1$ is strictly increasing
continuous in $\alpha$
\sn
\item "{${{}}$}"  $(vi) \quad$ for $\alpha < \beta,N_{i,\beta}$ is pseudo
$L^0_\beta(N_{i,\alpha})$-generic, see \scite{88r-4.3}(4) and take care 
\nl

\hskip20pt of $\bold Q$, i.e., if $\gamma < \beta,p(y,\bar x)$ a complete
$L^0_\gamma$-type, 
\nl

\hskip20pt $(N_{i,\beta},N_{i,\alpha})
\Vdash^{\aleph_1}_{\frak K} (\bold Q y) \wedge p(y,\bar a)$ 
\nl

\hskip20pt then for
some $b \in N_{i,\beta +1} \backslash N_{i,\beta}$ we have $(N_{i,\beta
+1},N_{i,\alpha}) \Vdash^{\aleph_1}_{\frak K} \wedge p(b,\bar a)$
\sn
\item "{${{}}$}"  $(vii) \quad$ if 
$\alpha < \beta$ and $\bar a,\bar b \in N_{\beta-1}$
materialize different $\Bbb L^0_{\omega_1,\omega}(N_{i,\alpha})$-types in
\nl

\hskip20pt $N_{i,\beta}$, \ub{then} $\bar a,\bar b$ realize different
$(\Bbb L_{\omega_1,\omega}(\tau^{+0}) \cap L_{\xi_{i,\beta +1}})
(N_\alpha)$-types in 
\nl

\hskip20pt $N_{i,\beta}$
\sn
\item "{${{}}$}"  $(viii) \quad N_i = \cup\{N_{i,\alpha}:\alpha < \omega_1\}$
\sn
\item "{${{}}$}"  $(ix) \quad$ for $\alpha < \beta$, if $\bar a \in {}^{\omega
>}(N_{i,\beta})$ then for some complete $L^0_{\xi_{i,\beta}}$-type
$p,\bar a$ 
\nl

\hskip20pt materialize $p$ in $(N_{i,\beta},N_{i,\alpha})$ and $p$ has a
unique extension to a 
\nl

\hskip20pt $L^0_{< \omega_1}$-type
\sn
\item "{${{}}$}"  $(x) \quad$ if $\alpha_\ell < \beta$ for
$\ell=1,2,\gamma < \beta,n < \omega$ and $\bar a_1 \in
{}^n(N_{i,\beta})$ then for some
\nl

\hskip20pt $\bar a_2 \in {}^n(N_{i,\beta})$ we
have gtp$_{L^0_\gamma}(\bar a_1;N_{i,\alpha_1};N_{i,\beta}) =$
\nl

\hskip20pt $\text{\rm gtp}_{L^0_\gamma}(\bar a_2;N_{i,\alpha_2};N_{i,\beta})$.
\ermn
This is possible by the earlier claims.  By clause (e) clearly
\mr
\item "{$\boxtimes_2$}"  $(N_i,N_0)$ is $L^{-1}_{<
\omega_1}(\tau^{+0})$-homogeneous.
\ermn
We could below use $D_i$ a set of complete $L^0_{\delta(i)}$-types,
the only problem is that the countable $(D_i,\aleph_0)$-homogeneous
models have to be redefined using ``materialized" instead ``realized".
As it is we need to use clause (e) to translate the results on
$L^0_{\delta(i)}$ to $L^{-1}_{\delta(i)}$. 
\nl
Let $\tau^* = \{\in,Q_1,Q_2\} \cup \{c_\ell:\ell < 5\},c_\ell$ an
individual constant and ${\frak A}^*_i$ be $({\Cal H}(\aleph_2),\in)$ 
expanded to a $\tau^*$-model, by predicates for $K,\le_{\frak K}$
with $Q^{{\frak A}^*_i}_1 = K \cap {\Cal H}(\aleph_2),
Q^{{\frak A}^*}_2 = \{(M,N):M \le_{\frak
K} N$ both in ${\Cal H}(\aleph_2)\},c^{{\frak A}^*}_0,\dotsc,c^{{\frak
A}^*_i}_4$ being $\{\langle N_{i,\alpha}:\alpha < \omega_1 \rangle\},
\langle \xi_{i,\alpha}:\alpha < \omega_1 \rangle,\{\bar a_i\},N_i$ and
$\{i\}$ respectively.

Let ${\frak A}_i$ be a countable elementary submodel of ${\frak A}^*_i$
so $|{\frak A}_i| \cap \omega_1$ is an ordinal $\delta(i) < \omega_1$.
It is also clear that $c^{{\frak A}_i}_3$ is 
$N_{i,\delta(i)}$ as $c^{{\frak A}^*_i}_3 = N_i$.  As
${\frak A}_i$ is defined for $i < \omega_2$, for some unbounded $S
\subseteq \omega_2$ and $\delta < \omega_1$, for every $i \in S,
\delta(i) = \delta$ and for $i,j \in S$, some sequence from $N_j$
materializes $p$ in the pair $(N_{i,0},N_{i,\delta,i})$ iff $i=j$.  For $i
\in S$ let $D_i = \{p:p$ is a complete 
$L^{-1}_{\delta(i)}$-type materialized in
$(N_{i,\delta(i)},N_{i,0})\}$.  Because of the $\xi_{i,\alpha}$'s choice
and $\boxtimes_2$ the pair 
$(N_{i,\delta},N_0)$ is $(D_i,\aleph_0)$-homogeneous and $D_i$ is a countable 
 set of complete $L^{-1}_\delta$-types.  
Note that $(N_{i,\delta},N_{i,0},\bar a_{i,0})
\ncong (N_{j,\delta},N_{j,0},\bar a_{j,0})$ for $i \ne j (\in S)$ by the
choice of $S$, hence
$|\{j \in S:D_j = D_i\}| \le \aleph_0$, hence \wilog \, 
$i \ne j (\in S) \Rightarrow D_i \ne D_j$.

Let $\Gamma = \{D:D$ a countable set of complete 
$L^{-1}_\delta$-types, such that for some model ${\frak A} = {\frak A}_D$ of
$\dbca_{i \in S} \text{ Th}_{{\Bbb L}_{\omega,\omega}}
({\frak A}_i)$, with
$\{a:{\frak A}_D \models$ ``a countable ordinal$\} = \delta$ (and the
usual order) we have $D = \{\{\varphi(\bar x):\varphi(\bar x) \in
L^{-1}_\delta$ and ${\frak A}_D \models
``(N_{i,\alpha},N_0) \Vdash^{\aleph_1}_{\frak K} 
\varphi[\bar a]"\}:\bar a \in \dbcu_{i < \delta} N^{\frak A}_{i,\alpha}\}\}$.

So $D_i \in \Gamma$ for $i < \omega_2$, hence $\Gamma$ is uncountable.

By standard descriptive set theory $\Gamma$ (is an analytic set hence)
has cardinality continuum.  So
let $D(\zeta) \in \Gamma$ be distinct for $\zeta < 2^{\aleph_0}$.  For
each $\zeta$, let ${\frak A}^0_{D(\zeta)}$ be as in the definition of
$\Gamma$.  We define by induction on $\alpha < \omega_1,
{\frak A}^\alpha_{D(\zeta)}$ such that
\mr
\item "{$(\alpha)$}"  ${\frak A}^\alpha_{D(\zeta)}$ is countable
\sn
\item "{$(\beta)$}"  $\alpha < \beta \Rightarrow 
{\frak A}^\alpha_{D(\zeta)} \prec_{{\Bbb L}_{\omega,\omega}} 
{\frak A}^\beta_{D(\zeta)}$
\sn
\item "{$(\gamma)$}"  for limit $\alpha$ we have ${\frak A}^\alpha_{D(\zeta)} =
\dbcu_{\beta < \alpha} {\frak A}^\beta_{D(\zeta)}$
\sn
\item "{$(\delta)$}"  if $d \in {\frak A}^{\alpha + 1}_{D(\zeta)}
\backslash {\frak A}^\alpha_{D(\zeta)},{\frak A}^{\alpha
+1}_{D(\zeta)} \models ``d$ a countable ordinal" \ub{then} for 
$a \in {\frak A}^\alpha_{D(\zeta)}$ we have 
${\frak A}^{\alpha +1}_{D(\zeta)} \models$ 
``if $a$ is a countable ordinal then $a < d$"
\sn
\item "{$(\varepsilon)$}"  for $\alpha = 0$ in clause $(\delta)$ 
there is no minimal such $d$
\sn
\item "{$(\zeta)$}"  for every $\alpha$ there is $d_{\zeta,\alpha} \in
{\frak A}^{\alpha +1}_{D(\zeta)} \backslash 
{\frak A}^\alpha_{D(\zeta)}$ satisfying 
${\frak A}^{\alpha +1}_{D(\zeta)} \models 
``d_{\zeta,\alpha}$ a countable ordinal" and for 
$\alpha \ne 0$ it is minimal.
\ermn
Let $M_{\zeta,\alpha}$ be the $d_{\zeta,\alpha}$-th member of the
$\omega_1$-sequence of models in ${\frak A}^\beta_{D(\zeta)}$ for
$\beta > \alpha$ (remember $c^{{\frak A}^*_i}_0 = 
\langle N_{i,\alpha}:\alpha < \omega_1 \rangle$).
Let $M_\zeta = \dbcu_{\alpha < \omega_1} M_{\zeta,\alpha}$.  
By absoluteness from ${\frak A}^\beta_{D(\zeta)}$ we have
$M_{\zeta,\alpha} \le_{\frak K} M_{\zeta,\beta} \in K_{\aleph_0}$.
Now
\mr
\item "{$(*)$}"   $0 < \alpha < \beta,(M_{\zeta,\beta},M_{\zeta,\alpha})$ is
$(D(\zeta),\aleph_0)$-homogeneous.
\ermn
[Why?  Assume ${\frak A}^\alpha_{D(\zeta)} \models ``d_1 < d_2$ are
countable ordinals $> \gamma$" when $\gamma < \delta$.  Now if $\bar
a,\bar b \in {}^{\omega >}(N^{{\frak A}^\alpha_{D(\zeta)}}_{d_2})$ and
$[\gamma < \delta \Rightarrow \text{\rm gtp}_{L^0_\gamma}(\bar
a,N^{{\frak A}^\alpha_{D(\zeta)}}_{d_1};N^{{\frak
A}^\alpha_{D(\zeta)}}_{d_2}) = \text{\rm gtp}_{L^0_\gamma}(\bar b,
N^{{\frak A}^\alpha_{D(\zeta)}}_{d_1};N^{{\frak
A}^\alpha_{D(\zeta)}}_{d_2})]$ and also ${\frak A}^\alpha_{D(\zeta)}$
satisfies this but ``${\frak A}^\alpha_{D(\zeta)}$ thinks that the
countable ordinals are well ordered" \ub{hence} for some $d,{\frak
A}^\alpha_{D(\zeta)} \models ``d$ is a countable ordinal $> \gamma$"
for each $\gamma < \delta$ and we have ${\frak A}^\alpha_{D(\zeta)}
\models ``\text{\rm gtp}_{L^0_d}(\bar a,N_{d_1};N_{d_2}) = \text{\rm
gtp}_{L^0_d}(\bar a,N_{d_1},N_{d_2})$.  Hence if ${\frak
A}^\alpha_{D(\zeta)} \models ``d' < d"$ then for every $a \in
N^{{\frak A}^\alpha_{D(\zeta)}}_{d_2}$ for some $b \in N^{{\frak
A}^\alpha_{D(\zeta)}}_{d_2}$ we have

$$
{\frak A}^\alpha_{D(\zeta)} \models ``\text{gtp}_{L^0_d}(\bar a
{}^\frown \langle a \rangle,N_{d_1};N_{d_2}) = \text{
gtp}_{L^0_d}(\bar b {}^\frown \langle b \rangle;N_{d_1};N_{d_2})"
$$
\mn
hence 
gtp$_{L^0_\gamma}(\bar a {}^\frown \langle a \rangle;
N^{{\frak A}^\alpha_{D(\zeta)}}_{d_2}) = \text{ gtp}(\bar b {}^\frown
\langle b \rangle;;N^{{\frak A}^\alpha_{D(\zeta)}}_{d_1}; 
N^{{\frak A}^\alpha_{D(\zeta)}}_{d_2})$.

Also we can replace $L^0_-$ by $L^{-1}_-$.  By clause (x) of
$\boxtimes_1$ the set $\{\text{gtp}_{L^0_\delta}(\bar a,
N^{{\frak A}^\alpha_{D(\zeta)}}_{d_1};N^{{\frak
A}^\alpha_{D(\zeta)}}_{d_2}):a \in {}^{\omega >}(N^{{\frak
A}^\alpha_{D(\zeta)}}_{d_2})\}$ is $D_i$.

So $(N^{{\frak A}^\alpha_{D(\zeta)}}_{d_2},N^{{\frak
A}^\alpha_{D(\zeta)}}_{d_2})$ is $(D_i,\aleph_0)$-homogenous.    

So from the isomorphism type of $M_\zeta$ we can compute $D(\zeta)$.
So $\zeta \ne \xi \Rightarrow M_\zeta \ncong M_\xi$.  As $M_\zeta \in
K_{\aleph_1}$ we finish.
\bn
\ub{Case B}:  $2^{\aleph_0} < 2^{\aleph_1}$.

By \scite{88r-3.5}, ${\frak K}$ has the $\aleph_0$-amalgamation property.  So
clearly if $N \le_{\frak K} M \in K_{\aleph_0},\bar a \in M$, \ub{then} $\bar
a$ materializes in $(M,N)$ a complete 
$\Bbb L^0_{\omega_1,\omega}(\tau^{+0})$-type.  We would now like 
to use descriptive set theory.

We represent a complete $\Bbb L^0_{\omega_1,\omega}(\tau^{+0})$-type 
materialized in some $(N,M)$ by a real, by representing the isomorphism
type of some $(N,M,\bar a),
N \le_{\frak K} M \in K_{\aleph_0},\bar a \in M$.  The set of
representatives is analytic recalling ${\frak K}$ is PC$_{\aleph_0}$, 
and the equivalence relation is $\Sigma^1_1$.
[As $(N_1,M_1,\bar a_1),(N_2,M_2,\bar a_2)$ represents the same type
if and only if for some $(N,M),N \le_{\frak K} M \in K_{\aleph_0}$, there are
$\le_{\frak K}$-embeddings 
$f_1:M_1 \rightarrow M,f_2:M_2 \rightarrow M$ such that
$f_1(N_1) = f_2(N_2) = N,f_1(\bar a) = f_2(\bar a)$.]

By Burgess \cite{Bg} (or see \cite{Sh:202}) as there are $> \aleph_1$
equivalence classes, there is a perfect set of representation,
pairwise representing different types.

\relax From this we easily get that \wilog \, that their restriction to some
$L^0_\alpha$ are distinct, contradicting part
(a).
\mn
\ub{Clause $(g)$}:   Easy by the proof of clause (f), Case A above. 
\mn
\ub{Clause $(h)$}:  As in the proof of clause (e).
\hfill$\square_{\scite{88r-4.8}}$ 

\remark{\stag{88r-4.8B} Remark}   1) Note that in the proof of Clause (f)
of \scite{88r-4.8}, in Case (A) we get many types too but it
was not clear whether we can make the $N_\zeta$ to be generic enough,
to get the contradiction we got in Case (B) but is not crucial here. \nl
2)  We may like to replace
$\Bbb L^0_{\omega_1,\omega}$ by $\Bbb L^1_{\omega_1,\omega}$ in \scite{88r-4.6},
\scite{88r-4.7} and \scite{88r-4.8} (except that, for our benefit, in
\scite{88r-4.8}(e); we may retain the definition of $L^1(N))$.  We 
lose the ability to build $L$-generic models in $K_{\aleph_1}$ 
(as the number of (even unary)
relations on $N \in K_{\aleph_0}$ is $2^{\aleph_0}$, which may be $>
\aleph_1$).  However, we can say ``$\bar a$ materializes in $N \in
K_{\aleph_0}$ the type $p=p(\bar x)$ which is a complete type in
$\Bbb L^1_{\omega_1,\omega}(N_n,N_{n-1},\dotsc,N_0)$; where $N_0
\le_{\frak K} \ldots \le_{\frak K} N_n \le_{\frak K} N,N_\ell$ 
countable)". \nl
[Why?  Let some $N^1,\bar a^1$ be as above, 
$\bar a^1$ materialize $p$ in $(N^1,N_n,\dotsc,N_0)$ \ub{then} 
this holds for $(N,\bar a)$ \ub{iff}
for some $N',f$ we have $N \le_{\frak K} N' \in K_{\aleph_1}$ and
$f$ is an isomorphism from $N^1$ onto $N''$ mapping $\bar a^1$ to $\bar
a$ and $N_\ell$ to $N_\ell$ for $\ell \le n$.  If there is no such
pair $(N^1,\bar a^1)$ this is trivial.] \nl
We can get something on formulas.  

This suffices for \scite{88r-4.6}.
\endremark
\bigskip

\remark{\stag{88r-4.9} Concluding remarks for \S4}  We can get more
information on the case $1 \le \dot I(\aleph_1,K) < 2^{\aleph_1}$ (and
the case $1 \le \dot I(\aleph_1,K^{\bold F}_{\aleph_1})$, etc.).
\nl
1) As in \scite{88r-3.5}, there is no difficulty in getting the results of
this section
for the class of models of $\psi \in \Bbb L_{\omega_1,\omega}(\bold Q)$ as
using $(K,\le_{\frak K})$ 
from the proof of \scite{88r-3.9}(2) in all constructions
we get many non-isomorphic models for appropriate $\bold F$, as in
\scite{88r-4.5.3}(2).
\nl
2) For generic enough $N \in K_{\aleph_1}$ with $\le_{\frak
K}$-representation $\langle 
N_\alpha:\alpha < \omega_1 \rangle$, we have determined the
$N_\alpha$'s (by having that \wilog \, $K$ is categorical in
$\aleph_0$).  In this section we have shown that for some club $E$ of
$\omega_1$, for all $\alpha < \beta$ from $E$ the isomorphism
type of $(N_\beta,N_\alpha)$ essentially
\footnote{why only essentially? as the number of relevant complete
types can be $\aleph_1$; we can get rid of this shrinking ${\frak K}$}
is unique.
We can continue the analysis, e.g., deal with sequences $N_0
\le_{\frak K}  N_1 \le_{\frak K} \ldots \le_{\frak K} N_k \in 
K_{\aleph_0}$ such that $N_{\ell +1}$ is pseudo
$L^0_\alpha(N_\ell,N_{\ell -1},\dotsc,N_0)$-generic.  
We can prove by induction on $k$ that for any countable $L \subseteq \Bbb
L^0_{\omega_1,\omega}
(\tau^{+k})$ for some $\alpha$, any strong $L$-generic $N \in
K_{\aleph_1}$ is $L$-determined. 
That is, for any $\langle N_\alpha:\alpha < \omega_i\rangle,N_\alpha
\le_{\frak K} N$ countable $\le_{\frak K}$-increasing continuous with
union $N$, for some club $E$ for all $\alpha_0 < \ldots < \alpha_k$
from $N$ the isomorphic type of $\langle
N_{\alpha_k},N_{\alpha_k},\dotsc,N_{\alpha_0}\rangle$ is the same;
i.e., determining for $\Bbb L_{\infty,\omega}(aa)$.
\nl
3) We can do the same for stronger logics, let us elaborate.  

Let us define a logic ${\Cal L}^*$.  It has as variable

variables for elements $x_1,x_2$...and

variables for filters ${\Cal Y}_1,{\Cal Y}_2$...
\sn
The atomic formulas are:
\mr
\item "{$(i)$}"  the usual ones
\sn
\item "{$(ii)$}"  $x \in \text{ Dom}({\Cal Y})$.
\ermn
The logical operations are:
\mr
\item "{$(a)$}"  $\wedge$ conjunction, $\neg$ negation
\sn
\item "{$(b)$}"  $(\exists x)$ existential quantification where $x$ is
individual variable
\sn
\item "{$(c)$}"  the quantifier aa acting on variables ${\Cal Y}$ so
we can form $(aa\, {\Cal Y})\varphi$
\sn
\item "{$(d)$}"  the quantification 
$(\exists x \in \text{ Dom}({\Cal Y})) \varphi$
\sn
\item "{$(e)$}"  the quantification $(\exists^f x \in \text{\rm
Dom}({\Cal Y}))\varphi$.
\ermn
It should be clear what are the free variables of a formula $\varphi$.
The variable ${\Cal Y}$ vary on pairs (a countable set, a filter on
the set).
Now in $\exists x[\varphi,{\Cal Y}],(\exists x \in \text{ Dom}({\Cal
Y}))\varphi,(\exists^fx \in \text{ Dom}({\Cal Y}))\varphi,x$ is 
bounded but not ${\Cal Y}$ and in
$aa{\Cal Y},{\Cal Y}$ is bounded.  
The satisfaction relation is defined as usual plus
\mr
\item "{$(\alpha)$}"  $M \models (\exists x \in \text{ Dom}
({\Cal Y})\varphi(x,{\Cal Y},\bar a)$ if and only if for 
some $b$ from the domain of ${\Cal Y},M \models \varphi[b,{\Cal Y},\bar a]$
\sn
\item "{$(\beta)$}"  $M \models \exists^f x \in \text{ Dom}
({\Cal Y})\varphi(x,{\Cal Y}_{\bar a})$ 
if and only if $\{x \in \text{ Dom}({\Cal Y}):\models \varphi
(x,{\Cal Y},\bar a)\} \in {\Cal Y}$
\sn
\item "{$(\gamma)$}"  $M \models (aa\,{\Cal Y},\bar a) 
\varphi({\Cal Y})$ if and 
only if there is a function $\bold F$ from ${}^{\omega >}([M]^{< \aleph_1})
\rightarrow [M]^{< \aleph_1}$ such that:
if $A_n \subseteq M,|A_n| \le \aleph_0,A_n \subseteq A_{n+1}$ and
$\bold F(A_0,\dotsc,A_n) \subseteq A_{n+1}$ \ub{then} $M \models
\varphi[{\Cal Y}_{\langle A_n:n < \omega \rangle},\bar a]$ where 
${\Cal Y}_{\langle A_n:n < \omega \rangle}$ is the filter 
on $\dbcu_{n < \omega} A_n$, generated
by $\{\cup\{A_n:n < \omega\} \backslash A_\ell:\ell < \omega\}$. 
\ermn
4) We, of course, can define ${\Cal L}^*_{\mu,\kappa}$ (extending $\Bbb
L_{\mu,\ell}$).  As we like to analyze models in $\aleph_1$, it is most
natural to deal with ${\Cal L}^*_{\omega_1,\omega}$.
\sn 
We can prove that (if $1 \le \dot I(\aleph_1,{\frak K}) < 2^{\aleph_1})$
the quantifier $aa\,{\Cal Y}$ is determined on $K_{\aleph_1}$ (i.e.,
for almost all ${\Cal Y},\varphi({\Cal Y})$ iff not for 
almost all ${\Cal Y},\neg \varphi({\Cal Y})$. \nl
5) The logic from (3) strengthens the stationary logic $\Bbb L(aa)$, see
\cite{Sh:43}, \cite{BKM78}.
\nl
Not so strongly: looking at PC class for $\Bbb L_{\omega_1,\omega}(aa)$
(i.e., $\{M \restriction \tau:M$ a model of $\psi$ of cardinal
$\aleph_1\}$), we can assume that $\psi \vdash ``<$ is an
$\aleph_1$-like order".  Now we can express $\varphi \in 
{\Cal L}^*_{\omega_1,\omega}$, but the determinacy tells us more.  Also we
can continue to define higher variables ${\Cal Y}$.
\endremark
\newpage

\head {\S5 There is a superlimit model in $\aleph_1$} \endhead  \resetall \sectno=5
 \spuriousreset
\bn
Here we make
\demo{\stag{88r-5.0} Hypothesis}  Like \scite{88r-4.5}, but also $2^{\aleph_0}
< 2^{\aleph_1}$.

This section is the deepest (of this paper = chapter).  
The main difficulties are proving the
facts which are obvious in the context of \cite{Sh:48}.  So while it
was easy to show that every $p \in \bold D^*(N)$ is definable over a finite
set ($\bold D^*(N)$ defined below), it was not clear to me how to
prove that if you extend the type $p$ to
$q \in \bold D^*(M),(N \le_{\frak K} M \in 
K_{\aleph_0})$ by the same definition, then
$q \models p$ (remember $p,q$ are types materialized not realized, and
at this point in the paper we still do not have the tools to replace
the models by uncountable generic enough models).  So we rather have
to show that failure is a non-structure property, i.e., implies
existence of many models.

Also symmetry of stable amalgamation becomes much more complicated.
We prove existence of stable amalgamation by four stages 
(\scite{88r-5.15},\scite{88r-5.17}(3),\scite{88r-5.20},\scite{88r-5.22}).  The
symmetry is proved as a consequence of uniqueness of one sided
amalgamation, (so it cannot be used in its proof).  The culmination of
the section is the existence of a superlimit models in $\aleph_1$
(\scite{88r-5.24}).  This seems a natural stopping point as it seems
reasonable to expect that the next step
should be phrasing the induction on $n$, i.e., dealing with $\aleph_n$
and ${\Cal P}(n - \ell)$-diagrams of models of power $\aleph_\ell$ as
in \cite{Sh:87a}, \cite{Sh:87b}; (so this is done in \chaptercite{705}).
\enddemo
\bigskip

\definition{\stag{88r-5.1} Definition}   We define functions $\bold D,\bold D^*$
with domain $K_{\aleph_0}$. \nl
1) For $N \in K_{\aleph_0}$ let $\bold D(N) 
= \{p:p$ is a complete $\Bbb L^0_{\omega_1,\omega}(N)$-type over $N$ such
that for some $\bar a \in M \in K_{\aleph_0},N \le_{\frak K} M$ and $\bar a$
materializes $p$ in $(M,N)\}$, (i.e. the members of $p$ have the form
$\varphi(\bar x,\bar a)$, ($\bar x$ finite and fixed for $p$) $\bar a$ a
finite sequence from $N$ and $\varphi \in \Bbb L^0_{\omega_1,\omega}(N))$. \nl
2) For $N \in K_{\aleph_0}$ let $\bold D^*(N) = \{p:p$ a complete
$\Bbb L^0_{\omega_1,\omega}(N;N)$-type such that for some $\bar a \in M \in
K_{\aleph_0},N \le_{\frak K} M$ and $\bar a$ materializes $p$ in $(M,N;N)\}$.
\enddefinition
\bn
\margintag{88r-5.1A}\ub{\stag{88r-5.1A} Explanation}:  0) Recall that any formula in $\Bbb
L^0_{\omega_1,\omega}(N)$ has finitely many free variables.
\nl
1) So for every finite $\bar b \in N$ 
and $\varphi(\bar x,\bar y) \in \Bbb L^0_{\omega_1,\omega}(N)$, 
if $p \in \bold D(N)$, then $\varphi(\bar x,\bar b)
\in p$ or $\neg \varphi(\bar x,\bar b) \in p$. \nl
2) But a formula from $p \in \bold D^*(N)$ 
may have all $c \in N$ as parameters whereas a formula from $p \in
\bold D(N)$ can mention only finitely many members of $N$.
\bigskip

\proclaim{\stag{88r-5.2} Lemma}  1) ${\frak K}$ has the $\aleph_0$-amalgamation
property. \nl
2) If $N_* \le_{\frak K} N \in K_{\aleph_0}$, 
$A_i \subseteq N_*$ for $i \le n$
\ub{then} for every sentence 
\nl
$\psi \in \Bbb L^1_{\infty,\omega}(N_*,A_n,\dotsc,A_1;A_0)$ we have

$$
N \Vdash^{\aleph_1}_{\frak K} \psi \text{ or } N 
\Vdash^{\aleph_1}_{\frak K} \neg \psi.
$$
\mn
3) If $N \le_{\frak K} M \in K_{\aleph_0}$, \ub{then} 
every $\bar a \in M$ materializes in $(M;N)$ one and only 
one type from $\bold D^*(N)$ and also materializes in $(M,N)$ one and
only one type from $\bold D(N)$.  Also
for every $N \le_{\frak K} M \in K_{\aleph_0}$ and $q \in \bold D^*(N)$ 
for some $M',M \le_{\frak K} M' \in  K_{\aleph_0}$ and some
$\bar b \in M'$ materializes $q$ in $(M;N)$. \nl
4) For every $N \in K_{\aleph_0}$ and countable $L \subseteq 
\Bbb L^0_{\omega_1,\omega}(N;N)$ the 
number of complete $L(N;N)$-types $p$ 
such that $N \Vdash^{\aleph_1}_{\frak K} ``(\exists \bar x) \wedge p"$
is countable; note that pedantically $L \subseteq \Bbb
L_{\omega_1,\omega}(\tau^+ \cup\{c:c \in N\})$ and we
restrict ourselves to models $M$ such that $P^M = |N|,c^M =c$. \nl
5) For $N \in K_{\aleph_0}$ there are 
countable $L^0_\alpha \subseteq \Bbb L^0_{\omega_1,\omega}(N;N)$ 
for $\alpha < \omega_1$ increasing
continuous in $\alpha$, closed under finitary operations (and
subformulas) such that:
\mr
\item "{$(*)$}"  for each complete $L^0_\alpha(N;N)$-type $p$ we have

$$
[N \Vdash^{\aleph_1}_{\frak K} \exists \bar x \wedge p \Rightarrow \wedge p
\in L^0_{\alpha +1}].
$$
\ermn
Hence for every $\Bbb L^0_{\omega_1,\omega}(N;N)$ formula $\psi(\bar x)$ for
some $\varphi_n(\bar x) \in \dbcu_{\alpha < \omega} L^0_\alpha$ 
for $n < \omega$ for every $N \in K_{\aleph_0}$

$$
(N,N) \Vdash^{\aleph_1}_{\frak K} (\forall \bar x)[\psi(\bar x) \equiv \dsize
\bigvee_{n < \omega} \varphi_n(\bar x)].
$$
\mn
6) For $N \in K_{\aleph_0}$ we have $|\bold D^*(N)| \le \aleph_1$ and
$|\bold D(N)| \le \aleph_1$. \nl
7) If $p \in \bold D^*(N)$ \ub{then} there is $q$ such that if $N \le_{\frak
K} M \in K_\lambda,\bar a \in M$ materializes $p$ in $(M;N)$ 
\ub{then} the complete $\Bbb L^0_{\infty,\omega}(N,N)$-type which $\bar a$
realizes in $M$ over $N$ is $q$; also $q$ belongs to $\bold D(N)$ and
is unique.  Moreover, we can replace $p$ by the complete
$\Bbb L^{-1}_{\omega_1,\omega}(N,N)$-type which $\bar a$ materializes
in $M$.
Similarly for $\bold D(N),\Bbb L^0_{\infty,\omega}(N),\Bbb
L^{-1}_{\omega_1,\omega}(N)$. \nl
8) If $n < \omega$ and $\bar b,\bar c \in {}^n N$ realize the same $\Bbb
L_{\omega_1,\omega}(\tau)$-type in $N$ then they materialize
the same $\Bbb L^1_{\omega_1,\omega}(\tau^{+0})$-type in $(N,N)$.
\nl
9) If $f$ is an isomorphism from $N_1 \in K_{\aleph_0}$ onto $N_2 \in
K_{\aleph_0}$ \ub{then} $f$ induces a one to one function from 
$\bold D(N_1)$ onto $\bold D(N_2)$ and from $\bold D^*(N_1)$ onto
$\bold D^*(N_2)$.
\endproclaim
\bigskip

\demo{Proof}  1) By \scite{88r-3.5}. \nl
2) By 1). \nl
3) By 2) and 1). \nl
4) Like the proof of \scite{88r-4.7} (just easier). \nl
5) Like the proof of \scite{88r-4.8}(a). \nl
6) Like the proof of \scite{88r-4.8}(f) (and \scite{88r-0.9}). \nl
7) Clear as in $p \in \bold D^*(N)$ we allow formulas than for $q \in
\bold D(N)$.
\nl
8),9)  Easy, too.    \hfill$\square_{\scite{88r-5.2}}$
\enddemo
\bn
We shall use from now on a variant of gtp (from Definition \scite{88r-4.2}(4)).
\definition{\stag{88r-5.2.7} Definition}  1) If 
$N_0 \le_{\frak K} N_1 \in K_{\aleph_0},\bar a \in N_1$, 
gtp$(\bar a,N_0,N_1)$ is the $p \in \bold D(N_0)$ such that 
$(N_1,N_0) \Vdash^{\aleph_1}_{\frak K} \wedge p[\bar a]$.  
So $\bar a$ materializes (but does not necessarily realize)
gtp$(\bar a,N_0,N_1)$.  We may omit $N_1$ when clear from context. 
We define gtp$^*(\bar a,N_0,N_1) \in \bold D^*(N_0)$ similarly.
\nl
2) We say $p = \text{ gtp}^*(\bar b,N_0,N_1)$ is definable over $\bar
a \in N_0$ \ub{if} gtp$(\bar b,N_0,N_1) = p^-  =: p \cap 
\Bbb L^0_{\omega_1,\omega}(N_0)$ is definable over
$\bar a$ (see Definition \scite{88r-5.4.1} below,
note that $p \mapsto p^-$ is a one-to-one mapping from
$\bold D^*(N_0)$ onto $\bold D(N_0)$ by \scite{88r-5.5}(1) below).  So 
stationarization is defined for $p \in \bold D^*(N_0)$, too, after we
know \scite{88r-5.5}(1).
\enddefinition
\bigskip

\proclaim{\stag{88r-5.4} Claim}  1) Each $p \in \bold D(N)$ does not $\Bbb
L^0_{\omega_1,\omega}(\tau^{+0})$-split (see Definition \scite{88r-5.4.1}
below; also see more below) over 
some finite subset $C$ of $N$, hence $p$ is definable over it. \nl
Moreover, letting $\bar c$ list $C$ there is a function $g_p$
satisfying $g_p(\varphi(\bar x,\bar y))$ is $\psi_{p,\varphi}(\bar
y,\bar z) \in \Bbb L_{\omega_1,\omega}(\tau)$ such
that for each $\varphi(\bar x,\bar y) \in \Bbb
L^0_{\omega_1,\omega}(N)$ and $\bar a \in N$ we have
$[\varphi(\bar x,\bar a) \in p
\Leftrightarrow N \models \psi_{p,\varphi}(\bar a,\bar c)]$, 
(in particular, $\bold Q$ is ``not necessary"). \nl
2) Every automorphism of $N$ maps $\bold D(N)$ onto itself and each $p \in
\bold D(N)$ has at most $\aleph_0$ possible images; we may also call
them conjugates.  So if $g$ is an
isomorphism from $N_0 \in K_{\aleph_0}$ onto $N_1 \in K_{\aleph_0}$
then $g(\bold D(N_0)) = \bold D(N_1)$.
\nl
3) If $N_0 \le_{\frak K} N_1 \le_{\frak K} N_2 \in 
K_{\aleph_0}$ and $\bar a \in N_1$ then
{\rm gtp}$(\bar a,N_0,N_1) = \text{\rm gtp}(\bar a,N_0,N_2)$.
\endproclaim
\bn
Before we prove \scite{88r-5.4}:
\definition{\stag{88r-5.4.1} Definition}  Assume
\mr
\item "{$(a)$}"  $N$ is a model
\sn
\item "{$(b)$}"  $\Delta_1$ is a set of formulas (possibly in a
vocabulary $\nsubseteq \tau_N$) closed under negation
\sn
\item "{$(c)$}"  $\Delta_2$ is a set of formulas in the vocabulary 
$\tau = \tau_N$
\sn
\item "{$(d)$}"  $p$ is a $(\Delta_1,n)$-type over $N$ (i.e., each member
has the form $\varphi(\bar x,\bar a),\bar a$ from $N,\varphi(\bar
x,\bar y)$ from $\Delta_1,\bar x = \langle x_\ell:\ell < n \rangle$;
no more is required
(we may allow other formulas but they are irrelevant)
\sn
\item "{$(e)$}"  $A \subseteq N$.
\ermn
0) We say $p$ is a complete $\Delta_1$-type over $B$ when:
\mr
\widestnumber\item{$(iii)$}
\item "{$(i)$}"  $B \subseteq N$
\sn
\item "{$(ii)$}"  $\varphi(\bar x,\bar b) \in p \Rightarrow \bar b
\subseteq A \wedge \varphi(\bar x,\bar y) \in \Delta_1$
\sn
\item "{$(iii)$}"  if $\varphi(\bar x,\bar y) \in \Delta_1$ and $\bar b
\in {}^{\ell g(\bar y)} A$ then $\varphi(\bar x,\bar b) \in p$ or $\neg
\varphi(\bar x,\bar b) \in p$.
\ermn
The default value here for $\Delta_1$ is $\Bbb
L_{\omega_1,\omega}(\tau_{\frak K})$.
\nl    
1) We say that $p$ does $(\Delta_1,\Delta_2)$-split over $A$ when there are
$\varphi(\bar x,\bar y) \in \Delta_1$ and $\bar b,\bar c \in 
{}^{\ell g(\bar y)} N$ such that
\mr
\item "{$(\alpha)$}"   $\varphi(\bar x,\bar b),\neg \varphi(\bar
x,\bar c) \in p$
\sn
\item "{$(\beta)$}"  $\bar b,\bar c$ realize the same $\Delta_2$-type
over $A$.
\ermn
2) We say that $p$ is $(\Delta_1,\Delta_2)$-definable over $A$
\ub{when}: for every formula $\varphi(\bar x,\bar y) \in \Delta_1$
there is a formula $\psi(\bar y,\bar z) \in \Delta_2$ and $\bar c \in
{}^{\ell g(\bar z)}A$ such that

$$
\varphi(\bar x,\bar b) \in p \Rightarrow N \models \psi[\bar b,\bar c]
$$

$$
\neg \varphi(\bar x,\bar b) \in p \Rightarrow N \models \neg \psi[\bar
b,\bar c]
$$
\mn
(in the case $p$ is complete over $B,\bar b \subseteq B$ we get ``iff")
\enddefinition
\bigskip

\demo{\stag{88r-5.4.2} Observation}  Assume
\medskip

$\quad (a),(b),(c),(d),(e) \qquad$ as in \scite{88r-5.4.1} and in addition

$\quad (d)^+ \qquad \qquad \qquad \qquad p$ is a complete
$(\Delta_1,n)$-type over $N$, i.e., if
\nl

\hskip100pt $\varphi(\bar x,\bar y) \in \Delta_1,\bar d \in {}^{\ell g(\bar y)} N,\bar x = 
\langle x_\ell:\ell < n \rangle$ \nl

\hskip100pt then $\varphi(\bar x,\bar d) 
\in p$ or $\neg \varphi(\bar x,\bar d) \in d$.
\sn
\ub{Then} the following conditions are equivalent:
\mr
\item "{$(\alpha)$}"    $p$ does $(\Delta_1,\Delta_2)$-splits over $A$  
\sn
\item "{$(\beta)$}"   there is a sequence of $\langle
g_{\varphi(\bar x,\bar y)}:\varphi(\bar x,\bar y) \in \Delta_1
\rangle$ of functions such that:
\sn
\item "{${{}}$}"   $(i) \quad g_{\varphi(\bar x,\bar y)}$ is a
function with domain including $\{\text{tp}_{\Delta_2}(\bar
b,A,N)$:
\nl

\hskip25pt $\bar b \in {}^{\ell g(\bar y)} N\}$
\sn
\item "{${{}}$}"    $(ii) \quad$ the values of 
$g_{\varphi(\bar x,\bar y)}$ are truth values
\sn
\item "{${{}}$}"    $(iii) \quad$ if $\varphi(\bar x,\bar y) \in
\Delta_1,\bar b \in {}^{\ell g(\bar y)} N$ and $q = \text{ tp}_{\Delta_2}(\bar
b,A,N)$ \ub{then}:
\nl

\hskip25pt  $\varphi(\bar x,\bar b) \in p \Rightarrow
g_{\varphi(\bar x,\bar y)}(q) =$ true, and
\nl

\hskip25pt  $\neg \varphi(\bar x,\bar b)
\in p \Rightarrow g_{\varphi(\bar x,\bar y)}(q) =$ false.
\endroster
\enddemo
\bigskip

\demo{Proof of \scite{88r-5.4.2}}  Reflect on the definitions.
\enddemo
\bigskip

\demo{Proof of \scite{88r-5.4}}  1) Assume this fails, 
$N \le_{\frak K} M,\bar a \in M$
materializes $p$ and for every $\bar b \in M,(M,N) \Vdash \wedge q[\bar
b]$ for some $q(\bar x) \in \bold D(N)$ and let 
$\langle b^*_\ell:\ell < \omega \rangle$ list $N$.
We choose by induction on $n,\langle C^0_\eta,C^1_\eta,f_\eta,\bar
a^0_\eta,\bar a^1_\eta:\eta \in {}^n 2 \rangle$ such that
\mr
\item "{$(a)$}"  $C^\ell_\eta$ is a finite subset of $N$ for $\ell <
2,\eta \in {}^n 2$
\sn
\item "{$(b)$}"  $f_\eta$ is an automorphism of $N$ mapping $C^0_\eta$
onto $C^1_\eta$
\sn
\item "{$(c)$}"  $\{b^*_{\ell g(\eta)}\} \cup C^0_\eta \cup C^1_\eta
\subseteq C^0_{\eta \char 94 <\ell>} \cap C^1_{\eta \char 94 <\ell>}$
for $\ell =0,1$
\sn
\item "{$(d)$}"  $\bar a^0_\eta,\bar a^1_\eta \in N$ realize in $N$ the same
$\Bbb L_{\omega_1,\omega}(\tau)$-type over $C^0_\eta \cup C^1_\eta \cup
\{b^*_{\ell g(\eta)}\}$ in $(M,N)$ but $\bar a \char 94 \bar a^0_\eta,\bar a
\char 94 \bar a^1_\eta$ do not materialize the same 
$\Bbb L^0_{\omega_1,\omega}(\tau^{+0})$ in $(M,N)$ (this exemplifies splitting)
\sn
\item "{$(e)$}"   $f_{\eta \char 94 <0>}(\bar a^0_\eta) = 
\bar a^1_\eta,f_{\eta \char 94 <1>}(\bar a^1_\eta) = \bar a^1_\eta$
\sn
\item "{$(f)$}"  $f_\eta \restriction C^0_\eta \subseteq f_{\eta \char
94 <\ell>}$ for $\ell=0,1$
\sn
\item "{$(g)$}"  $\bar a^0_\eta \char 94 \bar a^1_\eta \subseteq
C^0_{\eta \char 94 <\ell>} \cap C^1_{\eta \char 94 <\ell>}$.
\ermn
For $n=0$ let $C^0_\eta,C^1_\eta = \emptyset,f_\eta = \text{ id}_N$.
Recall that $K_{\aleph_0}$ is categorical in $\aleph_0$ and $N$ is
countable, hence if $n < \omega,\bar b',\bar b'' \in {}^n N$ realize the
same $\Bbb L_{\omega_1,\omega}(\tau)$-type over a finite subset $B$ of
$N$, then some automorphism of $N$ over $B$ maps $\bar b'$ to $\bar b''$.
If $(C^0_\eta,C^1_\eta,f_\eta)$ are defined and satisfies clauses (a), (b)
we recall that by our assumption toward contradiction as
$C^0_\eta \cup C^1_\eta \cup \{b^*_{\ell g(\eta)}\}$ is a finite subset
of $N$, there are $\bar a^0_\eta,\bar a^1_\eta \in {}^{\omega >}
N$ as required in clause (d).  So clearly there are automorphisms
$f_{\eta \char 94 <0>},f_{\eta \char 94 <1>}$ extending $f_\eta
\restriction C^0_\eta$ such that $f_{\eta \char 94 <0>}(\bar a^0_\eta)
= \bar a^1_\eta,f_{\eta \char 94 <1>}(\bar a^1_\eta) = \bar a^1_\eta$
as required in clause (e), (f).

Lastly, choose $C^0_{\eta \char 94 <\ell>} = C^0_\eta \cup \{b^*_{\ell
g(\eta)},f^{-1}_{\eta \char 94 <\ell>}(b^*_{\ell g(\eta)}),\bar a^0_\eta
\char 94 \bar a^1_\eta,f^{-1}_{\eta \char 94 <\ell>}(\bar a^0_\eta
\char 94 \bar a^1_\eta)\}$ and $C^1_{\eta \char 94 <\ell>} = f_{\eta
\char 94 <\ell>}(C^0_{\eta \char 94 <\ell>})$.

Having carried the induction, for every $\eta \in {}^\omega 2$ clearly
$f_\eta = \cup\{f_{\eta \restriction n} \restriction C^0_\eta:n <
\omega\}$ is an automorphism of $N$.
\nl
[Why?  As $\langle f_{\eta \restriction n} \restriction C^0_{\eta
\restriction n}:n < \omega \rangle$ is an increasing sequence of
functions, by clauses (c) + (f) the union $f_\eta$ is
a partial automorphism of $N$ by clause (b), and $f_\eta$ have domain 
$N$ by clause (c) and as $f_{\eta \restriction
n}(C^0_{\eta \restriction n}) = C^1_{\eta \restriction n}$ the union $f_\eta$
has range $N$ by clause (c).]  Hence for some $M_\eta \in
K_{\aleph_0}$ there is an isomorphism $f^+_\eta$ from $M$ onto
$M_\eta$ extending $f$.  Now for some $p_\eta \in \bold D(N),f_\eta(\bar a)$
materialize $p_\eta$ in $(M_\eta,N)$, but $\eta \ne \nu \in {}^\omega 2
\Rightarrow p_\eta \ne p_\nu$ by clauses (d) + (e), 
contradiction to \scite{88r-5.2}(4) as we
can use $\le \aleph_0$ formulas to distinguish. \nl
2) Follows. 
\nl
3) Trivial.   \hfill$\square_{\scite{88r-5.4}}$
\enddemo
\bigskip

\proclaim{\stag{88r-5.5} Claim}  1) Suppose $N_0 \le_{\frak K} N_1 \in 
K_{\aleph_0}$ and $N_1$ forces that $\bar a,\bar b$ (in $N_1$)
realize the same $\Bbb L^0_{\omega_1,\omega}(N_0)$-type 
over $N_0$, \ub{then} $N_1$ forces that they realize 
the same $\Bbb L^0_{\omega_1,\omega}(N_0;N_0)$-type; (the
inverse is trivial). \nl
2) If $N_0 \le_{\frak K} N_1 \le_{\frak
K} N_2 \in K_{\aleph_0}$
and $\bar a,\bar b \in N_2$ (remember $N_2$ determines the complete
$\Bbb L^0_{\omega_1,\omega}(N_1)$-generic types of $\bar a,\bar b$)
\ub{then} from the $\Bbb L^0_{\omega_1,\omega}(N_1)$-generic type of $\bar
a$ over $N_1$ we can compute the $\Bbb L^0_{\omega_1,\omega}(N_0)$-generic
type of $\bar a$ over $N_0$ (hence if the
$\Bbb L^0_{\omega_1,\omega}(N_1)$-generic types of $\bar a,\bar b$ over
$N_1$ are equal, \ub{then} so are the $\Bbb L^0_{\omega_1,\omega}(N_0)$-generic
types of $\bar a,\bar b$ over $N_0$). 
\endproclaim
\bigskip

\remark{Remark}  1) So there is no essential difference 
between $\bold D(N)$ and $\bold D^*(N)$. \nl
2) Recall that in a formula of
$\Bbb L^0_{\omega_1,\omega}(N_0;N_0)$ all $c \in N_0$ may appear as
individual constants.
\endremark
\bigskip

\demo{Proof}  1) We shall prove there are $N_2$ 
such that $N_1 \le_{\frak K} N_2 \in K_{\aleph_0}$ and an 
automorphism of $N_2$ over $N_0$ taking
$\bar a$ to $\bar b$; this clearly suffices; and we prove the existence
of such $N_2$, of course, by hence and forth
arguments.  So by renaming and symmetry, it suffices to prove that
\mr
\item "{$(*)$}"  if $m < \omega,N_0 \le_{\frak K} N_0$ and $\bar a,\bar
b \in {}^m(N_1)$ materialize the same $\Bbb
L^0_{\infty,\omega}(N_0)$-type over $N_0$ \ub{then} for
every $c \in N_1$, there are $N_2$ and $d \in N_2$ such 
that $\bar a \char 94 <c>,\bar b \char 94 <d>$ materialize the same
$\Bbb L^0_{\omega_1,\omega}(N_0)$-type over $N_0$.
\ermn  
However, by the previous claim \scite{88r-5.2} for some $\bar a^* \in
{}^{\omega >}(N_0)$ the $\Bbb L^0_{\omega_1,\omega}(N_0)$-type over 
$N_0$ that $\bar a \char 94
<c>$ materialize in $(N_1,N_0)$ does not split over $\bar a^*$.
Now $\bar a,\bar b$ materialize in $(N_1,N_0)$ the same $\Bbb
L^0_{\omega_1,\omega}(N_0)$-type over $N_0$ hence $\bar a^* \char 94 \bar
a,\bar a^* \char 94 \bar b$ materialize in $(N_1,N_0)$ the same $\Bbb
L^0_{\omega_1,\omega}(N_0)$-type.  Hence there is $N_2,N_1 \le_{\frak
K} N_2 \in K_0$ and an automorphism $f$ of $N_2$ mapping $N_0$ onto
$N_1$ and mapping $\bar a^* \char 94 \bar a$ to $\bar a^* \char 94 \bar b$
(but possibly $f \restriction N_0 \ne \text{\rm id}_{N_0}$).  Let
$d=f(c)$, hence if 
$\bar a \char 94 <c>,\bar b \char 94 <d>$ materialize the same
$\Bbb L^0_{\omega_1,\omega}(N_0)$-type in $(N_2,N_0)$ then 
they materialize the same 
$\Bbb L^0_{\omega_1,\omega}(N_0)$-type over $N_0$ in $(N_2,N_0)$. \nl
2) Clearly it suffices to prove the ``hence " part.
By the assumption and proof of \scite{88r-5.5}(1) there 
are $N_3$ satisfying $N_2 \le_{\frak K} N_3 \in K_{\aleph_0}$ 
and $f$ an automorphism of $N_3$ over $N_1$ taking $\bar a$ 
to $\bar b$.  Now the conclusion follows.  ${{}}$ \hfill$\square_{\scite{88r-5.5}}$
\enddemo
\bigskip

\definition{\stag{88r-5.6.1} Definition}  1) We say that $\bold D_*$ is a
${\frak K}$-diagram function when
\mr
\item "{$(a)$}"  $\bold D_*$ is a function with domain $K_{\aleph_0}$
(later we shall lift it to $K$)
\sn
\item "{$(b)$}"  $\bold D_*(N) \subseteq \bold D(N)$ and has at least
one non-algebraic member for $N \in K_{\aleph_0}$
\sn
\item "{$(c)$}"  if $N_1,N_2 \in K_{\aleph_0}$ and $f$ is an
isomorphism from $N_1$ onto $N_2$ then $f$ maps $\bold D_*(N_1)$ onto
$\bold D_*(N_2)$, this applies in particular to an automorphism of $N
\in K_{\aleph_0}$.
\ermn
1A) Such $\bold D_*$ is called weakly good if:
\mr
\item "{$(d)$}"  $\bold D_*(N)$ is closed under subtypes, that is: if $p(\bar
x) \in \bold D_*(N),\bar x = \langle x_\ell:\ell < m\rangle,\pi$ is a
function from $\{0,\dotsc,m-1\}$ into $\{0,\dotsc,n-1\}$ then some
(necessarily unique) $\bar q(\langle x_0,\dotsc,x_{n-1}\rangle) \in
\bold D_*(N)$ is equal to $\{\varphi(\langle x_0,\dotsc,\bar
x_{n-1}\rangle):\varphi(x_{\pi(0)},\dotsc,x_{\pi(m-1)}) \in p(\bar
x)\}$.
\ermn
2) Such $\bold D_*$ is called countable if $N \in K_{\aleph_0}
\Rightarrow |\bold D_*(N)| \le \aleph_0$.
\nl
3) Such $\bold D_*$ is called good when it is weakly good (i.e.,
clause (d) holds) and
\mr
\item "{$(e)$}"  $\bold D_*(N)$ has amalgamation (i.e., if $p_0(\bar
x),p_1(\bar x,\bar y),p_2(\bar x,\bar z) \in \bold D_*(N)$ and $p_0
\subseteq p_1 \cap p_2$ then there is $q(\bar x,\bar y,\bar z) \in
\bold D_*(N)$ which includes $p_1(\bar x,\bar y) \cup p_2(\bar x,\bar z))$.
\ermn
4) Such $\bold D_*$ is called very good if it is good and:
\mr
\item "{$(f)$}"  $N_0 \le_{\frak K} N_1 \le_{\frak K} N_2 \in
K_{\aleph_0},\bar a_0 \subseteq \bar a_1 \subseteq \bar a_2$ and $\bar
a_\ell \subseteq N_\ell$ for $\ell=0,1,2$ and gtp$(\bar a_{\ell
+1},N_\ell,N_{\ell +1})$ is definable over $\bar a_\ell$ and belongs
to $\bold D_*(N_\ell)$ for $\ell=0,1$ then gtp$(\bar a_2,N_0,N_2)$
belongs to $\bold D_*(N_0)$ and is definable over $\bar a_0$.
\endroster
\enddefinition
\bigskip

\demo{\stag{88r-5.6} Fact}  1) There are $\bold D_\alpha,
\bold D^*_\alpha$ for $\alpha <
\omega_1$, functions with domain $K_{\aleph_0}$ such that:
\mr
\item "{$(a)$}"  for $N \in K_{\aleph_0},\bold D_\alpha(N),
\bold D^*_\alpha(N)$
is a countable subset of $\bold D(N),\bold D^*(N)$ respetively
\sn
\item "{$(b)$}"  for each $N \in K_{\aleph_0},\langle
\bold D_\alpha(N):\alpha < \omega_1 \rangle$ as well as $\langle
\bold D^*_\alpha(N):\alpha < \omega_1 \rangle$ are increasing continuous
\sn
\item "{$(c)$}"  $\bold D(N) = \dbcu_{\alpha < \omega_1} \bold D_\alpha(N)$
and $\bold D^*(N) = \dbcu_{\alpha < \omega_1} \bold D^*_\alpha(N)$ 
\sn
\item "{$(d)$}"  if $N_1,N_2 \in K_{\aleph_0},f$ is an isomorphism from
$N_1$ onto $N_2$ \ub{then} $f$ maps $\bold D_\alpha(N_1)$ onto 
$\bold D_\alpha(N_2)$ and $\bold D^*_\alpha(N_1)$ onto 
$\bold D^*_\alpha(N_2)$ for $\alpha < \omega_1$
\sn
\item "{$(e)$}"  for every $\alpha < \omega_1$ and 
$N \in K_{\aleph_0}$ there is a
$(\bold D_\alpha(N),\aleph_0)^*$-homogeneous model (see below
Definition \scite{88r-5.7}(1))
(obviously it is unique up to isomorphism over $N$)
\sn
\item "{$(f)$}"  if $N_0 \le_{\frak K} N_1 \le_{\frak K} N_2 \in 
K_{\aleph_0},N_2$ is $(\bold D_\alpha(N_1),\aleph_0)^*$-homogeneous 
(see Definition \scite{88r-5.7}(1) below) and
$N_1$ is $(\bold D_\alpha(N_0),\aleph_0)^*$-homogeneous or just 
$(\bold D_\beta(N_0),\aleph_0)^*$-homogeneous for some $\beta \le
\alpha$
\ub{then} $N_2$ is $(\bold D_\alpha(N_0),\aleph_0)^*$-homogeneous
\sn
\item "{$(f)^+$}"  if $\langle \alpha_\varepsilon:\varepsilon \le
\zeta \rangle$ is increasing continuous sequence of countable
ordinals, $\zeta > 0$ and $\langle N_\varepsilon:\varepsilon \le \zeta
\rangle$ is $\le_{\frak K}$-increasing continuous, $N_\varepsilon \in
{\frak K}_{\aleph_0}$, for every $\bar a \in N_{\varepsilon +1}$,
gtp$(\bar a,N_\varepsilon,N_{\varepsilon +1}) \in \bold
D_\alpha(N_\varepsilon)$ and for every $\xi < \zeta$ for some
$\varepsilon \in [\xi,\zeta),N_{\varepsilon +1}$ is $(\bold
D_{\alpha_\varepsilon}(N_\varepsilon),\aleph_0)^*$-homogeneous
\ub{then} $N_\zeta$ is $(\bold D_{\alpha_\zeta}(N_0),\aleph_0)^*$-homogeneous
\sn
\item "{$(g)$}"  $N_1$ is $(\bold D_\alpha(N_0),
\aleph_0)^*$-homogeneous \ub{if and only if} $N_1$ is 
$(\bold D^*_\alpha(N_0),\aleph_0)^*$-homogeneous 
where $N_0 \le_{\frak K} N_1 \in K_{\aleph_0}$
\sn
\item "{$(h)$}"  $\bold D_\alpha$ is a very good ${\frak K}$-diagram function.
\ermn
2) If $\bold D$ is very good then clauses (d),(e),(f),(f)$^+$ hold.
\enddemo
\bigskip

\remark{\stag{88r-5.6.4} Remark}  1) We can add
\mr
\item "{$(h)$}"   if ${\frak K},<^*$ are as derived from the $\psi \in
\Bbb L_{\omega_1,\omega}(\bold Q)$ in the proof of \scite{88r-3.9} 
then we can add: if $N_0
\le_{\frak K} N_1 \in K_{\aleph_0}$ and every $p \in \bold D_0(N_0)$
is materialized in $N_1$ then $N_0 <^* N_1$.
\ermn
2) So our results apply to $\psi \in \Bbb L_{\omega_1,\omega}(\bold Q)$, too.
\nl
3) So it follows that if $\langle N_i:i \le \alpha \rangle$ is
$\le_{\frak K}$-increasing in $K_{\aleph_0},N_{i+1}$ is $(\bold
D_{\beta_i}(N_0),\aleph_0)^*$-homogeneous 
and $\langle \beta_i:i < \alpha \rangle$
is non-decreasing with supremum $\beta$ then $N_\alpha$ is $(\bold
D_\beta,\aleph_0)^*$-homogeneous. 
\nl
4) So each $\bold D_\alpha$ is very good and countable.
\endremark
\bigskip

\demo{Proof of \scite{88r-5.6}}   First, $\bold D$ is a ${\frak
K}$-diagram function by Definition \scite{88r-5.1} and \scite{88r-5.2}(9).  As
$\bold D(N)$ has cardinality $\le \aleph_1$ by \scite{88r-5.2}(6) we can
find a sequence $\langle \bold D_\alpha:\alpha < \omega_1\rangle$ such
that
\mr
\item "{$\circledast$}"  $(a) \quad \bold D_\alpha$ is a countable
${\frak K}$-diagram function
\sn
\item "{${{}}$}"  $(b) \quad$ for every $N \in K_{\aleph_0},\langle
\bold D_\alpha(N):\alpha < \omega_1\rangle$ is increasing continuous
with union $\bold D(N)$.
\ermn
Second, $\bold D$ is very good (clause (f) by \scite{88r-5.16} the other -
easier). Third, note that for each of the demands (d),(e),(f) from
Definition \scite{88r-5.6.1}, for a club of $\delta < \omega_1,\bold
D_\delta$ satisfies it.  So \wilog \, each $\bold D_\alpha$ is very
good.

The parts on $\bold D^*_\alpha$ follow (using also
\scite{88r-5.8}(1) below which does not rely on \scite{88r-5.6}-\scite{88r-5.7A}
(and see proof of \scite{88r-5.10}).  \hfill$\square_{\scite{88r-5.6}}$
\enddemo
\bigskip

\definition{\stag{88r-5.7} Definition}  Assume $N_0 \le_{\frak K} N_1 \in 
K_{\aleph_0}$.
\nl
1) We say that $(N_1,N_0)$ or just
$N_1$ is $(\bold D_\alpha(N_0),\aleph_0)^*$-homogeneous if:
\mr
\item "{$(a)$}"  every $\bar a \in N_1$ materializes in $(N_1,N_0)$
over $N_0$ some $p \in \bold D_\alpha(N_0)$ and every $q \in 
\bold D_\alpha(N_0)$ is materialized in $(N_0,N_1)$ by some $\bar b \in N_1$
\sn
\item "{$(b)$}"  if $\bar a,\bar b \in N_1,\bar a,\bar b$ materialize
in $(N_1,N_0)$ the same type over $N_0$ and $c \in N_1$ \ub{then} for
some $d \in N_1$ sequence $\bar a \char 94 <c>,\bar b \char 94 <d>$ materialize
in $(N_1,N_0)$ the same type from $\bold D_\alpha(N_0)$.
\ermn
2) Similarly for $(\bold D^*_\alpha(N_0),\aleph_0)^*$-homogeneity.
\enddefinition
\bigskip

\remark{\stag{88r-5.7A} Remark}   1) Now this is meaningful only for $N
\le_{\frak K} M \in K_{\aleph_0}$, but later it becomes meaningful for
any $N \le_{\frak K} M \in K$. \nl
2) Uniqueness for such countable models hold in this context too.
\endremark
\bn
Now by \scite{88r-5.5}.
\demo{\stag{88r-5.8} Conclusion}   If 
$(N_1,N_0)$ is $(\bold D_\alpha(N_0),\aleph_0)^*$-homogeneous
\ub{then} $(N_1,N_0,c)_{c \in N_0}$ is
$(\bold D^*_\alpha(N_0),\aleph_0)^*$-homogeneous. 
\enddemo
\bigskip

\demo{Proof}  This is easy by \scite{88r-5.5}(1) and clause (g) of \scite{88r-5.6}.
  \hfill$\square_{\scite{88r-5.8}}$
\enddemo
\bigskip

\proclaim{\stag{88r-5.9} Lemma}  There is $N^* \in K_{\aleph_1}$ such that
$N^* = \dbcu_{\alpha <\omega_1} N_\alpha$ and $N_\alpha \in K_{\aleph_0}$ is
$\le_{\frak K}$-increasing continuous with $\alpha$ and $N_{\alpha +1}$ is 
$(\bold D_{\alpha +1}(N_\alpha),\aleph_0)^*$-homogeneous for 
$\alpha < \omega_1$.
\endproclaim
\bigskip

\demo{Proof}  Should be clear.  \hfill$\square_{\scite{88r-5.9}}$
\enddemo
\bigskip

\proclaim{\stag{88r-5.10} Theorem}  The $N^* \in K_{\aleph_1}$ from
\scite{88r-5.9} is unique (even not depending on the choice of
$\bold D_\alpha(N)$'s), is universal and is 
model-homogeneous (for ${\frak K}$).
\endproclaim
\bigskip

\demo{Proof} 
\enddemo
\bn
\ub{Uniqueness}:  For $\ell =0,1$ let $N^\ell_\alpha,\bold D^\ell_\alpha \,
(\alpha < \omega_1)$ be as in \scite{88r-5.6}, \scite{88r-5.9} and we should
prove $\dbcu_{\alpha < \omega_1} N^0_\alpha \cong \dbcu_{\alpha <
\omega_1} N^1_\alpha$; because of clause (g) of \scite{88r-5.6} it does
not matter if we use the $\bold D$ or $\bold D^*$ version.  
As $\bold D^\ell_\alpha \, (\alpha < \omega_1)$ is
increasing and continuous, $|\bold D^\ell_\alpha(N)| \le \aleph_0$ and
$\dbcu_{\alpha< \omega_1} \bold D^\ell_\alpha(N) = \bold D(N)$ 
for every $N \in K_{\aleph_0}$ and the $\bold D^\ell_\alpha$'s 
commute with isomorphisms,  clearly there is a
closed unbounded $E \subseteq \omega_1$, such that $\alpha \in E
\Rightarrow \bold D^0_\alpha = \bold D^1_\alpha$.  Let $E = \{\alpha(i):i <
\omega_1\},\alpha(i)$ increasing and continuous.  Now we define by
induction on $i < \omega_1$, an isomorphism $f_i$ from
$N^0_{\alpha(i)}$ on $N^1_{\alpha(i)}$, increasing with $i$.  For
$i=0$ use the $\aleph_0$-categoricity of $K$ and for limit $i$ let
$f_i = \dbcu_{j<i} f_j$.  Suppose $f_i$ is defined, then by clause (d) of
\scite{88r-5.6} the function
$f_i$ maps $\bold D^0_{\alpha(i+1)}(N^0_{\alpha(i)})$
onto $\bold D^0_{\alpha(i+1)}(N^1_{\alpha(i)})$ and by the choice of $E,
\bold D^0_{\alpha(i+1)} = \bold D^1_{\alpha(i+1)}$.  By the assumption on the
$N^\ell_\alpha$ and clause (f)$^+$ of \scite{88r-5.6}, $N^\ell_{\alpha(i+1)}$ is
$(\bold D^\ell_{\alpha(i+1)}(N^\ell_{\alpha(i)}),\aleph_0)^*$-homogeneous.
Summing up those facts and \scite{88r-5.6}(e) we see that we can extend
$f_i$ to an isomorphism $f_{i+1}$ from $N^0_{\alpha(i+1)}$ onto
$N^1_{\alpha(i+1)}$.

Now $\dbcu_{i < \omega_1} f_i$ is the required isomorphism.
\bn
\ub{Universality}: Let $M \in K_{\aleph_1}$, so $M = \dbcu_{\alpha <
\omega_1} M_\alpha,M_\alpha$ increasing continuous and $\|M_\alpha\|
\le \aleph_0$.  We now define $f_\alpha,N_\alpha,\beta_\alpha$ by induction on
$\alpha < \omega_1$ such that: $\gamma_\alpha \in [\alpha,\omega_1)$ is
increasing continuous with $\alpha,f_\alpha$ is a $\le_{\frak K}$-embedding of
$M_\alpha$ into $N_\alpha \in K_{\aleph_0},N_\alpha$ is $\le_{\frak
K}$-increasing continuous, $f_\alpha$ is 
increasing and continuous, and for $\beta < \alpha,N_{\beta +1}$ is
$(\bold D_{\gamma_\beta +1}(N_\beta),\aleph_0)^*$-homogeneous.  
The only novelty over \scite{88r-5.9} is the use of the
$\aleph_0$-amalgmation property (which holds by \scite{88r-3.5},
\scite{88r-4.5}).  So $f = \cup\{f_\alpha:\alpha < \omega_1\}$ embeds $M$
into $N = \cup\{N_\alpha:\alpha < \omega_1\}$ which
is isomorphic to $N^*$ by the uniqueness.
So the universality follows from the
uniqueness.
\bn
\ub{Model-homogeneity}:  So let $\langle N_\alpha:\alpha < \omega_1
\rangle,\bold D_\alpha,N^*$ be as in
\scite{88r-5.6}, \scite{88r-5.9} and $M_\ell \le_{\frak K} N^*$ (for $\ell = 0,1$) are
countable, $f$ an isomorphism from $M_0$ onto $M_1$.  For some $\gamma
< \omega_1$ we have $M_0,M_1 \le_{\frak K} N_\gamma$. 
Every type in $\bold D(N_\gamma)$ is
realized in $N^*$ and $N_\beta$ is $(\bold D_\beta(N_\gamma),
\aleph_0)^*$-homogeneous for $\beta > \gamma$.
\nl
Now
\mr
\item "{$(*)$}"  if $M_\ell \le_{\frak K} M^+_\ell \in K_{\aleph_0},\bar
a {}^\frown \bar b \in M^+_\ell$ and $\beta \in (\gamma,\omega_1),\bar
a' \in N_\beta$ materializes gtp$(\bar a,M_\ell,M^+_\ell)$ then for some
$\beta_1 \in (\beta,\omega_1)$ and $\bar b' \in N_{\beta_1}$ the
sequence $\bar a' {}^\frown \bar b'$ materializes in $N_{\beta_1}$ the
type gtp$(\bar a {}^\frown \bar b,M_\ell,M^+_\ell)$.
\nl
[Why?  By amalgamation and equality of types \wilog \, $N_\beta
\le_{\frak K} M^+_\ell$ and $\bar a = \bar a'$ hence gtp$(\bar
b,N_\beta,M^+_\ell) \in \bold D(N_\beta)$ hence is materialized by
some $\bar b'$ from some $N_{\beta_1},\beta_1 \in [\beta,\omega_1)$.]
\ermn
So for a club of $\alpha \in (\gamma,\omega_1)$ the
model $N_\alpha$ is
$(\bold D_\alpha(M_\ell),\aleph_0)^*$-homogeneous for $\ell = 0,1$, so $f$
can be extended to an automorphism of $N_\alpha$, hence as in the
uniqueness part, to an automorphism of $N^*$.  \nl
${{}}$  \hfill$\square_{\scite{88r-5.10}}$
\bigskip

\definition{\stag{88r-5.11} Definition}  1) If $N_0 \le_{\frak K} N_1 \in
K_{\aleph_0}$ and $p_\ell \in \bold D(N_\ell)$ for $\ell =1,2$ and they are
definable in the same way (see Definition \scite{88r-5.4.1} (and
\scite{88r-5.4}), so in particular both do not split
over the same finite subset of $N_0$).  \ub{Then} we call $p_1$ the
stationarization of $p_0$ over $N_1$. \nl
2) For $p_\ell \in \bold D(N_\ell)$ for 
$\ell=0,1$ let $p_1 \models p_0$ mean that if $N_1
\le_{\frak K} N_2 \in K_{\aleph_0}$ and $\bar a \in N_2$ materializes
$p_1$ then it materializes $p_0$.
\enddefinition
\bigskip

\remark{\stag{88r-5.11A} Remark}   It is easy to justify the uniqueness
implied by ``\ub{the} stationarization". 
\endremark
\bn
Observe
\proclaim{\stag{88r-5.12.1} Claim}  If 
$p_\ell = \text{\rm gtp}(\bar a,N_\ell,N_2)$ for $\ell=0,1$ and $N_0
\le_{\frak K} N_1 \le_{\frak K} N_2 \in K_{\aleph_0}$ \ub{then}
$p_1 \models p_0$.
\endproclaim
\bigskip

\demo{Proof}  Easy.  \hfill$\square_{\scite{88r-5.12.1}}$
\enddemo
\bigskip

\proclaim{\stag{88r-5.16} Claim}   Suppose $N_0 \le_{\frak K} N_1 \le_{\frak
K} N_2 \in K_{\aleph_0},\bar a_i \in N_i$, (for $i = 0,1,2$), 
$\bar a_0 \subseteq \bar a_1
\subseteq \bar a_2,{\text{\rm gtp\/}}
(\bar a_1,N_0,N_1)$ is definable over $\bar a_0$ and 
${\text{\rm gtp\/}}(\bar a_2,N_1,N_2)$ is definable over $\bar a_1$.
\ub{Then} ${\text{\rm gtp\/}}(\bar a_2,N_0,N_2)$ is definable over $\bar a_0$.
Moreover, the definition depends only on the definition mentioned previously.
\endproclaim
\bigskip

\demo{Proof}  So we have to prove that gtp$(\bar a_2,N_0,N_2)$ does
not split over $\bar a_0$.
Let $n < \omega$ and $\bar b,\bar c \in
{}^n N_0$ realize the same type in $N_0$ over $\bar a_0$ (in the logic
$\Bbb L_{\omega_1,\omega}(\tau_{\frak K})$, or even first order logic
as every $N \in
K_{\aleph_0}$ is atomic).  Now also $\bar b \char 94 \bar a_1,\bar c
\char 94 \bar a_1$ materialize the same $\Bbb L_{\omega_1,\omega}(N_0)$-type
in $N_1$ hence they realize the same 
$\Bbb L_{\omega_1,\omega}(\tau_{\frak K})$-type
(recall \scite{88r-5.2}(8)).  Hence $\bar b,\bar c$ realize the same
$\Bbb L_{\omega_1,\omega}(\tau_{\frak K})$-type over $\bar a_1$ in 
$N_1$.  But gtp$(\bar a_2,N_0,N_2)$ does not split over $\bar a_1$, so
by the previous sentence we get that
$\bar b \char 94 \bar a_2,\bar c \char 94 \bar a_2$ materializes the same
$\Bbb L_{\omega_1,\omega}(N_0)$-type in $N_2$.  
\hfill$\square_{\scite{88r-5.16}}$
\enddemo
\bigskip

\proclaim{\stag{88r-5.12} Lemma}  1) Suppose $N_0 \le_{\frak K} N_1 \in
K_{\aleph_0},p_\ell \in \bold D(N_\ell)$ and $p_1$ is a stationarization
of $p_0$ over $N_1$, \ub{then} $p_1 \models p_0$, i.e., every sequence
materializing $p_1$ materializes $p_0$ in any $N_2$ such that $N_1
\le_{\frak K} N_2$.
\endproclaim
\bigskip

\remark{Remark}  1) In \cite{Sh:48}, \cite{Sh:87a}, \cite{Sh:87b} 
and \cite{Sh:c} the
parallel proof of the claims were totally trivial, but here we need to
invoke $\dot I(\aleph_1,K) < 2^{\aleph_1}$. \nl
2) A particular case can be proved in the context of \S4.
\endremark
\bn

\demo{Proof}  So suppose $N_0,N_1,p_0,p_1$ contradict the claim and let
$\bar a^* \in N_0$ be such that $p_0$ is definable over $\bar a^*$ so
$p_1$, too.
By \scite{88r-5.6}(e)+(f) there are $\delta < \omega_1$ and 
$N_2 \in K_{\aleph_0}$ satisfying $N_1 \le_{\frak K} N_2$
such that $N_2$ is $(\bold D^*_\delta(N_\ell),\aleph_0)^*$-homogeneous for
$\ell=0,1$.  We can
find $p_2 \in \bold D(N_2)$ which is the stationarization of $p_0,p_1$.
It is enough to prove that $p_2 \models p_1$. \nl
[Why?  First, note that there is an automorphism $f$ of $N_2$ which
maps $N_1$ onto $N_0$ and $f(\bar a^*) = \bar a^*$ hence $f(p_2) =
p_2,f(p_1) = p_0$ hence $p_2 \models p_0$.  
Now assume that $N_1 \le_{\frak K} N^+_1 \in K_{\aleph_0}$ and $\bar a_1
\in N^+_1$ materializes $p_1$ clearly 
we can find $N^+_2,\bar a_2$ such that $N_2
\le_{\frak K} N^+_2 \in K_{\aleph_0}$ and $\bar a_2 \in N^+_2$ which
materializes $p_2$, as we
are assuming $p_2 \models p_1$ there are $N_3,f$ such that $N^+_1
\le_{\frak K} N_3 \in K_{\aleph_0}$ and $f$ is a $\le_{\frak K}$-embedding
of $N^+_2$ into $N_3$ over $N_1$ mapping $\bar a_2$ to $\bar a_1$.
But $p_2 \models p_0$ (see above) hence $f(\bar a_2) = \bar a_1$ 
materializes $p_0$ and $p_1$, too.]

So \wilog \, for some $\delta$
\mr
\item "{$\circledast$}"  $N_1$ is
$(\bold D^*_\delta(N_0),\aleph_0)^*$-homogeneous.
\ermn
For $N \in K_{\aleph_0},N_0 \le_{\frak K} N$, let 
$p_N$ be the stationarization of $p$ over $N$; without loss of
generality the universes
of $N_0,N_1$ are $\omega,\omega \times 2$ respectively.
Now we choose by induction on $\alpha$ a model $N_\alpha \in
K_{\aleph_0} \, (\alpha < \omega_1),|N_\alpha| = \omega (1 +
\alpha),[\beta < \alpha \Rightarrow N_\beta \le_{\frak K} N_\alpha]$;
$N_0,N_1$ are the ones mentioned in the claim and $\bar a_\alpha \in
N_{\alpha +1}$ materializes the stationarization $p_\alpha \in
\bold D^*_\delta(N_\alpha)$ of $p_0$ over $N_\alpha$ and for $\alpha <
\beta,N_\beta$ is $(\bold D^*_\delta(N_\alpha),\aleph_0)$-homogeneous (see
\scite{88r-5.6}(f)).  Recalling that ${\frak K}$ is categorical in
$\aleph_0$ (and uniqueness of $(\bold D_\delta(N_0),\aleph_0)^*$-homogeneous)
we have  $\alpha > \beta \Rightarrow (N_\alpha,N_\beta) \cong
(N_1,N_0)$ so clearly $\bar a_\alpha$ does not materialize $p_{N_\beta}$ (in
$N_{\alpha +1}$).  Let $N = \cup\{N_\alpha:\alpha < \omega_1\}$.
Let ${\frak B}$ be $({\Cal H}(\aleph_2),\in)$ expanded by $N,K \cap
{\Cal H}(\aleph_2),\le_{\frak K} \restriction {\Cal H}(\aleph_2)$ and 
anything else which is necessary.
Let ${\frak B}^-$ be a countable elementary submodel of ${\frak B}$ to
which $\langle N_\alpha:\alpha < \omega_1 \rangle,N$ belong and let
$\delta(*) = {\frak B}^- \cap \omega_1$.
For any stationary co-stationary 
$S \subseteq \omega_1$, let ${\frak B}_S$ be a model which is
\mr
\item "{$(\alpha)$}"  ${\frak B}_S$ an elementary extension of ${\frak B}^-$
\sn
\item "{$(\beta)$}"  ${\frak B}_S$ is an end extension of ${\frak
B}^-$ for $\omega_1$,
that is, if ${\frak B}_S \models ``s < t$ are countable ordinals" and $t
\in {\frak B}^-$ then $s \in {\frak B}^-$
\sn
\item "{$(\gamma)$}"  among the ${\frak B}_S$-countable ordinals not
in ${\frak B}^-$ there is no first one
\sn
\item "{$(\delta)$}"  ``the set of countable 
ordinals" of ${\frak B}_S$ is $I_S,I_S = \dbcu_{\alpha <
\omega_1} I^S_\alpha$, even $I^S_0$ is not well ordered, each $I_\alpha$ a
countable initial segment of $I_S,\alpha < \beta \Rightarrow I^S_\alpha
\subseteq I^S_\beta \wedge I^S_\alpha \ne I^S_\beta$ 
\sn
\item "{$(\varepsilon)$}"  $I_S \backslash I^S_\alpha$ has a first 
element if and only if $\alpha \in S$ and then we call it
$s(\alpha)$.  
\ermn
In particular $\omega$ and finite sets are standard in
${\frak B}_S$.  For $s \in I_S,N_s[{\frak B}_s] =: 
N^{{\frak B}_S}_s$ is defined
naturally, and so is $N^S = N^{{\frak B}_S}$; clearly $N^S_s \in
K_{\aleph_0}$ is $\le_{\frak K}$-increasing with $s \in I$ as those
definitions are $\Sigma^1_1$ (as ${\frak K}$ is PC$_{\aleph_0}$).  
Let $N^S_\alpha = \dbcu_{s \in I_\alpha} N^{{\frak B}_S}_s$ and let
$s+1$ be the successor of $s$ in $I_S$.
So if ${\frak B}_S \models ``s < t$ are countable ordinals", 
then $(N^{{\frak B}_S}_t,N^{{\frak B}_S}_s)$ is $(\bold
D^*_\delta(N^{{\frak B}_S}_s),\aleph_0)^*$-homogeneous. 

If $\alpha \in S$ then clearly the type $p = p_{N^S_\alpha}$ satisfies
(using absoluteness from ${\frak B}_S$ because $N^S_\alpha$ is
definable in ${\frak B}_S$ as $N^{{\frak B}_S}_{s(\alpha)}$):
\mr
\item "{$(a)$}"  $p$ is materialized in $N^S$ (i.e. in $N^S_\beta$ for
a club of $\beta \in S$) 
\ermn
but by the assumption toward contradiction
\mr
\item "{$(b)$}"  for a closed unbounded $E \subseteq \omega_1$ for no
$\beta \in E \cap S,\beta > \alpha$ does a sequence from $N^S$ materialize
both $p$ and its stationarization over $N^S_\beta$ (again remember
$N^S_\alpha = N^{{\frak B}_S}_{s(\alpha)}$ because $\alpha \in S$) \nl
and similarly
\sn
\item "{$(c)$}"   for a closed unbounded set of $\beta >
\alpha,N^S_\beta$ is $(\bold D^*_\delta(N^S_\alpha),\aleph_0)^*$-homogeneous.
\ermn
We shall prove that every $\alpha < \omega_1$,
\mr
\item "{$\boxdot$}"   if $\alpha \notin S$
then $\alpha$ cannot satisfy the statements (a),(b),(c) above.
\ermn
This is sufficient because if $S(1),S(2) \subseteq \omega_1$,f is an 
isomorphism from $N^{S(1)}$ onto $N^{S(2)}$ mapping $\bar a^*$ to
itself, then for a closed
unbounded set $E \subseteq \omega_1$, for each $\alpha < \omega_1$ 
the function $f$ maps $N^{S(1)}_\alpha$
onto $N^{S(2)}_\alpha$, hence the property above is preserved, hence
$S(1) \cap E = S(2) \cap E$.  But there is a sequence
$\langle S_i:i < 2^{\aleph_1} \rangle$ of subsets of $\omega_1$ 
such that for $i \ne j$ the set $S_i \backslash S_j$ is stationary.
So by \scite{88r-0.9} we have $\dot I(\aleph_1,K) = 2^{\aleph_1}$, contradiction.

So suppose $\alpha \in \omega_1 \backslash S,p = p_{N^S_\alpha}$ and
clauses (a), (b), (c) above hold; let $\bar a \in N^S$ materialize
$p$ in $N^S$ and we shall get a contradiction.

There are elements $0 = t(0) < t(1) < \ldots < t(k)$ of $I^S$ and 
$\bar a_0 \in N_0 = N^{{\frak B}_S}_{t(0)},\bar a_{\ell +1} \in 
N^{{\frak B}_S}_{t(\ell)+1}$ such that
$\bar a \subseteq \bar a_k,\bar a^* \subseteq \bar a_0,
\bar a_\ell \subseteq \bar a_{\ell +1}$ and
gtp$(\bar a_{\ell +1},N^{{\frak B}_S}_{t(\ell)},
N^{{\frak B}_S}_{t(\ell +1)})$ is definable over $\bar a_\ell$ and if
$t(\ell +1)$ is a successor (in $I_S$) then it is the successor of
$t(\ell)$ and if limit in $I^S$ then $\bar a_\ell = \bar a_{\ell+1}$. \nl
[Why do they exist?  Because of the sentence saying that
for every $\bar a$ we can find such $k,t(\ell)(\ell \le k),\bar
a_\ell(\ell \le k)$ as above is satisfied by ${\frak B}$ and involve
parameters which belong to ${\frak B}^-$ hence to ${\frak B}_S$, etc.,
so ${\frak B}_S$ inherits it (and finiteness is absolute from 
${\frak B}_S$)].
It follows that gtp$(\bar a,N^{{\frak B}_S}_{t(\ell)},
N^{{\frak B}_S}_{t(k)})$ is definable over $\bar a_\ell$.

We would like to show that $N^{{\frak B}_S}_{t(\ell+1)}$ is
$(D^*_\delta(N^{{\frak B}_S}_\alpha),\aleph_0)^*$-homogeneous; now if
we replace $N^{{\frak B}_S}_\alpha$ by $N^{{\frak B}_S}_s,t(\ell)
\le_I s \in I_\alpha$ we know it by the choice of ${\frak B}_S$; and
we shall use it.

Clearly $t(0) = 0 \in I_\alpha$ but $t(k) \notin I_\alpha$
(otherwise $t(k) +1 \in I_\alpha$ hence $\bar a \in
N^{{\frak B}_S}_{t(k)+1} \le_{\frak K} N^S_\alpha$, impossible as $p$
is a non-algebraic type over $N^{{\frak B}_S}_\alpha$).  Hence for 
some $\ell,t(\ell) \in I_\alpha,t(\ell +1) \notin I_\alpha$.  
By the construction $t(\ell +1)$ is limit hence
$\bar a_{\ell +1} = \bar a_\ell$.
As $\alpha \notin S$ we can choose $t(*) \in I_S \backslash I^S_\alpha,
t(*) < t(\ell +1)$.  
As we are assuming (toward contradiction) that 
$\alpha,p$ satisfy clause (c), for some $\beta \in S,s(\beta)$ is
well defined $s(\beta) > t(k)$ (on the definition of $s(\gamma)$ for
$\gamma \in S$ see clause $(\varepsilon)$ above) and $N^S_\beta$ is
$(\bold D^*_\delta(N^S_\alpha),\aleph_0)^*$-homogeneous.  Now
$N^{{\frak B}_S}_\beta,N^{{\frak B}_S}_{t(\ell +1)}$ are 
isomorphic over $N_{t(*)}$ (being $(\bold D^*_\delta
(N^{{\frak B}_S}_{t(*)}),\aleph_0)^*$-homogeneous by the choice of
${\frak B}_S$).  Hence also $N^{{\frak B}_S}_{t(\ell +1)}$ is
$(\bold D^*_\delta(N^S_\alpha),\aleph_0)^*$-homogeneous, too, hence
$(N^{{\frak B}_S}_{t(\ell +1)},N^S_\alpha,\bar a^*) \cong (N_1,N_0,\bar a^*)$.

As clearly $N^S_\alpha,N^{{\frak B}_S}_{t(*)}$ are
$(\bold D^*_\delta(N^{{\frak B}_S}_{t(\ell)+1)}),
\aleph_0)^*$-homogeneous  there is an
isomorphism $f_0$ from $N^S_\alpha$ onto $N^{{\frak B}_S}_{t(*)}$ over
$N^{{\frak B}_S}_{t(\ell)+1}$.  As $N^{{\frak B}_S}_{t(\ell +1)}$ is
$(\bold D^*_\delta(N^{{\frak B}_S}_{t(*)}),\aleph_0)^*$-homogeneous and
$(\bold D^*_\delta(N^S_\alpha),\aleph_0)^*$-homogeneous by the
previous paragraph (where we use $\beta$) we can extend $f_0$
to an automorphism $f_1$ of $N^{{\frak B}_S}_{t(\ell +1)}$.  Let $\gamma \in S
\cap E$ satisfy
$s(\gamma) \ge t(k) +1$.  As gtp$(\bar a_k,N^{{\frak B}_S}_{t(\ell+1)},
N^S_\gamma)$ is
definable over $\bar a_\ell$ and $\bar a_\ell = f_0(\bar a_\ell) =
f_1(\bar a_\ell)$ (as $\bar a_\ell \in N^{{\frak B}_S}_{t(\ell)+1})$ and
$N^S_{\gamma +1}$ is $(\bold D^*_\delta(N^{{\frak B}_S}_{t(\ell +1)}),
\aleph_0)^*$-homogeneous, we can extend $f_1$ to an automorphism
$f_2$ of $N^S_\gamma$ satisfying $f_2(\bar a_k) = \bar a_k$.

Notice that by the choice of $\langle \bar a_\ell:\ell \le k \rangle$
and $\langle t(\ell):\ell \le k \rangle$ it follows that for any
$m<k$, gtp$(\bar a_k,N_{t(m)},N_{t(k)+1})$ does not split over $\bar
a_m$ hence is definable over it by \scite{88r-5.16}, and recall that we
know that $\bar a_\ell = \bar a_{\ell+1}$.

So there is in $N^S$ a sequence materializing both gtp$(\bar
a,N^S_\alpha,N^S_\gamma) = p_{N^S_\alpha}$ and its stationarization over
$N^S_{t(\ell +1)}$: just $\bar a (\subseteq \bar a_k)$ (so use
$f_2$). \nl
This contradicts the assumption as $(N_1,N_0,\bar a^*) \cong 
(N^{{\frak B}_S}_{t(\ell +1)},N^S_\alpha,\bar a^*)$.  \hfill$\square_{\scite{88r-5.12}}$
\enddemo
\bigskip

\remark{Remark}   Imitating the proof, we can show that (c) holds 
for any $\alpha < \omega_1$.
\endremark
\bigskip

\proclaim{\stag{88r-5.13} Claim}  1) If $\bar a \in N_0 \le_{\frak K} 
N_1 \le_{\frak K} N_2 \in K_{\aleph_0},
\bar b \in N_2,p_1 = { \text{\rm gtp\/}}(\bar b,N_1,N_2)$ is
definable over $\bar a \in N_0$, \ub{then} 
$p_0 = { \text{\rm gtp\/}}(\bar b,N_0,N_2)$ is
definable in the same way over $\bar a$, hence ${\text{\rm gtp\/}}
(\bar b,N_1,N_2)$ is its stationarization. \nl
2) For a fixed countable $M \in K_{\aleph_0}$ to have a common
stationarization in $\bold D(N')$ for some $N'$ satisfying 
$M \le_{\frak K} N'$ or $N' \le_{\frak K} M$ is an 
equivalence relation over $\{p$: for some $N
\le_{\frak K} M,p \in \bold D(N)\}$ (and we can choose the common 
stationarization in $\bold D(M)$ as a representative). 
So if $N_0 \le_{\frak K} N_1 \le_{\frak K} N_2 \in K_{\aleph_0},p_\ell
\in \bold D(N_\ell)$ for $\ell=0,1,2$ and $p_1,p_2$ are
stationarizations of $p_0$ then $p_2 \models p_1$.
\nl
3) If $N_\alpha \in K_{\aleph_0} \, (\alpha \le \omega +1)$ is
$\le_{\frak K}$-increasing and continuous and $\bar a \in 
N_{\omega +1}$ \ub{then} for
some $n < \omega$, for every $k$ we have: $n < k \le \alpha \le \omega$ implies
${\text{\rm gtp\/}}(\bar a,N_\alpha,N_{\omega +1})$ is the stationarization of
${\text{\rm gtp\/}}(\bar a,N_k,N_{\omega +1})$. \nl
4) If $N \le_{\frak K} M \in K,N \in K_{\aleph_0},\bar a \in M$ 
\ub{then} for all $M' \in K_{\aleph_0}$, satisfying $\bar a \in M',
N \le_{\frak K} M' \le_{\frak K} M,{\text{\rm gtp\/}}
(\bar a,N,M')$ is the same, we call it ${\text{\rm gtp\/}}
(\bar a,N,M)$ (the new point is
that $M$ is not necessarily countable).
\nl
5) Suppose $N_0 \le_{\frak K} N_1$ (in $K$), $\bar a \in N_1$,
\ub{then} there is a countable $M \le_{\frak K} N_0$, 
such that for every countable $M'$ satisfying 
$M \le_{\frak K} M' \le_{\frak K} N_0$ we have
${\text{\rm gtp\/}}(\bar a,M',N_1)$ is the stationarization of 
${\text{\rm gtp\/}}(\bar a,M,N_1)$.
Moreover there is a finite $A \subseteq N_0$ such that any countable
$M \le_{\frak K} N_0$ which includes $A$ is O.K.
\nl
6) The parallel of Part 
(3) holds for $N_\alpha \in K$, too, and any limit ordinal
instead of $\omega$.  That is if $\langle N_\alpha:\alpha \le \delta
+1 \rangle$ is $\le_{\frak K}$-increasing continuous and $\bar a \in
N_{\delta +1}$, \ub{then} for some $\alpha < \delta$ and countable $M
\le_{\frak K} N_\alpha$ we have: $M \le_{\frak K} M' \le_{\frak K}
M_\delta \Rightarrow { \text{\rm gtp\/}}(\bar a,M',M_\delta)$ is the
stationarization of ${\text{\rm gtp\/}}(\bar a,M,M_\delta)$. \nl
7) If $N_0 \le_{\frak K} N_1 \le_{\frak K} N_2 \le_{\frak K} N_3
\le_{\frak K} N_4$ and $\bar a \in N_4$ and 
{\rm gtp}$(\bar a,N_3,N_4)$ is the stationarization of
${\text{\rm gtp\/}}(\bar a,N_0,N_4)$ then ${\text{\rm gtp\/}}
(\bar a,N_2,N_4)$ is the stationarization of 
${\text{\rm gtp\/}}(\bar a,N_1,N_3)$.
Also if $\bar b'$ satisfies {\rm Rang}$(\bar b') 
\subseteq \text{\rm Rang}(\bar a)$ and {\rm gtp}$(\bar a,N_2,N_4)$ is 
the stationarization of {\rm gtp}$(\bar a,N_1,N_4)$ 
\ub{then} this holds for $\bar b$.
\nl
8) If $N_0 \le_{\frak K} N_1 \le_{\frak K} N_2 \in K_{\aleph_0}$ and
$p_\ell \in \bold D(N_\ell)$ for $\ell =0,1,2$ and $p_{\ell +1}$ is
the stationarization of $p_\ell$ for $\ell =0,1$ \ub{then} $p_2$ is
the stationarization of $p_0$.
\nl
9) If $\langle M_\alpha:\alpha \le \delta +1\rangle$ is $\le_{\frak
K}$-increasing continuous and $\bar a \in {}^{\omega >}(M_{\delta
+1})$ \ub{then} 
\mr
\item "{$(a)$}"  for some $\alpha < \delta$ we have {\rm gtp}$(\bar
a,M_\beta,M_{\delta +1})$ is the stationarization of {\rm gtp}$(\bar
a,M_\alpha,M_{\delta +1})$ whenever $\beta \in [\alpha,\delta)$
\sn
\item "{$(b)$}"  if {\rm gtp}$(\bar a,M_\alpha,M_{\delta +1})$ is the
stationarization of {\rm gtp}$(\bar a,M_0,M_{\delta +1})$ for every $\alpha
< \delta$ then this holds for $\alpha = \delta$, too.
\endroster
\endproclaim
\bigskip

\demo{Proof}  1) As we can replace $N_2$ by any $N'_2$ satisfying $N_2
\le_{\frak K} N'_2 \in K_{\aleph_0}$,.
without loss of generality for some $\alpha,N_2$ is
$(\bold D^*_\alpha(N_0),\aleph_0)^*$-homogeneous and
$(\bold D^*_\alpha(N_1),\aleph_0)^*$-homogeneous.  
Let $p_2 \in \bold D(N_2)$ be the stationarization of $p_1$ over $N_2$.

So by
\scite{88r-5.12} we get $p_2 \models p_1$.  On the other hand, clearly there is an
isomorphism $f_0$ from $N_0$ onto $N_1$ such that 
$f_0(\bar a) = \bar a$; and by
the assumption above on $N_2,f_0$ can be extended to an automorphism
$f_1$ of $N_2$.

Note that $f_1$ maps $p_0 = \text{\rm gtp}(\bar b,N_0,N_2)$ to 
$p'_0 =: \text{ gtp}
(f_1(\bar b),f_1(N_0),N_2)$ and maps $p_2$ to itself as $f_0(\bar a) =
\bar a$.

Now $p_1 \models p_0$ (by the choices of $p_1,p_0$) and $p_2 \models
p_1$ by \scite{88r-5.5}(1), together $p_2 \models p_0$.  As $f_1(p_2) =
p_2,f_1(p_0) = p'_0$ it follows that $p_2 \models p'_0$.  As also $p_2
\models p_1$ and $p'_0,p_1 \in \bold D(N_1)$ it follows that $p'_0 =
p_1$ hence $p_1,p'_0$ have the same definition over $\bar a$, but now
also $p_0 \in \bold D(N_0),p'_0 \in \bold D(N_1)$ have the same
definition over $\bar a$ (using $f_1$), together also $p_1,p_0$ have
the same definition over $\bar a$, which means that $p_1$ is the
stationarization of $p_0$ over $N_1$ and we are done.
\nl
2) Trivial. \nl
3) By part (1). \nl
4) Easy. \nl
5) By (3) and (4). \nl
6),7),8),9)  Easy by now.  \hfill$\square_{\scite{88r-5.13}}$
\enddemo
\bigskip

\definition{\stag{88r-5.14} Definition}  By \scite{88r-5.13}(5) gtp$(\bar a,M,N)$
can be reasonably 
defined when $M \le_{\frak K} N,\bar a \in {}^{\omega >} N$ and 
we can define
$\bold D(N)$, gtp$(\bar a,N,M)$ and stationarization for not necessarily
countable $N$ and $N \le_{\frak K} M \in K$.  
Everything still holds, except that maybe some
$p$'s are not materialized in any $\le_{\frak K}$-extension of $N$.

More formally
\mr
\item "{$(a)$}"  if $N \le_{\frak K} M$ and $\bar a \in {}^{\omega >}
M$ then gtp$(\bar a,N,M)$ is defined as $\cup\{\text{\rm gtp}(\bar a,N',M'):N_0
\le_{\frak K} N' \le_{\frak K} M' \in K_{\aleph_0},M' \le_{\frak K}
M,N' \le_{\frak K} N\}$ for every countable large enough $N_0
\le_{\frak K} N$; it is well defined by \scite{88r-5.13}(5) and we say $\bar a$
materializes gtp$(\bar a,N,M)$ in $M$
\sn 
\item "{$(b)$}"  if $N \in {\frak K},N \le_{\frak K} M$ and $p \in \bold
D(N)$ is definable over the countable $N_0 \le_{\frak K} N$, \ub{then}
the stationarization of $p$ over $M$ is $\cup\{p_{M_0}:N_0 \le_{\frak K} 
M_0 \le_{\frak K} M,M_0$ is countable$\}$ where
$p_{M_0}$ is the stationarization of $p \in \bold D(M)$ over $M_0$
\sn
\item "{$(c)$}"  if $N \in {\frak K}$ let $\bold D(N) = \{q$: for some
countable $M \le_{\frak K} N$ and $p \in \bold D(M),q$ is well defined
by \scite{88r-5.15}(2) and is the stationarization of $p$ over $N$ so is
definable over some $\bar a \subseteq M\}$
\sn
\item "{$(d)$}"  if $N \in K_{\aleph_0},M \in K,N \le_{\frak K} M,p
\in \bold D(N)$ \ub{then} the stationarization of $p$ over $M$ is 
defined as in \scite{88r-5.11}
\sn
\item "{$(e)$}"  those definitions are compatible with the ones for
countable model 
\sn
\item "{$(f)$}"  gtp$(\bar a,N,M)$ (where $\bar a \in M,N \le_{\frak K} M$ are both is
$K$) is the stationarization over $N$ of gtp$(\bar a,N',M)$ for every
large enough countable $N' \le_{\frak K} N$, see \scite{88r-5.13}(5).
\endroster
\enddefinition
\bigskip

\remark{Remark}  1) So easily the parallel of \scite{88r-5.13} holds for
not necessarily countable models. \nl
2) Claim \scite{88r-5.15} below strengthens \scite{88r-3.5}, it is a step 
toward non-forking amalgamation.
\endremark
\bigskip

\proclaim{\stag{88r-5.15} Claim}  Suppose $N_0 \le_{\frak K} N_1 \in
K_{\aleph_0},N_0 \le_{\frak K} N_2 \in K_{\aleph_0},\bar a \in N_1$.  
\ub{Then} we can find $M,N_0 \le_{\frak K} M \in K_{\aleph_0}$ and 
$\le_{\frak K}$-embeddings $f_\ell$ of $N_\ell$ into $M$ over
$N_0$ (for $\ell = 1,2$)  such that ${\text{\rm gtp\/}}
(f_1(\bar a),f_2(N_2),M)$ is a stationarization of $p_0 = 
{ \text{\rm gtp\/}}(\bar a,N_0,N_1)$ (so $f_1(\bar a) \notin f_2(N_2))$.
\endproclaim
\bigskip

\demo{Proof}  Let $p_2 \in \bold D(N_2)$ be the stationarization of $p_0$.
Clearly we can find an $\alpha < \omega_1$ (in fact, a closed
unbounded set of $\alpha$'s) and $N'_1,N'_2$ from $K_{\aleph_0}$ which
are $(D^*_\alpha(N_0),\aleph_0)^*$-homogeneous and $N_\ell \le_{\frak K} 
N'_\ell$ (for $\ell = 1,2)$ and some $\bar b \in N'_2$ materializing 
$p_2$.  But by
\scite{88r-5.12}, $\bar b$ materializes $p_0$ hence there is an
isomorphism $f$ from $N'_1$ onto $N'_2$ over $N_0$ satisfying
$f(\bar a) = \bar b$.
Now let $M = N'_2,f_1 = f \restriction N_1,f_2 = \text{ id}$. \nl
${{}}$  \hfill$\square_{\scite{88r-5.15}}$
\enddemo
\bigskip

\proclaim{\stag{88r-5.17} Claim}  1) For any $N_0 \le_{\frak K} 
N_1 \in K_{\aleph_1}$ so $N_0 \in K_{\le \aleph_1}$,
there is $N_2$ such that $N_1 \le_{\frak K} N_2 \in K_{\aleph_1}$ and $N_2$ is
$(\bold D(N_0),\aleph_0)^*$-homogeneous. \nl
2) Also \scite{88r-5.15} holds for $N_2 \in K_{\aleph_1}$ (but still $N_0,N_1
\in K_{\aleph_0}$). \nl
3) If $N_0 \le_{\frak K} N_1 \in K_{\aleph_0}$ and $N_0 \le_{\frak K}
 N_2 \in K_{\le \aleph_1}$ \ub{then} we can find $M \in K_{\le \aleph_1}$
 and $\le_{\frak K}$-embeddings $f_1,f_2$ of $N_1,N_2$ into $M$ over
 $N_0$ respectively such that 
{\rm gtp}$(f_1(\bar c),f_2(N_2),M)$ is a stationarization of 
{\rm gtp}$(\bar c,N_0,N_1)$ for every $\bar c \in N_1$, hence 
$f_1(N_1) \cap f_2(N_2) = N_0$. \nl
4) $K_{\aleph_2} \ne \emptyset$.
\endproclaim
\bigskip

\remark{Remark}  1) Note that \scite{88r-5.17}(3) is another step toward
stable amalgamation.
\nl
2) Note that \scite{88r-5.17}(3) strengthen \scite{88r-5.17}(2) hence \scite{88r-5.15}.
\endremark
\bigskip

\demo{Proof}  1) As we can iterate, it is 
enough to prove that: if $p(\bar x,\bar y) \in
\bold D(N_0)$ and 
$\bar a \in N_1$ materializes $p(\bar x,\bar y) \restriction \bar
x$ in $(N_1,N_0)$ then for some $N_2 \in K_{\aleph_1},N_1 \le_{\frak
K} N_2$ and
for some $\bar b \in N_2$ the sequence 
$\bar a \char 94 \bar b$ materializes $p(\bar
x,\bar y)$ in $(N_2,N_0)$.  Let $M_0 \le_{\frak K} N_0$ be countable and $q \in
\bold D(M_0)$ be such that $p(\bar x,\bar y)$ a stationarization of
$q$.  Without loss of generality if $N_0$ is countable then $M_0 =
N_0$.  Note that the case $N_0 = M_0$ is easier.
Choose $M_i(0 < i < \omega_1)$ such that $M_i \le_{\frak K} N_1,N_1 
= \dbcu_{i < \omega_1} M_i,\langle M_i:i < \omega_1 \rangle$ is 
$\le_{\frak K}$-increasing continuous sequence of countable 
models, $M_0 \cup \bar a \subseteq M_1$.
As $\langle M_i \cap N_0:i < \omega_1 \rangle$ is an increasing
continuous sequence of countable sets with union $N_0$ clearly for a
club of $i < \omega_1,M_i \cap N_0 \le_{\frak K} N_0$ hence $M_i \cap
N_0 \le_{\frak K} M_i$.  So \wilog \, $i < \omega_1 \Rightarrow M_i \cap N_0
\le_{\frak K} N_0,M_i$.  For every $\bar c \in N_1$ there is a
countable $N_{0,\bar c}$ such that $M_0 \le_{\frak K} N_{0,\bar c}
\le_{\frak K} N_0$ and: if $N_{0,\bar c} \le_{\frak K} N' \le_{\frak
K} N_0$ and $N' \in K_{\aleph_0}$ then  gtp$(\bar c,N',N_1)$ is the
stationarization of gtp$(\bar c,N_{0,\bar c},N_1)$.  Without loss of
generality $\bar c \in M_i \Rightarrow N_{0,\bar c} \subseteq M_i$ hence
\mr
\item "{$(*)$}"  for every $\bar c \in M_i$, gtp$(\bar c,N_0,N_1)$ is
a stationarization of gtp$(\bar c,N_0 \cap M_i,M_i)$.
\ermn
We can find $M^*_1 \in K_{\aleph_0}$ satisfying 
$M_1 \le_{\frak K} M^*_1$ and $\bar b \in M^*_1$ such that $q = 
\text{ gtp}(\bar a \char 94 \bar b,M_0,M^*_1)$.
We can find $\bar a_2,\bar a_1,\bar a_0$ such that 
$\bar a_0 \in M_1 \cap N_0,\bar
a_1 \in M_1,\bar a_2 \in M^*_1,\bar b \subseteq \bar a_2,\bar a
\subseteq \bar a_1$ and $\bar a_0 \trianglelefteq \bar a_1
\trianglelefteq \bar a_2$ and gtp$(\bar a_2,M_1,M^*_1)$, gtp$(\bar a_1,M_1 \cap
N_0,M_1)$ are definable over $\bar a_1,\bar a_0$, respectively.  Now
we define $f_i,M^*_i,1 < i < \omega_1$ by induction on $i$ such that:
\mr
\widestnumber\item{$(iii)$}
\item "{$(i)$}"   $\langle M^*_j:1 \le j \le i \rangle$ is
$\le_{\frak K}$-increasing continuous
\sn
\item "{$(ii)$}"  $M^*_j$ is countable
\sn
\item "{$(iii)$}"   $f_j$ is a $\le_{\frak K}$-embedding of $M_j$ into
$M^*_j$
\sn
\item "{$(iv)$}"   $f_j$ is the identity on $M_1$
\sn
\item "{$(v)$}"  $f_j$ is increasing continuous with $j$
\sn
\item "{$(vi)$}"  gtp$(\bar a_2,M_j,M^*_j)$ is the stationarization of
gtp$(\bar a_2,M_1,M^*_1)$ (so definable over $\bar a_1$).
\ermn
For $j=1$ we have it letting $f^*_j = \text{ id}_{M_1}$. \nl
For $j>1$ successor, use \scite{88r-5.15} to define $(M_j,f_j)$ such that
grp$(\bar a_2,f_j(M_j),M^*_j)$ is the stationarization of gtp$(\bar
a_2,f_{j-1}(M_{j-1}),M^*_{j-1})$.  So clauses (i)-(v) clearly holds.
Clause (vi) follows by \scite{88r-5.13}(8).
\nl
For $j$ limit: let $M^*_j = \dbcu_{1 \le i < j} M^*_i$
and $f_j = \cup\{f_i:1 \le i < j\}$, condition (v)
holds by \scite{88r-5.13}(3).

By renaming \wilog \, $f_j = \text{ id}_{M_j}$ for $j \in
[1,\omega_1)$. \nl
By $(*)$ we get that gtp$(\bar a_1,N_0 \cap M_j,M^*_j) = \text{\rm
gtp}(\bar a_1,N_0 \cap M_j,M_j)$ is definable over $\bar a_0$ (as this
holds for $j=1$).  Combining this and clause (iv) by \scite{88r-5.16} we get
for every $j \ge 1$, that gtp$(\bar a_2,N_0 \cap
M_j,M^*_j)$ is the stationarization of gtp$(\bar a_2,N_0 \cap
M_1,M^*_1)$.  Hence by the choice of $\bar a_2,\bar a_1,a_0$ and
\scite{88r-5.13}(7), easily gtp$(\bar a \char 94 \bar b,N_0 \cap
M_j,M^*_j)$ is the stationarization of gtp$(\bar a \char 94 \bar b,N_0
\cap M_1,M^*_1)$.

So by \scite{88r-5.14}, clause (c) 
and the first sentence in the proof, we finish. \nl
2) Similar proof \footnote{here $N_1 \in K_{\aleph_1}$ is O.K.;
similar to \scite{88r-2.7}(1)} (or use the proof of part (3)). \nl
3) Without loss of generality $N_2 \cong N^*$ from \scite{88r-5.9} (as we
can replace $N_2$ by an extension so use \scite{88r-5.10}).

Also (by \scite{88r-5.17}(1)) there is $M,N_2 \le_{\frak K} M 
\in K_{\aleph_1}$ such that $M$ is $(\bold
D(N_2),\aleph_0)^*$-homogeneous.  
As $N_1$ is countable there is 
$\alpha < \omega_1$ such that for every $\bar c \in N_1$, gtp$(\bar
c,N_0,N_1) \in \bold D_\alpha(N_0)$.  Let $M = \dbcu_{i < \omega_1} M_i$
with $M_i \in K_{\aleph_0}$ being
$\le_{\frak K}$-increasing continuous.  So for some $i$ we have
$\alpha < i < \omega_1,M_i \cap N_2 \le_{\frak K} M$ 
and (recalling \scite{88r-5.13}(6)) for every $\bar c \in M_i$, gtp$(\bar
c,N_2,M)$ is stationarization of gtp$(\bar c,N_2 \cap M_i,M_i)$ and
$M_i$ is $(D_i(N_2 \cap M_i),\aleph_0)^*$-homogeneous.  Now we can
find an isomorphism $f_0$ from $N_0$ onto $N_2 \cap M_i$ (as $K$ is
$\aleph_0$-categorical) and extend it to an automorphism $f_2$ of
$N_2$ (by \scite{88r-5.10}-model homogeneity).  Also there is
$N'_1$ such that $N_1 \le_{\frak K} N'_1 \in K_{\aleph_0}$ and $N'_1$ is
$(\bold D_i(N_1),\aleph_0)^*$-homogeneous, hence is
$(\bold D_i(N_0),\aleph_0)^*$-homogeneous (by the choice of $\alpha$ as
$\alpha < i$ see \scite{88r-5.6}(f)), hence there
is an isomorphism $f'_1$ from $N'_1$ onto $M_i$ extending $f_0$.  Now
$f_0,f'_1 \restriction N_1,f_2,M$ show that amalgamation as required
exists (we just change names). \nl
4) Immediate, use 1) or 2) or 3) $\omega_2$-times.
\hfill$\square_{\scite{88r-5.17}}$ 
\enddemo
\bigskip

\definition{\stag{88r-5.18} Definition}  For any $\bold D_* = 
\bold D_\alpha$ for some $\alpha <
\omega_1$ (or just any reasonable such $\bold D_*$, i.e., satisfies the
demands on each $\bold D_\alpha$ in \scite{88r-5.6}) we define: \nl
1) $M \le_{\bold D_*} N$ if $M \le_{\frak K} N$ and for every $\bar a \in N$

$$
\text{gtp}(\bar a,M,N) \in \bold D_* (M).
$$
\mn
2) $K_{\bold D_*}$ is the class of $M \in K$ which are the union of a family
of countable submodels, which is directed by $\le_{{\bold D}_*}$. \nl
3) ${\frak K}_{\bold D_*} = (K_{{\bold D}_*} \le_{{\bold D}_*})$, or
pedantically $(K_{\bold D_*},\le_{\bold D_*} \restriction K_{\bold D_*})$.
\enddefinition
\bigskip

\proclaim{\stag{88r-5.19} Claim}  1) The pair $(K_{{\bold D}_*},
\le_{{\bold D}_*})$ is an
$\aleph_0$-presentable a.e.c., that is it 
satisfies all the axioms from \scite{88r-1.2}(1) and is {\rm PC}$_{\aleph_0}$. \nl
2) Also for $(K_{{\bold D}_*},\le_{{\bold D}_*})$, we get $\bold D(N)$
countable and equal to $\bold D_*(N)$ for every countable 
$N \in K_{{\bold D}_*}$.
\endproclaim
\bigskip

\demo{Proof}  1) Obviously $K_{{\bold D}_*}$ is a class of $\tau$-models and
$\le_{{\bold D}_*}$ is a two-place relation on $K_{D_*}$; also they are
preserved by isomorphisms.  About being PC$_{\aleph_0}$ note that
\mr
\item "{$\circledast_1$}"  $M \in K_{{\bold D}_*}$ iff $M \in K$ and for some
model ${\frak B}$ with universe $|M|$ and countable vocabulary, for
every countable ${\frak B}_1 \subseteq {\frak B}_2 \subseteq {\frak
B}$ we have $M \restriction {\frak B}_1 \le_{{\bold D}_*} M \restriction
{\frak B}_2$ \ub{iff} there is a directed partial order and $\langle
M_t:t \in I\rangle$ such that $M_t \in K_{\aleph_0}$ and $s <_I t
\Rightarrow M_s \le_{\frak K} M_t$ and $\bar a \subseteq M_t
\Rightarrow \text{\rm gtp}(\bar a,M_s,M_t) \in \bold D_*(M_s)$
\sn
\item "{$\circledast_2$}"  similarly for $M \le_{{\bold D}_*} N$.
\endroster
\enddemo
\bn
\ub{Ax I}:  If $M \le_{{\bold D}_*} N$ then $M \le_{\frak K} N$ hence $M
\subseteq N$.
\bn
\ub{Ax II}:  The transitivity of $\le_{\bold D_*}$ holds by 
\scite{88r-5.16} + Definition \scite{88r-5.14} (works as $\bold D_*$ is closed
enough or use clause (f) of \scite{88r-5.6}).  
The $M \le_{\bold D_*} M$ is trivial \footnote{recall that $M
\restriction {\frak B} = M \restriction \{a \in M:a \in {\frak B}\}$}.
\bn
\ub{Ax III}:  Assume $\langle M_i:i < \lambda \rangle$ is 
$\le_{{\bold D}_*}$-increasing continuous and $M = \cup\{M_i:
i < \lambda\}$.  As
${\frak K}$ is an a.e.c. clearly $M \in K$ and $i < \lambda
\Rightarrow M_i \le_{\frak K} M$.  Also for each $i < \lambda$ and
$\bar a \in M$ for some $j \in (i,\lambda)$ we have $\bar a \in M_j$
hence gtp$(\bar a,M_i,M_j) \in \bold D_*(M_i)$ hence gtp$(\bar
a,M_i,M) = \text{ gtp}(\bar a,M_i,M_j) \in \bold D_*(M_i)$.  So $i <
\lambda \Rightarrow M_i \le_{{\bold D}_*} M$.  By applying
$\circledast_1$ to every $M_i$ and coding we can easily show that $M \in
K_{{\bold D}_*}$ thus finishing.
\bn
\ub{Ax IV}:  Assume $\langle M_i:i < \lambda \rangle,M$ are as above
and $i < \lambda \Rightarrow M_i \le_{{\bold D}_*} N$.  To prove
$M \le_{{\bold D}_*} N$ note that as ${\frak K}$ is an a.e.c., we have
$M \le_{\frak K} N$ and consider $\bar a \in N$.  By \scite{88r-5.13}(6)
for some $i < \lambda$, gtp$(\bar a,M,N)$ is the stationarization of
gtp$(\bar a,M_i,N)$ but the latter belongs to $\bold D_*(M_i)$ hence
$\text{ gtp}(\bar a,M,N) \in \bold D_*(M)$ as required.
\bn
\ub{Ax V}:  By $\circledast_2$ this is translated to the case
$N_0,N_1,M \in K_{\aleph_0}$ but then it holds easily.
\bn
\ub{Ax VI}:  By $\circledast_1 + \circledast_2$ + Ax VI for ${\frak K}$.
\nl
2) So we replace ${\frak K}$ by ${\frak K}' = {\frak K}_{\bold D_*}$
and easily all that we need for $\bold D$ for ${\frak K}'$ is
satisfied by $\bold D_*$ (actually repeating the works in \S5 till now
on ${\frak K}'$ we get it).   \hfill$\square_{\scite{88r-5.19}}$ 
\bigskip

\proclaim{\stag{88r-5.20} Claim}  Suppose $N_0 \le_{\frak K} N_\ell 
\in K_{\aleph_0} \, (\ell = 1,2)$ and 
$\bar c \in N_2$, \ub{then} there is $M,N_0 \le_{\frak K} M$ and
$\le_{\frak K}$-embeddings $f_\ell$ of $N_\ell$ into $M$ over $N_0$ such that
\mr
\item "{$(i)$}"  for every $\bar a \in N_1,{\text{\rm gtp\/}}(f_1(\bar
a),f_2(N_2),M)$ is a stationarization of ${\text{\rm gtp\/}}(\bar a,N_0,N_1)$
\sn
\item "{$(ii)$}"  ${\text{\rm gtp\/}}(f_2(\bar c),
f_1(N_1),M)$ is a stationarization of ${\text{\rm gtp\/}}(\bar c,N_0,N_2)$.
\endroster
\endproclaim
\bigskip

\remark{Remark}  This is one more step toward stable amalgamation: in
\scite{88r-5.15} we have gotten it for one $\bar a \in N_1$, in \scite{88r-5.17}(3)
for every $\bar a \in N_1$, which gives disjoint amalgamation.
\endremark
\bigskip

\demo{Proof}   Clearly we can for $\ell=1,2$ replace $N_\ell$ by any
$N'_\ell,N_\ell \le_{\frak K} N'_\ell \in K_{\aleph_0}$, 
and \wilog \, $N_0 = N_1
\cap N_2$.  By \scite{88r-5.17}(3) there is $N_3 \in K_{\aleph_0}$ such
that $N_\ell \le_{\frak K} N_3$ for $\ell < 3$ and $\bar a \in
{}^{\omega >}(N_1) \Rightarrow \text{ gtp}(\bar a,N_2,N_3)$ is the
stationarization of $\text{gtp}(\bar a,N_0,N_1)$.
So we can assume that for some $\bold D_\alpha$ as in
Definition \scite{88r-5.18} we have $N_\ell$ is 
$(\bold D_\alpha(N_0),\aleph_0)^*$-homogeneous for $\ell = 1,2$.
As in the proof of \scite{88r-5.12}, we can find a countable linear order $I$,
such that every element $s \in I$ has an immediate successor $s+1,0$
is first element and $I^*$ has a subset isomorphic to the
rationals (follow really) and 
models $M_s \in K_{\aleph_0}$, (for $s \in I$)  such that $s < t
\Rightarrow M_s \le_{\frak K} M_t$ and $M_t$ is $(\bold
D_\alpha(M_s),\aleph_0)$-homogeneous when $s <_I t$, etc.  
So by \scite{88r-5.13}(3) for every initial
segment $J$ of $I$ and $t \in I$ such that $J < t$, that is, $(\forall
s \in J)(s <_I t)$  if $J$ has no last
element and $I \backslash J$ has no first element then $M_t$ is
$(\bold D_\alpha(M_J),\aleph_0)^*$-homogeneous, where $M_J = \dbcu_{s \in J}
M_s = \dbca_{t \in I \backslash J} M_t$.  We let $N^J_0 = M_J,N^J_1 = M_I$ and
$N^J_2$ be a $(\bold D_\alpha(N^J_0),\aleph_0)^*$-homogeneous model satisfying
$N^J_0 \le_{\frak K} N^J_2$ and \wilog \, 
$N^J_1 \cap N^J_2 = N^J_0$.  Also easily there is $N'_0 <_{\frak K}
N_0$ such that gtp$(\bar c,N_0,N_1)$ is definable over some $\bar c_0
\subseteq N'_0$ and $N_0$ is $(\bold D_\alpha(N'_0),\aleph_0)$-homogeneous.
Clearly the triples
$(N_0,N_1,N_2),(N^J_0,N^J_1,N^J_2)$ are isomorphic and let
$f^J_0,f^J_1,f^J_2$ be appropriate isomorphisms 
such that $f^J_0 \subseteq f^J_1,f^J_2$ and \wilog \, $f^J_0(N'_0) = M_0$.  
Now by \scite{88r-5.17}(3), there is 
$M^J \in K_{\aleph_0}$ satisfying $N^J_\ell \le_{\frak K} M^J \, 
(\ell = 0,1,2)$ such that for 
every $\bar a \in N^J_1$, gtp$(\bar a,N^J_2,M^J)$ is the
stationarization of gtp$(\bar a,N^J_0,N^J_1)$ and there are $N_3 \in
K_{\aleph_0},N_\ell \le_{\frak K} N_3$ for $\ell=0,1,2$ and an
isomorphism $f^J_3 \supseteq f^J_1 \cup f^J_2$ from $N_3$ onto $M^J$.

Suppose our conclusion fails, then gtp$(f^J_2(\bar c),N^J_1,M^J)$ is not
the stationarization of gtp$(f^J_2(\bar c),N^J_0,M^J)$.  Moreover, as
in the proof of \scite{88r-5.12}, 
$t \in I \backslash J \Rightarrow M_I = N^J_1,M_t$ are
isomorphic over $N^J_0 = M_J$, hence we can
replace $N^J_1$ by $M_t$ for any $t \in I \backslash J$ so as
we assume that our conclusion fails, $t \in I \backslash J \Rightarrow
\text{ gtp}(f^J_2(\bar c),M_t,M^J)$ is not a stationarization of
gtp$(f^J_2(\bar c),N^J_0,M^J)$ and the latter is the stationarization of
gtp$(f^J_2(\bar c),M_0,M^J)$.  Let $p_J = \text{ gtp}
(f^J_2(\bar c),N^J_1,M^J) = \text{ gtp}(\bar c,M_I,M^J)$; all this was
done for any appropriate $J$.  So it is easy
to check that $J_1 \ne J_2 \Rightarrow p_{J_1} \ne p_{J_2}$, but as
$I^* \subseteq I \and |I| = \aleph_0$, we have 
continuum many such $J$'s hence such $p_J$.  If CH fails, we are
done.  Otherwise, note that moreover, we can ensure that for 
$J_1 \ne J_2$ as above there is an 
automorphism of
$M_I$ taking $p_{J_1}$ to $p_{J_2}$, hence for some $\beta <
\omega_1,\{p_J:J$ as above$\} \subseteq \bold D_\beta(M_I)$, i.e.,
$(f^{J_2}_1) \circ (f^{J_1}_1)^{-1}$ maps one to the other, 
contradiction by clause (d) of \scite{88r-5.6} (alternatively
repeat the proof of \scite{88r-5.12}.  More elaborately by the way $\bold
D_\alpha$ was chosen, Claim \scite{88r-5.17}(3) holds for ${\frak
K}_{\bold D_*}$ hence \wilog \, $M^J$ is $(\bold
D_\alpha(N_1),\aleph_0)$-homogeneous and so \wilog \, for some $t(*)
\in I \backslash J,N^J_1 = M_{t(*)}),N^J=M_{t(*)+1}$ and we get a
contradiction as in the proof of \scite{88r-5.12} (i.e., the choice of
$\langle \bar a_\ell:\ell \le \ell(*)\rangle$ there).  
\hfill$\square_{\scite{88r-5.20}}$
\enddemo
\bigskip

\definition{\stag{88r-5.21} Definition}  1) ${\frak K}$ has the symmetry property
when the following holds: if $N_0 \le_{\frak K} N_\ell \le_{\frak K}
N_3 \, (\ell = 1,2)$ and
for every $\bar a \in N_1$, gtp$(\bar a,N_2,N_3)$ is the
stationarization of gtp$(\bar a,N_0,N_3)$, \ub{then} for every $\bar b
\in N_2$, gtp$(\bar b,N_1,N_3)$ is the stationarization of 
gtp$(\bar b,N_0,N_3)$. \nl
2) If $N_0,N_1,N_2 \le_{\frak K} N_3$ 
satisfies the assumption and conclusion of
part (1) we say that $N_1,N_2$ are in stable amalgamation over $N_0$
inside $N_3$.  If only
the hypothesis of (1) holds we say they are in a one sided stable
amalgamation over $N_0$ inside $N_3$ 
(then the order of $(N_1,N_2)$ is important). \nl
3) We say that ${\frak K}$ has unique [one sided] amalgamation when: if
$N_0 \le_{\frak K} N_\ell \in K_{\aleph_0}$ for $\ell=1,2$ then
$N_1,N_2$ has unique [one sided] stable amalgamation, see part (4).
We say $N_1,N_2$ have a unique [one sided] stable amalgamation
over $N_0$ provided.
\nl
4) That \ub{if}: clauses (a)-(d) below holds then
$(*)$ below holds, where:
\mr
\item "{$(a)$}"  $N_1 \le_{\frak K} N_3,N_2 \le_{\frak K} N_3$ and
$(N_1,N_2)$ in [one sided] stable amalgamation inside $N_3$ over $N_0$
and $\|N_3\| \le \|N_1\| + \|N_2\|$
\sn
\item "{$(b)$}"  $M_0 \le_{\frak K} M_\ell \le_{\frak K} M_3$ for
$\ell=1,2$ and $(M_1,M_2)$ are in [one sided] stable amalgamation
inside $M_3$ over $M_0$
\sn
\item "{$(c)$}"  $f_\ell$ is an isomorphism from $N_\ell$ onto
$M_\ell$ for $\ell=0,1,2$
\sn
\item "{$(d)$}"  $f_0 \subseteq f_1$ and $f_0 \subseteq f_2$
{\roster
\itemitem{ $(*)$ }  we can find 
$M'_3,M_3 \le_{\frak K} M'_3$ and $f_3$, a $\le_{\frak
K}$-embedding of $N_3$ into $M'_3$ extending $f_1 \cup f_2$.
\endroster}
\endroster
\enddefinition
\bn
We at last get the existence of stable amalgamation (to which we
earlier get approximations).
\proclaim{\stag{88r-5.22} Claim}  For any $N_0 \le_{\frak K} N_1,N_2$, all from
$K_{\aleph_0}$, we can find $M,N_0 \le_{\frak K} M \in K_{\aleph_0}$ and
$\le_{\frak K}$-embeddings 
$f_1,f_2$ of $N_1,N_2$ respectively over $N_0$ into $N$
such that $N_0,f_1(N_1),f_2(N_1)$ are in stable amalgamation.
\endproclaim
\bigskip

\remark{Remark}  In the proof we could have ``inverted the tables" and
used $\bar c_\zeta$ in the $\omega_1$ direction.
\endremark
\bigskip

\demo{Proof}  We define by induction on $\zeta < \omega_1,\langle
M^\zeta_\alpha:\alpha < \omega_1 \rangle$ and $\bar c_\zeta$ such
that:
\mr
\widestnumber\item{$(iii)$}
\item "{$(i)$}"   $\langle M^\zeta_\alpha:\alpha < \omega_1 \rangle$
is $\le_{\frak K}$-increasing continuous and $M^\zeta_\alpha \in K_{\aleph_0}$
\sn
\item "{$(ii)$}"   for $\alpha < \zeta,M^\zeta_\alpha =
M^\alpha_\alpha$ and $\xi < \zeta \and \alpha < \omega_1 \Rightarrow
M^\xi_\alpha \le_{\frak K} M^\zeta_\alpha$
\sn
\item "{$(iii)$}"   for $\zeta$ limit, $M^\zeta_\alpha = \dbcu_{\xi <
\zeta} M^\xi_\alpha$
\sn
\item "{$(iv)$}"   for $\zeta \le \alpha < \omega_1,\zeta$ non-limit 
$M^\zeta_{\alpha +1}$ is
$(\bold D_{\alpha +1}(M^\zeta_\alpha),\aleph_0)^*$-homogeneous
\sn
\item "{$(v)$}"   for every $\bar c \in M^\zeta_{\alpha +1}$,
gtp$(\bar c,M^{\zeta +1}_\alpha,M^{\zeta +1}_{\alpha +1})$ is a
stationarization of gtp$(\bar c,M^\zeta_\alpha,M^\zeta_{\alpha +1})$
\sn
\item "{$(vi)$}"  $\bar c_\zeta \in M^{\zeta +1}_{\zeta +1}$ and for
$\zeta +1 < \alpha < \omega_1$, gtp$(\bar
c_\zeta,M^\zeta_\alpha,M^{\zeta +1}_\alpha)$ is the stationarization of
gtp$(\bar c_\zeta,M^\zeta_{\zeta +1},M^{\zeta +1}_{\zeta +1})$
\sn
\item "{$(vii)$}"   for every $p \in \bold D(M^\xi_\alpha)$ for some
$\zeta$ satisfying $\xi + \alpha < \zeta < \omega_1$ we have 
gtp$(\bar c_\zeta,M^\zeta_{\zeta
+1},M^{\zeta +1}_{\zeta +1})$ is a stationarization of $p$.
\ermn
There is no problem doing this (by \scite{88r-5.20} and as in earlier
constructions); in limit stages we use local character \scite{88r-5.13}(3)
and $\bold D_\alpha$ closed under stationarization. \nl
Now easily for a thin enough closed 
unbounded set of $E \subseteq \omega_1$, for
every $\zeta \in E$ we have
\mr
\item "{$(*)_\zeta(a)$}"  $M^\zeta_\zeta$ is
$(\bold D_\zeta(M^0_\zeta),\aleph_0)^*$-homogeneous
\sn
\item "{${{}}(b)$}"  for every $\bar c \in M^\zeta_\zeta$, gtp$(\bar
c,\dbcu_{\alpha < \omega_1} M^0_\alpha,\dbcu_{\xi < \omega_1}
M^\xi_\xi)$ is a stationarization of gtp$(\bar c,M^0_\zeta,M^\zeta_\zeta)$
\sn
\item "{${{}}(c)$}"  for every $\bar c \in M^0_{\zeta +1}$, 
gtp$(\bar c,M^{\zeta +1}_\zeta,M^{\zeta +1}_{\zeta +1})$ is a stationarization
of gtp $(\bar c,M^0_\zeta,M^0_{\zeta +1})$. 
\ermn
[Why?  Clause (c) holds by clause (v) of the construction (as $\langle
M^\zeta_\varepsilon:\varepsilon \le \zeta \rangle$ is $\le_{\frak
K}$-increasing continuous).  Clause (b) holds as $E$ is thin enough,
i.e., is proved as in earlier constructions (i.e., see the proof of
$(*)$ in the proof of \scite{88r-5.17}(1).  As for Clause (a) first
note that by clauses (i),(ii),(iii) the sequence
$\langle M^\zeta_\varepsilon:\varepsilon \le
\zeta \rangle$ is $\le_{\frak K}$-increasing continuous.  By clause
(vi) we have $\varepsilon < \zeta \Rightarrow \text{ gtp}(\bar
c_\varepsilon,M^\varepsilon_\zeta,M^{\varepsilon +1}_\zeta)$ does not
fork over $M^\varepsilon_\zeta$ and 
(vii) of the construction we have: if $p \in \bold
D_\zeta(M^\zeta_\varepsilon),\varepsilon < \zeta$ then for some $\xi
\in (\varepsilon,\zeta)$, tp$(\bar c_\xi,M^\zeta_\xi,M^\zeta_{\xi
+1})$ is a non-forking extension of $p$.  As $E$ is thin enough we
have $\bar d \in M^\zeta_\zeta \Rightarrow \text{ gtp}(\bar
d,M^\zeta_0,M^\zeta_\zeta) \in \bold D_\zeta(M^\zeta_0)$.  Together it is
easy to get clause (a), e.g., see \scite{88r-5.25}.]

So as in the proof of \scite{88r-5.17}(3) we can finish (choose $\zeta
\in E,f_0$ an isomorphism from $N_0$ onto $M^0_\zeta,f_1 \supseteq
f_0$ is an $\le_{\frak K}$-embedding of $N_1$ into $M^\zeta_\zeta$ and
$f_2 \supseteq f_0$ a $\le_{\frak K}$-embedding of $N_2$ into
$M^0_{\zeta +1}$).  \hfill$\square_{\scite{88r-5.22}}$
\enddemo
\bigskip

\remark{\stag{88r-5.22.3} Remark}  Note that in 
\chaptercite{600} we use only the results up to this point.
\endremark
\bigskip

\proclaim{\stag{88r-5.23} Theorem}  1) Suppose in addition to the
hypothesis of this section that $2^{\aleph_1} < 2^{\aleph_2}$ and the
weak diamond ideal on $\aleph_1$ is not $\aleph_2$-saturated and
$\dot I(\aleph_2,K) < 2^{\aleph_2}$ or just 
$\dot I(\aleph_2,K(\aleph_1$-saturated)) $< 2^{\aleph_2}$.
\ub{Then} ${\frak K}$ has the symmetry property.
\nl
2) If $\dot I(\aleph_2,K(\aleph_1$-saturated)) $<
\mu_{\text{unif}}(\aleph_2,2^{\aleph_1})$; this number is always $>
2^{\aleph_1}$, usually $2^{\aleph_2}$, see \scite{88r-0.wD}.  Then 
stable amalgamation in
$K_{\aleph_0}$ is unique (we know that it always exists and it follows
by (1) + (2) that one sided amalgamation is unique).
\endproclaim
\bn
\margintag{88r-5.23.7}\ub{\stag{88r-5.23.7} Discussion}:  1) This certainly gives a desirable
conclusion.  However, part (2) is not used so we shall return to it in
\cite{Sh:838}, \cite{Sh:849}.
\nl
2) As for part (1), we can avoid using it (except in \scite{88r-5.24} below).  
More fully, in
\sectioncite[\S3]{600} dealing with ${\frak K}$ as here by
\marginbf{!!}{\cprefix{600}.\scite{600-Ex.1}} for every $\alpha < \omega_1$ we derive a good
$\aleph_0$-frame ${\frak s}_\alpha$ with $K^{{\frak s}_\alpha} =
{\frak K}_{\bold D_\alpha}$ (if we would have liked to derive a good
$\aleph_1$-frame we would need \scite{88r-5.23}).

Then in \chaptercite{705} if ${\frak s}$ is
successful (holds if $2^{\aleph_0} < 2^{\aleph_1} < 2^{\aleph_2}$ and
$\dot I(\aleph_2,{\frak K}^{{\frak s}_\alpha}) < 2^{\aleph_2}$ and
WDmId$_{\aleph_1}$ is not $\aleph_2$-saturated)
then we derive the successor ${\frak s}^+_\alpha$, a
good $\aleph_1$-frame with $K^{{\frak s}^+_\alpha} \subseteq \{M \in
K^{{\frak s}_\alpha}_{\aleph_1}:M$ is $\aleph_1$-saturated for 
$K^{{\frak s}_\alpha}\}$, and ${\frak s}^+_\alpha$ is even good$^+$
(see Claim \marginbf{!!}{\cprefix{705}.\scite{705-stg.3}}(2) and Definition
\marginbf{!!}{\cprefix{705}.\scite{705-stg.1}}).  This suffices for the main conclusions of
\sectioncite[\S9]{600} and end of \sectioncite[\S12]{705}.  
\nl
3) Still we may wonder
is $\le_{{\frak s}^+_\alpha} = \le_{\frak K} \restriction {\frak
K}_{{\frak s}^+_\alpha}$?  If ${\frak s}_\alpha$ is good$^+$ then the
answer is yes (see \marginbf{!!}{\cprefix{705}.\scite{705-stg.3}}(1)).  That is, the 
present theorem \scite{88r-5.23} is used in \sectioncite[\S1]{705} to prove
${\frak s}$ is ``good$^+$", in fact proved in \scite{88r-5.24}.  In fact
part (1) of \scite{88r-5.23} is enough to prove that ${\frak s}_{\bold
D_*}$ is good$^+$; we shall return to this in \cite{Sh:838}, \cite{Sh:849}).
\nl
3) The proof of \scite{88r-5.23}(1) gives that if ${\frak K}$ fails the
symmetry property then $\dot I(\aleph_2,K) \ge 2^{\aleph_1}$ even
if $2^{\aleph_1} = 2^{\aleph_2}$ and do not use $2^{\aleph_0} =
2^{\aleph_1}$ directly (but use earlier results of \S5).  The case
``${\Cal D}_{\aleph_1}$ is $\aleph_2$-saturated, $2^{\aleph_0} <
2^{\aleph_1} < 2^{\aleph_2},\dot I(\aleph_2,\aleph_2) <
\mu_{\text{unif}}(\aleph_2,2^{\aleph_2})$" is covered in \cite{Sh:838},
\cite{Sh:849}. 
\bigskip

\demo{Proof}  1)  So in the first part 
toward contradiction we can assume that $K^4 \ne \emptyset$ where
$K^4$ is the class of quadruple $\bar N = (N_0,N_1,N_2,N_3)$ such that
$N_1,N_2$ are one sided stably amalgamated over 
$N_0$ inside $N_3$ but $N_2,N_1$
are not.  Hence there is $\bar c \in N_2$ such that gtp$(\bar c,
N_1,N_3)$ is not the stationarization of gtp$(\bar c,N_0,N_2) =
\text{ gtp}(\bar c,N_0,N_3)$.  We define a two-place relation $\le$ on
$K^4$ by $\bar N^1 \le \bar N^2$ iff $N^1_0 = N^2_0,N^1_\ell
\le_{\frak K} N^2_\ell$ for $\ell=0,1,2$ and $\bar a \in N^1_1
\Rightarrow \text{ gtp}(\bar a,N^2_2,N^2_3)$ is definable over some
$\bar b \in N^1_0$.  Easily this is a partial order and
$K^4$ is closed under union of
increasing countable sequences.  Hence \wilog \, for some $\bold
D_*,\bar N^*$
\mr
\item "{$(*)$}"  $(a) \quad \bold D_* \in \{\bold D_\alpha:\alpha <
\omega_1\}$
\sn
\item "{${{}}$}"  $(b) \quad \bar N^* \in K^4$
\sn
\item "{${{}}$}"  $(c) \quad N^*_\ell$ is $(\bold
D_*(N^*_0),\aleph_0)^*$-homogeneous over $N^*_0$ for $\ell=1,2$
\sn
\item "{${{}}$}"  $(d) \quad N^*_3$ is $(\bold D_*(N^*_\ell),
\aleph_0)^*$-homogeneous over $N^*_\ell$ for $\ell=1,2$
\endroster
\enddemo
\bn
We have proved
\demo{\stag{88r-5.23.8} Observation}  To prove \scite{88r-5.23}, we can assume
that $\bold D = \bold D_\alpha$ for $\alpha < \omega_1$, i.e., $\bold
D$ is countable.
\enddemo
\bigskip

\demo{Continuation of the proof}  A problem is that we still 
have not proven the existence of a superlimit
model of $K$ of cardinality $\aleph_1$ though we have a candidate
$N^*$ from \scite{88r-5.9}.  So we use $N^*$, but to ensure we get it at
limit ordinals (in the induction on $\alpha < \aleph_2$), we 
have to take a stationary $S_0 \subseteq \omega_1$
with $\omega_1 \backslash S_0$ not small, i.e., $\omega_1 \backslash
S_0$ does not belong to the ideal WDmId$_{\aleph_1}$ from Theorem \scite{88r-0.wD} 
and ``devote" it to ensure this, using
\scite{88r-5.22}. 
\enddemo
\bn
The point of using $S_0$ is as follows (this is supposed to help
to understand the quotation):
\definition{\stag{88r-5.23A} Definition}  1) Let $K^{\text{qt}} 
= \{\bar N:\bar N = \langle N_\alpha:
\alpha < \omega_1 \rangle$ be $\le_{\frak K}$-increasing
continuous, $N_\alpha \in K_{\aleph_0},N_{\alpha +1}$ is
$(\bold D_\alpha(N_\alpha),\aleph_0)^*$-homogeneous$\}$. \nl
2) On $K^{\text{qt}}$ we define a two-place relation
$<^a_S$ (for $S \subseteq \omega_1$) as follows: \nl
$\bar N^1 <^a_S \bar N^2$ if and only if
for some closed unbounded $E \subseteq \omega_1$
\mr
\item "{$(a)$}"  for every $\alpha \in C$ we 
have $N^1_\alpha \le_{\frak K} N^2_\alpha,N^1_\alpha$ is $(\bold
D_*(N^1_\alpha),\aleph_0)^*$-homogenous and $N^1_{\alpha +1}
\le_{\frak K} N^2_{\alpha +1}$
\sn
\item "{$(b)$}"    for every $\alpha < \beta$ from $E$ we have 
$N^2_\beta \cap \dbcu_{\alpha < \omega_1}
N^1_\alpha = N^1_\beta$ and $N^1_\beta,N^2_\alpha$ are in one sided 
stable amalgamation over $N^1_\alpha$ inside $N^2_\beta$
\sn
\item "{$(c)$}"    if $\alpha \in S \cap C$ \ub{then}
$N^2_\alpha,N^1_{\alpha +1}$ are in stable amalgamation over
$N^1_\alpha$ inside $N^2_{\alpha +1}$.
\endroster
\enddefinition
\bigskip

\demo{\stag{88r-5.23B} Fact}  0) The two-place relation $<^a_S$ defined in 
\scite{88r-5.23A} are partial orders on $K^{\text{qt}}$ for $n < \omega$. \nl
1)  If $\bar N^n \le^a_{S_0} \bar N^{n+1}$ and
let $E_n$ exemplify this (as in the Definition \scite{88r-5.23A}) and 
let $E_\omega = \dbca_{n < \omega} E_n,E'_\omega = 
\{\alpha,\alpha +1:\alpha \in
C_\omega\}$ and let $N^\omega_\alpha = \dbcu_{n < \omega} N^n_\beta$
when $\beta = \text{ Min}[E'_\omega \backslash \alpha]$.  \ub{Then} 
$\langle N^\omega_\alpha:\alpha < \omega_1 \rangle \in K_{< \aleph_1}$
and $\bar N^n \le^a_{S_0} \langle
N^\omega_\alpha:\alpha < \omega_1 \rangle$ for $n < \omega$. \nl
2) If $\langle \bar N^\varepsilon:\varepsilon < \omega_1\rangle$ is
$<^a_S$-increasing and $N^\varepsilon =
\cup\{N^\varepsilon_\alpha:\alpha < \omega_1\} \in K_{\aleph_1}$ 
is $\le_{\frak K}$-increasing continuous, the club
$E_{\varepsilon,\zeta}$ witness $\bar N^\varepsilon \le \bar N^\zeta$
for $\varepsilon < \zeta < \aleph_1$ 
and $\langle N_\alpha:\alpha < \omega_1\rangle$ a
$\le_{\frak K}$-representation of $N$, and for a club of $\alpha <
\aleph_1,N_\alpha = \cup\{N^\varepsilon_\alpha:\varepsilon <
\alpha\},N_{\alpha +1} = \cup\{N^\varepsilon_{\alpha +1}:\varepsilon <
\alpha\}$ \ub{then} $\varepsilon < \omega_1
\Rightarrow \bar N^\varepsilon \le^a_{S_0} \bar N$.
\enddemo
\bigskip

\demo{Proof}  Should be easy by now.  \hfill$\square_{\scite{88r-5.23B}}$

Returning to the proof of \scite{88r-5.23} it is done as in
\cite[\S3]{Sh:576} or better in \cite{Sh:838}, basically as follows.

There is $\langle S_\varepsilon:\varepsilon < \omega_1\rangle$ such
that $S_\varepsilon \subseteq \omega_1,\zeta < \varepsilon \Rightarrow
S_\zeta \cap S_\varepsilon$ countable and $S_0,S_{\varepsilon
+1} \backslash S_\varepsilon \in ({\Cal D}_{\omega_1})^+$,
possible by an assumption.

Now for any $u \subseteq \omega_2$ we choose
$N^u_\varepsilon,N^u_\varepsilon$ by induction on $\varepsilon <
\omega_1$ such that
\mr
\item "{$\circledast$}"  $(a) \quad \bar N^u_\varepsilon = \langle
N^u_{\varepsilon,\alpha}:\alpha < \omega_1\rangle \in K^{\text{qt}}$
\sn
\item "{${{}}$}"  $(b) \quad N^u_\varepsilon =
\cup\{N^u_{\varepsilon,\alpha}:\alpha < \omega_1\} \in K_{\aleph_1}$
\sn
\item "{${{}}$}"  $(c) \quad$ for $\zeta < \varepsilon$ we have $\bar
N^u_\zeta <^1_{S_\xi} \bar N^u_\varepsilon$ when $\xi \notin
[\zeta,\varepsilon) \cap u$ (we can use $S'_{[\zeta,\varepsilon)}$,
\nl

\hskip25pt the compliment of the diagonal union of 
$\{\langle S_\xi:\varepsilon \in [\zeta,\varepsilon)\rangle \cap u\}$
\sn
\item "{${{}}$}"  $(d) \quad$ we can demand continuity as defined
implicitly in Fact \scite{88r-5.23B}
\sn
\item "{${{}}$}"  $(e) \quad$ for each $\varepsilon \in u$ for a club
of $\alpha < \omega_1$ if $\alpha \in S_\varepsilon$ then
$N^u_{\varepsilon +1,\alpha},N^u_{\varepsilon,\alpha +1}$ 
\nl

\hskip25pt are not in
stable amalgamation over $N^u_{\varepsilon,\alpha}$ inside
$N^u_{\varepsilon +1,\alpha +1}$ 
\nl

\hskip25pt (though is in one side).
\ermn
Lastly, let $N^u = \cup\{N^u_\varepsilon:\varepsilon < \omega_1\} \in
K_{\aleph_2}$.  Now we can prove thatif $u,v \subseteq w_2$ and $N^u
\approx N^v$ then for some club $C$ of $\omega_2,u \cap C = v \cap
C$.  So we can easily get $\dot I(\aleph_2,{\frak K}) = 2^{\aleph_2}$.
 \hfill$\square_{\scite{88r-5.23}}$
\enddemo
\bigskip

\proclaim{\stag{88r-5.24} Theorem}  Suppose ${\frak K}$ has the symmetry
property (holds if the conclusion of \scite{88r-5.23}(1) hold).  
\ub{Then} ${\frak K}$ has a superlimit model in $\aleph_1$.
\endproclaim
\bigskip

\demo{Proof}  We have a candidate $N^*$ from \scite{88r-5.9}.  So let
$\langle N_i:i < \delta \rangle$ be $\le_{\frak K}$-increasing, 
$N_i \cong N^*$ and \wilog \, $\delta = \text{ cf}(\delta)$.
If $\delta = \omega_1$ this is very easy.  If $\delta = \omega$, let
$N_\omega = \dbcu_{i < \omega} N_i$ and for each $i \le \omega$ let
$\langle N^\alpha_i:\alpha < \omega_1 \rangle$ be $\le_{\frak
K}$-increasing continuous with union $N_i$ and $N^\alpha_i \in
K_{\aleph_0}$.  Now by restricting ourselves to a
club $E$ of $\alpha$'s and renaming it $E = \omega_1$, we get:
for $i < j \le \omega,N^\alpha_i = N_i \cap N^\alpha_j$,
and
\mr
\item "{$\circledast_1$}"   for any $\alpha < \beta < \omega_1,
\bar a \in N^\alpha_\omega$ and $i < \omega$, the type 
gtp$(\bar a,N^\beta_i,N^\beta_\omega)$ is a stationarization of 
gtp$(\bar a,N^\alpha_i,N^\alpha_\omega)$.
\ermn  
To prove $N_\omega \cong N^*$ it is enough to prove:
\mr
\item "{$\circledast_2$}"  if $\alpha < \omega_1,p \in \bold
D(N^\alpha_\omega)$ then some $\bar b \subseteq N_\omega$ realizes $p$
in $N_\omega$.
\ermn
By \scite{88r-5.13}(3) there is $i < \omega$ such that $p$ is the
stationarization of $q = p \restriction N^\alpha_i \in \bold
D(N^\alpha_i)$.  As $N_i \cong N^*$, there is $\bar b \subseteq N_i$
which realizes $q$ and we can find $\beta \in (\alpha,\omega_1)$ such
that $\bar b \subseteq N^\beta_i$.  By $\circledast_2$
we have $N^\alpha_\omega,N^\beta_i$ is in one sided
stable amalgamation over $N^\alpha_i$ inside $N^\beta_\omega$
(see \scite{88r-5.21}(2)).

As we assume the conclusion of Theorem \scite{88r-5.23}(1), also
$N^\beta_i,N^\alpha_\omega$ is in stable amalgamation over
$N^\alpha_i$ inside $N^\beta_\omega$.  
In particular, as $\bar b \subseteq N^\beta_i$, we have
gtp$(\bar b,N^\alpha_\omega,N^\beta_\omega)$ is the stationarization
of gtp$(\bar b,N^\alpha_i,N^\beta_i)$ but the latter is $p
\restriction N^\alpha_i$ so by uniqueness of stationarization, $p =
\text{ gtp}(\bar b,N^\alpha_\omega,N^\beta_\omega)$ which is gtp$(\bar
b,N^\alpha_\omega,N_\omega)$, so $p$ is realized in $N_\omega$ as
required.    \hfill$\square_{\scite{88r-5.24}}$
\enddemo
\bn
We have implicitly proved
\proclaim{\stag{88r-5.24.5} Claim}  Assume that $N_0 \le_{\frak K} N_1 \in
K_{\aleph_0}$ and $\bar a_\ell \in {}^{\omega >}(N_1)$ for
$\ell=1,2$.  \ub{Then} $(*)_1 \Leftrightarrow (*)_2$ where for
$\ell=1,2$
\mr
\item "{$(*)_\ell$}"  there are $M_1,M_2,\bar b_1,\bar b_2$ such that
{\roster
\itemitem{ $(a)$ }  $N_0 \le_{\frak K} M_1 \le_{\frak K} M_2 \in
K_{\aleph_1}$
\sn
\itemitem{ $(b)$ }  $\bar a_k\in {}^{\omega >}(M_k)$ for $k=1,2$
\sn
\itemitem{ $(c)$ }  {\rm gtp}$(\bar b_\ell,N_0,M_1) = \text{\rm
gtp}(\bar a_{3 - \ell},N_0,N_1)$
\sn
\itemitem{ $(d)$ }  {\rm gtp}$(\bar b_\ell,M_1,M_2)$ is the
stationarization of {\rm gtp}$(\bar a_\ell,N_0,N_1)$ from $\bold
D(M_1)$
\sn
\itemitem{ $(e)$ }  {\rm gtp}$(\bar b_1 \char 94 \bar b_2,N_0,M_2) =
\text{\rm gtp}(\bar a_1 \char 94 \bar a_2,N_0,N_1)$
\endroster}
\endroster
\endproclaim
\bigskip

\demo{Proof}  We can deduce it from \scite{88r-5.20} (or immitate the
proof of \scite{88r-5.12}).

In detail by symmetry it is enough to assume $(*)_2$ and prove
$(*)_1$.  So let $M_1,M_2,\bar b_1,\bar b_2$ witness $(*)_2$.

By \scite{88r-5.26} we can find $M'_2,f$ such that: $M_2 \le_{\frak K}
M'_2 \in K_{\aleph_0},f$ is a $\le_{\frak K}$-embedding of $N_1$ into
$N'_3$ over $N_0$ such that $M_1,f(M_2)$ is in stable amalgamation
over $N_0$ inside $M_2$.  Now, as $f(M_2),M_1$ are in one sided stable
amalgamation over $N_0$ inside $N'_3$ by the choie of $(M_1,M_2,\bar
b_1,\bar b_2)$ we get gtp$(f(\bar b_2),M_1,M'_2) = \text{\rm gtp}(\bar
b_2,M_1,M'_2)$ hence gtp$(\bar b_1 \char 94 \bar b_2,N_0,M'_2) =
\text{\rm gtp}(\bar b_1 \char 94 f(\bar b_2),N_0,M'_2)$.

By the choice of $M^2_1,f$, tp$(\bar b_1,f(M_2),M'_2)$ is the
stationarization of gtp$(\bar b_1,N_0,M_2) = \text{\rm gtp}(\bar
a_1,N_0,N_1)$.  Now $(*)_1$ holds as exempflified $(f(M_2),M'_2,f(\bar
b_2),\bar b_1)$. 
\enddemo
\bigskip

\demo{\stag{88r-5.25} Exercise}  Assume $\alpha \le \omega_1$ and
\mr
\item "{$(a)$}"  $\langle M_i:i \le \delta \rangle$ is $\le_{\frak
K}$-increasing continuous $\delta$ a limit ordinal
\sn
\item "{$(b)$}"  if $p \in \bold D(M_i)$ is realized in $M_{i+1}$ then it
$\in \bold D_\alpha(M_i)$ or just $p \restriction M_0 \in \bold D(M_0)$
\sn
\item "{$(c)$}"  if $i <\delta,p \in \bold D_\alpha(M_i)$ then $p$ is
materialized in $M_j$ for some $j \in (i,\delta)$.
\ermn
\ub{Then} $M_\delta$ is $(\bold D_\alpha(M_0),\aleph_0)^*$-homogeneous.
\enddemo
\bigskip

\demo{Proof}  Easy.
\enddemo
\bn
\margintag{88r-5.26}\ub{\stag{88r-5.26} Discussion}:  1) Consider 
$\psi \in \Bbb L_{\omega_1,\omega}(\bold
Q),|\tau_\psi| \le \aleph_0,1 \le \dot I(\aleph_1,\psi) <
2^{\aleph_0}$.  We translate it to ${\frak K}$ and $<^{**}$ as
earlier, see \scite{88r-3.9}.
\nl
2) What if we waive categoricity in $\aleph_0$?  The adoption of this
was O.K. as we shrink ${\frak K}$ but not too much.  But without
shrinking probably we still can say something on the model in ${\frak
K}^* = \{M \in {\frak K}_{\ge \aleph_0}$: if $N_0 \le_{\frak K} M,N_0
\in K_{\aleph_0}$ then for some $N_1,N_0 <^* N_1 \le_{\frak K} M\}$
as there are good enough approximations.
\newpage

\head {\S6 Counterexamples} \endhead  \resetall \sectno=6
 \spuriousreset
\bigskip

In \cite{Sh:48} the statement of Conclusion \scite{88r-3.5} was proved for
the first time where $K$ is the class of atomic models of a first
order theory assuming Jensen's diamond $\diamondsuit_{\aleph_1}$
(taking $\lambda = \aleph_0$).  In \cite{Sh:87a} and \cite{Sh:87b} the
same theorem was
proved using $2^{\aleph_0} < 2^{\aleph_1}$ only (using \scite{88r-0.wD}).
Let us now concentrate on the
case $\lambda = \aleph_0$.  We asked whether the assumption
$2^{\aleph_0} < 2^{\aleph_1}$ is necessary to get Conclusion \scite{88r-3.5}.
In this section we construct three classes of models $K^1,K^2,K^3,K^4$
failing amalgamation, i.e., failing the conclusion of \scite{88r-3.5},
$K^2,K^3,K^4$ are a.e.c. with LS-number $\aleph_0$
while $K^1$ satisfy all the axioms needed in the proof of
Conclusion \scite{88r-3.5} (but it is not an abstract elementary class -
fails to satisfy Ax.IV).

$K^2$ is PC$_{\aleph_0}$ and is axiomatizable in
$\Bbb L_{\omega_1,\omega}(\bold Q)$.

$K^3$ is PC$_{\aleph_0}$ and is axiomatizable in $\Bbb L(\bold Q)$.
\sn
Now the common phenomena to $K^1,K^2,K^3,K^4$ are that all of them satisfy
the hypothesis of Conclusion \scite{88r-3.5}, i.e., for $\ell =
1,2,3$ we have 
$\dot I(\aleph_0,K^\ell) = 1$ and the $\aleph_0$-amalgamation property fails in
$K^\ell$, but assuming $\aleph_1 < 2^{\aleph_0}$ and MA$_{\aleph_1}$
for $\ell = 1,2,3$ we have $\dot I(\aleph_1,K^\ell)=1$.
\bigskip

\definition{\stag{88r-6.1} Definition}  Let $Y$ be an infinite set.  A
family ${\Cal P}$ of infinite subsets of $Y$ is called independent 
if for every $\eta \in {}^{\omega >} 2$ and pairwise distinct
$X_0,X_1,\dotsc,X,\ell g(\eta)-1$ (notation:
for $X \in {\Cal P}$ denote $X^0 = X$ and $X^1 = Y \backslash X$) 
the following set $\dbca_{k < \ell g(\eta)} X^{\eta[k]}_k$ is infinite.
\enddefinition
\bigskip

\definition{\stag{88r-6.2} Definition}  1) The class of models $K^0$ is
defined by

$$
\align
K^0 = \{M:&M = \langle|M|,P^M,Q^M,R^M \rangle,|M| = P^M \cup Q^M, \\
  &P^M \cap Q^M = \emptyset,|P^M| = \aleph_0 \le |Q^M| \text{ and} \\
  &R \subseteq P^M \times Q^M\}.
\endalign
$$
\mn
2) For $M \in K^0$, let 
$A^M_y = \{x \in P^M:x R^M y\}$ for every $y \in Q^M$. \nl
3) Let $K^1$ be the class of $M \in K^0$ such that
\mr
\item "{$(a)$}"   the family 
$\{A^M_y:y \in Q^M\}$ is independent, which means that if $m<n$ and
$y_0,\dotsc,y_{n-1}$ are pairwise distinct members of $Q^M$ 
then the set $\{x \in P^M:x R^M y_\ell \equiv \ell < m$ 
for every $\ell<n\}$ is infinite
\sn
\item "{$(b)$}"   for every disjoint finite subsets $u,w$ of $P^M$ we have 
$\|M\| = |A^M_{u,w}|$ where $A^M_{u,w} := \{y \in Q^M:a \in u
\Rightarrow (a R^M y)$ and $b \in w \Rightarrow \neg(b R^My)\}$. 
\ermn
4) The notion of (strict) substructure $\le_{{\frak K}^1}$ is 
defined by: for $M_1,M_2 \in K^1,M_1 \le_{{\frak K}^1} M_2$ \ub{iff} 
$M_1 \subseteq M_2,P^{M_1}= P^{M_2}$ 
and for any finite disjoint $u,w \subseteq P^{M_2}$ the set $A^{M_2}_{u,w}
\backslash M_1$ is infinite when $M_1 \ne M_2$ (equivalently - non-empty). \nl
5) ${\frak K}^1 = (K^1,\le_{{\frak K}^1})$.
\enddefinition
\bigskip

\proclaim{\stag{88r-6.3} Lemma}  The class $(K^1,<_{{\frak K}^1})$ satisfies \nl
0) Ax 0. \nl
1) Ax I. \nl
2) Ax II. \nl
3) Ax III. \nl
4) Ax IV fails even for $\lambda = \aleph_0$. \nl
5) Ax V fails for countable models. \nl
6) Ax VI holds with {\rm LS}$({\frak K}^1) = \aleph_0$. \nl
7) For every $M \in K^1,\|M\| \le 2^{\aleph_0}$.
\endproclaim
\bigskip

\demo{Proof}  0), 1), 2) follows trivially from the definition. \nl
3) To prove that $M = \dbcu_{i < \lambda} M_i \in K^1$, it is enough
to verify that for every finite disjoint $u,w \subseteq
P^M,|A^M_{u,w}|  = \|M\|$.   
If $\langle M_i:i < \lambda \rangle$ is eventually constant we are
done hence \wilog \, $\langle M_i:i < \lambda\rangle$ is
$<_{{\frak K}^1}$-increasing; from the definition of 
$<_{{\frak K}^1}$ it follows that for each $i,M_{i+1}$
has a new $y = y_i$ as above, i.e., $y_i \in A^{M_{i+1}}_{u,w} \backslash M_i$
for every $i < \lambda$.   
Also for each $i$ there are at least $\|M_i\|$ many members in
$A^{M_i}_{u,w} \subseteq A^M_{u,w}$.
Together there are at least $\|M\|$ members in $A^M_{u,w}$. \nl
4) Let $\{M_n:n < \omega\} \subseteq K^1_{\aleph_0}$ be an 
$<_{{\frak K}^1}$-increasing
chain, let $M = \dbcu_{n < \omega} M_n$; by part 3) we have 
$M \in K^1_{\aleph_0}$.
Since $|Q^M| = \aleph_0$ by Claim \scite{88r-6.5}(a) below there
exists $A \subseteq P^M \backslash \{A^M_y:y \in Q^M\}$ 
infinite such that $\{A_y:y \in Q^M\}
\cup \{A\}$ is independent.  Now define $N \in K^1$ by $P^N = P^M$, let $y_0
\notin M,Q^N = Q^M \cup \{y_0\}$ and finally
let $R^N  = R^M \cup \{\langle a,y_0 \rangle:a \in P^N \and a \in A\}$.
Clearly for every $n < \omega,M_n \le_{{\frak K}^1} N$ but 
$N$ is not an $\le_{{\frak K}^1}$-extension of $M = \dbcu_{n < \omega}
M_n$ because the second part in Definition \scite{88r-6.2}(4) is
violated. \nl
5) Let $N_0 <_{{\frak K}^1} N \in K^1$ be given; as
in 4) define $N_1 \subseteq N,|N_0| \subseteq |N_1|$ by adding a single
element to $Q^{N_0}$ (from the elements of $Q^N \backslash Q^{N_0}$) it is
obvious that $N_0 \le_{{\frak K}^1} N,N_1 \le_{{\frak K}^1} N$ 
but $N_0 \nleq_{{\frak K}^1} N_1$. \nl
6) By closing the set under the second requirement in Definition
\scite{88r-6.2}(3).
\nl
7) Let $y_1 \ne y_2 \in Q^M$, we show that $A^M_{y_1} \ne A^M_{y_2}$;
if $A^M_{y_1} \subseteq A^M_{y_2}$ then $A^M_{y_1} \cap (P^M
\backslash A^M_{y_2}) =
\emptyset$ contradiction to the requirement that $\{A_y:y \in Q\}$ 
is independent hence $|Q^M| \le 2^{|P^M|} = 2^{\aleph_0}$ and as
$|P^M| = \aleph_0$ we are done.  \hfill$\square_{\scite{88r-6.3}}$
\enddemo
\bigskip

\proclaim{\stag{88r-6.4} Theorem}  ${\frak K}^1 = (K^1,<_{{\frak K}^1})$
satisfies the hypothesis of Conclusion \scite{88r-3.5}.  Namely
\mr
\item  $\dot I(\aleph_0,K^1) = 1$
\sn
\item  every $M \in K^1_{\aleph_0}$ has a proper $\le_{{\frak K}^1}$-extension
in $K^1_{\aleph_0}$ 
\sn
\item  ${\frak K}^1$ is closed under chains of length $\le \omega_1$
\sn
\item  ${\frak K}^1$ fails the $\aleph_0$-amalgamation property.
\endroster
\endproclaim
\bigskip

\demo{Proof}  1) Let $M_1,M_2 \in K^1_{\aleph_0}$, pick the following
enumerations $|M_1| = \{a_n:n < \omega\}$ and $|M_2| = \{b_n:n <
\omega\}$.  It is enough to define an increasing sequence of finite
partial isomorphisms $\langle f_n:n < \omega \rangle$ from 
$M_1$ to $M_2$ such that for every $k < \omega$ for some 
$n(k) < \omega$ satisfy $a_k \in \text{ Dom}(f_{n(k)})$ and 
$b_k \in \text{ Range}(f_{n(k)})$ (finally take $f =
\dbcu_{n < \omega} f_n$ and this will be an isomorphism from $M_1$
onto $M_2$).

Define the sequence $\langle f_n:n < \omega \rangle$ by 
induction on $n < \omega$:
let $f_0 = \emptyset$, if $n = 2m$ denote $k = \text{ min}\{k <
\omega:a_k \notin \text{ Dom}(f_n)\}$.  Distinguish between the
following two alternatives:
\mr
\item "{$(A)$}"    if $a_k \in P^{M_1}$ let $\{a'_0,\dotsc,a'_{j-1}\} = 
Q^{M_1} \cap \text{ Dom}(f_n)$.  Without loss of generality there
exists $i \le j-1$ such that for all $\ell < i,a_k R^{M_1} a'_\ell$
and for all $i \le \ell \le j-1,\neg a_k R a'_\ell$.  By \scite{88r-6.2}(1),
$P^{M_\ell}$ is infinite, hence by clause (b) of \scite{88r-6.2}(2) also
$Q^{M_\ell}$ is infinite.  Hence by clause (a) of \scite{88r-6.2}(2) there
are infinitely many $y \in P^{M_2}$ such that 
$y R^{M_2} f_n(a'_\ell)$ for all $\ell < i$ and
for all $i \le \ell < j -1,\neg y R^{M_2} f_n(a'_\ell)$.
But Rang$(f_n)$ is finite.   Hence there is such
$y \in P^{M_2} \backslash \text{ Rang}(f_n)$.  Finally 
$f_{n+1} = f_n \cup \{\langle a_k,y \rangle\}$
\sn
\item "{$(B)$}"    if $a_k \in Q^{M_1}$ let $\{a'_0,\dotsc,a'_{j-1}\} =
P^{M_1} \cap \text{ Dom}(f_n)$ and as before we may assume that there
exists $i \le j-1$ such that for all $\ell < i,a'_\ell R^{M_1} 
a_k$ and for all $i \le \ell
<j - 1$ we have $\neg(a'_\ell) R^{M_1} a_k$.  By the second requirement
in Definition \scite{88r-6.2}(3) there exists $y \in Q^{M_2}$ such that 
$(\forall \ell < i)[f_n(a'_\ell) R^{M_2} y]$
and $(\forall \ell)[i \le \ell < j -1 \Rightarrow \neg f_n(a'_\ell) 
R^{M_2} y]$.  Now define
$f_{n+1} = f_n \cup \{\langle a_k,y \rangle\}$.
\ermn
When $n=2m+1$ act similarly on $b_{\text{min}\{k < \omega:b_k \notin
\text{ Rang}(f_n)\}}$. \nl
2) First we prove the following.
\enddemo
\bigskip

\demo{\stag{88r-6.5} Observation}
\mr
\item "{$(a)$}"  Let $P$ be a countable set.  For every countable family
${\Cal P}$ of infinite subsets of $P$ \ub{if} ${\Cal P}$ is 
independent \ub{then} there exists an infinite $A \subseteq P$ 
such that ${\Cal P} \cup \{A\}$ is independent and $A \notin {\Cal
P}$, of course
\sn
\item "{$(b)$}"  if $A,{\Cal P}$ are as in (a) \ub{then} for every infinite $B
\subseteq P$ satisfying $|A \Delta B| < \aleph_0$ also ${\Cal P} \cup \{B\}$
is independent (and $B \notin {\Cal P}$)
\sn
\item "{$(c)$}"  moreover in clause (a) we can require in addition that: for
any disjoint finite $u,w \subseteq P$ there exists $A \subseteq P$
as in (a) satisfying $u \subseteq A$ and $A \cap w = \emptyset$.
\endroster
\enddemo
\bigskip

\demo{Proof of Claim \scite{88r-6.5}}
\sn
\ub{Clause $(a)$}:    Let ${\Cal P}^* = \{X \subseteq P:(\exists n <
\omega)(\exists X_0 \in {\Cal P}) \ldots (\exists X_{n-1} \in {\Cal P})
(\exists k < n)$ [$X$ or $P \backslash X$ is equal to 
$\cap\{X_i:i < k\} \cap \cap\{P \backslash
X_i:k \le i < n\}\}$.

Clearly $|{\Cal P}^*| = \aleph_0$ hence we can find a sequence
$\langle A_n:n < \omega \rangle$ such that  $\{A_n:n < \omega\} =
{\Cal P}^*$ and such that for every $k < \omega$ there exists $n >k$ satisfying
$A_n = A_k$ hence for some $n > k,A_n = P \backslash A_k$.  Let 
$P = \{a_n:n < \omega\}$ without repetition.

Now define by induction $i:\omega \rightarrow \omega$.  \nl
\sn
Let $i(0) = 0$.  \nl
If $n = k+1$, let $i(n) = \text{ Min}\{\ell < \omega:i(n-1) < \ell$ and

$$
a_\ell \in (A_k \backslash \{a_{i(0)},\dotsc,a_{i(n-1)})\}.
$$
\mn
It is easy to verify that the construction is possible.  
Directly from the construction it follows
that $A = \{a_{i(n)}:n < \omega\}$ is a set as required.
\mn
\ub{Clause $(b)$}:  Easy.
\mn
\ub{Clause $(c)$}:  Let $u,w \subseteq P$ be finite disjoint and ${\Cal P}$ a
countable family of subsets of $P$ which is independent.

Let $A' \subseteq P$ be as proved in 
clause (a).  According to (b) also $A = (A' \cup u) \backslash w$ 
satisfies: the family ${\Cal P} \cup \{A\}$ is independent.
\enddemo
\bigskip

\demo{Return to the proof of Theorem \scite{88r-6.4}(2)}  Let ${\Cal P} = 
\{A_y^M \subseteq P^M:y \in Q^M\}$.  Let $\langle s_n:n < \omega\rangle$ be an
enumeration of $[P^M]^{< \aleph_0}$ with repetitions such that for
every finite disjoint $u,w \subseteq P^M$ there exists $n < \omega$ such
that $s_{2n} = u,s_{2n+1} = w$ and for all $k < \omega,s_{2k} \cap
s_{2k+1} = \emptyset$.

It is enough to define $\{{\Cal P}_n:n < \omega\}$ increasing chain of
countable independent families of subsets of $P^M$ 
such that ${\Cal P}_0 = {\Cal P}$ and for all 
$k < \omega$ and every finite
disjoint $u,w \subseteq P^M,(\exists n < \omega)(\exists A \in
{\Cal P}_n \backslash {\Cal P}_k)[u \subseteq A \wedge A \cap 
w = \emptyset]$ because $\dbcu_{n < \omega} {\Cal P}_n$ enables us 
to define $N \in K^1_{\aleph_0}$ such that $M 
\le_{{\frak K}^1} N$ as required.  
Assume ${\Cal P}_n$ is defined;  apply Claim
\scite{88r-6.5}(c) on ${\Cal P}_n$ when substituting $u = s_{2n},w = s_{2n+1}$
let $A \subseteq P$ be supplied by the Claim and define ${\Cal
P}_{n+1} = {\Cal P}_n \cup \{A\}$.  It is easy to check that 
$\{{\Cal P}_n:n < \omega\}$ satisfies
our requirements. \nl
3) This is a private case of Ax III which we 
checked in Lemma \scite{88r-6.3}(3). \nl
4) Let $M \in K^1_{\aleph_0}$ and we shall find $M_\ell \in K^1_{\aleph_0}
(\ell=0,1),M \le_{{\frak K}^1} M_\ell$, which cannot be amalgamated
over $M$.
By part (2) a model $M_1$ such that $M <_{{\frak K}^1} M_1 \in
K^1_{\aleph_0}$ and choose $y \in Q^{M_1} \backslash Q^M$.  Define
$M_2 \in K^1_{\aleph_0}$; its universe is $|M_1|,P^{M_2} =
P^{M_1},Q^{M_2} = Q^{M_1}$ and $R^{M_2} = \{(a,b):aR^{M_1} b \and b
\ne y$ or $a \in P^M \and b = y \and \neg(aRy)\}$.
Clearly $M_1,M_2$ cannot be amalgamated over $M$ 
(since the amalgamation must contain a set and 
its complement).  \hfill$\square_{\scite{88r-6.4}}$
\enddemo
\bigskip

\proclaim{\stag{88r-6.6} Theorem}  Assume ${\text{\rm MA\/}}_{\aleph_1}$ (hence
$2^{\aleph_0} > \aleph_1$).  The class $(K^1,<_{{\frak K}^1})$ is categorical in
$\aleph_1$. 
\endproclaim
\bigskip

\demo{Proof}  Let $M,N \in K^1_{\aleph_1}$ and we shall prove that they
are isomorphic.  By repeated use of the
idea in the proof of Lemma \scite{88r-6.3}(6) for Ax.VI we get (strictly)
increasingly continuous chains $\{M_\alpha:\alpha <
\omega_1\},\{N_\alpha:\alpha < \omega_1\} \subseteq K^1_{\aleph_0}$
such that $M = \dbcu_{\alpha < \omega_1} M_\alpha$ and $N =
\dbcu_{\alpha < \omega_1} N_\alpha$ such that for $\alpha <
\beta,M_\alpha <_{{\frak K}^1} M_\beta,N_\alpha <_{{\frak K}^1} N_\beta$.

Now define a forcing notion which supplies an isomorphism $g:M
\rightarrow N$.

$$ 
\align
\Bbb P = \{f:&f \text{ is a partial finite isomorphism from } M 
\text{ into } N \text{ satisfying} \\
  &(\forall \alpha < \omega_1)(\forall a \in \text{ Dom}(f))
[a \in M_\alpha \Leftrightarrow f(a) \in N_\alpha]\},
\endalign
$$
\mn
the order is inclusion.  It is trivial to check that if $G \subseteq
\Bbb P$ is a directed subset 
then $g = \cup G$ is a partial isomorphism from $M$ to
$N$, we show that Dom$(g) = |M|$ if $G$ is generic enough.  
For every $a \in |M|$ define ${\Cal J}_a = \{f \in \Bbb P:a 
\in \text{ Dom}(f)\}$, and we shall show that for all $a \in |M|$ the 
set ${\Cal J}_a$ is dense.  For
$a \in M$ let $\alpha(a) = \text{ Min}\{\alpha < \omega_1:
a \in M_\alpha\}$, clearly it is zero or a successor ordinal.  
Let $f \in \Bbb P$ be a given condition,
it is enough to find $h \in {\Cal J}_a$ such that $f \subseteq h$ and $a
\in \text{ Dom}(h)$.  Let $A = \text{ Dom}(f)$, let $B,C \subseteq
A$ be disjoint sets such that $B \cup C = A$ and $B = \text{ Dom}(f) \cap P^M,C
= \text{ Dom}(f) \cap Q^M$.  Without loss of generality $a
\notin B \cup C$.  If $a \in P^M$ let $\varphi(x,\bar c) = \wedge\{\pm xRc:c
\in C$ and $M \models \pm aRc\}$.
\relax From the definition of $K^1$ there exists $b \in P^N \backslash \text{
Rang}(f)$ such that $N \models \varphi[b,f(\bar c)]$.  
If $a \in Q^M$ let $\varphi(x,\bar b) = \wedge\{ \pm bRx:b \in
B,M \models \pm bRa\}$, we can find infinitely many 
$b \in Q^{N_{f(a)}}
\backslash \dbcu_{\beta < f(a)} N_\beta$, satisfying
$\varphi(x,f(\bar b))$.  
\nl
Why?  This is as $\cup\{N_\beta:\beta < \alpha(a)\} <_{{\frak K}^1}
N_{\alpha(a)}$ as $C$ is finite \wilog \, $b \notin f(C)$.    
\nl
Finally, let $h = f \cup \{\langle a,b \rangle\}$.  

The proof that Range$(g) = |N|$ is analogous to the proof that Dom$(g)
= |M|$.  
In order to use ${\text{\rm MA\/}}$ we just have to
show that $R$ has the c.c.c.  Let $\{f_\alpha:\alpha < \omega_1\}
\subseteq R$ be given.  It is enough to find $\alpha,\beta < \omega_1$
such that $f_\alpha,f_\beta$ have a common extension.  Without loss of
generality we may assume $|M| \cap |N| = \emptyset$.  By the finitary
$\Delta$-system lemma there exists $S \subseteq \omega_1,|S| =
\aleph_1$ such that $\{\text{Dom}(f_\alpha) \cup \text{
Range}(f_\alpha):\alpha \in S\}$ is a $\Delta$-system with heart $A$.
Let $B \subseteq |M|,C \subseteq |N|$ be such that $A = B \cup C$, now
\wilog \, for every $\alpha \in S,f_\alpha$ maps $B$ into $C$.
\nl
[Why?  If not, $S_1 = \{\alpha \in S$: for some $b=b_\alpha \in
B,f_\alpha(b_\alpha) \notin C\}$ is uncountable hence for some $b \in
B,S_2 = \{\alpha \in S_1:b_\alpha =b\}$ is uncountable; so $\langle
f_\alpha(b):\alpha \in S_2 \rangle$ is without repetitions hence is
uncountable.  But
$\{f(b):f \in \Bbb P$ and $b \in \text{ Dom}(f) \cap B\}$ is countable
because $f \in \Bbb P \and b \in \text{ Dom}(f) \and \alpha < \omega_1
\Rightarrow [b \in M_\alpha \equiv f(b) \in N_\alpha]$.  
Similarly, $f^{-1}_\alpha$ maps $C$ into $B$, so necessarily
$f_\alpha$ maps $B$ onto $C$; but the
number of possible functions from $B$ to $C$ is $|C|^{|B|} <
\aleph_0$.  Hence there exists $S_1 \subseteq S,|S_1| = \aleph_1$ such
that for all $\alpha,\beta \in S_1,f_\alpha \restriction B = f_\beta
\restriction B$ and Dom$(f_\alpha) \cap M_0 \subseteq B$,
Rang$(f_\alpha) \cap N_0 \subseteq C$. As $P^{M_\alpha} = P^{M_0}
\subseteq M_0,P^{N_\alpha} = P^{N_0} \subseteq N_0$ for every
$\alpha \in S_1$ we have $P^M \cap \text{ Dom}(f_\alpha) \subseteq B,P^N
\cap \text{ Range}(f_\alpha) \subseteq C$, therefore for all $\alpha,\beta \in
S_1,f_\alpha \cup f_\beta \in \Bbb P$ and in particular there exists
$\alpha \ne \beta < \omega_1$ such that $f_\alpha \cup f_\beta \in
\Bbb P$.  \hfill$\square_{\scite{88r-6.6}}$
\enddemo
\bn
In the terminology of \cite{GrSh:174} Theorems \scite{88r-6.4} and
\scite{88r-6.6} give us together:
\demo{\stag{88r-6.7} Conclusion}  Assuming $2^{\aleph_0} > \aleph_1$ and
MA$_{\aleph_1},{\frak K}^1$ is a nice category which has a universal
object in $\aleph_1$, moreover it is categorical in $\aleph_1$.
\enddemo
\bigskip

\definition{\stag{88r-6.8} Definition}  1) $K^2$ is the class of $M \in
K^0$ (see Definition \scite{88r-6.2}) satisfying:
\mr
\item "{$(a)$}"  $(\forall x \in Q^M)(\forall u \in [P^M]^{< \aleph_0})
(\exists y \in Q)[A^M_x \Delta A^M_y = u]$
\sn
\item "{$(b)$}"  if $k < \omega$ and $y_0,\dotsc,y_{k-1} 
\in Q$ satisfies $|A_{y_\ell} \Delta A_{y_m}| \ge \aleph_0$
for $\ell < m < k$ \ub{then} the set $\{A^M_{y_\ell}:\ell < k\}$ is an
independent family of subsets of $P^M$
\sn
\item "{$(c)$}"  $Q(y) \wedge Q(z) 
\wedge (\forall x \in P)[x R y \leftrightarrow x R z] \rightarrow y =
z$,
\sn
\item "{$(d)$}"  for every $k < \omega$ for some $y_0,\dotsc,y_k \in
Q^M$ we have $\dsize \bigwedge_{\ell < m \le k} |A_{y_\ell} \Delta A_{y_m}| \ge
\aleph_0$. 
\ermn
2) For $M_1,M_2 \in K^2$

$$
M_1 \le_{{\frak K}^2} M_2 \Leftrightarrow^{\text{df}} 
M_1 \subseteq M_2,P^{M_1} = P^{M_2}.
$$
\mn
3) ${\frak K}^2 = (K^2,\le_{{\frak K}^2})$. \nl
4) $K^3$ is the class of models $M = (|M|,P^M,Q^M,R^M,E^M)$ such that
\mr
\item "{$(a)$}"  $(|M|,P^M,Q^M,R^M) \in K^1$
\sn
\item "{$(b)$}"  $E^M$ is an equivalence relation on $Q^M$
\sn
\item "{$(c)$}"  $E^M$ has infinitely many equivalence classes
\sn
\item "{$(d)$}"   each equivalence class of $E^M$ is countable
\sn
\item "{$(e)$}"  if $u,w \subseteq P^M$ are finite disjoint and $y \in
Q^M$ then for some $y' \in y/E^M$ we have $a \in u \Rightarrow aR^M
y'$ and $b \in w \Rightarrow \neg(b R^M y')$. 
\ermn
5) We define $\le_{{\frak K}^3}:M_1 \le_{{\frak K}^3} M_2 
\Leftrightarrow^{\text{df}} M_1 \subseteq M_2$ and $a \in M_1
\Rightarrow a/E^{M_2} = a/E^{M_1}$. 
\nl
6) ${\frak K}^3 = (K^3,\le_{{\frak K}^3})$.
\enddefinition
\bn
If we like to have a class defined by a sentence from $\Bbb
L_{\omega_1,\omega}$ (rather than $\Bbb L_{\omega_1,\omega}(\bold Q)$)
we can use:
\definition{\stag{88r-6.8A} Definition}  1) ${\frak K}^4$ is defined as
follows:
\mr
\item "{$(A)$}"  $\tau({\frak K}^4) = \{P,Q,R\} \cup \{P_n:n <
\omega\}$, $R$ two-place predicates, $P,Q,P_n$ are unary predicates
\sn
\item "{$(B)$}"  $M \in K^4$ \ub{iff} $M$ is a 
$\tau({\frak K}^4)$-model such that $M \restriction \{P,Q,R\} \in K^2$ and
{\roster
\itemitem{ $(a)$ }  $\langle P^M_n:n < \omega \rangle$ is a partition
of $P^M$
\sn
\itemitem{ $(b)$ }  $P^M_n$ has exactly $2^n$ elements
\sn
\itemitem{ $(c)$ }  $(\forall x \in Q)(\forall u \in [P^M]^{< \aleph_0})
(\exists y \in Q^M)[A^M_x \Delta A^M_y = u]$
\sn
\itemitem{ $(d)$ }  if $k < \omega$ and $y_0,\dotsc,y_{k-1} 
\in Q$ satisfies $|A_{y_\ell} \Delta A_{y_m}| \ge \aleph_0$
for $\ell < m < k$ \ub{then} the set $\{A^M_{y_\ell}:\ell < k\}$ is an
independent family of subsets of $P^M$; moreover for any $n$ large
enough for any $\eta \in {}^k 2$ the set $P^M_n \cap
\cap\{A^M_{y_\ell}:\eta(\ell)=1\} \backslash \cup
\{A^M_{y_\ell}:\eta(\ell) =0\}$ has exactly $2^{n-k}$ elements
\sn
\itemitem{ $(e)$ }  $Q^M(y) \wedge Q^M(z) 
\wedge (\forall x \in P^M)[x R^M y \leftrightarrow x R^M z] \rightarrow y = z$,
\sn
\itemitem{ $(f)$ }  for every 
$k < \omega$ for some $y_0,\dotsc,y_k \in Q^M$ we have $\dsize 
\bigwedge_{\ell < m \le k} |A_{y_\ell} \Delta A_{y_m}| \ge \aleph_0$
\endroster}
\item "{$(C)$}"  $M \le_{{\frak K}^4} N$ iff $M,N \in K^4$ and $M
\subseteq N$ and $P^M = P^N$.
\endroster
\enddefinition
\bigskip

\proclaim{\stag{88r-6.9} Theorem}  1) $(K^2,<_{{\frak K}^2})$ is an
$\aleph_0$-presentable abstract elementary class which is categorical 
in $\aleph_0$. \nl
2)  Also ${\frak K}^3$ and ${\frak K}^4$  are
$\aleph_0$-presentable a.e.c. categorical in $\aleph_0$.
\endproclaim
\bigskip

\demo{Proof}  Similar to the proof for ${\frak K}^1$.
\enddemo
\bigskip

\proclaim{\stag{88r-6.11} Theorem}  1) ${\frak K}^1_{\aleph_1}$ has an
axiomatization in $\Bbb L(\bold Q)$ and $\le_{{\frak K}^1}$ is 
$<^{**}$ from the proof of \scite{88r-3.9} (this is $<^{**}$ from
\cite{Sh:87a} and \cite{Sh:87b}). \nl
2) ${\frak K}^2$ has an axiomatization in 
$\Bbb L_{\omega_1,\omega}(\bold Q)$ and $\le_{{\frak K}^2}$
is $\le^*$ from the proof of \scite{88r-3.9} (this is
$<^*_{\omega_1,\omega}$ from \cite{Sh:87a} and \cite{Sh:87b}). \nl
3) ${\frak K}^3$ has an axiomatization in $\Bbb L(\bold Q)$ 
and $\le_{{\frak K}^3}$ is $<^*$
from \cite{Sh:87a} and \cite{Sh:87b}. \nl
4) ${\frak K}^4$ has an axiomatization in $\Bbb L_{\omega_1,\omega}$ and
$\le_{{\frak K}^4}$ is just being a submodel. \nl
5) $(\forall \ell \in \{1,2,3,4\})[K^\ell$ is {\rm PC}$_{\aleph_0}]$. 
\endproclaim
\bigskip

\proclaim{\stag{88r-6.11A} Theorem}  If ${\text{\rm MA\/}}_{\aleph_1}$ 
\ub{then} $K^\ell$ is categorical in $\aleph_1$ for $\ell = 2,3,4$.
\endproclaim
\bigskip

\demo{Proof}  Easy.
\enddemo
\bigskip

\demo{\stag{88r-6.12} Conclusion}  Assuming MA$_{\aleph_1}$ there exists an
abstract elementary class, which is PC$_{\aleph_0}$, categorical in
$\aleph_0,\aleph_1$ but without the $\aleph_0$-amalgamation property;
moreover, there is such class axiomatize by some $\psi \in 
\Bbb L_{\omega_1,\omega}$.
\enddemo
\bn
\margintag{88r-6.13}\ub{\stag{88r-6.13} Example}:  There is a complete first order which has a
model-homomorphic model which is not sequence-homomorphic $T$, countable
theory $T$ with (i.e., $N \prec M^* \Rightarrow N = M$).
\bigskip

\demo{Proof}  Let $\Bbb P$  be the set of pairs $p = (u,{\Cal
F}_\alpha)_{\alpha \le \omega} = (u^{\Cal P},{\Cal
F}^P_\alpha)_{\alpha \le\omega}$
\mr
\item "{$\circledast_1$}"  $(a) \quad u$ is a finite subset of
$\omega$
\sn
\item "{${{}}$}"  $(b) \quad {\Cal F}_\alpha$ is a set of partial one
to one functions from $u$ to $u$
\sn
\item "{${{}}$}"  $(c) \quad {\Cal F}_\alpha$ is closed under
composition and inverse 0 and contains id$_v$ for any $v \subset u$
\sn
\item "{${{}}$}"  $(d) \quad {\Cal F}_\alpha$ decreases with $\alpha$
\sn
\item "{${{}}$}"  $(e) \quad$ there are $f \in {\Cal F}_\omega$ and $a
\in \text{ Dom}(f)$ such that $f(a) \ne a$.
\ermn
On $\Bbb P$ we define an order $\le_{\Bbb P}$
\mr
\item "{$\circledast_2$}"  $p \le_{\Bbb P} q$ iff: $(a) \quad u^p
\subseteq u^q$
\sn
\item "{${{}}$}"  $\qquad (b) \quad {\Cal F}^p_\alpha \subseteq {\Cal
F}^q_\alpha$
\sn 
\item "{${{}}$}"  $\qquad (c) \quad$ if $\beta \le \alpha,f \in
F^q_\beta \backslash F^p_\beta$ then for some $g \in F^p_\beta$ we
have $f \subseteq g$ moreover $g(a) = b \wedge a \in u^p \wedge b \in
u^p \Rightarrow f(a) = b$.
\ermn
Clearly we can find $\langle p_n:n < \omega \rangle$ such that
\mr
\item "{$\circledast_3$}"  $(a) \quad p_n \in \Bbb P$
\sn
\item "{${{}}$}"  $(b) \quad p_n \le p_{n+1}$
\sn
\item "{${{}}$}"  $(c) \quad$ if $f \in {\Cal F}^{p_n}_\alpha$ and $a
\in u^{p_n}$ and $\beta < \alpha$ then for some $m>m$ and $g \in {\Cal
F}^{p_n}_\alpha$ we have $f \subseteq g \wedge a \in \text{ Dom}(g)$
\sn
\item "{${{}}$}"  $(d) \quad \omega = \cup\{u^{p_n}:n < \omega\}$.
\ermn
Let $E_{n,m}$ be the folloiwng equivalence relation on $m_\omega:\bar
a E_{n,m} \bar b$ iff $\bar a,\bar b \in E {}^m \omega$ and for some
$f \in \cup\{{\Cal F}^{p_k}_n:k < \omega\}$ maps $\bar a$ to $\bar b$.

For any $\lambda$ we define a model $M_\lambda$:

its universe: $\lambda \times \omega$

relations: if $n,m < \omega$ and $e$ is an $E_{n,m}$-equivalence class
then 

$$
\align
R_e = \{\langle(\gamma,\ell_v),\dots,(\gamma,\ell_{m-1})\rangle:&\gamma <
\lambda \\
  &\text{ and } \langle \ell_0,\dotsc,\ell_{m-1}\rangle \in e\}
\endalign
$$

$$
E = \{(\gamma,n_1(\gamma,n_2)):\gamma < \lambda \text{ and } n_1,n_2 <
\omega\}.
$$
\mn
Note
\mr
\item "{$\circledast_4$}"  $(a) \quad$ if $f \in {\Cal
F}^{p_0}_\omega,f(a) = b \ne a$ (this occurs!) \ub{then}
$(\alpha,a),(\alpha,b) \in M_\lambda$ realizes the same type in
$M_\lambda$
\sn
\item "{${{}}$}"  $(b) \quad$ no automorphism of $M_\lambda$ maps any
$a \in M_\lambda$ to $b \in a/E^{M_\lambda} \backslash \{a\}$;
moreover, $M_\lambda \restriction (a) E^{M_\lambda}$ has no
automorphism excep the identity
\sn
\item "{${{}}$}"  $(c) \quad$ if $N \prec M_\lambda$ then $a \in N
\Leftrightarrow a/E^{M_\lambda} \subseteq N$ and $\{a/E:a \in N\}$ is
infinite so $(N)$ has the form $A \times \omega,A \subseteq \lambda$
infinite
\sn
\item "{${{}}$}"  $(d) \quad$ if $\pi$ is a partial one to one
function from $A \subseteq \lambda$ onto $B \subseteq \lambda$ and
$\hat \pi$ is the function mapping $(\alpha,n) \in A \times \omega$ to
$(\pi(\alpha),n) \in B \times \omega$ then $\hat \pi$ is an
isomorphism from $M_\lambda \restriction (A \times \omega)$ onto
$M_\lambda \restriction (B \times \lambda)$
\sn
\item "{${{}}$}"  $(e) \quad$ if $N_\ell \prec M_\lambda$ and
$|N_\ell| = A_\ell \times \lambda$ and $f$ is an isomorphism from
$N_1$ onto $N_2$ then $f = \hat \pi$ for some one to one function from
$A_1$ onto $A_2$.
\endroster
\enddemo 
\newpage

\nocite{ignore-this-bibtex-warning} 
\newpage
    
REFERENCES.  
\bibliographystyle{lit-plain}
\bibliography{lista,listb,listx,listf,liste}

\enddocument